\documentclass[smallextended,natbib,runningheads]{svjour3}
\journalname{Annals of the Institute of Statistical Mathematics}
\smartqed  
\usepackage{graphicx}
\usepackage{hyperref}
\usepackage{verbatim}
\usepackage{placeins}
\begin{document}

\title{Building A Theoretical Foundation for Combining Negative Controls and Replicates
\thanks{TPS would like to thank the Institute of Mathematical Sciences of the National University of Singapore for the invitation to speak at their January 2022 Workshop on Statistical Methods in Genetic/Genomic Studies at which a rudimentary version of the ideas in this paper were presented. 
Thanks also to Professor Bin Yu for her very helpful comments on an earlier version of this paper. 
Jiming Jiang's research was supported by the National Science Foundation of the United States grant DMS-1713120. The authors are grateful to a referee, whose insightful comments have helped improve the manuscript.}
}


\author{Jiming Jiang         \and
        Johann A. Gagnon-Bartsch \and
        Terence P. Speed
}


\institute{J. Jiang \at
              Department of Statistics\\
              University of California\\
              1 Shields Avenue,\\
              Davis, CA 95656, U. S. A.\\
              \email{jimjiang@ucdavis.edu}           
           \and
           J.A. Gagnon-Bartsch \at
              Department of Statistics, University of Michigan\\
              323 West Hall, 1085 South University,\\
              Ann Arbor, MI 48109-1107, U.S.A.\\
              \email{johanngb@umich.edu}
           \and
           T.P. Speed \at
              Walter and Eliza Hall Institute of Medical Research\\
              1G Royal Parade\\
              Parkville, VIC 3052, Australia\\
              \email{terry@wehi.edu.au}
}

\date{}

\maketitle

\begin{abstract}
Studies using assays to quantify the expression of thousands of genes on tens to thousands of cell samples have been carried out for over 20 years. Such assays are based on microarrays, DNA sequencing or other molecular technologies. All such studies involve unwanted variation, often called batch effects, associated with the cell samples and the assay process.  Removing this unwanted variation is essential before the measurements can be used to address the questions that motivated the studies. Combining the results of replicate assays with measurements on negative control genes to estimate the unwanted variation and remove it has proved to be effective at this task. The main goal of this paper is to present asymptotic theory that explains this effectiveness. The approach can be widened by using pseudo-replicate sets of pseudo-samples, for use with studies having no replicate assays. Theory covering this case is also presented. The established theory is supported by results of empirical investigations, including simulation studies and a real-data example.
\end{abstract}



\section{Introduction}

\subsection{Background and some history}

Reproducibility of results has long been a concern in the field of omics research, with unidentified batch effects being a major contributor to irreproducibility. These issues have been extensively discussed, with 
\citet{shi2006microarray} and \citet{leek2010tackling} being early and influential contributors. There have been many papers on both topics in the years since those papers, and it seems fair to say that the last word on dealing with irreproducibility due to batch effects will never be said. 
Nevertheless, we believe that real progress is being made, and that our removing unwanted variation (RUV) methods \citep{gagnonbartsch2012using, jacob2016correcting}, are contributing to this progress.

The 3-step method RUV-III  we discuss below was introduced in \citet{molania2019new}, where we saw clear empirical evidence of the value of technical replicates and pseudo-replicates and negative controls for library size normalization and batch effect removal with Nanostring nCounter mRNA abundance estimation. More recently \citet{molania2022removing} showed the power of the approach for normalizing large RNA-seq gene expression studies when no suitable technical replicates or pseudo-replicates could be defined, a not uncommon situation. In such cases suitable pseudo-samples and pseudo-replicate sets of pseudo-samples were defined, and the normalization proceeded as though these were genuine technical replicate sets of genuine samples. Not only was this strategy effective, it worked well in cases where genuine technical replicates existed but were not used. It also worked in cases where technical replicates are not adequate for removing certain unwanted variation, such as is the case with tumour purity in cancer studies. More broadly, the approach has proved to be effective for normalizing and removing batch effects from large proteomic datasets \citep{poulos2020strategies} and single-cell gene expression datasets \citep{lin2023atlas}. It is very effective in the omics quite broadly.

What is missing? To date the theory underlying the approach in its different forms has largely been heuristic. The method has an underlying linear model and there is a compelling logic for the steps in its implementation. However, its credentials have been entirely scientific, embodied in the individual cases to which it has been applied. Arguments for accepting the results make use of positive controls and other known biological facts about the system under study, together with a variety of formal and informal summaries and plots. The latter show the extent to which the known biology in the data has been preserved or enhanced, and the known unwanted variation has been reduced or removed. What has been lacking is exact or asymptotic theory which proves, under a framework embodying realistic assumptions on the model, that our use of replicates and negative controls does indeed capture and remove the unwanted variation term in the model that we seek to remove. In brief, having shown the approach works in practice, we want to show that it works in theory. That is the main contribution of this paper. It is worth pointing out that we know of no other normalization and batch effect removal procedure which comes with such theory. The published work nearest to the present paper is \citet{wang2017confounder}, which provides asymptotic theory for several methods, including RUV-4 \citep{gagnonbartsch2013removing, jacob2016correcting}, but in that work the focus is on consistent asymptotically normal estimates of regression parameters in the presence of unwanted variation, not on estimating and removing the unwanted variation itself.

With that background, we turn to the present paper. The approach was developed for normalizing datasets consisting of thousands of gene expression measurements from tens to thousands of assays done on collections of cells obtained from tumour biopsies or brain or other tissue samples. Although more generally applicable, we will keep this context in mind in the discussion which follows. Our use of the term unwanted variation is broad, encompassing what are known as batch effects due to sample collection, duration and conditions of sample storage, assay details, including reagents, equipment, operators; see \citet{stephensstatistica} for more background on this type of application. We can also include unwanted biological features such as the extent of normal cell contamination in tumour tissue. The term normalization refers to adjusting the measurements in such gene expression datasets with a view to removing the effects of unwanted variation; this is the goal of RUV-III. Before we introduce the linear model which underlies the approach, we will give an informal description of the method and summarize some relevant history. It combines two notions that were used for removing unwanted variation in agronomy in the late 19th and early 20th century: replicates, and what we now call negative control measurements.

In the late 1800s, plant breeders began interspersing plots containing a single check variety among the plots containing un-replicated test varieties, and adjusting the yields of the test varieties using the yields from nearby check plots, see \citet{holtsmark1905om}, \citet{thorne1907interpretation}, and 
\citet{kempton1984design}. This practice was almost killed by \citet{yates1936new}, who showed the advantages of Fisher’s randomized block experiments and Yates' incomplete block experiments; see also \citet{atiqullah1962use}. However, the need for large un-replicated trials has remained and the use of check plots with a replicated variety has continued in agronomy to this day \citep{haines2021augmented}. More recently the same idea has been used in other contexts where replication is difficult or expensive: proteomics \citep{neve2006collection}, microarrays \citep{walker2008empirical}, and mass cytometry \citep{van2020cytonorm, schuyler2019minimizing}. In addition to being called check varieties, other names for the replicated units include standard, reference, control, and anchor samples.

By the term negative controls we mean variables that should have no true association with the factors of interest in a study. For example, they could be bacterial DNA probes spiked-in on a mammalian cell microarray or so-called housekeeping genes. For such variables, all variation is unwanted. An agronomic analogue is the field uniformity trial, where an area is divided into plots and the same variety of a crop is grown under the same conditions on every plot. The yield of each plot is then recorded at harvest. That yield may be regarded as a negative control, as there is no factor of interest. Such trials played a valuable role in agronomic experimentation in the early part of the 20th century \citep{cochran1937catalogue}, and they have continued to be carried out \citep{jones2021automating}. The main use of uniformity trials is to learn the spatial nature and magnitude of unwanted variation in order to determine the best size and shape of plots \citep{smith1938empirical}. However, \citet{sanders1930note} analyzed the published results of uniformity trials carried out on two fields in Denmark between 1906 and 1911 to see whether soil variations were sufficiently stable from year to year to the permit removing unwanted variation from the yields of experimental plots using their yields from previous uniformity trials. In doing so Sanders used the analysis of covariance. His paper was the first publication of this technique, which appeared two years later as §49.1 of \citet{fisher1932statistical}.  Sanders’ results for annual crops such as oats, barley and rye were mixed, but within a few years his approach was found to be very effective for perennials: tea trees in Ceylon (now Sri Lanka) \citep{eden1931studies}, cacao in the West Indies \citep{cheesman1932uniformity}, and rubber trees in Sumatra \citep{murray1934value}. Although Sanders’ method does not seem to be used today in agronomy, it is striking that his rationale is identical to that adopted in \citep{lucas2006sparse} over 70 years later, who wrote ``The model in [1] utilizes multiple principal components of the set of normalisation control probes, and a selection of the housekeeping probes, as covariates.''

\subsection{Overview of the approach and our contribution}

As we will see, our approach takes the deviations (residuals) of replicate gene expression measurements (on the log scale) from their averages and combines them across replicate sets to give estimates, denoted by $\hat{\alpha}$, of the gene-level responses to unwanted variation. Centered negative control measurements are then regressed against these estimates to obtain estimates, denoted by $\hat{W}$, of unobserved sample-level covariates for unwanted variation. These are not refitted but simply multiplied by the estimated gene-level responses, and the resulting product $\hat{W}\hat{\alpha}$ is subtracted from the data. This idea is executed below via a linear model, and some linear algebra.

From now on we will call the method we discuss RUV-III, to be broadly consistent with the usage in our earlier publications on the topic \citep{molania2019new, molania2022removing}. Before we can describe the linear model underlying RUV-III we need to introduce an $m \times s$ mapping matrix, $M=[M(i,h)]_{1\leq i\leq m,1\leq h\leq s}$, connecting assays to distinct samples, which captures the pattern of replication in our assays.  Here, $m$ is the number of assays and $s$ the number of distinct samples being assayed. Let $M(i,h)=1$ if assay $i$ is on sample $h$, and $M(i,h)=0$ otherwise. There is exactly one $1$ in each row of $M$. Each column of $M$ sums to a distinct sample replication number. By rearranging the order of the assays according to the distinct samples, $h=1,\dots,s$, one can rewrite $M$ as
\begin{eqnarray}
M={\rm diag}(1_{m_{h}}, 1\leq h\leq s),\label{eq:M_expression}
\end{eqnarray}
where $1_{d}$ is the $d\times 1$ vector of $1$s, $m_{h}$ is the total number of assays that are on sample $h, 1\leq h\leq s$ (so that $\sum_{h=1}^{s}m_{h}=m$), and ${\rm diag}(M_{h}, 1\leq h\leq s)$ denote the block-diagonal matrix with $M_{h}, 1\leq h\leq s$ on the diagonal. We also define an $s\times p$ matrix, $X$, to capture the biological factor(s) of interest indexed by sample rather than assay. There is no restriction on $p$; indeed $X$ could be $I_{s}$, the $s\times s$ identity matrix. More generally, a linear mixed model [LMM; e.g., \citep{jiang2021linear}] can be expressed as
\begin{eqnarray}
Y=1_{m}\mu+MX\beta+W\alpha+\epsilon,\label{eq:LMM}
\end{eqnarray}
where $Y=(y_{ij})_{1\leq i\leq m, 1\leq j\leq n}$ and $\epsilon=(\epsilon_{ij})_{1\leq i\leq m, 1\leq j\leq n}$ are the data matrix and matrix of unobserved errors, respectively; $\mu$ is the $1\times n$ row vector of gene means; $\beta$ is an $p\times n$ matrix of unknown regression parameters; $W$ is an $m\times k$ matrix whose columns capture the unwanted variation, and $\alpha$ is a $k\times n$ random matrix whose $n$ columns are i.i.d. with mean $0$ and covariance matrix $\Sigma$. It is assumed that $W\perp 1_{m}$, that is, $1_{m}'W=0$.

Also, we suppose that we have a subset of $n_{\rm c}$ negative control genes whose corresponding $m\times n_{\rm c}$ submatrix of $Y$, $Y_{\rm c}$, satisfies (\ref{eq:LMM}) with the corresponding submatrix of $\beta$, $\beta_{\rm c}=0$, that is, we have
\begin{eqnarray}
Y_{\rm c}=1_{m}\mu_{\rm c}+W\alpha_{\rm c}+\epsilon_{\rm c},\label{eq:LMM_nc}
\end{eqnarray}
where $\alpha_{\rm c}$ and $\epsilon_{\rm c}$ are the corresponding submatrices of $\alpha$ and $\epsilon$, respectively; in other words, there is no true association between the negative control genes and the biology of interest.

The projection matrix, $P_{M}=M(M'M)^{-1}M'$, replaces the entries of $Y$, $y_{ij}$, by the averages of the assays on the same distinct sample, that is, the average of $y_{i'j}$ such that $y_{i'j}$ and $y_{ij}$ are replicate assays of the same distinct sample. Note that $M'M$ is non-singular by construction (we are not considering the possibility of samples that are never assayed at all).
Let $P_{M^{\perp}}=I_{m}-P_{M}$ be the orthogonal projection complementary to $P_{M}$. If the replication is technical at some level, then $P_{M^{\perp}}Y$ mainly contains information about unwanted variation in the system after the technical replicates were created. Depending on the study details, technical replicates could be created immediately before the assay was run, in parallel with or immediately after sample was collected, or somewhere in between. The earlier the creation of technical replicates the more unwanted variation will be captured in their differences.

Now consider the eigenvalue decomposition
\begin{eqnarray}
P_{M^{\perp}}YY'P_{M^{\perp}}=UDU',\label{eq:eigen_decom}
\end{eqnarray}
where $U$ is an $m\times m$ orthogonal matrix and $D$ is an $m\times m$ diagonal matrix with the diagonal elements being the eigenvalues of the left side of (\ref{eq:eigen_decom}), ordered from the largest to the smallest. Let $P_{1}=P_{1_{m}}=1_{m}1_{m}'/m$ and $P_{1^{\perp}}=I_{m}-P_{1}$. For any $1\leq k\leq m-s$, define $\hat{\alpha}_{(k)}=U_{(k)}'Y$, where $U_{(k)}$ is the submatrix consisting of the first $k$ columns of $U$, and
\begin{eqnarray}
\hat{W}_{(k)}&=&P_{1^{\perp}}Y_{\rm c}\{U_{(k)}'Y_{\rm c}\}'\{(U_{(k)}'Y_{\rm c})(U_{(k)}'Y_{\rm c})'\}^{-1}
\nonumber\\
&=&P_{1^{\perp}}Y_{\rm c}Y_{\rm c}'U_{(k)}\{U_{(k)}'Y_{\rm c}Y_{\rm c}'U_{(k)}\}^{-1}.\label{eq:Wkhat}
\end{eqnarray}
$\hat{W}_{(k)}$ may be viewed as an estimator of $W$ obtained by regressing the centered negative controls, $P_{1^{\perp}}Y_{\rm c}$, on $\hat{\alpha}_{(k){\rm c}}'$, where $\hat{\alpha}_{(k){\rm c}}$ is the c part of $\hat{\alpha}_{(k)}$.

Note that the definition in the last paragraph does not require that $k$ is the true dimension of $\alpha$, which may be unknown (see Section 3). The unwanted variation is estimated, for the given $k$, by $\hat{W}_{(k)} \hat{\alpha}_{(k)}$ so that, after RUV, the adjusted $Y$ is $Y_{(k)}=Y-\hat{W}_{(k)}\hat{\alpha}_{(k)}$. We anticipate that, with a sufficiently large $k$, we will have $\hat{W}_{(k)}\hat{\alpha}_{(k)}\approx W\alpha$, hence the adjustment, $Y_{(k)}$, achieves the goal of RUV. A couple of points are noteworthy. First, the dimension of the unwanted variation, $W\alpha$, is $m\times n$, and both $m$ and $n$ are increasing under our (realistic) asymptotic setting, so the approximation $\approx$ here needs to be defined in a suitable sense. Second, $k$ is typically unknown in practice. In Section 2, we first consider a simpler case in which $k$, the true dimensionality of $\alpha$, is known. The more realistic situation where $k$ is unknown is considered in Section 3. Specifically, here we show that RUV-III has very different asymptotic behaviors when $k<k_{0}$, and when $k\geq k_{0}$, where $k_{0}$ is the true dimension of $\alpha$. Although the goal of RUV is achieved for any $k\geq k_{0}$, because the latter is unknown, we propose to choose $k=K$, a sufficiently large upper bound, such as $K=m-s$, as long as it is computationally feasible. In fact, the asymptotic theory can justify the use of any given sequence of $k$, say $k_{N}$, with the property that $k_{N}\geq k_{0}$ with probability tending to one.

Section 4 concerns what we call pseudo-replicates of pseudo-samples, abbreviated by PRPS. It is a general practice that significantly improves computational efficiency of RUV-III. We show that similar asymptotic results hold for PRPS as well. Some empirical results, including simulation studies and a real-data example are provided in Section 5. Proofs of the theoretical results are provided in Section \ref{app:proofs}.

The biggest technical challenge has to do with having a vision, combined with skills in matrix algebra. The vision provided the big-picture ideas, and the skills execute those ideas. In terms of the insights and practical implications of the theoretical results, a highly interpretable result shows that there is a change point for the candidate values of $k$, the dimension of the unwanted variation, and the change point is equal to the true dimension of the unwanted variation. In a way, this is similar to the consistency property of some information criteria, such as the BIC, in model selection. The change point separates the under-fitting models and over-fitting models, whose impacts on the information criterion are of different orders. Thus, as in the information criteria, our asymptotic theory may provide guidance to the practice.
\section{Asymptotic analysis when $k$ is known}

Let us first consider a simpler case, in which the dimension of the random effects representing the
unwanted variation, $\alpha$, is known, that is, $k={\rm dim}(\alpha_{j})$ ($1\leq j\leq n$) is known, where
$\alpha_{j}$ is the $j$th column of $\alpha$.

As noted, the estimated unwanted variation, $\hat{W}_{(k)}\hat{\alpha}_{(k)}$, is an $m\times n$ matrix,
and $m, n\rightarrow\infty$. Therefore, the approximation, $\hat{W}_{(k)}\hat{\alpha}_{(k)}\approx W\alpha$
needs to be defined in a suitable way.

{\bf Definition 1.} We say $\hat{W}_{(k)}\hat{\alpha}_{(k)}$ approximates $W\alpha$, or
$\hat{W}_{(k)}\hat{\alpha}_{(k)}\approx W\alpha$, if $\hat{W}_{(k)}\hat{\alpha}_{(k)}$ can be expressed as
$\hat{W}_{(k)}\hat{\alpha}_{(k)}=W\alpha+L_{1}$ such that $\|L_{1}\|=o(\|W\alpha\|)$, where and throughout this paper, $\|A\|=\sqrt{\lambda_{\max}(A'A)}$ denotes the spectral norm of matrix $A$ ($\lambda_{\max}$ denotes the largest eigenvalue).

When the above approximation holds, we say the RUV-III under such a $k$ achieves the goal of RUV. Note
that the definition of RUV-III [see (\ref{eq:Wkhat})] naturally requires that $U_{(k)}'Y_{\rm c}Y_{\rm c}'
U_{(k)}$, which is $k\times k$, to be non-singular. This is guaranteed by the conditions of the following
theorem. Without loss of generality, let $Y_{\rm c}$ be the first $n_{\rm c}$ columns of $Y$ such that $Y=(Y_{\rm
c}\;Y_{\rm d})$, where $Y_{\rm d}$ represents the remaining columns; similar notation also applies to $\alpha$
and $\epsilon$. Let $w_{i}'$ denote the $i$th row of $W$, and ${\cal I}_{h}=\{1\leq i\leq m, M(i,h)=1\}$, that is,
the subset of indexes of the replicate assays corresponding to sample $h$. Let $W_{h}=(w_{ij})_{i\in{\cal
I}_{h},1\leq j\leq k}$ be the $m_{h}\times k$ submatrix so that $W=(W_{h})_{1\leq h\leq s}$, and $\bar{W}_{h}=
m_{h}^{-1}W_{h}'1_{m_{h}}$. Let $\lambda_{m}=\lambda_{\min}[\sum_{h=1}^{s}(W_{h}-1_{m_{h}}\bar{W}_{h}')'(W_{h}-1_{m_{h}}\bar{W}_{h}')]$ and 
$\lambda_{M}=\lambda_{\max}(W'W)$. We assume that the following regularity conditions hold:
\begin{eqnarray}
\lambda_{M}=\lambda_{m}O(1),\label{eq:Wcond1}\\
\lim\lambda_{m}=\infty.\label{eq:Wcond2}
\end{eqnarray}
Intuitively, condition (\ref{eq:Wcond2}) means that asymptotically there is variation of the unwanted variation within the samples;
condition (\ref{eq:Wcond1}) means that the unwanted variation is of the same order as the variation of the unwanted variation within the samples. The latter holds quite typically. For example, in the case of i.i.d. random variables $X_{1},\dots,X_{n}$, this means that $S_{1}=\sum_{i=1}^{n}X_{i}^{2}$ is of the same order as $S_{2}=\sum_{i=1}^{n}(X_{i}-\bar{X})^{2}$, where $\bar{X}=n^{-1}\sum_{i=1}^{n}X_{i}$. Note that $n^{-1}\sum_{i=1}^{n}X_{i}^{2}\stackrel{\rm P}{\longrightarrow}{\rm E}(X_{1}^{2})$ and $n^{-1}\sum_{i=1}^{n}(X_{i}-\bar{X})^{2}\stackrel{\rm P}{\longrightarrow}{\rm var}(X_{1})$, provided that ${\rm E}(X_{1}^{2})<\infty$. Thus, if further ${\rm var}(X_{1})>0$, $S_{1}$ and $S_{2}$ are of the same order.

We now state the main result of this section. We use the customary notation that, for a symmetric matrix
$A$, $A>0$ means that $A$ is positive definite.

{\bf Theorem 1.} $\hat{W}_{(k)}\hat{\alpha}_{(k)}\approx W\alpha$ holds provided that
\renewcommand{\labelenumi}{(\Roman{enumi})}
\begin{enumerate}
    \item $n, m-s\rightarrow\infty$, and $\liminf(n_{\rm c}/n)>0$;
    \item $\alpha_{1},\dots,\alpha_{n}$ are i.i.d. with ${\rm E}(\alpha_{j})=0$, ${\rm Var}(\alpha_{j})=\Sigma>0$, and ${\rm E}(|\alpha_{j}|^{4})<\infty$;
    \item the entries of $\epsilon$ are i.i.d. with mean $0$ and finite fourth moment; and
    \item conditions (\ref{eq:Wcond1}) and (\ref{eq:Wcond2}) hold.
\end{enumerate}
The proof of Theorem 1 is given in Section \ref{app:proofs}. A basic idea of the proof is to establish the following expression:
\begin{eqnarray}
\hat{W}_{(k)}\hat{\alpha}_{(k)}=W\alpha+\sum_{s=1}^{4}T_{s},\label{eq:expr1}
\end{eqnarray}
where $T_{s}=L_{1}, s=1,2,3,4$. Hereafter, $L_{1}$ denotes a generic term satisfying $\|L_{1}\|=o(\|W\alpha\|)$
(see Definition 1). Note that the specific expression of an $L_{1}$ term may be different at different places. We
refer the details to the proof.

It is possible to weaken the i.i.d. assumption about the entries of $\epsilon$. Let $\epsilon_{ij}$
denote the entries of $\epsilon$. It is reasonable to assume that ${\rm var}(\epsilon_{ij})=\tau_{j}^{2}$, which
may depend on $j$. In such a case, the entries of $\epsilon$ are not i.i.d., because the variance depends
on $j$. To deal with this complication in a way that one is still able to utilize results from the random matrix theory (RMT; e.g., Jiang 2022, ch. 16), which is largely built on the i.i.d. assumption, we replace the previous $\epsilon$ by $\epsilon\tau$, where $\tau={\rm diag}(\tau_{j}, 1\leq j\leq n)$. Here,
the $\epsilon$ is still assumed to satisfy assumption (III) above; it is further assumed that ${\rm
E}(\epsilon_{ij}^{2})=1$. It is seen that, provided that the $\tau_{j}$'s are constants, the $j$ column of
$\epsilon\tau$ now has variance $\tau_{j}^{2}$, $1\leq j\leq n$; this is in line with the note above regarding
the unequal variance. The result of Theorem 1 can be extended to this situation without essential difficulty
(details omitted).

\section{Asymptotic analysis when $k$ is unknown}

In practice, the true dimensionality of $\alpha$, $k$, is unknown. Typically, one would choose a relatively
large number, $K$, and assume that it is an overestimate of $k$. One then uses $K$ in place of $k$ in
RUV-III. Is there a rationale for such a practice? Let $k_{0}$ denote the true and unknown dimension of
$\alpha$. In this section, we study the consequences of choosing $k<k_{0}$, or $k>k_{0}$, in RUV-III. We
do not need to study the consequence of choosing $k=k_{0}$ (i.e., we happen to be right about $k$)
because it was already studied in the previous section.

To this end, let us revisit expression (\ref{eq:expr1}). This expression was derived under the assumption
that one knows the true $k$; hence, the matrix $B=U_{(k)}'W$ is $k\times k$ [this is because $U_{(k)}$,
by definition, is $m\times k$, and $W$ is $m\times k_{0}$; under the assumption that $k_{0}=k$, $B$ is
$k\times k$]. Without assuming $k=k_{0}$, however, $B$ is $k\times k_{0}$; in this case, expression
(\ref{eq:expr1}) becomes
\begin{eqnarray}
\hat{W}_{(k)}\hat{\alpha}_{(k)}=HG^{-1}B\alpha+\sum_{s=1}^{4}T_{s},\label{eq:expr2}
\end{eqnarray}
where $G=BA_{\rm c}B'$ with $A_{\rm c}=\alpha_{\rm c}\alpha_{\rm c}'$, $H=WA_{\rm c}B'$, and $T_{s},
s=1,2,3,4$ are the same as in (\ref{eq:expr1}); see the proof of Theorem 1. It is seen that the key difference is
$$HG^{-1}B\alpha=WA_{\rm c}B'(BA_{\rm c}B')^{-1}B\alpha,$$
which equals $W\alpha$ when $k=k_{0}$.
Note that, however, when $k>k_{0}$, the matrix $BA_{\rm c}B'$ is not invertible. This is because $B$ is
$k\times k_{0}$, hence ${\rm rank}(B)\leq k_{0}$; therefore, ${\rm rank}(BA_{\rm c}B')\leq k_{0}<k$ [e.g., (5.37) of \citep{jiang2022large2}], but $BA_{\rm c}B'$ is $k\times k$. We will deal with this case later. For now, let us
first consider $k<k_{0}$.

\subsection{$k<k_{0}$}

As is for doing virtually everything, an idea is most important. Intuitively, the left side of (\ref{eq:expr2}) is
intended to capture the unwanted variation; thus, the larger $k$ is the more variation is captured (let us
ignore, for now, some potential biology, and noise, that are also captured). This seems to suggest that
one should choose $k$ that maximizes $\|\hat{W}_{(k)}\hat{\alpha}_{(k)}\|_{2}^{2}$. Then, by
(\ref{eq:expr1}) and focusing on the leading term, we should be looking at
\begin{eqnarray}
\|HG^{-1}B\alpha\|_{2}^{2}&=&{\rm tr}(\alpha'B'G^{-1}H'HG^{-1}B\alpha)\nonumber\\
&=&{\rm tr}\{(BA_{\rm c}B')^{-1}BA_{\rm c}W'WA_{\rm c}B'(BA_{\rm c}B')^{-1}B\alpha\alpha'B'\}.\label{eq:expr3}
\end{eqnarray}
Note that $\alpha\alpha'=\sum_{j=1}^{n}\alpha_{j}\alpha_{j}'\approx n\Sigma$ and $A_{\rm c}=\sum_{j=1}^{n_{\rm
c}}\alpha_{j}\alpha_{j}'\approx n_{\rm c}\Sigma$, by the law of large numbers. Thus, if we replace the $A_{\rm c}$
and $\alpha\alpha'$ on the right side of (\ref{eq:expr3}) by $n_{\rm c}\Sigma$ and $n\Sigma$, respectively, the
right side of (\ref{eq:expr3}) is approximately equal to
\begin{eqnarray}
n{\rm tr}\{(B\Sigma B')^{-1}B\Sigma W'W\Sigma B'\}=n{\rm tr}\{W\Sigma B'(B\Sigma B')^{-1}B\Sigma W'\}.
\label{eq:expr4}
\end{eqnarray}
Note that $\Sigma^{1/2}B'(B\Sigma B')^{-1}B\Sigma^{1/2}$ is a projection matrix; thus, we have
\begin{eqnarray}
 \Sigma^{1/2}B'(B\Sigma B')^{-1}B\Sigma^{1/2}\leq I_{k_{0}}\label{eq:ineq1}
 \end{eqnarray}
 for any $k\leq k_{0}$. It follows that the right side of (\ref{eq:expr4}) is bounded by $n{\rm tr}(W\Sigma W')$, and
 this upper bound is achieved when $k=k_{0}$.
 Essentially, this is the idea for dealing with $k<k_{0}$. Namely, we are going to show that, asymptotically, there
 is a gap in $\|\hat{W}_{(k)}\hat{\alpha}_{(k)}\|_{2}^{2}$ between any $k<k_{0}$ and $k=k_{0}$. Hereafter, $\|A\|_{2}=\sqrt{{\rm tr}(A'A)}$ denotes the Euclidean norm of matrix $A$. We now formally execute this idea by evaluating the asymptotic difference in $\|\hat{W}_{(k)}\hat{\alpha}_{(k)}\|_{2}^{2}$.
 
 {\bf Theorem 2.} Under the conditions of Theorem 1 (with $k$ replaced by $k_{0}$), there is a constant $\delta_{1}>0$ such that, with probability tending to one, we have
 \begin{eqnarray}
\min_{1\leq k<k_{0}}\left\{\|\hat{W}_{(k_{0})}\hat{\alpha}_{(k_{0})}\|_{2}^{2}-\|\hat{W}_{(k)}\hat{\alpha}_{(k)}\|_{2}^{2}
\right\}\geq\delta_{1}\lambda_{\min}(\Sigma)\lambda_{m}n.\label{eq:thm2}
\end{eqnarray}

The proof of Theorem 2 is given in Section \ref{app:proofs}. A direct implication of the theorem is the following. Let $\hat{k}$ be
the maximum of $\|\hat{W}_{(k)}\hat{\alpha}_{(k)}\|_{2}^{2}$ over $1\leq k\leq K$, where $K$ is a known upper
bound such that $k_{0}\leq K$. It is fairly convenient to choose such a $K$ in practice.

{\bf Corollary 1.} Under the conditions of Theorem 2, we have ${\rm P}(\hat{k}\geq k_{0})\rightarrow 1$.

Corollary 1 ensures that, asymptotically, $\hat{k}$ is an overestimate of $k_{0}$. It turns out that (see below)
this is sufficient for approximating $W\alpha$ via RUV-III, which is our ultimate goal.

\subsection{$k>k_{0}$}

As noted (see the paragraph before sec. 3.1), in this case, $BA_{\rm c}B'$ is not invertible. Among the key
techniques used in analyzing RUV-III in this case are block-matrix algebra, and projection decomposition. Write (\ref{eq:eigen_decom}) as
\begin{eqnarray}
P_{M^{\perp}}YY'P_{M^{\perp}}=\sum_{l=1}^{m-s}\lambda_{l}u_{l}u_{l}',\label{eq:eigen_decom1}
\end{eqnarray}
where $\lambda_{1}\geq\cdots\geq\lambda_{m-s}\geq 0=\cdots=0$ are the eigenvalues of the left-side matrix,
and $u_{l}$ is the $l$th column of $U$. Let us first consider a special case. Suppose that $n=n_{\rm c}$, that
is, all of the genes are negative controls. Then, we have $Y=Y_{\rm c}$ and $\hat{\alpha}_{(k)}=U_{(k)}'Y_{\rm c}$.
Write $U_{(k)}=[U_{(k_{0})}\;U_{(+)}]$, where $U_{(+)}=[u_{k_{0}+1}\;\dots\;u_{k}]$. It follows that
$$U_{(k)}'Y_{\rm c}Y_{\rm c}'U_{(k)}=\left[\begin{array}{cc}U_{(k_{0})}'Y_{\rm c}Y_{\rm c}'U_{(k_{0})}&
U_{(k_{0})}'Y_{\rm c}Y_{\rm c}'U_{(+)}\\
U_{(+)}'Y_{\rm c}Y_{\rm c}'U_{(k_{0})}&U_{(+)}'Y_{\rm c}Y_{\rm c}'U_{(+)}\end{array}\right].$$
By (\ref{eq:eigen_decom1}) and the fact that $u_{k_{1}}'u_{l}u_{l}'u_{k_{2}}=0$ for any $k_{1}\neq k_{2}$, we have
\begin{eqnarray}
&&U_{(k_{0})}'Y_{\rm c}Y_{\rm c}'U_{(+)}=U_{(k_{0})}'P_{M^{\perp}}YY'P_{M^{\perp}}U_{(+)}
=\sum_{l=1}^{m-s}\lambda_{l}U_{(k_{0})}'u_{l}u_{l}'U_{(+)}\nonumber\\
&&=\sum_{l=1}^{m-s}\lambda_{l}[u_{k_{1}}'u_{l}u_{l}'u_{k_{2}}]_{1\leq k_{1}\leq k_{0},k_{0}+1\leq k_{2}\leq k}=0.
\label{eq:orthog1}
\end{eqnarray}
Recall, for any matrix $Z$, $P_{Z}=Z(Z'Z)^{-1}Z'$. Then, (\ref{eq:orthog1}) implies
\begin{eqnarray}
P=P_{0}+P_{+},\label{eq:proj_decom}
\end{eqnarray}
where $P=P_{Y_{\rm c}'U_{(k)}}$, $P_{0}=P_{Y_{\rm c}'U_{(k_{0})}}$, and $P_{+}=P_{Y_{\rm c}'U_{(+)}}$.
(\ref{eq:proj_decom}) is what we call a projection (orthogonal) decomposition, from which it can be shown
that $\hat{W}_{(k)}\hat{\alpha}_{(k)}\approx W\alpha$ for every $k\geq k_{0}$ (see the proof of Theorem 3 in Section \ref{app:proofs}).

In the general case, that is, without assuming $n=n_{\rm c}$, the nice decomposition of (\ref{eq:proj_decom})
no longer holds. However, some deeper block-matrix algebra still carries through the final approximation result
(again, see the proof of Theorem 3 given in Section \ref{app:proofs}). We state the result formally as follows.

Recall that $Y=(Y_{\rm c}\;Y_{\rm d})$. Define $\xi_{\rm c}=W\alpha_{\rm c}+\epsilon_{\rm
c}$ and $\xi_{\rm d}=W\alpha_{\rm d}+\epsilon_{\rm d}$, where $\alpha_{\rm d}, \epsilon_{\rm d}$ are
defined similarly. We assume that
\begin{eqnarray}
\|(\xi_{\rm c}\xi_{\rm c}')^{-1/2}\xi_{\rm d}\xi_{\rm d}'(\xi_{\rm c}\xi_{\rm c}')^{-1/2}\|=O_{\rm P}(1).\label{eq:add_cond}
\end{eqnarray}
To see that (\ref{eq:add_cond}) is a reasonable assumption, note that the expression inside the spectral norm on the left side
can be written as
$$\frac{n_{\rm d}}{n_{\rm c}}\left(\frac{1}{n_{\rm c}}\sum_{j=1}^{n_{\rm c}}\xi_{j}\xi_{j}'\right)^{-1/2}\left(\frac{1}{n_{\rm
d}}\sum_{j=n_{\rm c}+1}^{n_{\rm c}+n_{\rm d}}\xi_{j}\xi_{j}'\right)\left(\frac{1}{n_{\rm c}}\sum_{j=1}^{n_{\rm c}}\xi_{j}\xi_{j}'\right)^{-1/2},$$
where $\xi_{j}=W\alpha_{j}+\epsilon_{j}, 1\leq j\leq n=n_{\rm c}+n_{\rm d}$. Thus, assuming that $n_{\rm d}/n_{\rm
c}$ is bounded, and some sort of weak law of large numbers holds for the averages of $\xi_{j}\xi_{j}'$ so that the
limit is bounded from above, and away from zero (in a suitable sense), (\ref{eq:add_cond}) is expected to hold. In fact,
(\ref{eq:add_cond}) can be weakened to that with the right side replaced by $o_{\rm P}(\|W\alpha\|^{2}/L_{2})$ for
a certain $L_{2}$, which can be specified, such that $\|L_{2}\|=o(\|W\alpha\|^{2})$.

{\bf Theorem 3.} Suppose that, in addition to the conditions of Theorem 1 (with $k=k_{0}$),  (\ref{eq:add_cond})
holds. Then,  we have $\hat{W}_{(k)}\hat{\alpha}_{(k)}\approx W\alpha$ for every $k>k_{0}$. In fact, the approximation holds for any fixed sequence $k=k_{N}$, where $N$ denotes the combination of all the ``sample sizes'' such as $m, n, n_{\rm c}$ that may increase, such that $k_{N}\leq m-s$ and $k_{N}\rightarrow\infty$.

{\bf Corollary 2.} Under the conditions of Theorem 3, we have $\hat{W}_{(k)}\hat{\alpha}_{(k)}\approx W\alpha$ for $k=m-s$.

The result of Theorem 3 can also cover certain random  (i.e., data-dependent) sequences of $k$. For example, recall that $\hat{k}$ is the maximizer of $\|\hat{W}_{(k)}\hat{\alpha}_{(k)}\|_{2}$ over $1\leq k\leq K$. Combining
Theorems 1, Corollary 1, and Theorem 3, we obtain the following result.

{\bf Theorem 4.} Under the conditions of Theorem 3 we have $\hat{W}_{(\hat{k})}\hat{\alpha}_{(\hat{k})}\approx W\alpha$.

Note that the difference between the conclusion of Theorem 1 and that of Theorem 4 is that the (true) $k$ in
Theorem 1, which is unknown in practice, is replaced by $\hat{k}$ in Theorem 4, which is an estimator of $k$;
yet, the same approximation still holds. The proof of Theorem 3 is given in Section \ref{app:proofs}; the proof of Theorem
4 is straightforward, and therefore omitted.

{\bf Note.} The asymptotic theory shows that $\hat{k}\geq k_{0}$ with probability tending to one (Corollary 1), and that $\hat{W}_{(\hat{k})}\hat{\alpha}_{(\hat{k})}\approx W\alpha$; however, the theory does not state that $\hat{k}$ is a consistent estimator of $k_{0}$, or that $\hat{k}$ is the best choice in some sense. In fact, for $\hat{k}$ to be consistent, some penalty on larger $k$ needs to be imposed in a way similar to the information criteria [e.g., \citep{jn2016}, sec. 1.1]. As noted in the latter reference, there are many choices of the penalty, which may depend on the sample size, that all satisfy the asymptotic requirement for consistency. In practice, it is unclear what is the best choice of the penalty. In this regard, the adaptive fence method [e.g., \citep{jn2016}, sec. 3.2] is potentially applicable.

Nevertheless, here, our main interest is not finding a consistent or optimal choice of $\hat{k}$; our main interest is RUV. As long as $k$ is chosen such that it removes asymptotically the unwanted variation, it serves our main objective of interest. Our practical recommendation, based on the established theory and our empirical investigation (see Section 5), is to use $k=K$, a (relatively) large number, such as $K=m-s$, as long as it is computationally feasible.
\section{Pseudo-replicates of pseudo-samples}

Suppose that we have study data consisting of the results of a large number of assays labelled $i=1,\dots,m$
on cells from tissue samples from patients. Some of these assays may be technical replicates, that is, some may
be carried out on (different) sets of cells from the same tissue sample. They can be used in RUV-III in the usual
way, but we will not be mentioning them in the discussion which follows, as the notion of pseudo-replicates of
pseudo-samples [PRPS;\citep{molania2022removing}] is specifically designed for the common case (e.g. TCGA),
where the study involved no technical replicates.

We require some knowledge of the (in general unknown) biology which we designate by $\{x_{i}, i=1,\dots,m\}$
embodied in some of the tissue samples. For example, if the tissue samples are from a breast cancer biopsy, we
may know the subtype (e.g. basal, luminal A, luminal B, Her2, normal-like) of the breast cancer. This is partial
knowledge of the samples' $\{x_{i}, i=1,\dots,m\}$ and is only required for the samples used to create PRPS.
Similarly we require some knowledge of the (in general unknown) unwanted variation designated by $\{w_{i},
i=1,\dots,m\}$ associated with some of the tissue samples and assays. For example, we may know the date on
which the assay was performed and the plate number in which the assay was carried out. This is partial
knowledge of the samples' $\{w_{i}, i=1,\dots,m\}$ and is only required for the assays used to create PRPS.
See \citep{molania2022removing}, pp. 39--41 for more detail.

The underlying idea is to create {\it in silico} pseudo-samples, that is, the assay results for samples that do not
actually exist, which concentrate the partial knowledge of the biology and the partial knowledge of the unwanted
variation. To do this we select a set of assay results labeled by $i\in{\cal I}_{\rm ps}$ that share the same (partially
known) biology $\{x_{i}, i\in{\cal I}_{\rm ps}\}$ and the same (partially known) unwanted variation $\{w_{i}, i\in{\cal
I}_{\rm ps}\}$. We then form the component-wise average (across genes) of the (log scale) gene expression measurements
from the assays in ${\cal I}_{\rm ps}$, giving the pseudo-sample (PS) assay results. The hope here is that the
pseudo-sample will share the same (partially known) biology and share the same (partially known) unwanted
variation as the samples contributing to it. Indeed, we anticipate that, it being an average, the biology of the
pseudo-sample is more typical of its kind and its unwanted variation is also more typical and more concentrated
of its kind. The assay results for these pseudo-samples created are now added to the assay results that were originally
labelled $i=1,\dots,m$, expanding the rows of our data matrix.
Now suppose that we have created two or more pseudo-samples ${\rm ps1}, {\rm ps2}, \dots$ with the same
shared biology but different levels of shared unwanted variation. We can then form the pseudo-replicate (PR)
set $\{{\rm ps1}, {\rm ps2}, \dots\}$ of two or more PS. The PS in this PR set will have the same shared biology
but different levels of unwanted variation. Declaring them to be PR (i.e., adding a new column) in the
$M$ matrix for the enlarged set of sample plus pseudo-samples will, upon implementation of RUV-III, preserve
the biology embodied in them but consign their differences to unwanted variation by permitting these differences
to contribute to the definition of $\hat{\alpha}$.

In practice we will create as many PS as we can to allow us to remove all the partially known unwanted variation
we have within all the partially known sample biology, with size of the ${\rm PS}$ ranging from 5 to 20 or more,
but at times going as low as 3. In a sense we are building up an $M$ matrix from what was in the case of zero
technical replicates originally the identity matrix, by adding new rows corresponding to the PS we create and new
columns corresponding to the PR sets we define, until we feel that we have covered all the unwanted variation we
know about to the extent we can within all the biology we know about. It has been found that when applying
RUV-III-PRPS with all the unwanted variation that is (partially) known and examining the results, 
visible unwanted variation that was originally unknown is found; one is then in a position to define yet more
PS and PR sets and remove more unwanted variation. See \citep{molania2022removing} for an illustration with the
TCGA breast cancer data.

To formulate PRPS mathematically in a way similar to (\ref{eq:LMM}), suppose that we initially have
$m_{0}$ assays. Some of these are non-replicated, that is, each assay corresponds to a different sample. Let
$s_{\rm nr}$ be the total number of such samples/assays (nr stands for non-replicated). The $M$ matrix
corresponding to these assays is the $s_{\rm nr}$-dimensional identity matrix.

The rest of the $m_{0}$ assays are real replicated assays, that is, multiple assays correspond to the same
sample. For these $m_{\rm rr}=m_{0}-s_{\rm nr}$ assays (rr stands for real replicated), the $M$ matrix has
$m_{\rm rr}$ rows and $s_{\rm rr}$ columns so that each column has multiple $1$s. Without loss of generality,
let this $M$ be expressed as $M_{\rm rr}={\rm diag}(1_{m_{h}}, s_{\rm pr}+1\leq h\leq s_{\rm pr}+s_{\rm rr})$,
where $s_{\rm pr}$ is the total number of pseudo replicated samples to be produced below.

The model for the initial assays, similar to (\ref{eq:LMM}), can be written as
\begin{eqnarray}
Y_{\rm 0}=1_{m_{0}}\mu+M_{0}X_{0}\beta+W_{0}\alpha+\epsilon_{0}\label{eq:LMM0}
\end{eqnarray}
with $M_{0}={\rm diag}(M_{\rm rr},I_{s_{\rm nr}})$, where the entries of $\epsilon_{0}$ satisfy the previous assumptions for the entries of $\epsilon$ [see condition (III) of Theorem 1]. We then produce the pseudo assays, which satisfy
\begin{eqnarray}
Y_{\rm pa}=A_{\rm a}Y_{0}=1_{m_{\rm pa}}\mu+A_{\rm a}M_{0}X_{0}\beta+A_{\rm a}W_{0}\alpha+A_{\rm
a}\epsilon_{0},\label{eq:Ypa}
\end{eqnarray}
where $A_{\rm a}$ is the matrix corresponding to taking the averages in producing the pseudo assays (the
subscript $a$ refers to average), and $m_{\rm pa}$ is the number of pseudo assays such produced (i.e., the
number of rows of $A_{\rm a}$). Note that we have $A_{\rm a}1_{m_{0}}=1_{m_{\rm pa}}$.

Now suppose there exist, at least, approximate relationships so that
\begin{eqnarray}
A_{\rm a}M_{0}X_{0}=M_{\rm pr}X_{\rm pr}\label{eq:Aa_equ}
\end{eqnarray}
for some matrix $m_{\rm pa}\times s_{\rm pr}$ matrix $M_{\rm pr}$, whose entries are $0$ or $1$ with exactly
one $1$ in each row and multiple $1$s in each column, and possibly unknown sub-matrix, $X_{\rm pr}$.
Without loss of generality, let $M_{\rm pr}={\rm diag}(1_{m_{h}}, 1\leq h\leq s_{\rm pr})$. Then, we can rewrite
(\ref{eq:Ypa}) as
\begin{eqnarray}
Y_{\rm pa}=1_{m_{\rm pa}}\mu+M_{\rm pr}X_{\rm pr}\beta+A_{\rm a}W_{0}\alpha+A_{\rm a}\epsilon_{0},
\label{eq:Ypa1}
\end{eqnarray}
so that, when combining (\ref{eq:LMM0}) and (\ref{eq:Ypa1}), we have
\begin{eqnarray}
&&Y=\left(\begin{array}{c}Y_{\rm pa}\\
Y_{\rm 0}\end{array}\right)=1_{m}\mu+\left(\begin{array}{cc}M_{\rm pr}&0\\
0&M_{0}\end{array}\right)\left(\begin{array}{c}X_{\rm pr}\\
X_{0}\end{array}\right)\beta+\left(\begin{array}{c}A_{\rm a}\\
I_{m_{0}}\end{array}\right)W_{0}\alpha+\left(\begin{array}{c}A_{\rm a}\\
I_{m_{0}}\end{array}\right)\epsilon_{0}\nonumber\\
&&=1\mu+MX\beta+W\alpha+\epsilon=1\mu+M\eta+W\alpha+\epsilon,\label{eq:LMM1}
\end{eqnarray}
with notations defined in obvious ways, that is, $1=1_{m}$ with $m=m_{0}+m_{\rm pa}$,
\begin{eqnarray}
M=\left(\begin{array}{cc}M_{\rm pr}&0\\
0&M_{0}\end{array}\right), X=\left(\begin{array}{c}X_{\rm pr}\\
X_{0}\end{array}\right), W=\left(\begin{array}{c}A_{\rm a}\\
I_{m_{0}}\end{array}\right)W_{0}, \epsilon=\left(\begin{array}{c}A_{\rm a}\\
I_{m_{0}}\end{array}\right)\epsilon_{0},\label{eq:newMXW}
\end{eqnarray}
and $\eta=X\beta$. What is important about expressions (\ref{eq:newMXW}) is that $M$ is known; it does not really
matter whether $X$ is known or not because, in the end, $X\beta$ will be treated together as an unknown matrix,
which we denote by $\eta$.

To proceed further, we need to separate the non-diagonal part of $M$ from the diagonal part, that is, write [see
below (\ref{eq:LMM0})] $M={\rm diag}(M_{\rm pr},M_{\rm rr},I_{s_{\rm nr}})={\rm diag}(M_{\rm r},I_{s_{\rm nr}})$
with
\begin{eqnarray}
M_{\rm r}={\rm diag}(M_{\rm pr},M_{\rm rr})={\rm diag}(1_{m_{h}}, 1\leq h\leq s_{\rm r}=s_{\rm pr}+s_{\rm rr}).
\label{eq:Mr}
\end{eqnarray}
Because, in practice, $m_{\rm r}<<m$, and always $s_{\rm r}\leq s$, an important contribution to the
speed-up of computation via $M_{\rm r}$ comes from the fact that the eigenvectors involved in defining
$\hat{\alpha}$ can be obtained by averaging and taking residuals using the $m_{\rm r}\times s_{\rm r}$ matrix
$M_{\rm r}$, with $m_{\rm r}=m_{\rm pa}+m_{\rm rr}$, and just the $m_{\rm r}$ rows of the $m\times n$ data
matrix $Y$ corresponding to assays or pseudo-assays involved in the replicate or pseudo-replicate sets.

We are now ready to extend the previous theoretical results to PRPS. First, some regularity conditions are
needed for the matrix $A_{\rm a}$ of operation via averages. We assume that\\
(a) each average (corresponding to a row of $A_{\rm a}$) involves a bounded number, say, $b_{1}$, of original
assays (i.e., rows of $Y_0$); and\\
(b) each original assay (i.e., row of $Y_{0}$) is involved in a bounded number, say, $b_{2}$, of averages (i.e.,
rows of $A_{\rm a}$).\\
These conditions intend to simplify the essential degrees of freedom, or effective sample sizes, involved in the asymptotic theory. It is well known [e.g., \citep{jiang2017}] that the effective sample size is complicated in the asymptotic analysis of mixed effects models, and this has particularly to do with the so-called replicates, that is, observations sharing the same random effects. Because the replicates are correlated for sharing the same random effects, their contribution to the effective sample size is not the same as the number of replicates. However, in case the number of replicates are bounded, the effective sample size is proportional to the number of clusters, which corresponds to the random effects to which the replicates are associated. Conditions (a) and (b) are similar to assuming that the number of replicates is bounded for each cluster in the case of clustered observations, with the number of original assays, or the number of measurements involved in the average, corresponding to the number of the replicates mentioned above. In case that the number of replicates also increase, the effective sample size is more complicated [e.g.,\citep{jwb2022}]. We avoid such a complexity in this work, which also seems reasonable in the context of the current applications.

{\bf Theorem 5.} Suppose the following  hold, where conditions (II) and (III) refer to those in Theorem 1:
\renewcommand{\labelenumi}{(\roman{enumi})}
\begin{enumerate}
    \item $s, n, m-s, m_{\rm r}-s_{\rm r}\rightarrow\infty$, $\liminf(n_{\rm c}/n)>0$;
    \item condition (II) for $\alpha$, and condition (III) for the entries of $\epsilon_{0}$;
    \item the averaging matrix, $A_{\rm a}$, satisfies the assumptions (a), (b) above;
    \item (\ref{eq:Wcond1}), (\ref{eq:Wcond2}) hold with $W$ replaced by $W_{\rm r}$, which is the first $m_{\rm r}$ rows of the $W$ in (\ref{eq:newMXW});
    \item (\ref{eq:add_cond}) holds with $\xi_{\rm c}=W_{\rm r}\alpha_{\rm c}+\epsilon_{\rm rc}$, $\xi_{\rm d}=W_{\rm r}\alpha_{\rm d}+\epsilon_{\rm rd}$, where $\epsilon_{\rm r}=(\epsilon_{\rm rc}\;\epsilon_{\rm rd})$ is the first $m_{\rm r}$ rows of the $\epsilon$ in (\ref{eq:newMXW}).
\end{enumerate}
Then, the conclusions of Theorems 2, Theorem 3, and Theorem 4 hold for PRPS, that is, RUV-III based on the combined data matrix satisfying (\ref{eq:LMM1}), with the $k_{N}$ being any sequence $k_{N}$ satisfying $k_{N}\leq m_{\rm r}-s_{\rm r}$ and $k_{N}\rightarrow\infty$.
\section{Some empirical results}

\subsection{Simulation study}

We conduct a series of simulation studies to provide empirical support for the theorems, to investigate rates of convergence, and to explore the impact of different choices of $k$. 

Our simulation study consists of multiple ``runs''.  In each simulation run, we allow $m$ to consecutively take on the values $2^4, 2^5, ..., 2^9$, and for each value of $m$ we generate 100 datasets (details below).  We fit each of the 100 datasets by RUV-III, then calculate the relative norm of approximation error,
$$
q = \frac{||\widehat{W\alpha} - W\alpha||}{||W\alpha||}.
$$
We then compute the average of the 100 values of $q$ for each value of $m$.  In this way we are able to see empirically how $q$ varies with $m$, that is, to check that it goes to 0 and to see how quickly it does so.  

In all simulation runs we set $n = m^2$.  In addition, we always fix the true value of $k$ to be $k_0 = 3$.  What we vary between runs are (1) how $n_c$ scales with $m$; (2) the replication structure; (3) the distributions of $W$, $\alpha$, and $\epsilon$; and (4) the choice of $k$ when fitting RUV-III.

More specifically, we consider three possible scalings for $n_c$.  In the first, we set $n_c = m^2/8 = n/8$ so that $n_c/n$ is a fixed constant.  In the second, we set $n_c = m^{3/2}/2$ so that $n_c/n \to 0$, but also $m/n_c \to 0$.  In the third, we set $n_c = 2m$.  Note that in all three cases, we have that $n_c = 32$ at the smallest sample size, that is, $n_c = 32$ when $m = 2^4 = 16$.

We consider two replication structures, \textit{samples increasing} and \textit{replicates increasing}.  In  \textit{samples increasing}, we let the number of samples $s$ grow with $m$ and set $s = m/4$; each sample has four replicates.  In \textit{replicates increasing} we do the opposite --- we fix $s = 4$ but let each sample have $m/4$ replicates.

When generating $W$, $\alpha$, and $\epsilon$ we consider two options: \textit{normal} and \textit{Pareto}.  In the first, all entries of $W$, $\alpha$, and $\epsilon$ are i.i.d. standard normal.  In the second, we replace the standard normal distribution with a rescaled, recentered Pareto distribution with the scale parameter equal to $1$ and shape parameter equal to $5$. The rescaling and recentering standardizes the distribution to have mean 0 and variance 1 so that it is comparable to the $\textit{normal}$ setting.  Note that we do not generate a $\mu$ term nor an $MX\beta$ term, and simply set $Y = W\alpha + \epsilon$.  Importantly, the RUV-III estimate $\widehat{W\alpha}$ is not a function of $MX\beta$ (assuming $\beta_c = 0$ and $k \le m-s$), and $q$ is therefore also not a function of $MX\beta$.  We do however compute $M$ itself, which is used by RUV-III.

Finally, we consider three cases for $k$.  In the first, we set $k$ equal to its true value, i.e., $k = k_0 = 3$.  In the second, we slightly overestimate $k$ and set $k = 10$.  In the third, we completely overestimate $k$ and set it equal to its maximum possible value, $k = k_{\mathrm{max}} \equiv m - s$.

\begin{figure}
    \centering
    \begin{tabular}{cccc}
    & $k=3$ & $k=10$ & $k=k_{\mathrm{max}}$ \\
    \rotatebox{90}{\hspace{0.4cm}Number of samples increasing} &
    \includegraphics[scale=0.25]{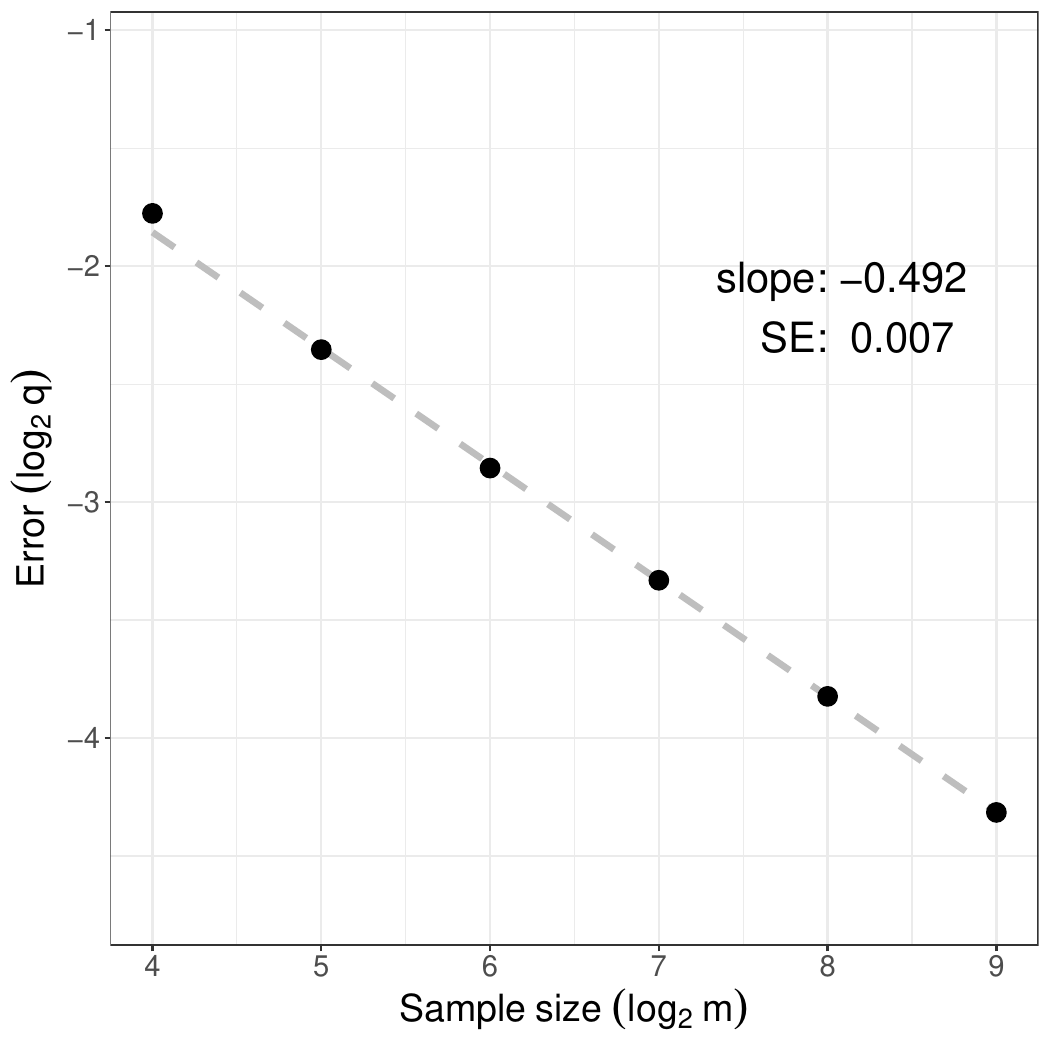} &
    \includegraphics[scale=0.25]{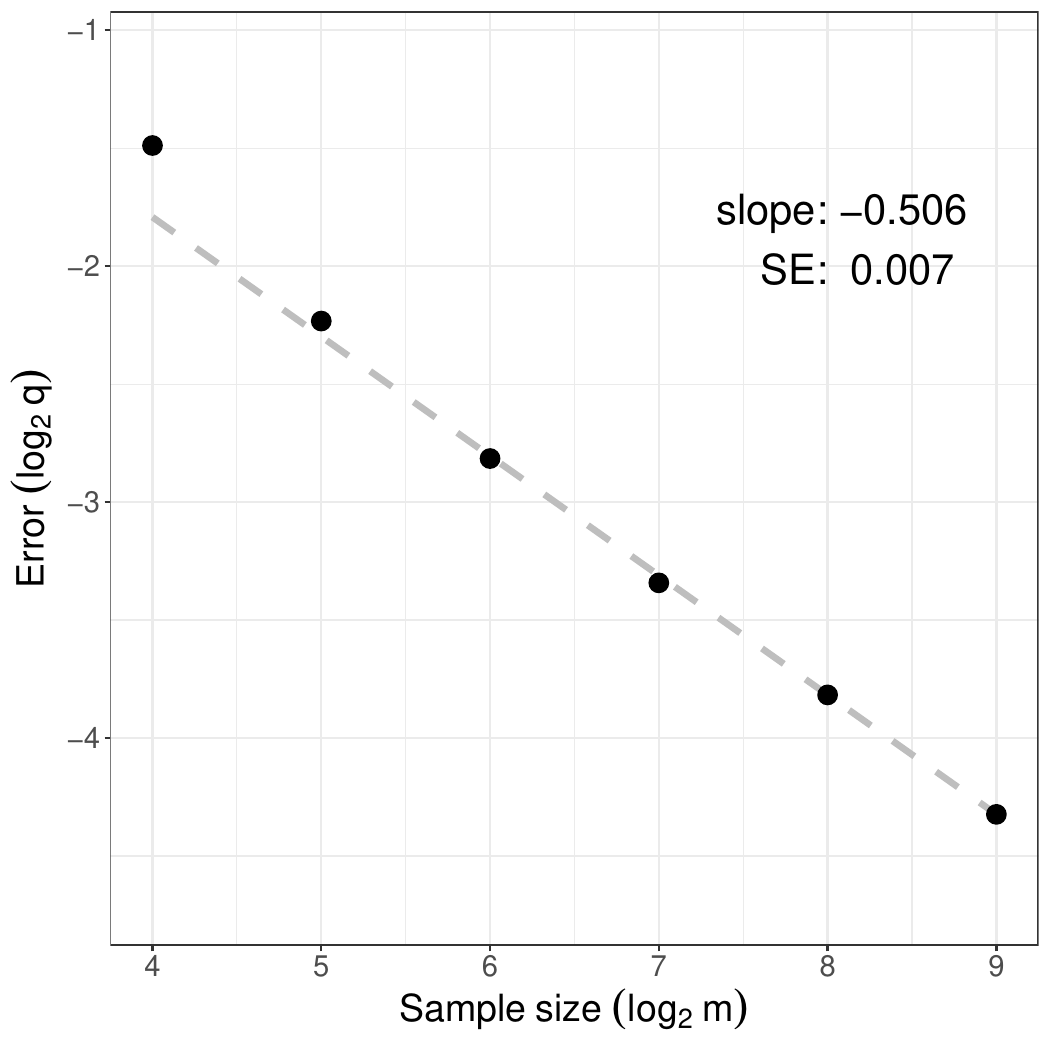} &
    \includegraphics[scale=0.25]{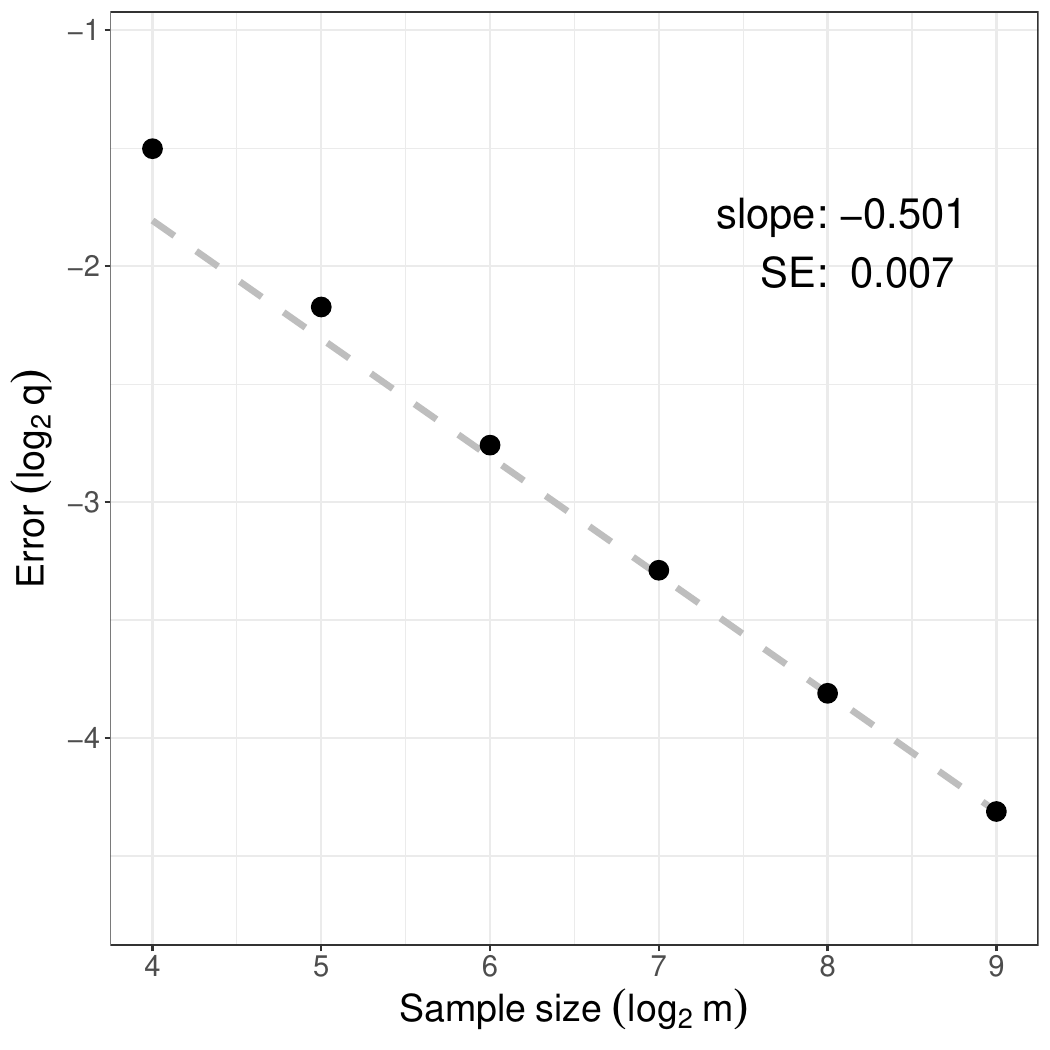} \\    
    \rotatebox{90}{\hspace{0.4cm}Number of replicates increasing} &
    \includegraphics[scale=0.25]{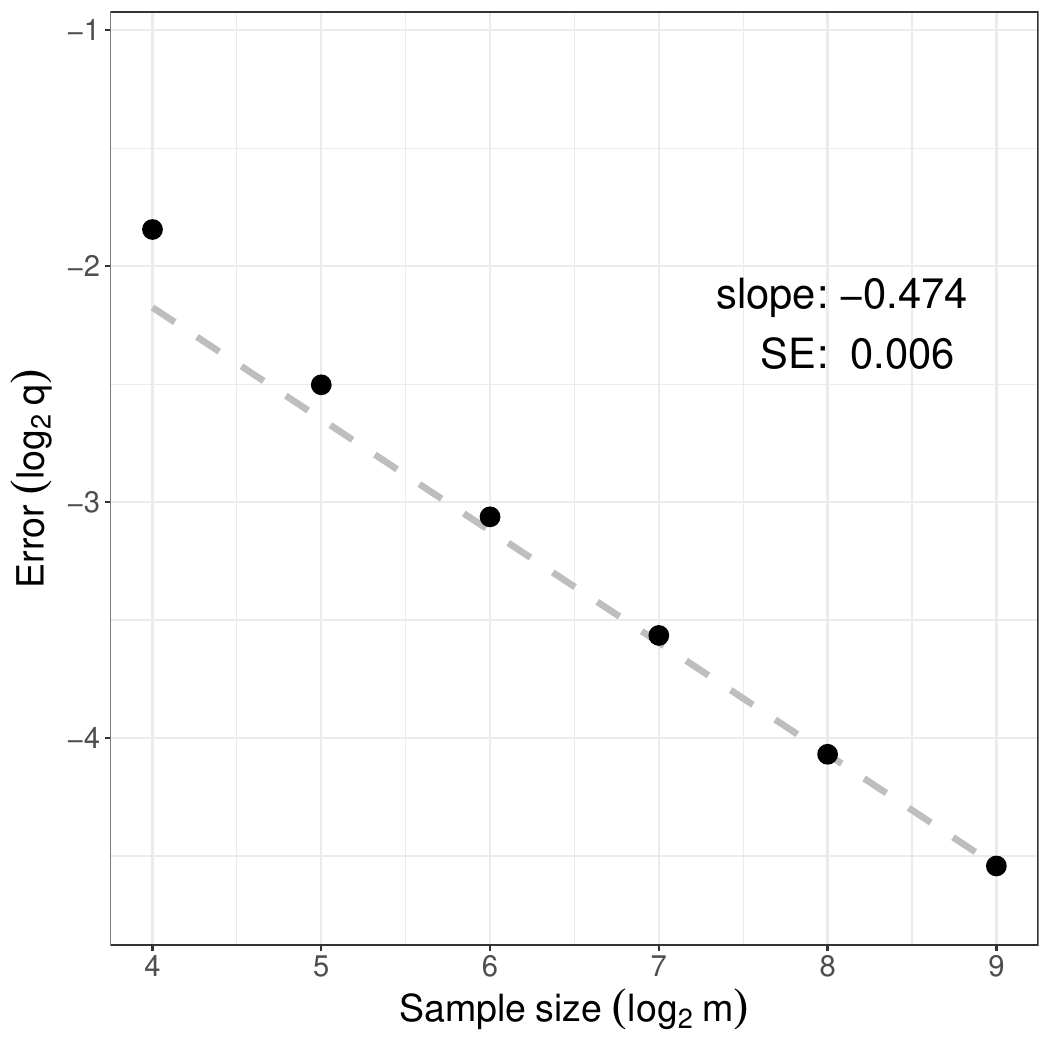} &
    \includegraphics[scale=0.25]{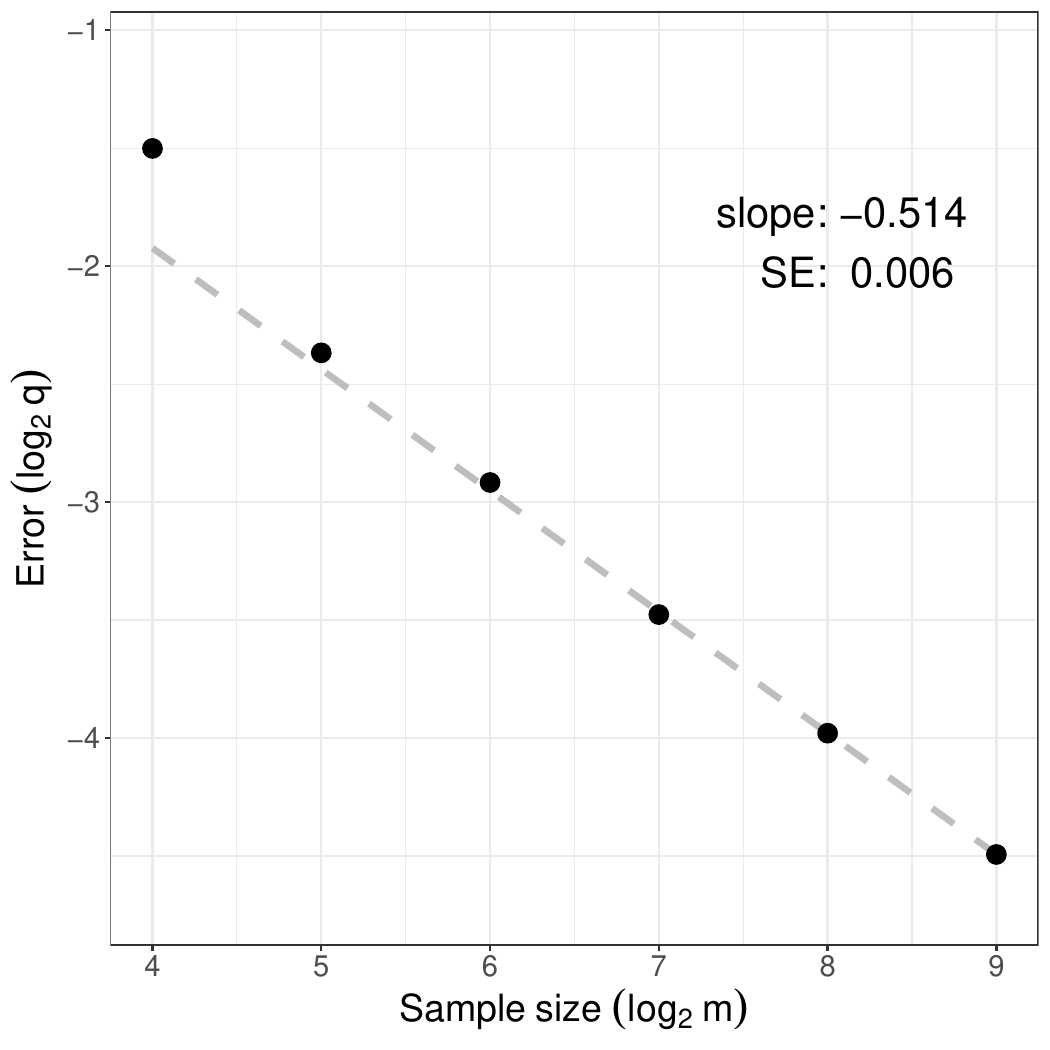} &
    \includegraphics[scale=0.25]{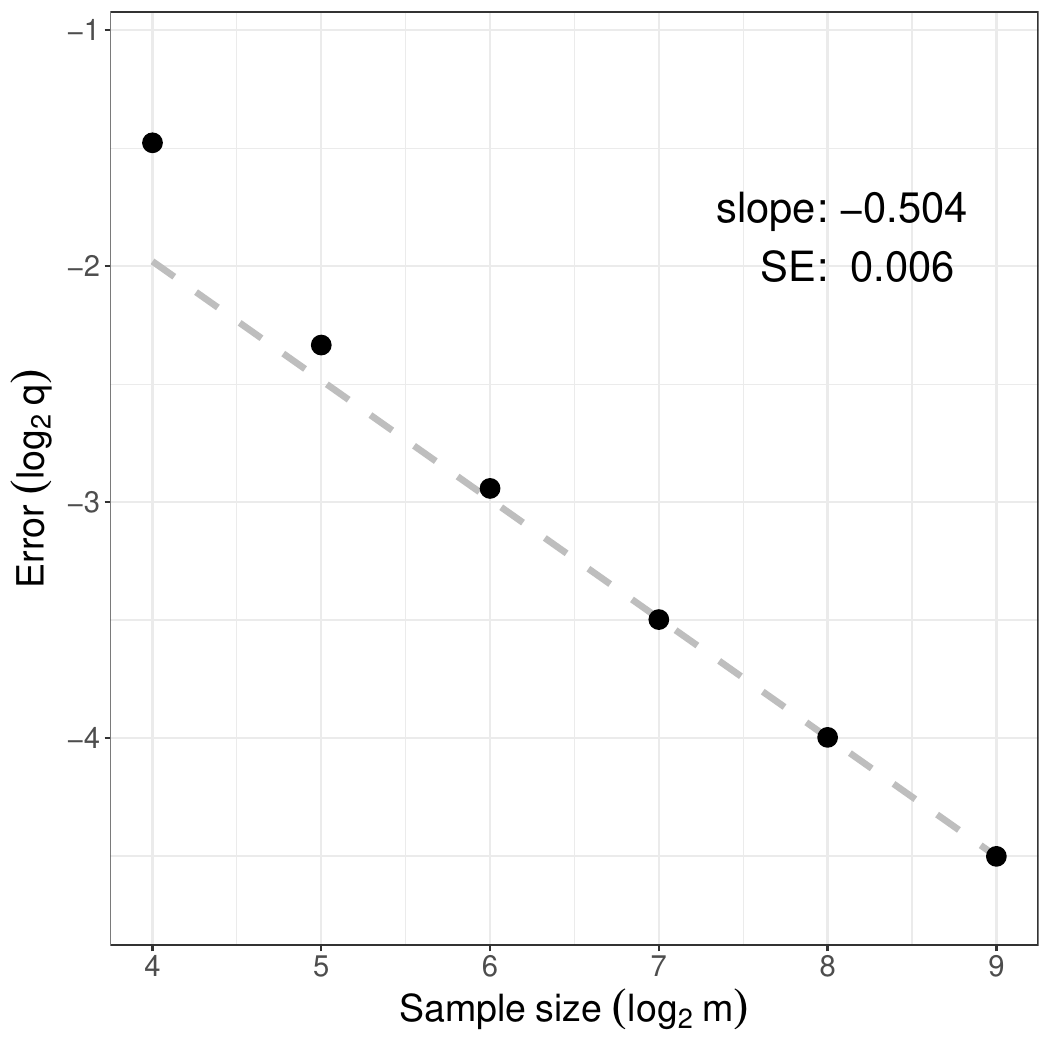}     
    \end{tabular}
    \caption{Simulation results for the case $n_c = m^2/8$ with normally distributed data.}
    \label{fig:simbasic}
\end{figure}

Figure \ref{fig:simbasic} shows results for the case $n_c = m^2/8$ with normally distributed data.  The horizontal axis is $m$ on a $\log_2$ scale and the vertical axis is the average value of $q$, also on a $\log_2$ scale.  The dotted line is fitted to the final two points (when $m$ is $2^8$ and $2^9$), and its slope therefore gives an indication of how $q$ decays with $m$ asymptotically.  This slope and its standard error are provided in each plot.  We see that in all cases, $q$ tends to 0, and appears to do so at a $m^{-1/2}$ rate.  This is true for both replication structures and even when $k$ is overestimated.  
A similar plot is provided in the supplmentary material for the Pareto case and the results are essentially the same.

\begin{figure}
    \centering
    \begin{tabular}{ccc}
    $k=3$ & $k=10$ & $k=k_{\mathrm{max}}$ \\
    \includegraphics[scale=0.25]{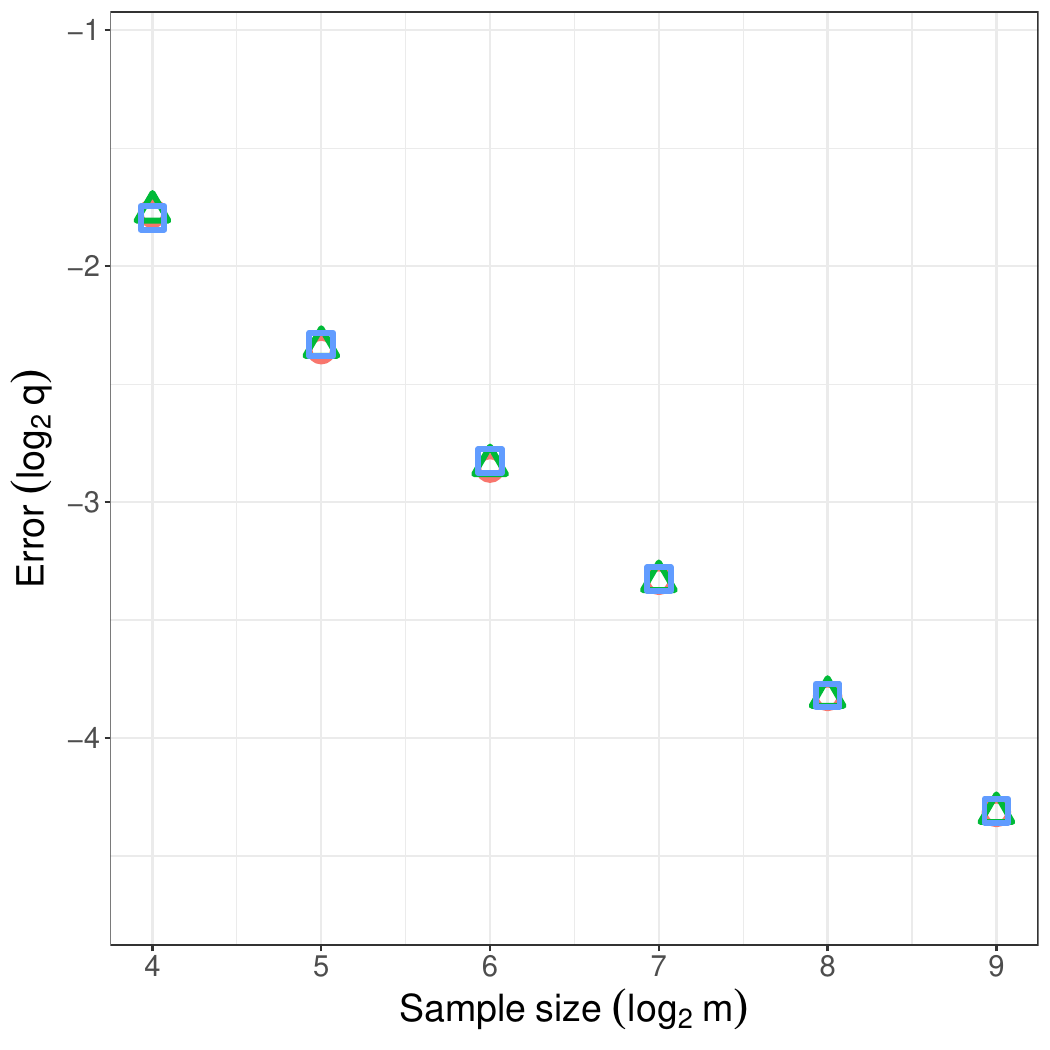} &
    \includegraphics[scale=0.25]{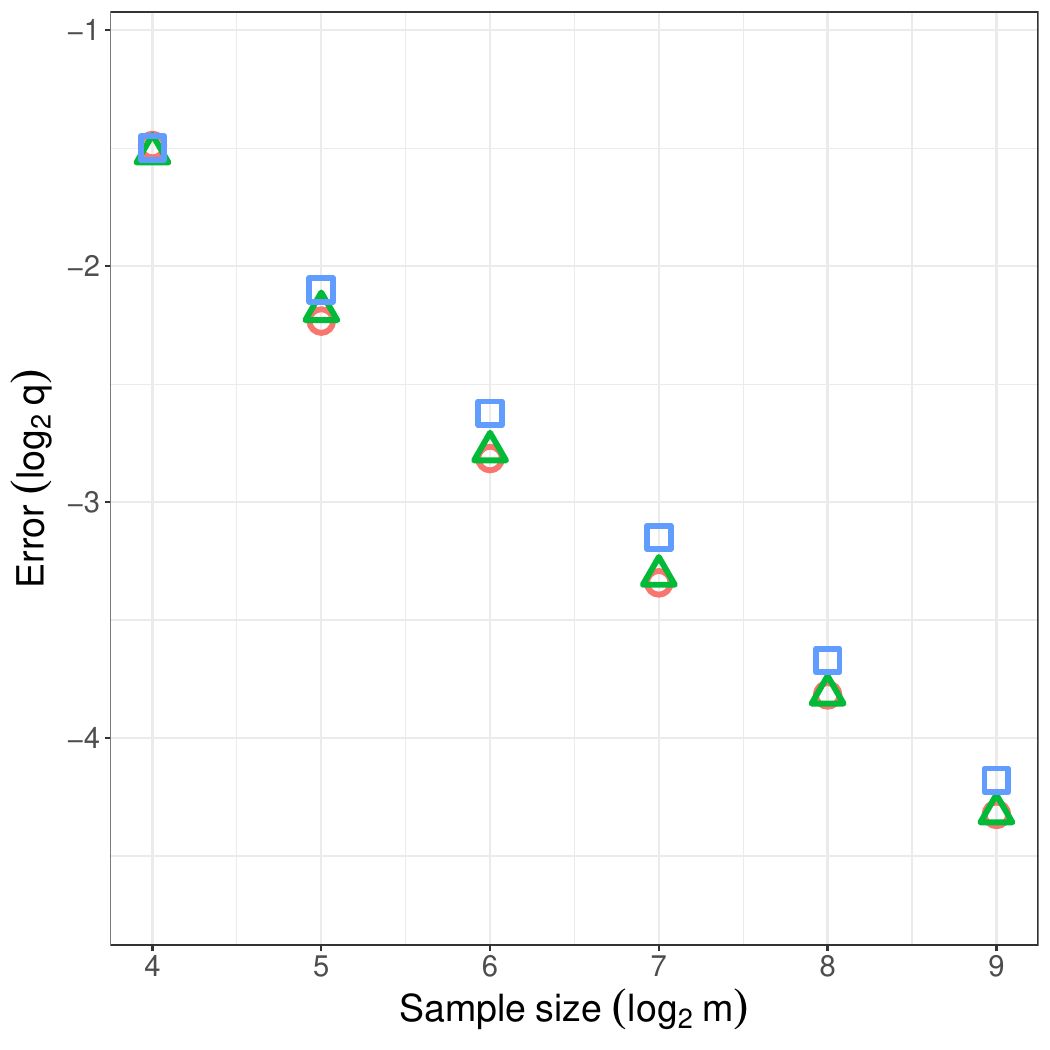} &
    \includegraphics[scale=0.25]{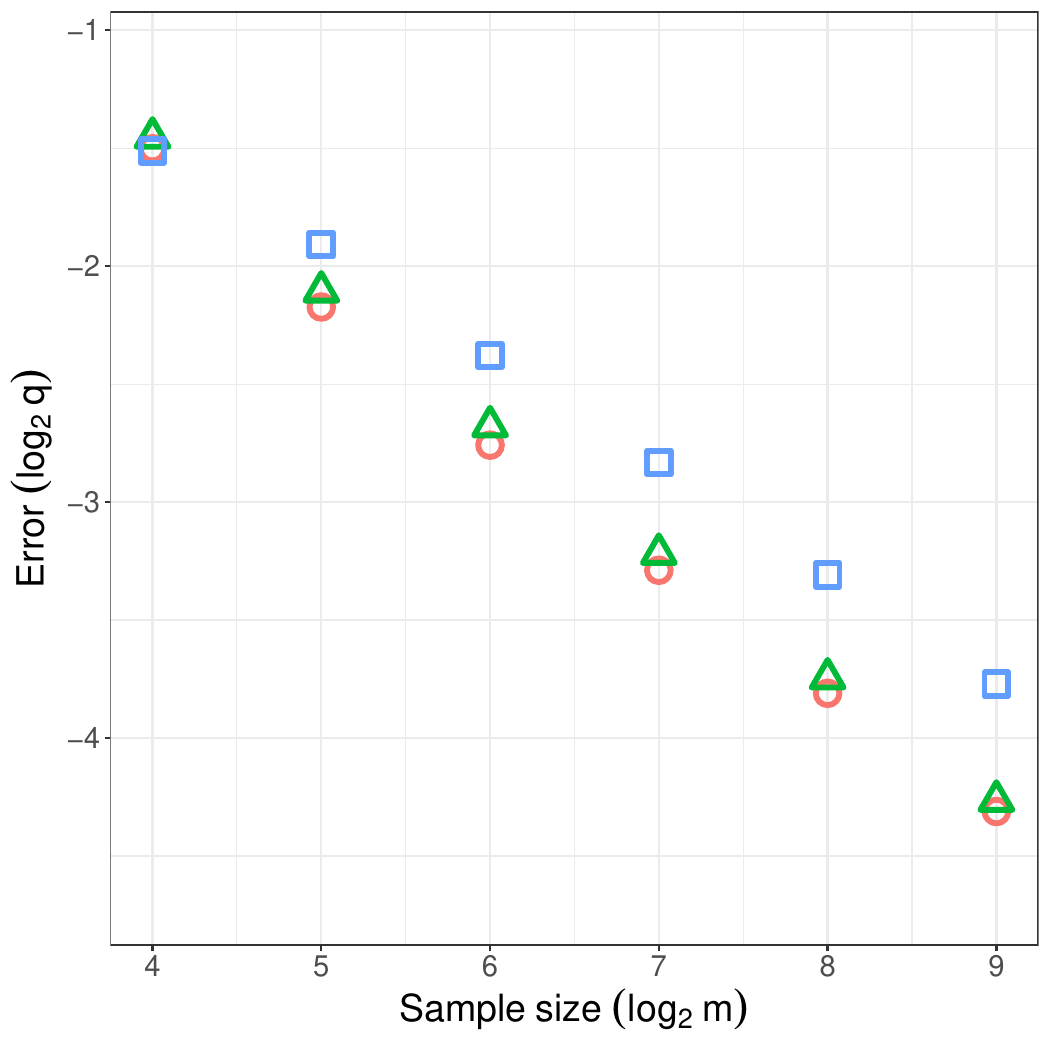} 
    \end{tabular}
    \caption{Simulation results comparing the cases $n_c = m^2/8$ (red circle), $n_c = m^{3/2}/2$ (green triangle), and $n_c = 2m$ (blue square).  Data is normally distributed, and \textit{samples increasing}.}
    \label{fig:simcomparison1}
\end{figure}

Figure \ref{fig:simcomparison1} compares the three cases of $n_c$.  Interestingly, when $k = 3$, the results are almost identical for all three cases, suggesting that we may obtain a good adjustment even with a fairly limited number of negative controls.  However, when $k$ is overestimated, especially when $k = k_{\mathrm{max}}$, the quality of the adjustment is more sensitive to the number of negative controls.  From a practical point of view, this suggests that when the number of negative controls is limited, extra attention should be given to the selection of $k$ \citep{gagnonbartsch2012using, molania2019new, molania2022removing}. Nonetheless, we note that even when $n_c = 2m$ and $k = k_{\mathrm{max}}$, it still appears that $q \to 0$, although not at a $m^{-1/2}$ rate.

\begin{figure}
    \centering
    \begin{tabular}{ccc}
    $n_c = \frac{1}{8}m^2$ & $n_c = \frac{1}{2}m^\frac{3}{2}$ & $n_c = 2m$ \\
    \includegraphics[scale=0.25]{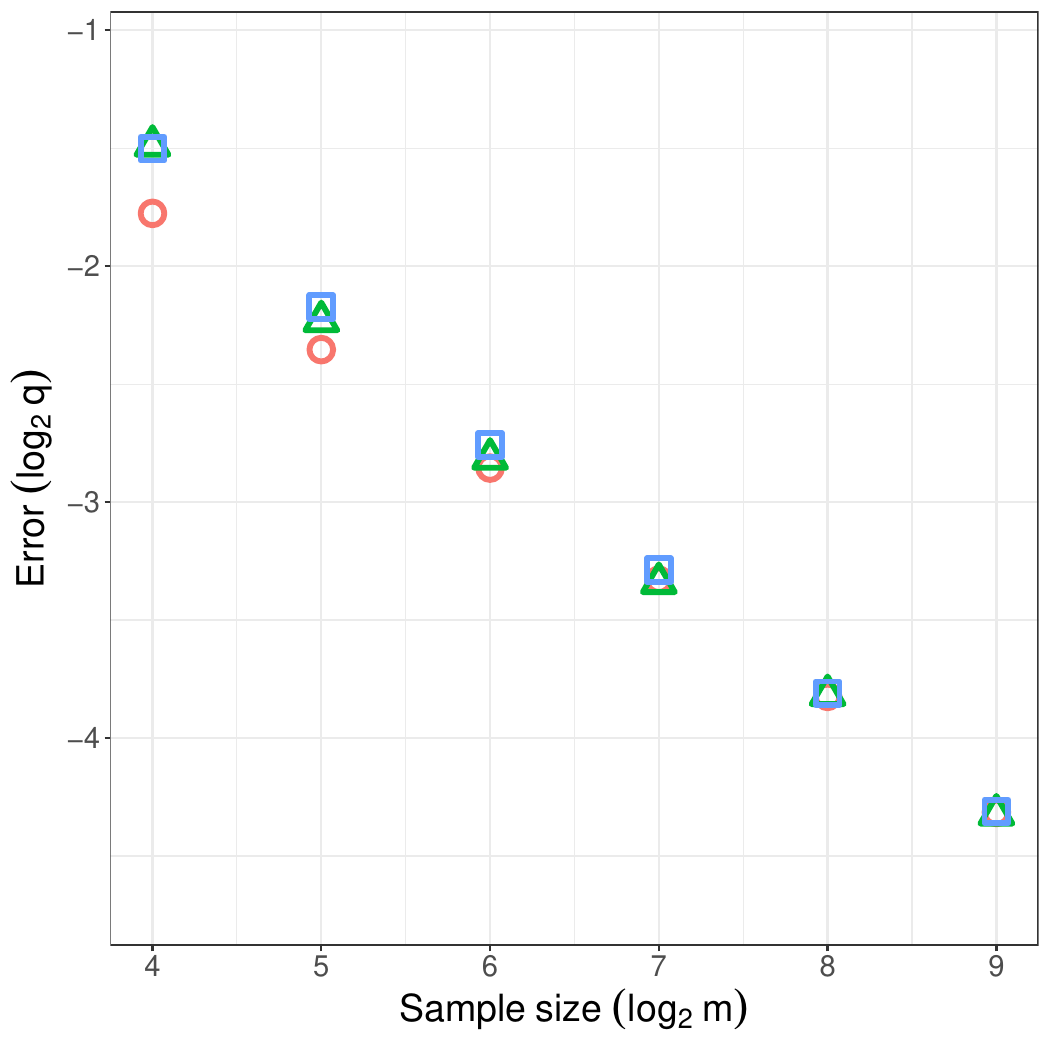} &
    \includegraphics[scale=0.25]{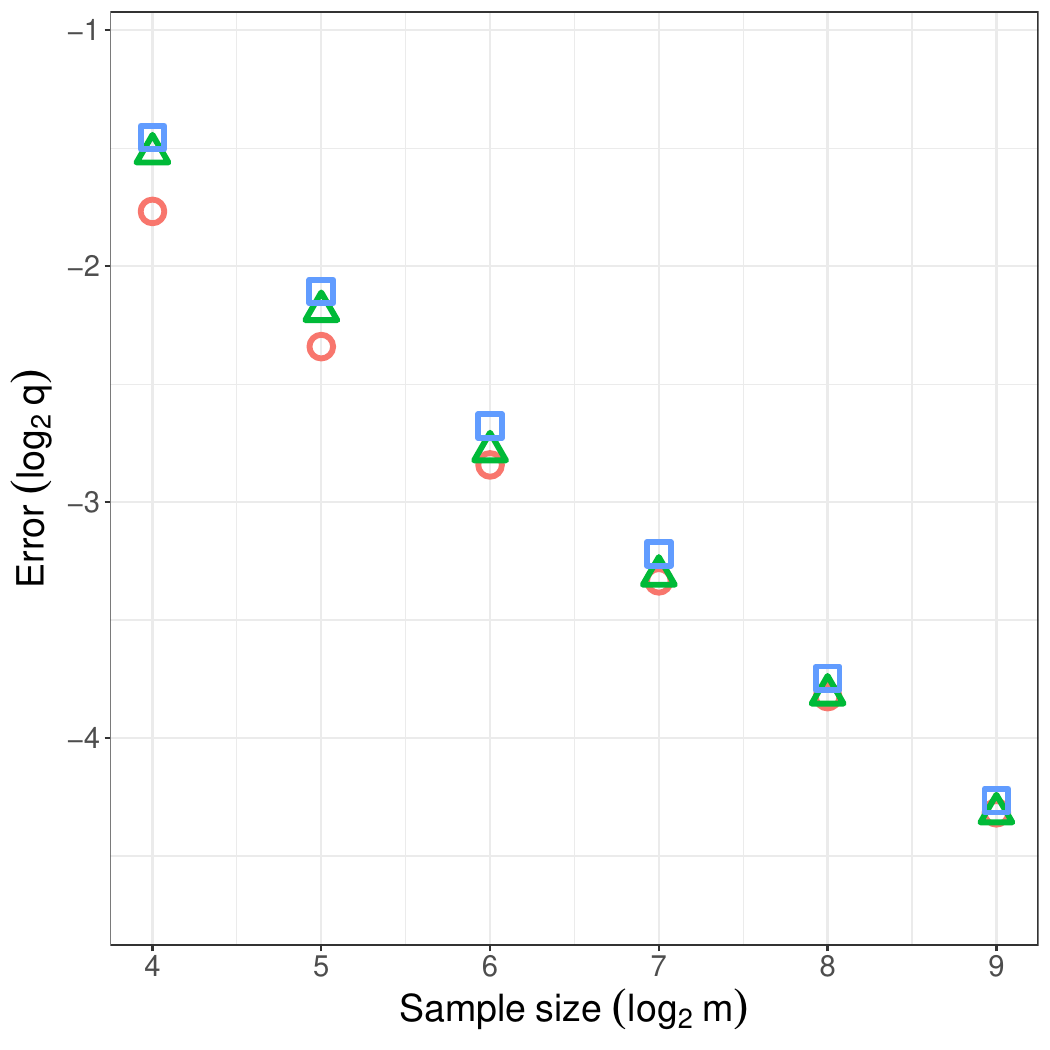} &
    \includegraphics[scale=0.25]{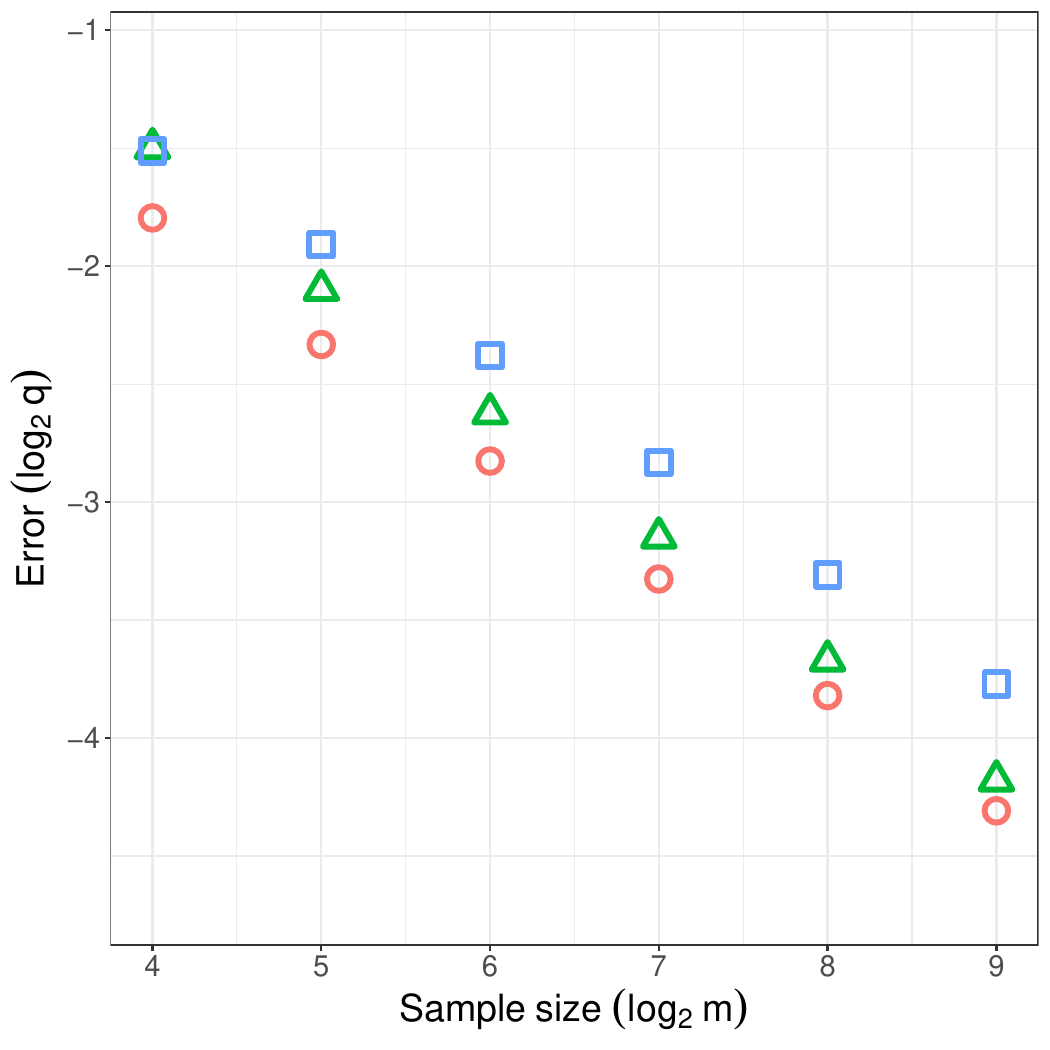} 
    \end{tabular}
    \caption{Simulation results comparing the cases $k=3$ (red circle), $k=10$ (green triangle), and $k = k_{\mathrm{max}}$ (blue square).  Data is normally distributed, and \textit{samples increasing}.}
    \label{fig:simcomparison2}
\end{figure}

Figure \ref{fig:simcomparison2} shows the same results as Figure \ref{fig:simcomparison1}, but arranged differently to more directly compare the impact of $k$.  Of particular note are the left and center panels, which suggest that asymptotically there is effectively no penalty for overestimating $k$ when the number of negative controls is moderately large, relative to the sample size.  However, at smaller sample sizes or with a limited number of negative controls the choice of $k$ has a modest impact.

The supplementary material provides a figure similar to Figure \ref{fig:simbasic} but showing results for $k=1$ and $k=2$, and show that $q$ does not tend to 0.  The supplementary material also provides figures similar to Figures \ref{fig:simcomparison1} and \ref{fig:simcomparison2} but with the Pareto distribution, and these look essentially the same as the Normal.

Finally, Figure S5 in the Supplementary Material shows the results of a simulation in which we implement PRPS.  The data is generated similarly as above, with some modifications.  Above, the $m$ observations are partitioned into $s$ groups of replicates.  Here, we consider the case that $m_0 = s$, so that there are no real replicates, but instead partition the $m_0$ observations into $s$ groups of known biology.  Here, $s$ may stand for \textit{subtype} --- e.g., basal, luminal A, luminal B, Her2, normal-like --- rather than \textit{sample}, but conveniently $s$ plays the same role in the data generating process here as in the previous simulations.  In addition, we introduce a systematic trend component to $W$, in order to simulate temporal drift corresponding to the run order of the assays.  More specifically, rather than let $W_{ij} \sim N(0, 1)$ as we did above, we now let $W_{ij} \sim N(t_{ij}, 1)$ where $t_{ij}$ is the temporal trend (we focus on Normal, not Pareto).  The trend $t_{ij}$ is piecewise linear in $i$, with piecewise components that vary by $j$ and by subtype.  The range of each individual component is generated independently as $[-U, U]$, where $U \sim U(0,1)$.  Otherwise, data generation proceeds exactly as above.

Having generated the $m_0$ original observations, we then create $s\times 2$ pseudo assays.  We create two pseudo assays for each subtype; one from the first half of the assays from that subtype, and the other from the second half.  The final $Y$ matrix therefore contains $m_0 + 2s$ rows.  We consider the case $s = m_0/4$, i.e., \textit{subtypes increasing} (analogous to \textit{samples increasing} above), so that the number of pseudo assays is proportional to $m$.  In Figure S5 we see that, similarly to Figure 1, $q$ tends to 0 and appears to do so at a $m^{-1/2}$ rate.

\subsection{A real-data example}

We apply RUV-III to a microarray dataset both to demonstrate the efficacy of RUV-III, as well as to highlight the role that negative controls and replicates can play in essentially defining what variation is wanted and what is unwanted.  

We use the dataset of \citet{vawter2004gender}, which is the same dataset analyzed in \citet{gagnonbartsch2012using}.  Briefly, the goal of the original study was to find genes differentially expressed in the brain with respect to sex.  Tissue samples were collected post-mortem from 10 individuals, 5 male and 5 female.  Three tissue samples were taken from each individual --- one from the anterior cingulate cortex (ACC), one from the dorsolateral prefrontal cortex (DLPFC), and one from the cerebellum (CB).  Thus, $10\times3=30$ tissue samples were collected in total.  Three aliquots were derived from each sample, and sent to three different laboratories for analysis (UC Davis, UC Irivine, U. Michigan); this was for sake of replication, and the analysis at each laboratory was in principle the same.  In particular, each laboratory assayed their 30 aliquots on Affymetrix HG U95 microarrays, resulting in a total of $30\times3=90$ microarray assays.  Data from 6 of the 90 are unavailable, presumably due to quality control issues, so we only have a final sample size of $m = 84$.  Each microarray measures the gene expression of $n = 12,600$ probesets (genes).  33 of these probesets correspond to spike-in controls, exogenous RNA added at a fixed concentration to each sample for quality control.  For additional details, see \citet{vawter2004gender} and \citet{gagnonbartsch2012using}. 

\begin{figure}
    \centering
    \begin{tabular}{cc}
    \hspace{-1cm} No Adjustment & 
    \hspace{-1cm} ``Technical'' Adjustment \\
    \includegraphics[scale=0.3]{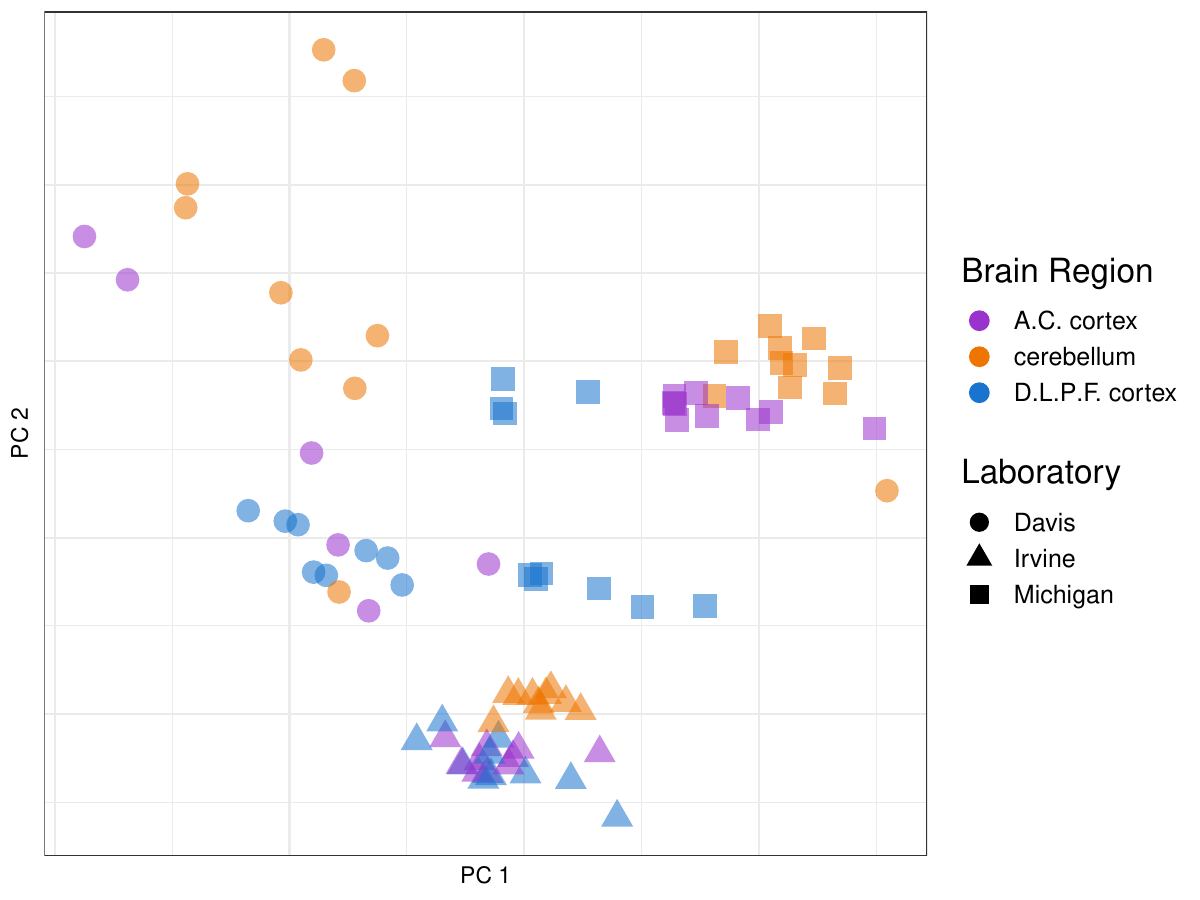} &
    \includegraphics[scale=0.3]{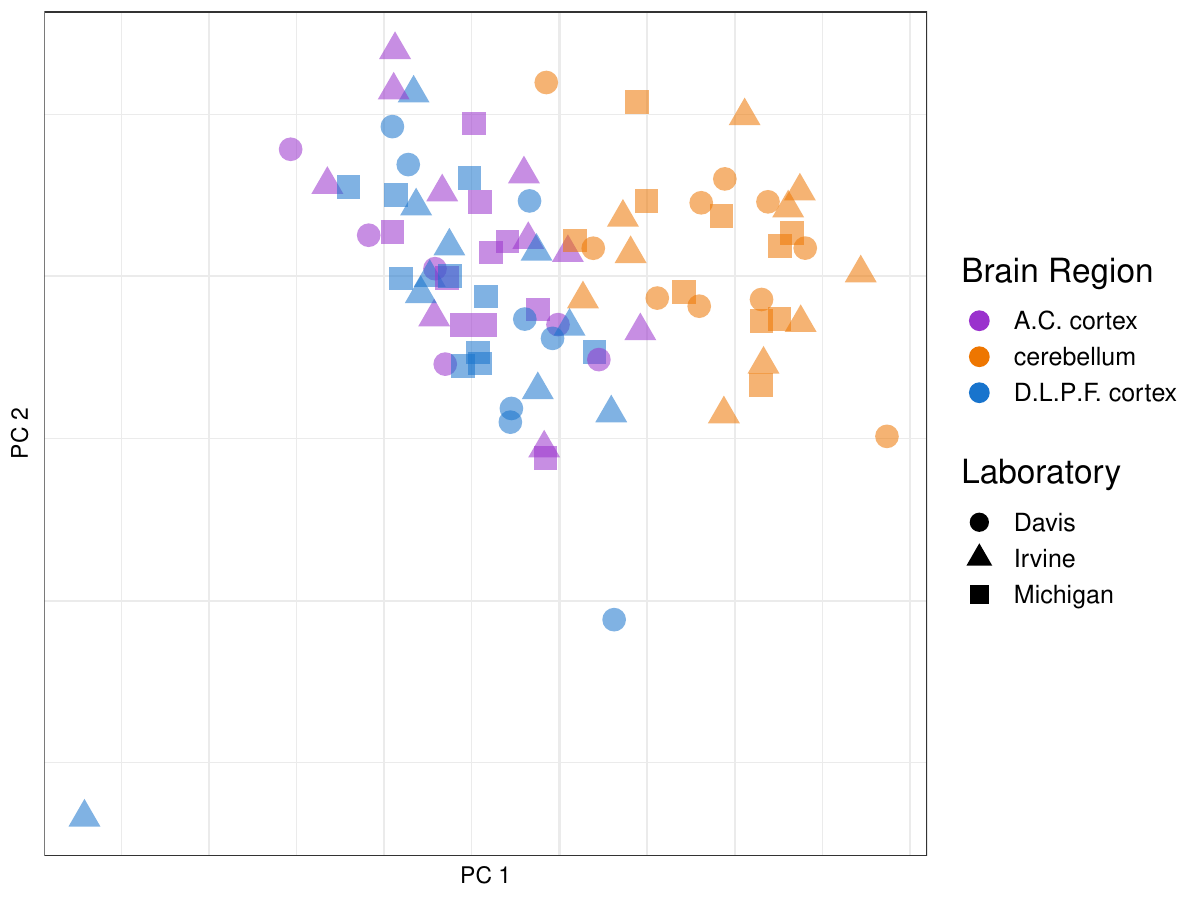} \\    
    \hspace{-1cm} ``Bio'' Adjustment & 
    \hspace{-1cm} ``Bio'' Adjustment (colored by patient) \\
    \includegraphics[scale=0.3]{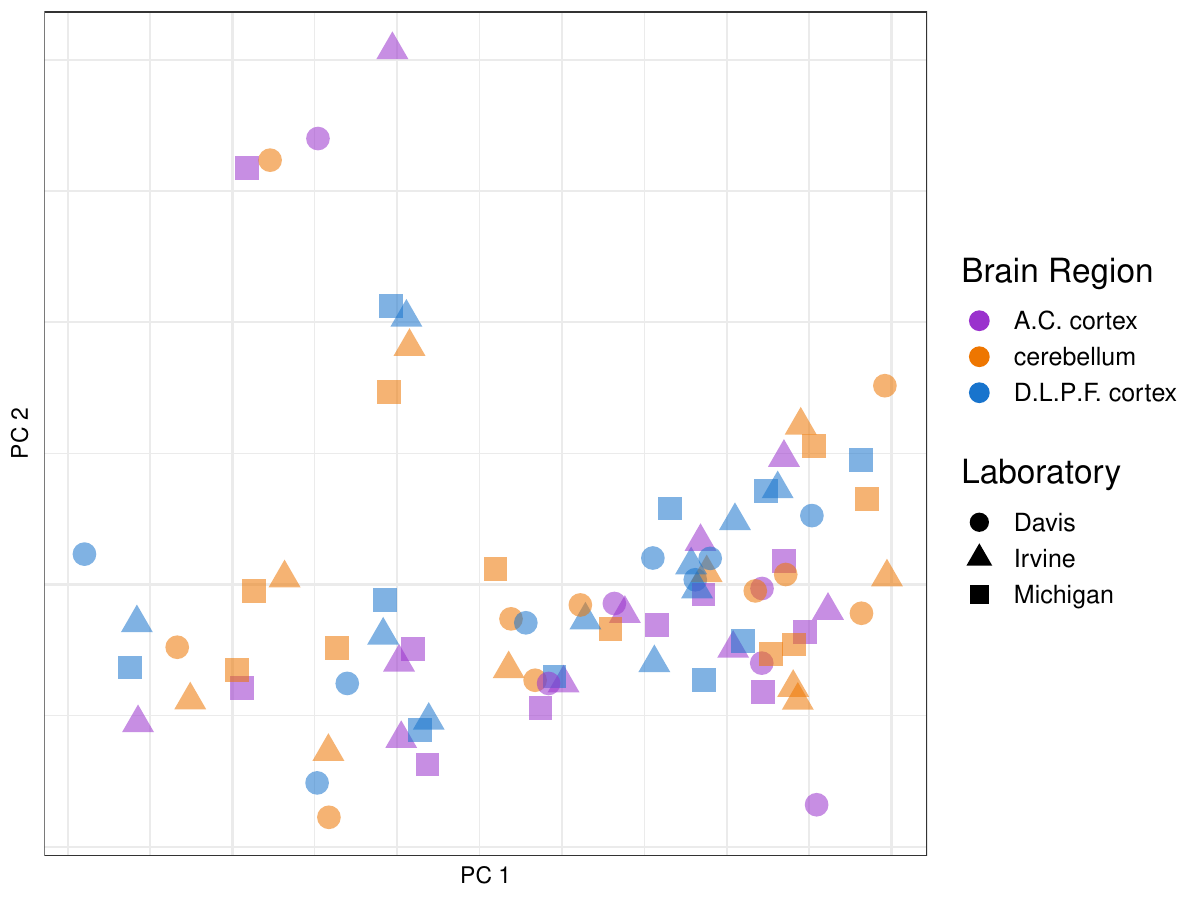}  &
    \includegraphics[scale=0.3]{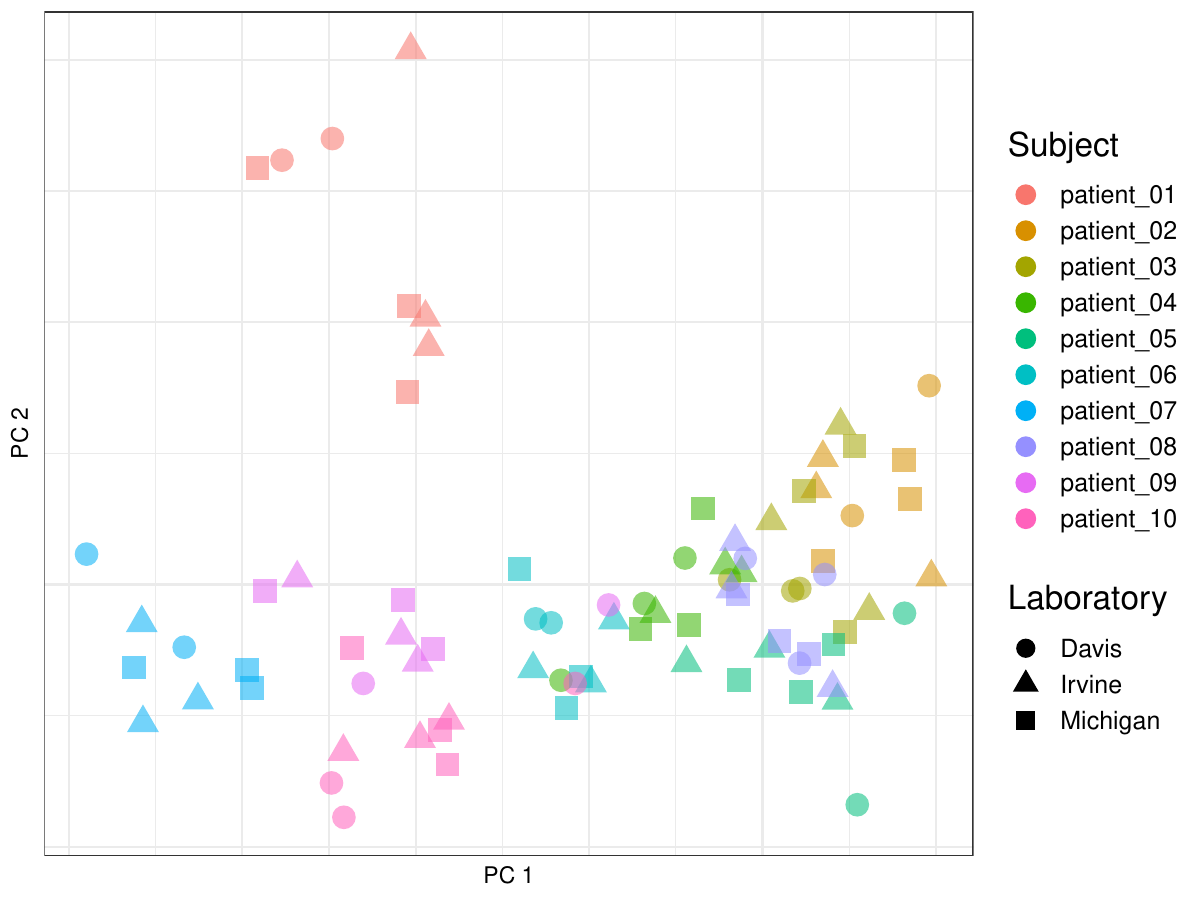}     
    \end{tabular}
    \caption{PC plots of the brain dataset before and after adjustment by RUV-III.  The ``technical'' adjustment uses spike-in controls and define replicates purely as technical replicates (same sample but different lab).  The ``bio'' adjustment uses housekeeping genes and defines to include all assays from a single individual.}
    \label{fig:svd.all}
\end{figure}

In Figure \ref{fig:svd.all}, top-left panel, we plot the first two PCs of the data.  We see that the primary clustering is by laboratory, indicating strong batch effects.  Notably, these batch effects persist even though the data have been quantile normalized by RMA \citep{bolstad2003comparison, irizarry2003summaries, irizarry2003exploration}.

We now apply RUV-III as follows.  As negative controls, we use the $n_c = 33$ spike-in probesets.  The mapping matrix is the $84\times30$ matrix in which each column corresponds to a single tissue sample.  We set $k=10$ as in \citet{gagnonbartsch2012using}, where a value of $k$ was determined via quality metrics such as RLE plots \citep{bolstad2005quality, brettschneider2008quality}.  The results are shown in Figure \ref{fig:svd.all}, top-right panel.  The laboratory batch effects have been removed, and the primary clustering is now by brain region; cerebellum forms one cluster, and the two types of cortex form the other cluster.  Thus RUV-III has succeeded in its primary goal --- it has removed unwanted technical variation while preserving biological variation.

\begin{figure}[h]
    \centering
    \begin{tabular}{cc}
    \multicolumn{2}{c}{No Adjustment}\\
    \includegraphics[scale=0.3]{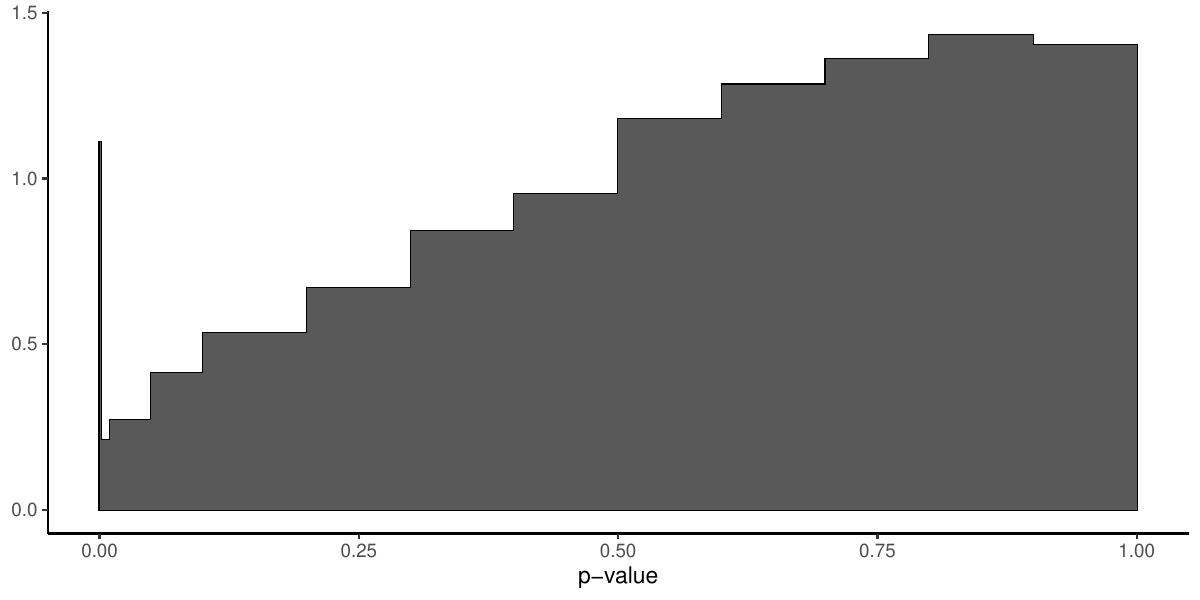} &
    \includegraphics[scale=0.3]{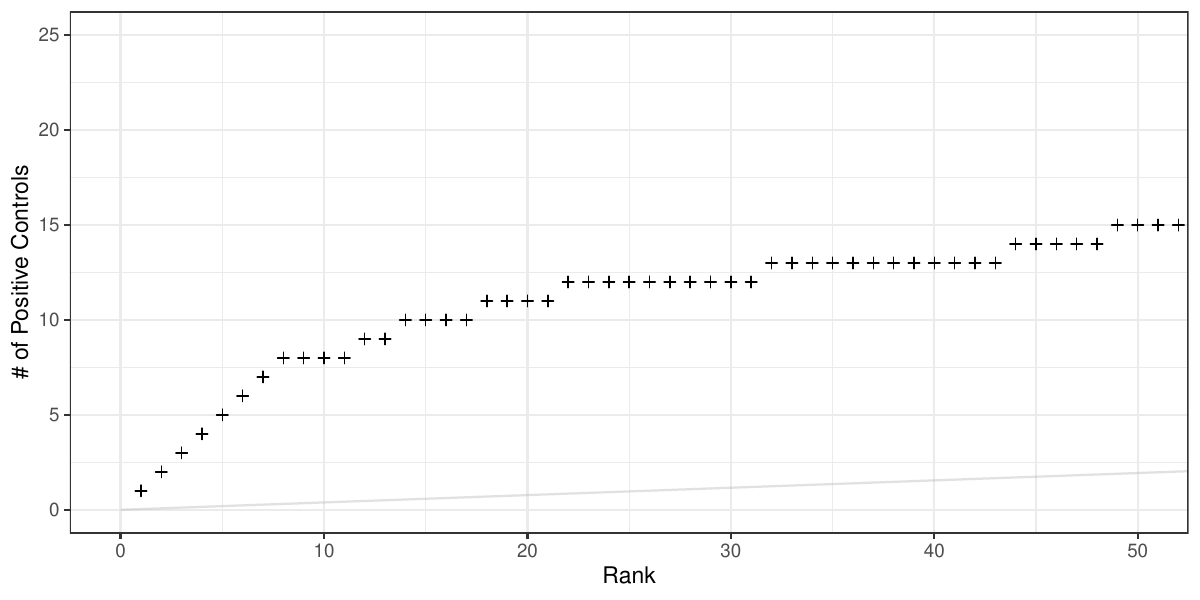} \\    
    \multicolumn{2}{c}{``Technical'' Adjustment}\\
    \includegraphics[scale=0.3]{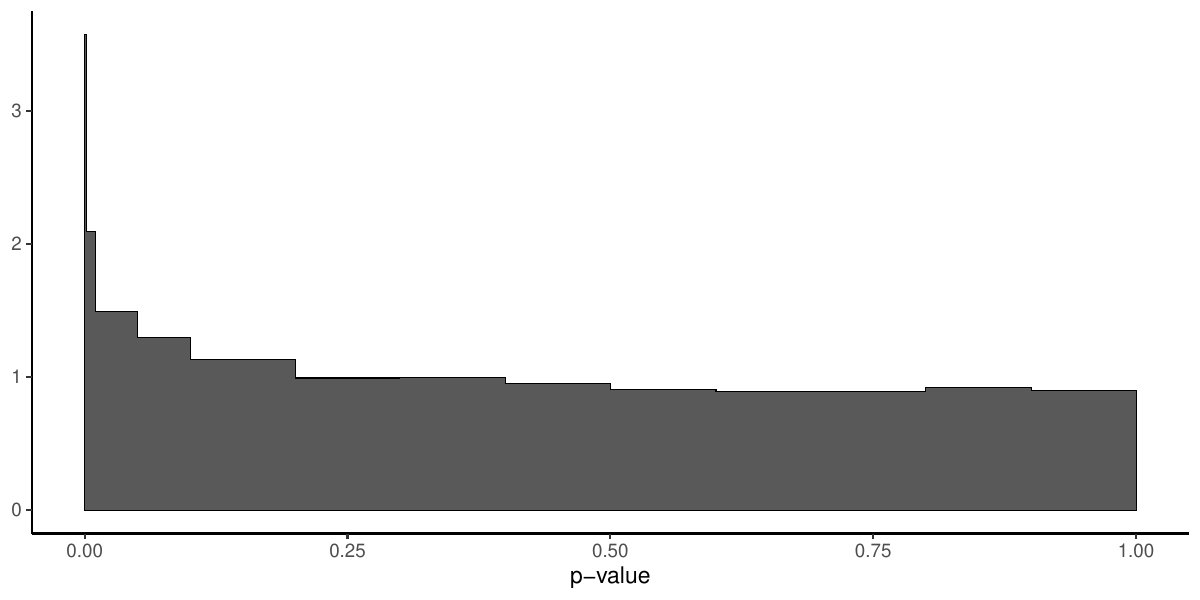} &
    \includegraphics[scale=0.3]{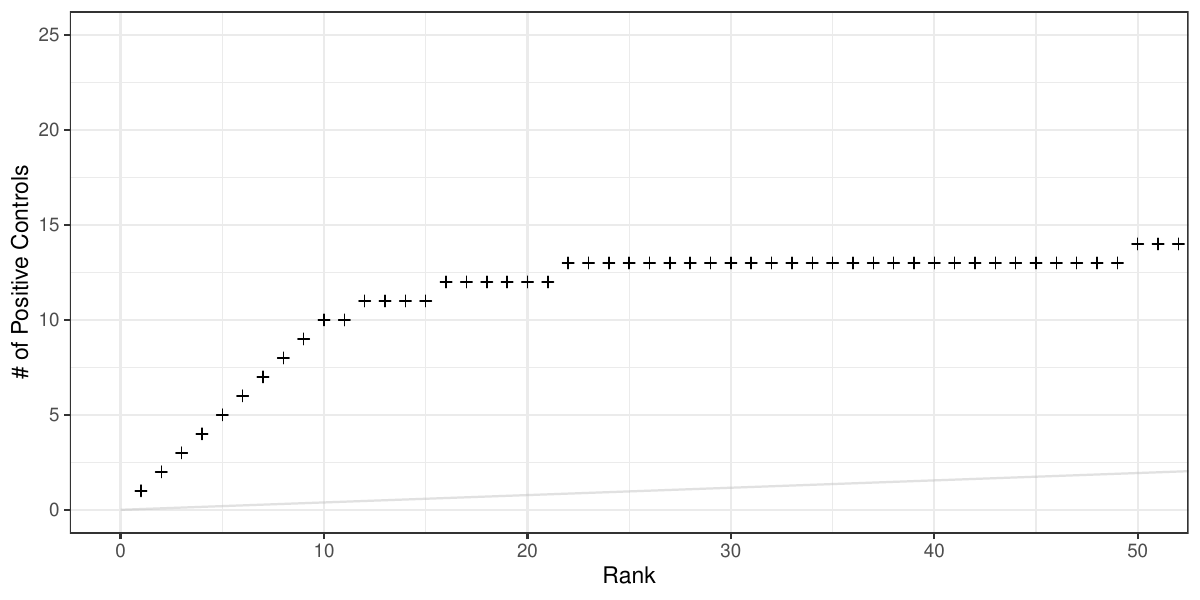} \\    
    \multicolumn{2}{c}{``Bio'' Adjustment}\\
    \includegraphics[scale=0.3]{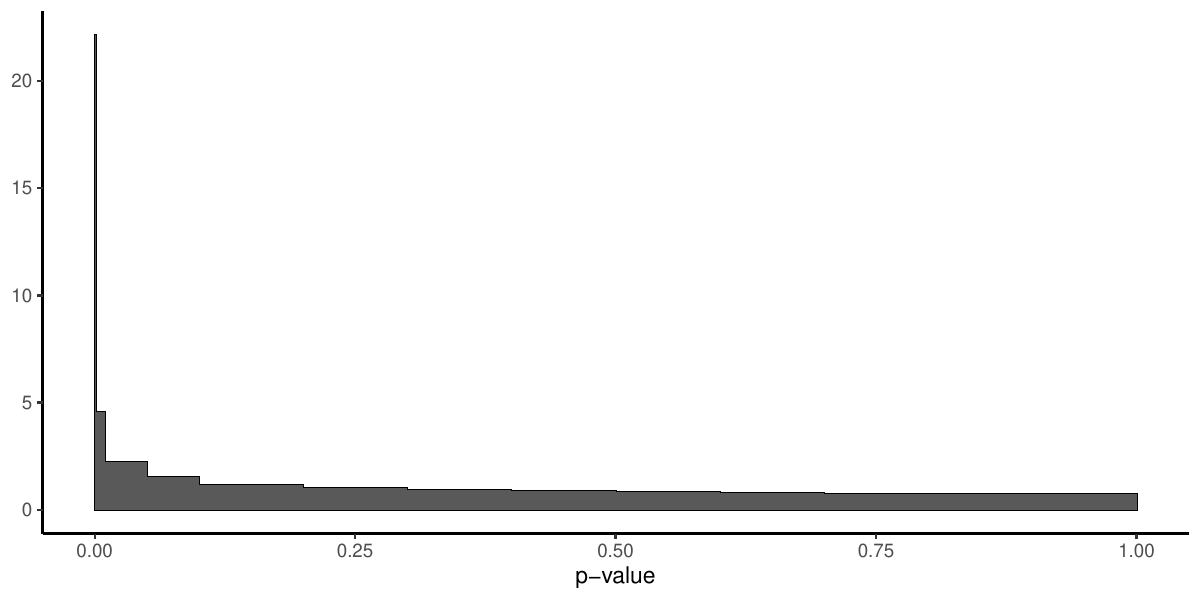} &
    \includegraphics[scale=0.3]{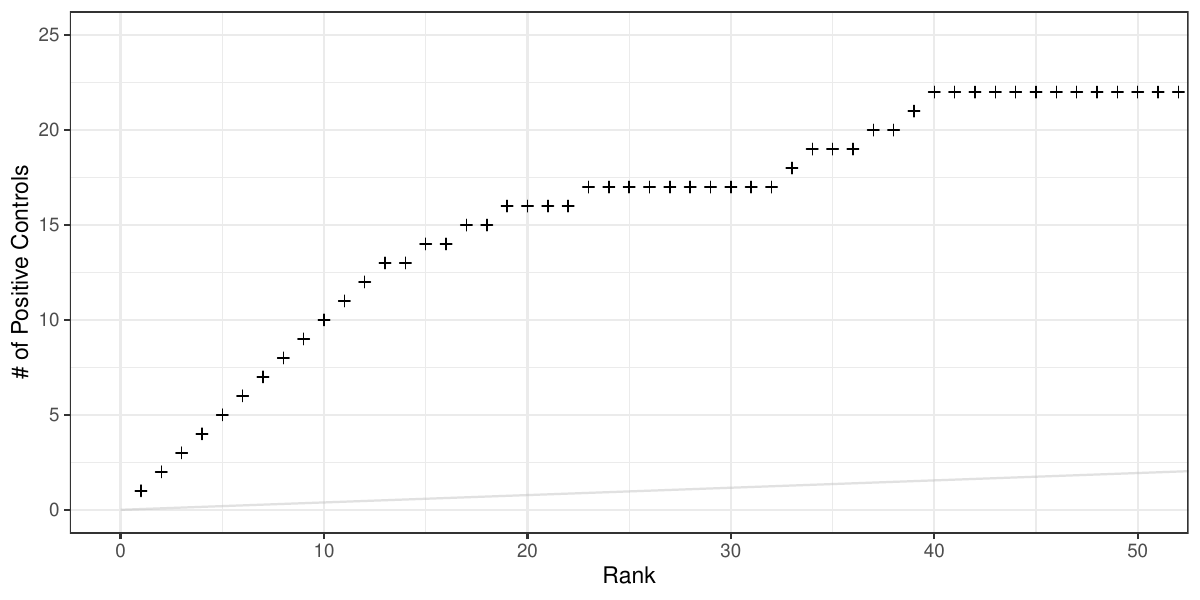} \\    
    \end{tabular}
    \caption{P-value histograms (left) and rank plots (right) before and after RUV-III.  Histogram break points are at 0.001, 0.01, 0.05, 0.1, 0.15, 0.2, ..., 0.95.}
    \label{fig:pvalrank.all}
\end{figure}

\begin{figure}[h]
    \centering
    \begin{tabular}{cc}
    \hspace{-1cm} (\textbf{a}) & 
    \hspace{-1cm} (\textbf{b}) \\
    \includegraphics[scale=0.3]{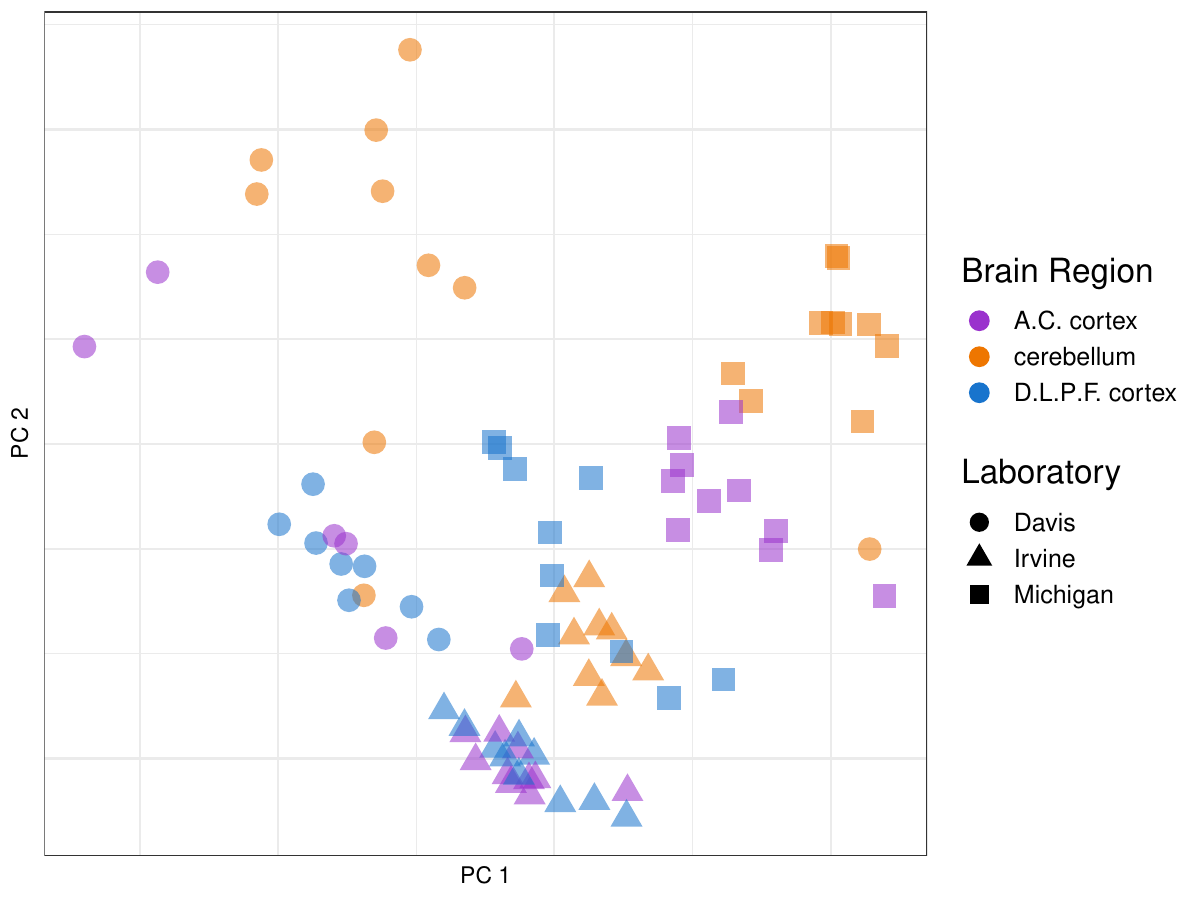} &
    \includegraphics[scale=0.3]{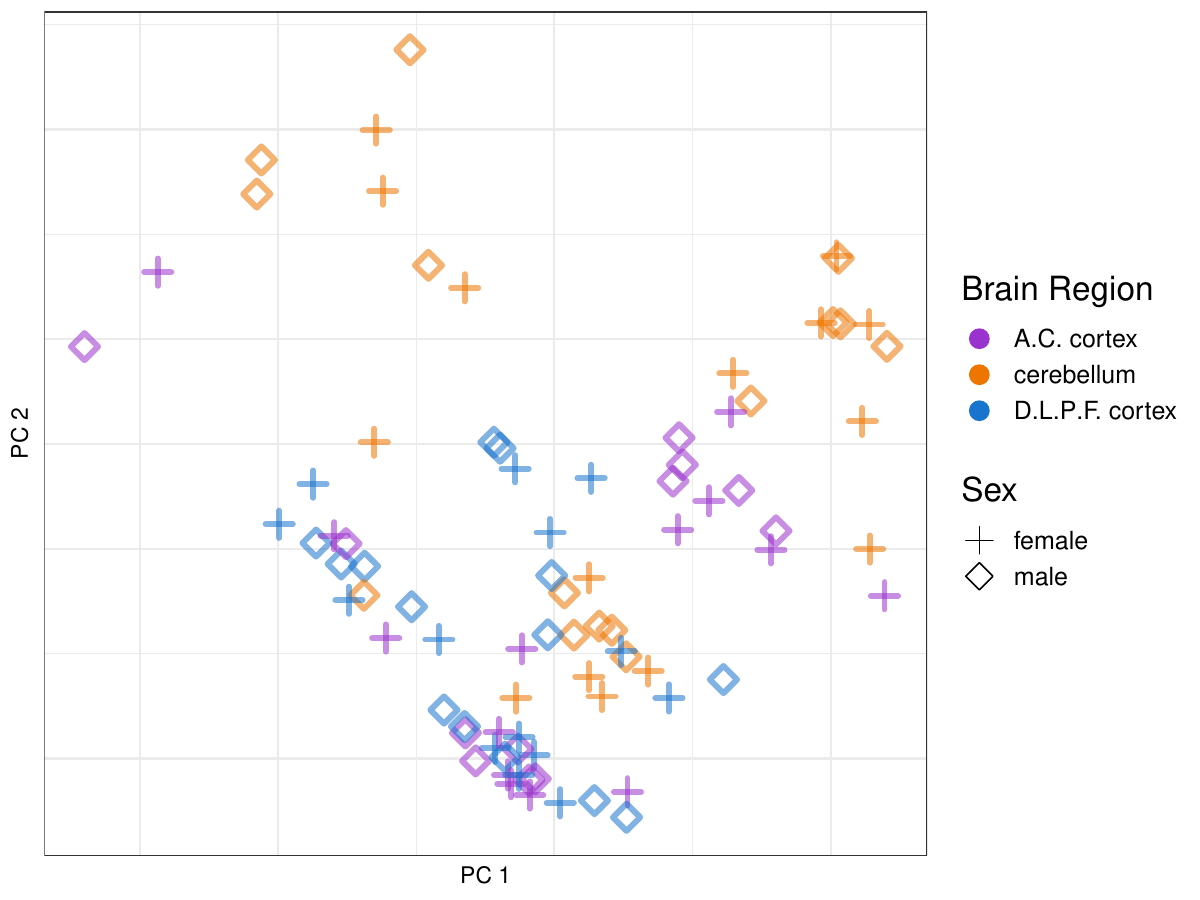} \\    
    \hspace{-1cm} (\textbf{c}) & 
    \hspace{-1cm} (\textbf{d}) \\
    \includegraphics[scale=0.3]{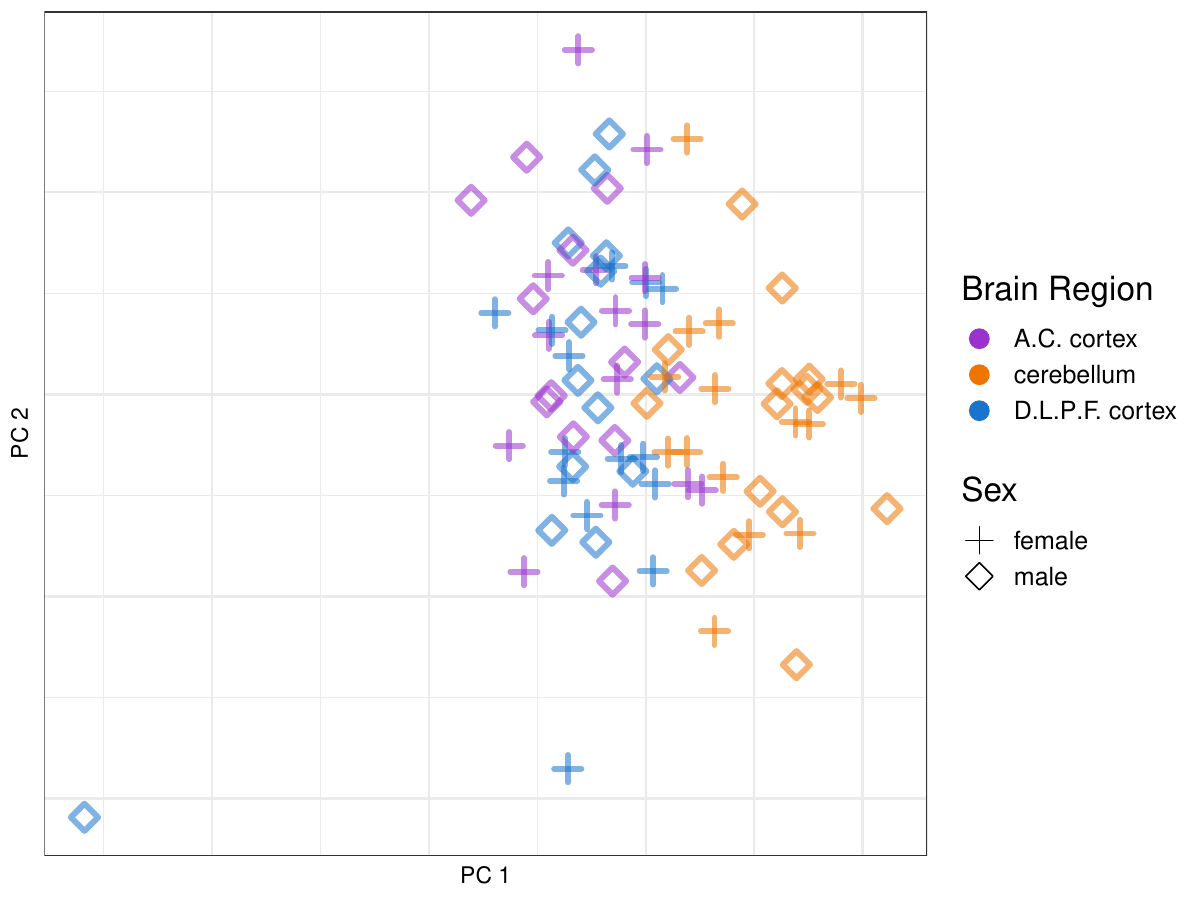}  &
    \includegraphics[scale=0.3]{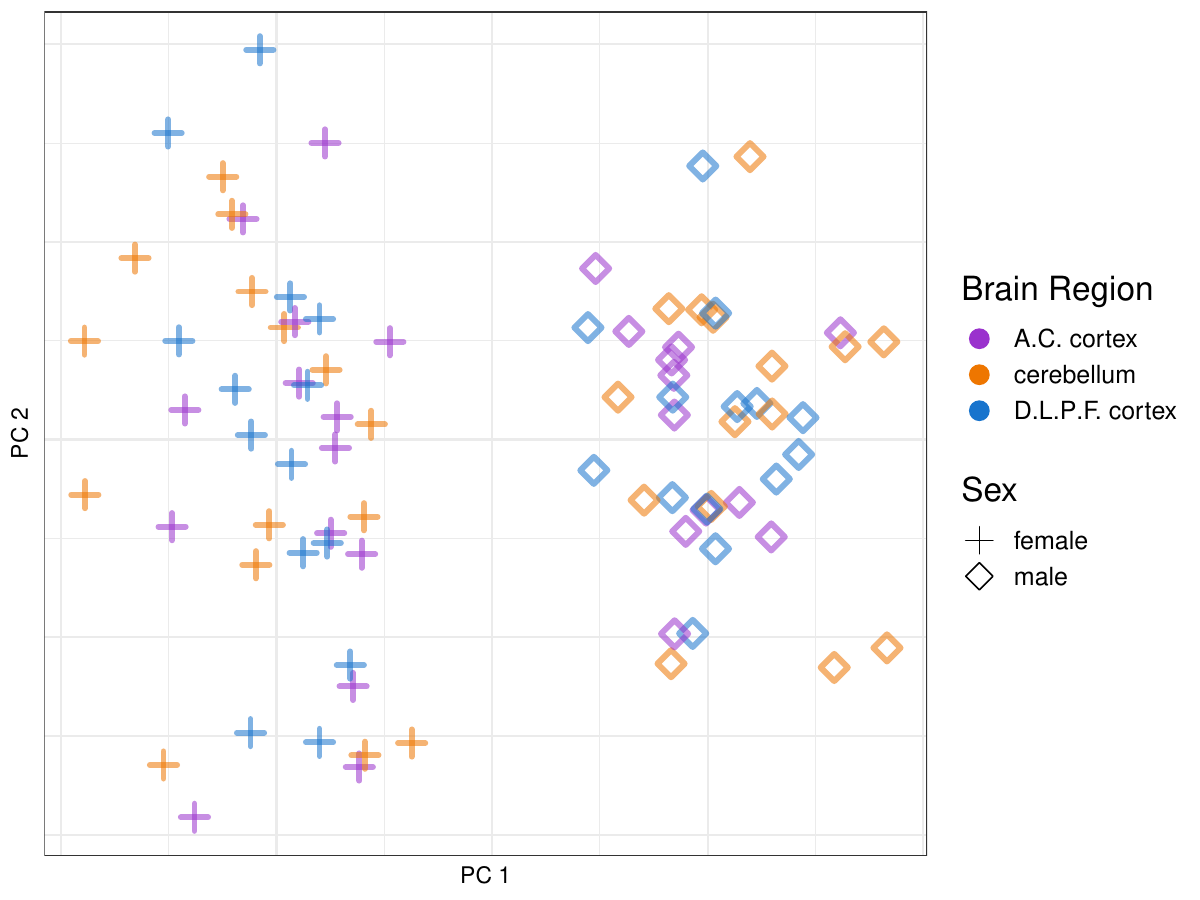}     
    \end{tabular}
    \caption{PC plots of X/Y genes before and after adjustment. (\textbf{a}) No adjustment.  Like top-left panel of Figure \ref{fig:svd.all}, but with PCs computed only on X/Y genes.  (\textbf{b}) No adjustment.  Like panel (a), but plot symbol denotes sex. (\textbf{c}) ``Technical'' Adjustment. Like panel (b), but after RUV-III technical adjustment. (\textbf{d}) ``Bio'' Adjustment. Like panel (b), but after RUV-III bio adjustment.}
    \label{fig:svd.xy}
\end{figure}

The Supplementary Material provides analogs of Figure \ref{fig:svd.all} using different values of $k$, specifically $k = 1, 2, 5, 15, 20, 30$  (S6--S11).  Notably, for larger values of $k$ (15, 20, 30) the results look qualitatively similar to those of Figure \ref{fig:svd.all}, which is consistent with Theorem 3. The Supplementary Material also provides analogs of Figures \ref{fig:pvalrank.all} and \ref{fig:svd.xy} (to be discussed below) using different values of $k$, and again the results are qualitatively stable for larger values of $k$.

One of the advantages of RUV-III is flexibility in the choice of negative controls and flexibility in the specification of replicates.  This flexibility can be used to effectively define which sources of variation are wanted and which are unwanted.  For example, although brain region effects might be of interest in some contexts, if the primary goal of the analysis is to find genes differentially expressed by sex, then we might wish to remove brain region effects in addition to batch effects.  

In this spirit, we consider an alternative application of RUV-III.  As negative controls, we use $n_c = 799$ housekeeping genes \citep{eisenberg2003human}.  Housekeeping genes tend to have relatively stable expression, and we expect them to be negative controls with respect to sex.  Nonetheless, even housekeeping genes may exhibit expression differences across tissue types, and we therefore expect some of these genes to capture brain region differences.  As the mapping matrix we use the $84\times10$ matrix in which each column corresponds to an individual person.  Thus the ``replicates'' span the three tissue types, allowing us to capture tissue-type variation in $W$.  On the other hand, since each set of replicates corresponds to one individual, the replicates should \emph{not} capture sex differences.

The results are shown in Figure \ref{fig:svd.all}, bottom-left panel.  Indeed, there is no longer clustering by brain region (nor by laboratory).  The bottom-right panel of Figure \ref{fig:svd.all} is the same plot but colored by individual.  We see that now the assays cluster by person.  Thus in this application of RUV-III we have preserved some sources of biological variation while removing others.

Figure \ref{fig:pvalrank.all} reveals how these adjustments impact a differential expression analysis.  We regress each gene on sex and compute a p-value.  The p-values were computed by applying Limma \citep{smyth2004linear} to the RUV-III adjusted data matrix.  A histogram of the 12,600 p-values should ideally be mostly flat, with a spike on the left corresponding to the genes that are truly differentially expressed.  As seen in Figure \ref{fig:pvalrank.all}, however, prior to RUV-III the p-value histogram is skewed to the right, with only a small spike on the left.  After the ``technical'' adjustment, the p-value distribution appears to be improved, and more so after the ``bio'' adjustment.  (We note, however, that the asymptotic theory developed in this paper does not address the validity of the p-values.)  In addition, we rank the genes by p-value, and check to see how many of the top-ranked genes lie on the sex chromosomes (X or Y).  More specifically, for each value of $N$ from 1 to 50 we count the number of genes located on X or Y that are ranked in the top $N$ genes.  The results are plotted on the right side of Figure \ref{fig:pvalrank.all}, and we see that the ``bio'' adjustment substantially increases the number of top-ranked genes on the sex chromosomes.  For example, before adjustment 15 of the 50 top-ranked genes are on the X or Y chromosome, but this increases to 24 after the ``bio'' adjustment.  (For context, 488 of the 12,600 genes --- fewer than 4\% --- lie on the X or Y chromosome.)

To further investigate the ability of RUV-III to remove unwanted variation while preserving variation of interest, in Figure 6 we create PC plots similar to those of Figure \ref{fig:svd.all}, but using only genes from the sex chromosomes.  (The application of RUV-III is exactly the same as in Figure \ref{fig:svd.all}; the only difference is which genes are used to calculate the PCs.)  We expect many of these genes to be differentially expressed with respect to sex, and therefore might naturally expect males and females to form clusters on a PC plot.   Panels (\textbf{a}) and (\textbf{b}) plot the first two PCs before adjustment.  Perhaps surprisingly, we see that even when we focus on just the sex chromosomes, there is no clustering by male/female apparent in the first two PCs.  In panel (\textbf{a}) we see that the strongest factor is still the laboratory batch effect, and in panel (\textbf{b}) we see that there is no discernible clustering by sex even within labs; the plusses and diamonds appear randomly distributed.  Likewise, in panel (\textbf{c}) we see that the primary signal after the ``technical'' adjustment is still brain region.  Again, although we have limited our PC analysis to X/Y genes, there is no discernible clustering by sex.  In panel (\textbf{d}), however, after applying the ``bio'' adjustment, sex differences appear prominently in the first PC.

We emphasize that this example is not intended to suggest intrinsic superiority of housekeeping genes over spike-in controls, or of using biological replicates over technical replicates.  Rather, we emphasize that the choice of negative controls and the specification of what constitutes a replicate plays a critical role in determining which variation is captured and removed by RUV-III and which variation is preserved.

\section{Appendix: Proofs} \label{app:proofs}

Throughout the proofs, $L_{r}$ denotes a generic matrix satisfying
$$\|L_{r}\|=o(\|W\alpha\|^{r}), \;\;r=1,2.$$
We use w. p. $\rightarrow 1$ for ``with probability tending to one''.


\subsubsection*{Proof of Theorem 1}

Assumption (I) implies $1\leq k\leq m-s$ w. p. $\rightarrow 1$. We have
\begin{eqnarray}
P_{1^{\perp}}Y_{\rm c}&=&W\alpha_{\rm c}+P_{1^{\perp}}\epsilon_{\rm c},\label{eq:equ1}\\
U_{(k)}'Y_{\rm c}&=&U_{(k)}'W\alpha_{\rm c}+U_{(k)}'\epsilon_{\rm c},\label{eq:equ2}\\
\hat{\alpha}_{(k)}&=&U_{(k)}'W\alpha+U_{(k)}'\epsilon,\label{eq:equ3}
\end{eqnarray}
For (\ref{eq:equ2}), note that, by (\ref{eq:M_expression}), we have
\begin{eqnarray}
P_{M}&=&{\rm diag}(1_{m_{h}}, 1\leq h\leq s){\rm diag}(m_{h}^{-1}, 1\leq h\leq s){\rm diag}(1_{m_{h}}', 1\leq h\leq s)
\nonumber\\
&=&{\rm diag}(P_{1_{m_{h}}}, 1\leq h\leq s),\label{eq:P_expression}
\end{eqnarray}
where $P_{1_{d}}=d^{-1}1_{d}1_{d}'$. It then follows that $P_{M^{\perp}}1=0$. Also, (\ref{eq:eigen_decom1})
implies
\begin{eqnarray}
0=P_{M^{\perp}}YY'P_{M^{\perp}}1=\sum_{l=1}^{m-s}\lambda_{l}u_{l}u_{l}'1.\label{eq:equ4}
\end{eqnarray}
The assumptions imply (see the arguments below) that $\lambda_{j}>0, 1\leq j\leq k$ w. p. $\rightarrow 1$.
For any $1\leq j\leq k$, by the orthonormality of the $u_{l}$'s, we have [by multiplying $u_{j}'$ from the left on
both sizes of (\ref{eq:equ4})], $\lambda_{j}u_{j}'1=0$; then, because $\lambda_{j}>0$, we have $u_{j}'1=0$, and
this holds for any $1\leq j\leq k$, implying $U_{(k)}'1=0$. Similarly, it can be shown that $U_{(k)}'M=0$, leading to
(\ref{eq:equ3}).

We now establish (\ref{eq:expr1}). We use the following notation: Define $A_{\rm c}$, $B$, $G$ and $H$ as below (\ref{eq:expr2}). Let $\xi=U_{(k)}'\epsilon$, $\xi_{\rm c}=U_{(k)}'\epsilon_{\rm c}$, $\eta_{\rm c}=P_{1^{\perp}}
\epsilon_{\rm c}$; then let $\Delta_{1}=B\alpha_{\rm c}\xi_{\rm c}'$, $\Delta_{2}=\xi_{\rm c}\xi_{\rm c}'$, $\Delta
=\Delta_{1}+\Delta_{1}'+\Delta_{2}$; $\Gamma_{1}=W\alpha_{\rm c}\xi_{\rm c}'$, $\Gamma_{2}=\eta_{\rm
c}\alpha_{\rm c}'B'$, $\Gamma_{3}=\eta_{\rm c}\xi_{\rm c}'$, $\Gamma=\Gamma_{1}+\Gamma_{2}+\Gamma_{3}$.
By (\ref{eq:Wkhat}), (\ref{eq:equ1})--(\ref{eq:equ3}), it can be shown that $U_{(k)}'Y_{\rm c}Y_{\rm c}'U_{(k)}=G
+\Delta$, $P_{1^{\perp}}Y_{\rm c}Y_{\rm c}'U_{(k)}=H+\Gamma$, thus, we have
$$\hat{W}_{(k)}=HG^{-1}(I_{k}+\Delta G^{-1})^{-1}+\Gamma G^{-1}(I_{k}+\Delta G^{-1})^{-1}.$$
Thus, it can be verified that (\ref{eq:expr1}) holds with $T_{1}=HG^{-1}(Q-I_{k})B\alpha$, $T_{2}=\Gamma G^{-1}Q
B\alpha$, $T_{3}=HG^{-1}Q\xi$, $T_{4}=\Gamma G^{-1}Q\xi$, where $Q=(I_{k}+\Delta G^{-1})^{-1}$.

It remains to show that all of the terms $T_{s}, s=1, 2, 3, 4$ are $L_{1}$ (see notation at the beginning of this
section). A main task is working on $G$. Write
\begin{eqnarray}
G&=&U_{(k)}'W\alpha_{\rm c}\alpha_{\rm c}'W'U_{(k)}\nonumber\\
&=&n_{\rm c}U_{(k)}'W\Sigma W'U_{(k)}+n_{\rm c}U_{(k)}'W(n_{\rm c}^{-1}\alpha_{\rm c}\alpha_{\rm c}'-\Sigma)
W'U_{(k)}.\label{eq:Gexpr1}
\end{eqnarray}
Assumption (II) implies $n_{\rm c}^{-1}\alpha_{\rm c}\alpha_{\rm c}'-\Sigma=n_{\rm c}^{-1}\sum_{j=1}^{n_{\rm
c}}(\alpha_{j}\alpha_{j}'-\Sigma)=O_{\rm P}(1)/\sqrt{n_{\rm c}}$. Furthermore, (\ref{eq:Wcond1}) implies $\|W\|^{2}
=\lambda_{m}O(1)$, hence
$$\|U_{(k)}'W(n_{\rm c}^{-1}\alpha_{\rm c}\alpha_{\rm c}'-\Sigma)W'U_{(k)}\| \leq(\lambda_{m}/\sqrt{n_{\rm
c}})O_{\rm P}(1).$$
It then follows, by (\ref{eq:Gexpr1}), that we can write
\begin{eqnarray}
G=n_{\rm c}U_{(k)}'W\Sigma W'U_{(k)}+\lambda_{m}\sqrt{n_{\rm c}}O_{\rm P}(1).\label{eq:Gexpr2}
\end{eqnarray}
On the other hand, because $P_{M^{\perp}}Y=P_{M^{\perp}}(W\alpha+\epsilon)$, we have
\begin{eqnarray}
P_{M^{\perp}}YY'P_{M^{\perp}}&=&P_{M^{\perp}}W\alpha\alpha'W'P_{M^{\perp}}+P_{M^{\perp}}(\epsilon\epsilon'
+\epsilon\alpha'W'+W\alpha\epsilon')P_{M^{\perp}}\nonumber\\
&=&nP_{M^{\perp}}W\Sigma W'P_{M^{\perp}}+nP_{M^{\perp}}W(n^{-1}\alpha\alpha'-\Sigma)W'P_{M^{\perp}}
\nonumber\\
&&+P_{M^{\perp}}(\epsilon\epsilon'+\epsilon\alpha'W'+W\alpha\epsilon')P_{M^{\perp}}.\label{eq:expan1}
\end{eqnarray}
By similar arguments, the second term on the right side of (\ref{eq:expan1}) is $\lambda_{m}\sqrt{n}O_{\rm P}(1)$. Furthermore,
note that $\|P_{M^{\perp}}\|\leq 1$, $\|\epsilon\|=\sqrt{n}O_{\rm P}(1)$ by the RMT (e.g., \citep{jiang2022large2}, ch. 16); also, it
can be shown that $\|\alpha\|=\sqrt{n}O_{\rm P}(1)$. It follows that the third term on the right side of (\ref{eq:expan1})
is $n\sqrt{\lambda_{m}}O_{\rm P}(1)$. Note that $U_{(k)}'P_{M^{\perp}}YY'P_{M^{\perp}}U_{(k)}=[u_{g}'P_{M^{\perp}}YY'
P_{M^{\perp}}u_{h}]_{1\leq g,h\leq k}$. Thus, assuming $k\leq m-s$, which holds w. p. $\rightarrow 1$, we have, by
(\ref{eq:eigen_decom1}),
$$u_{g}'P_{M^{\perp}}YY'P_{M^{\perp}}u_{h}=\sum_{l=1}^{m-s}\lambda_{l}u_{g}'u_{l}u_{l}'u_{h},$$
which is $0$ if $g\neq h$; if $g=h$, the expression is equal to $\lambda_{g}$. Thus, combining the above
results, we have
\begin{eqnarray}
\lambda_{k}I_{k}&\leq&{\rm diag}(\lambda_{1},\dots,\lambda_{k})\nonumber\\
&=&U_{(k)}'P_{M^{\perp}}YY'P_{M^{\perp}}U_{(k)}\nonumber\\
&=&nU_{(k)}'W\Sigma W'U_{(k)}+(\lambda_{m}\sqrt{n}+n\sqrt{\lambda_{m}})O_{\rm P}(1),\label{eq:ub1}
\end{eqnarray}
noting that $P_{M^{\perp}}U_{(k)}=U_{(k)}$ and $\|U_{(k)}\|=1$. By(\ref{eq:Gexpr2}) and (\ref{eq:ub1}), we have
\begin{eqnarray}
G&=&\left(\frac{n_{\rm c}}{n}\right)nU_{(k)}'W\Sigma W'U_{(k)}+\lambda_{m}\sqrt{n_{\rm c}}O_{\rm P}(1)\nonumber\\
&\geq&\left(\frac{n_{\rm c}}{n}\right)\{\lambda_{k}I_{k}-(\lambda_{m}\sqrt{n}+n\sqrt{\lambda_{m}})O_{\rm P}(1)\}+\lambda_{m}\sqrt{n_{\rm c}}O_{\rm
P}(1)\nonumber\\
&=&\lambda_{k}\left(\frac{n_{\rm c}}{n}\right)I_{k}-(\lambda_{m}\sqrt{n}+n\sqrt{\lambda_{m}})O_{\rm P}(1).\label{eq:G_lb1}
\end{eqnarray}

We now utilize a result from the eigenvalue perturbation theory. For any $m\times m$ symmetric matrix $A$, let
$\lambda_{1}(A)\geq\cdots\geq\lambda_{m}(A)$ denote its eigenvalues arranged in decreasing order. We have
the following inequality (e.g., \citep{jiang2022large2}, p. 156).

{\bf Lemma 1 (Weyl's eigenvalue perturbation theorem).} For any symmetric matrices $A, B$, we have
$\max_{1\leq i\leq m}|\lambda_{i}(A)-\lambda_{i}(B)|\leq\|A-B\|$.

Now $\lambda_{1},\dots,\lambda_{m}$ are the eigenvalues of $P_{M^{\perp}}YY'P_{M^{\perp}}$ arranged in
decreasing order. Let $\gamma_{1},\dots,\gamma_{m}$ denote the eigenvalues of $F_{M}=P_{M^{\perp}}W\alpha
\alpha'W'P_{M^{\perp}}$, also arranged in decreasing order. Then, by Lemma 1 and (\ref{eq:expan1}), we have
\begin{eqnarray}
\max_{1\leq i\leq k}|\lambda_{i}-\gamma_{i}|\leq\|P_{M^{\perp}}(\epsilon\epsilon'+\epsilon\alpha'W'
+W\alpha\epsilon')P_{M^{\perp}}\|=n\sqrt{\lambda_{m}}O_{\rm P}(1).\label{eq:eigen_dfb1}
\end{eqnarray}
Next, we have $F_{M}=nP_{M^{\perp}}W(n^{-1}\alpha\alpha')W'P_{M^{\perp}}=nJ_{M}+nP_{M^{\perp}}W(
n^{-1}\alpha\alpha'-\Sigma)W'P_{M^{\perp}}$, where $J_{M}=P_{M^{\perp}}W\Sigma W'P_{M^{\perp}}$. Let
$\rho_{1},\dots,\rho_{m}$ denote the eigenvalues of $J_{M}$, arranged in decreasing order. Then, by Lemma 1, we have
\begin{eqnarray}
\max_{1\leq i\leq k}|\gamma_{i}-n\rho_{i}|\leq n\|P_{M^{\perp}}W(n^{-1}\alpha\alpha'
-\Sigma)W'P_{M^{\perp}}\|=\lambda_{m}\sqrt{n}O_{\rm P}(1).\label{eq:eigen_dib2}
\end{eqnarray}
Note that $J_{M}$ is $m\times m$. Thus, assuming ${\rm rank}(J_{M})\geq k$, which holds w. p. $\rightarrow
1$ under the assumptions, we have $\rho_{k}\geq$ the smallest nonzero eigenvalue of $J_{M}$, which is the
same as the smallest eigenvalue of $\Sigma^{1/2}W'P_{M^{\perp}}W\Sigma^{1/2}$, provided that the latter,
which is $k\times k$, is positive definite, which again holds w. p. $\rightarrow 1$ under the assumptions. It
can be shown that
$W'P_{M^{\perp}}W=\sum_{h=1}^{s}(W_{h}-1_{m_{h}}\bar{W}_{h}')'(W_{h}-1_{m_{h}}\bar{W}_{h}')$.
Thus, we have $\rho_{k}\geq\lambda_{\min}(W'P_{M^{\perp}}W)\lambda_{\min}(\Sigma)
=\lambda_{\min}(\Sigma)\lambda_{m}$. Note that $\rho_{k}>0$ implies ${\rm rank}(J_{M})=k$. Combining the
above results, we have obtained a asymptotic lower bound for $\lambda_{k}$:
\begin{eqnarray}
\lambda_{k}&=&n\rho_{k}+\lambda_{k}-n\rho_{k}\nonumber\\
&\geq&n\lambda_{m}\lambda_{\min}(\Sigma)-\max_{1\leq i\leq k}|\lambda_{i}-n\rho_{i}|\nonumber\\
&\geq&\lambda_{m}n\left\{\lambda_{\min}(\Sigma)-\frac{O_{\rm P}(1)}{\sqrt{\lambda_{m}}}-\frac{O_{\rm P}(1)}{\sqrt{n}}\right\}.
\label{eq:ldk_lb1}
\end{eqnarray}

Combining (\ref{eq:G_lb1}) and (\ref{eq:ldk_lb1}), we have
\begin{eqnarray}
G&\geq&\left(\frac{n_{\rm c}}{n}\right)\lambda_{m}n\left\{\lambda_{\min}(\Sigma)-\frac{O_{\rm P}(1)}{\sqrt{\lambda_{m}}}-\frac{O_{\rm P}(1)}{\sqrt{n}}\right\}I_{k}-(\lambda_{m}\sqrt{n}+n\sqrt{\lambda_{m}})O_{\rm P}(1)\nonumber\\
&=&\lambda_{m}n\left[\left(\frac{n_{\rm c}}{n}\right)\left\{\lambda_{\min}(\Sigma)-\frac{O_{\rm P}(1)}{\sqrt{\lambda_{m}}}
-\frac{O_{\rm P}(1)}{\sqrt{n}}\right\}I_{k}-\frac{O_{\rm P}(1)}{\sqrt{n}}-\frac{O_{\rm P}(1)}{\sqrt{\lambda_{m}}}\right]\nonumber\\
&=&\lambda_{m}n\left[\left(\frac{n_{\rm c}}{n}\right)\lambda_{\min}(\Sigma)I_{k}-\frac{O_{\rm P}(1)}{\sqrt{\lambda_{m}}}-\frac{O_{\rm P}(1)}{\sqrt{n}}\right].\label{eq:G_lb2}
\end{eqnarray}
It then follows, by assumption (I), that we have the following key result:
\begin{eqnarray}
G^{-1}=(\lambda_{m}n)^{-1}O_{\rm P}(1).\label{eq:G_inv}
\end{eqnarray}

Also, we have $\Delta_{1}=U_{(k)}'W\alpha\epsilon'U_{(k)}=U_{(k)}'P_{M^{\perp}}W\alpha\epsilon'U_{(k)}$, hence $\|\Delta_{1}\|\leq\|P_{M^{\perp}}W\alpha\epsilon'\|$. We have $\|P_{M^{\perp}}\|=\sqrt{\lambda_{m}}O(1)$ by (\ref{eq:Wcond1}); $\|\alpha\|=\sqrt{n}O_{\rm P}(1)$ by assumption (II); and $\|\epsilon\|=\sqrt{n}O_{\rm P}(1)$
by assumption (III) and the RMT. It follows that $\|\Delta_{1}\|=\sqrt{\lambda_{m}}nO_{\rm P}(1)$. Similarly, it can be shown
that $\|\Delta_{2}\|=nO_{\rm P}(1)$. Thus, we have $\|\Delta\|=\sqrt{\lambda_{m}}nO_{\rm P}(1)$, hence $\Delta G^{-1}=
O_{\rm P}(1)/\sqrt{\lambda_{m}}=o_{\rm P}(1)$ by (\ref{eq:G_inv}), assumption (I) and (\ref{eq:Wcond2}), hence $Q-I_{k}=O_{\rm P}(1)/\sqrt{\lambda_{m}}$.

Note that $T_{1}=HG^{-1}(Q-I_{k})U_{(k)}'W\alpha$, and $H=W\alpha_{\rm c}\alpha_{\rm c}'W'P_{M^{\perp}}
U_{(k)}$. Thus, we have
$\|H\|\leq\|W\alpha_{\rm c}\|\sqrt{U_{(k)}'P_{M^{\perp}}W\alpha_{\rm c}\alpha_{\rm c}'W'P_{M^{\perp}}U_{(k)}}
=\lambda_{m}n_{\rm c}O_{\rm P}(1)$,
which is implied by the assumptions [in particular, (\ref{eq:Wcond1})]. It follows that
$$\|HG^{-1}(Q-I_{k})U_{(k)}'\|\leq\lambda_{m}n_{\rm c}O_{\rm P}(1)\left\{\frac{O_{\rm P}(1)}{\lambda_{m}n}\right\}\left\{\frac{O_{\rm
P}(1)}{\sqrt{\lambda_{m}}}\right\}\leq\frac{O_{\rm P}(1)}{\sqrt{\lambda_{m}}}.$$
Thus, by (\ref{eq:Wcond2}), we have $T_{1}=o_{\rm P}(1)W\alpha=L_{1}$ (see notation at the beginning of this section).

Next, by (\ref{eq:G_inv}), we have $T_{2}=\Gamma\{(\lambda_{m}n)^{-1}O_{\rm P}(1)\}U_{(k)}'W\alpha$. By similar
arguments, it can be shown that $\|\Gamma\|=\sqrt{\lambda_{m}}nO_{\rm P}(1)$, and $\|U_{(k)}\|=1$. Thus, we have
$T_{2}=o_{\rm P}(1)W\alpha=L_{1}$.

As for $T_{j}, j=3,4$, we first obtain an asymptotic lower bound for $\|W\alpha\|$. By the following
inequalities (e.g., Ex. 5.27 of \citep{jiang2022large2} 2022):
$$\lambda_{\max}(A)+\lambda_{\min}(B)\leq\lambda_{\max}(A+B)\leq\lambda_{\max}(A)
+\lambda_{\max}(B)$$
for matrices $A, B$ of same dimensions, and the arguments above, it can be shown that
$$\left|\|W\alpha\|^{2}-n\lambda_{\max}(W\Sigma W')\right|\leq\lambda_{m}\sqrt{n}O_{\rm P}(1).$$
Furthermore, using the fact that $W'W\geq\lambda_{\min}(W'W)I_{k}$ and some known matrix inequalities
[e.g., (ii), (iii) of \citep{jiang2022large2}, sec. 5.3.1], we have
\begin{eqnarray*}
&&\lambda_{\max}(W\Sigma W')=\lambda_{\max}(\Sigma^{1/2}W'W\Sigma^{1/2})\geq\lambda_{\min}(W'W)\lambda_{\max}(\Sigma)\\
&&\geq\lambda_{\min}(W'P_{M^{\perp}}W)\lambda_{\max}(\Sigma)=\lambda_{\max}(\Sigma)\lambda_{m}.
\end{eqnarray*}
It follows that, w. p. $\rightarrow 1$, $\|W\alpha\|/\sqrt{\lambda_{m}n}$ has a positive lower bound.
On the other hand, it can be shown that, under assumptions (II) and (III), and (\ref{eq:Wcond1}), one has
$\|T_{3}\|=\sqrt{n}O_{\rm P}(1)$, and $\|T_{4}\|=\sqrt{n/\lambda_{m}}O_{\rm P}(1)$. Therefore, by (\ref{eq:Wcond2}), both $T_{3}$ and $T_{4}$ are $L_{1}$.

The proof is complete by (\ref{eq:expr1}) and Definition 1 in Section 2.
\subsubsection*{Proof of Theorem 2}

First, let us obtain a lower bound for the difference between the right side of (\ref{eq:expr4}), denoted by $g(k)$,
and its maximum, $g(k_{0})$. We have
 \begin{eqnarray}
 g(k_{0})-g(k)=n{\rm tr}\{W\Sigma^{1/2}P_{\perp}(k)\Sigma^{1/2}W'\},\label{eq:diff_1}
 \end{eqnarray}
 where $P_{\perp}(k)=I_{k_{0}}-\Sigma^{1/2}B'(B\Sigma B')^{-1}B\Sigma^{1/2}$, which is a projection matrix. Therefore, we have
 \begin{eqnarray}
 {\rm tr}\{W\Sigma^{1/2}P_{\perp}(k)\Sigma^{1/2}W'\}&=&{\rm tr}\{P_{\perp}(k)\Sigma^{1/2}W'W\Sigma^{1/2}
 P_{\perp}(k)\}\nonumber\\
 &\geq&\lambda_{\min}(W'W){\rm tr}\{P_{\perp}(k)\Sigma P_{\perp}(k)\}\nonumber\\
 &\geq&\lambda_{\min}(W'P_{M^{\perp}}W)\lambda_{\min}(\Sigma){\rm tr}\{P_{\perp}(k)\}\nonumber\\
 &=&\lambda_{\min}(\Sigma)(k_{0}-k)\lambda_{m}.\label{eq:tr_lb1}
 \end{eqnarray}
 It follows that the right side of (\ref{eq:diff_1}) has an asymptotic lower bound, which  is $\lambda_{\min}(\Sigma)(k_{0}-k)\lambda_{m}n$.
 
 Next, consider the difference between the right sides of (\ref{eq:expr3}) and (\ref{eq:expr4}). It can be shown
 that the difference can be written as $\sum_{s=1}^{5}d_{s}$, where
 \begin{eqnarray*}
 d_{1}&=&{\rm tr}[\{(BA_{\rm c}B')^{-1}-(Bn_{\rm c}\Sigma B')^{-1}\}BA_{\rm c}W'WA_{\rm c}B'(BA_{\rm c}B')^{-1}
 B\alpha\alpha'B'],\\
 d_{2}&=&{\rm tr}\{(Bn_{\rm c}\Sigma B')^{-1}B(A_{\rm c}-n_{\rm c}\Sigma)W'WA_{\rm c}B'(BA_{\rm c}B')^{-1}B
 \alpha\alpha'B'\},\\
 d_{3}&=&{\rm tr}\{(Bn_{\rm c}\Sigma B')^{-1}Bn_{\rm c}\Sigma W'W(A_{\rm c}-n_{\rm c}\Sigma)B'(BA_{\rm
 c}B')^{-1}B\alpha\alpha'B'\},\\
 d_{4}&=&{\rm tr}[(Bn_{\rm c}\Sigma B')^{-1}Bn_{\rm c}\Sigma W'Wn_{\rm c}\Sigma B'\{(BA_{\rm c}B')^{-1}-
 (Bn_{\rm c}\Sigma B')^{-1}\}\\
 &&\times B\alpha\alpha'B'],\\
 d_{5}&=&{\rm tr}\{(Bn_{\rm c}\Sigma B')^{-1}Bn_{\rm c}\Sigma W'Wn_{\rm c}\Sigma B'(Bn_{\rm c}\Sigma
 B')^{-1}B(\alpha\alpha'-n\Sigma)B'\}.
 \end{eqnarray*}
 By the previous arguments [in particular, (\ref{eq:Gexpr2}) and  (\ref{eq:G_inv})], it can be shown that
 \begin{eqnarray}
 d_{j}=\lambda_{m}no_{\rm P}(1),\;\;j=1,2,3,4,5.\label{eq:ds_order}
 \end{eqnarray}
 By similar arguments, it can be also shown that
 \begin{eqnarray}
 \|T_{j}\|_{2}^{2}=\lambda_{m}no_{\rm P}(1),\;\;j=1,2,3,4\label{eq:Ts_order}
 \end{eqnarray}
 for the $T_{j}$ in (\ref{eq:expr2}). It follows, by (\ref{eq:expr2}) and the Cauchy-Schwarz inequality that
 \begin{eqnarray}
 \|\hat{W}_{(k)}\hat{\alpha}_{(k)}\|_{2}^{2}=\|HG^{-1}B\alpha\|_{2}^{2}+\lambda_{m}no_{\rm P}(1),\label{eq:cs_appr}
 \end{eqnarray}
 
 Combining (\ref{eq:expr3}), (\ref{eq:expr4}), and the above results, we have, w. p. $\rightarrow 1$,
 \begin{eqnarray}
\min_{1\leq k<k_{0}}\left\{\|\hat{W}_{(k_{0})}\hat{\alpha}_{(k_{0})}\|_{2}^{2}-\|\hat{W}_{(k)}\hat{\alpha}_{(k)}\|_{2}^{2}
\right\}\geq\frac{a_{1}}{2}\lambda_{\min}(\Sigma)\lambda_{m}n;
\end{eqnarray}
in other words, (\ref{eq:thm2}) holds with $\delta_{1}=a_{1}/2$.
This completes the proof.
\subsubsection*{Proof of Theorem 3}

First, let us continue on with the special case of $n=n_{\rm c}$, discussed in the early part of Section 3.2. By (\ref{eq:proj_decom}), we have
\begin{eqnarray}
\hat{W}_{(k)}\hat{\alpha}_{(k)}=P_{1^{\perp}}Y_{\rm c}P=P_{1^{\perp}}Y_{\rm c}P_{0}+P_{1^{\perp}}Y_{\rm c}P_{+}.
\label{eq:UV_decom1}
\end{eqnarray}
We have $P_{1^{\perp}}Y_{\rm c}=P_{1^{\perp}}W\alpha+P_{1^{\perp}}\epsilon=W\alpha+L_{1}$, because
$\|P_{1^{\perp}}\epsilon\|\leq\|\epsilon\|=\sqrt{n}O_{\rm P}(1)$ by the RMT. Thus (using the Cauchy-Schwarz
inequality for the cross products), we have
\begin{eqnarray}
&&P_{1^{\perp}}Y_{\rm c}PY_{\rm c}'P_{1^{\perp}}=(W\alpha+L_{1})P(\alpha'W'+L_{1}')\nonumber\\
&&=W\alpha P\alpha'W'+L_{2}\leq W\alpha\alpha'W'+L_{2}\label{eq:res1}
\end{eqnarray}
(see notation introduced at the beginning of Section 6). Furthermore, according to the proof of Theorem 1, we
have $P_{1^{\perp}}Y_{\rm c}P_{0}=W\alpha+L_{1}$. Thus, we have
\begin{eqnarray}
&&(P_{1^{\perp}}Y_{\rm c}P_{0}+P_{1^{\perp}}Y_{\rm c}P_{+})(P_{0}Y_{\rm c}'P_{1^{\perp}}+P_{+}Y_{\rm
c}'P_{1^{\perp}})\nonumber\\
&=&P_{1^{\perp}}Y_{\rm c}P_{0}Y_{\rm c}'P_{1^{\perp}}+P_{1^{\perp}}Y_{\rm c}P_{+}Y_{\rm c}'P_{1^{\perp}}
\nonumber\\
&=&(W\alpha+L_{1})(\alpha'W'+L_{1}')+P_{1^{\perp}}Y_{\rm c}P_{+}Y_{\rm c}'P_{1^{\perp}}\nonumber\\
&=&W\alpha\alpha'W'+L_{2}+P_{1^{\perp}}Y_{\rm c}P_{+}Y_{\rm c}'P_{1^{\perp}},\label{eq:res2}
\end{eqnarray}
using the fact that $P_{0}P_{+}=P_{+}P_{0}=0$ and, again, the Cauchy-Schwarz inequality for the cross
products. Combining the second equation in (\ref{eq:UV_decom1}), (\ref{eq:res1}), and (\ref{eq:res2}), we
have $W\alpha\alpha'W'+L_{2}+P_{1^{\perp}}Y_{\rm c}P_{+}Y_{\rm c}'P_{1^{\perp}}\leq W\alpha\alpha'W'
+L_{2}$, hence $P_{1^{\perp}}Y_{\rm c}P_{+}Y_{\rm c}'P_{1^{\perp}}\leq L_{2}$, which implies
$P_{1^{\perp}}Y_{\rm c}P_{+}=L_{1}$. Thus, by (\ref{eq:UV_decom1}) and combining the above results, we
have established that
\begin{eqnarray}
\hat{W}_{(k)}\hat{\alpha}_{(k)}=W\alpha+L_{1}
\end{eqnarray}
for every fixed $k_{0}<k\leq K$, where $K$ is the maximum number of $k$ under consideration.

It is worth noting that, in fact, the arguments of the above proof apply to any fixed sequence of $k$, say, $k=k_{N}$, which $N$ denotes the combination of $m, n, n_{\rm c}$, etc. that may increase, such that $k_{N}\leq m-s$ and $k_{N}\rightarrow\infty$. In particular, the arguments apply to $k_{N}=m-s$.

Now let us remove the $n=n_{\rm c}$ assumption. In this case, $U_{(k)}'Y_{\rm c}Y_{\rm c}'U_{(k)}$ is not
necessarily block-diagonal. However, by the inversion of blocked matrices (e.g., \citep{bernstein2009matrix} p.44), we have
\begin{eqnarray*}
&&\{U_{(k)}'Y_{\rm c}Y_{\rm c}'U_{(k)}\}^{-1}\\
&=&\left[\begin{array}{cc}U_{(k_{0})}'Y_{\rm c}Y_{\rm c}'U_{(k_{0})}&U_{(k_{0})}'Y_{\rm c}Y_{\rm c}'U_{(+)}\\
U_{(+)}'Y_{\rm c}Y_{\rm c}'U_{(k_{0})}&U_{(+)}'Y_{\rm c}Y_{\rm c}'U_{(+)}\end{array}\right]^{-1}\\
&=&\left[\begin{array}{c}Q_{00}^{-1}+Q_{00}^{-1}Q_{0+}\{U_{(+)}'Y_{\rm c}P_{0}^{\perp}Y_{\rm c}'U_{(+)}\}^{-1}Q_{+0}Q_{00}^{-1}\\
-\{U_{(+)}'Y_{\rm c}P_{0}^{\perp}Y_{\rm c}'U_{(+)}\}^{-1}Q_{+0}Q_{00}^{-1}\end{array}\right.\\
&&\left.\begin{array}{c}-Q_{00}^{-1}Q_{0+}\{U_{(+)}'Y_{\rm c}P_{0}^{\perp}Y_{\rm c}'U_{(+)}\}^{-1}\\
\{U_{(+)}'Y_{\rm c}P_{0}^{\perp}Y_{\rm c}'U_{(+)}\}^{-1}\end{array}\right],
\end{eqnarray*}
where $Q_{00}=U_{(k_{0})}'Y_{\rm c}Y_{\rm c}'U_{(k_{0})}$, $Q_{0+}=U_{(k_{0})}'Y_{\rm c}Y_{\rm c}'
U_{(+)}$, $Q_{+0}=Q_{0+}'$, and $P_{0}^{\perp}=I_{n_{\rm c}}-P_{0}$ with $P_{0}$ defined above. Also retain
the notation $P$ previously define. It follows that we have the following decomposition:
\begin{eqnarray}
P&=&P_{0}+Y_{\rm c}'U_{(k_{0})}Q_{00}^{-1}Q_{0+}\{U_{(+)}'Y_{\rm c}P_{0}^{\perp}Y_{\rm
c}'U_{(+)}\}^{-1}Q_{+0}Q_{00}^{-1}U_{(k_{0})}'Y_{\rm c}\nonumber\\
&&-Y_{\rm c}'U_{(+)}\{U_{(+)}'Y_{\rm c}P_{0}^{\perp}Y_{\rm c}'U_{(+)}\}^{-1}Q_{+0}Q_{00}^{-1}U_{(k_{0})}'Y_{\rm c}
\nonumber\\
&&-Y_{\rm c}'U_{(k_{0})}Q_{00}^{-1}Q_{0+}\{U_{(+)}'Y_{\rm c}P_{0}^{\perp}Y_{\rm c}'U_{(+)}\}^{-1}U_{(+)}'Y_{\rm c}
\nonumber\\
&&+Y_{\rm c}'U_{(+)}\{U_{(+)}'Y_{\rm c}P_{0}^{\perp}Y_{\rm c}'U_{(+)}\}^{-1}U_{(+)}'Y_{\rm c}
\nonumber\\
&=&P_{0}+Y_{\rm c}'\Lambda\{U_{(+)}'Y_{\rm c}P_{0}^{\perp}Y_{\rm c}'U_{(+)}\}^{-1}\Lambda'Y_{\rm c}\nonumber\\
&=&P_{0}+P_{\Lambda}\label{eq:proj_decom1}
\end{eqnarray}
where $\Lambda=U_{(+)}-U_{(k_{0})}Q_{00}^{-1}Q_{0+}$ and $P_{\Lambda}=P_{Y_{\rm c}'\Lambda}$. Note that
$U_{(+)}'Y_{\rm c}P_{0}^{\perp}Y_{\rm c}'U_{(+)}=\Lambda'Y_{\rm c}Y_{\rm c}'\Lambda$. Also note that
$U_{(k_{0})}'Y_{\rm c}Y_{\rm c}'\Lambda=0$, which implies $P_{0}P_{\Lambda}=0$. Furthermore, by the previous
argument, we have $P_{1^{\perp}}Y_{\rm c}P_{0}=W\alpha_{\rm c}+L_{1}$. Thus, by similar argument to (\ref{eq:res2}),
we have
\begin{eqnarray}
P_{1^{\perp}}Y_{\rm c}PY_{\rm c}'P_{1^{\perp}}=W\alpha_{\rm c}\alpha_{\rm c}'W'+L_{2}+P_{1^{\perp}}Y_{\rm
c}P_{\Lambda}Y_{\rm c}'P_{1^{\perp}}.\label{eq:res3}
\end{eqnarray}
On the other hand, by similar argument to (\ref{eq:res1}), we have
\begin{eqnarray}
P_{1^{\perp}}Y_{\rm c}PY_{\rm c}'P_{1^{\perp}}\leq W\alpha_{\rm c}\alpha_{\rm c}'W+L_{2}.\label{eq:res4}
\end{eqnarray}
Combining (\ref{eq:res3}) and (\ref{eq:res4}), we have $0\leq P_{1^{\perp}}Y_{\rm c}P_{\Lambda}Y_{\rm
c}'P_{1^{\perp}}\leq L_{2}$, from which, it follows that $P_{1^{\perp}}Y_{\rm c}P_{\Lambda}=L_{1}$. Thus, combining
the above results, we have
\begin{eqnarray}
P_{1^{\perp}}Y_{\rm c}P=W\alpha_{\rm c}+L_{1}.\label{eq:P_res}
\end{eqnarray}

Next, similar to (\ref{eq:proj_decom1}), we have the following decomposition:
\begin{eqnarray}
Q=Q_{0}+Q_{\Lambda},\label{eq:Q_decom}
\end{eqnarray}
where the $Q$'s are the corresponding $P$'s with the last multiplication factor, $Y_{\rm c}$, replaced by $Y_{\rm
d}$, and $Y=(Y_{\rm c}\;Y_{\rm d})$; namely, we have
\begin{eqnarray*}
&&Q=Y_{\rm c}'U_{(k)}\{U_{(k)}'Y_{\rm c}Y_{\rm c}'U_{(k)}\}^{-1}U_{(k)}'Y_{\rm d},\\
&&Q_{0}=Y_{\rm c}'U_{(k_{0})}\{U_{(k_{0})}'Y_{\rm c}Y_{\rm c}'U_{(k_{0})}\}^{-1}U_{(k_{0})}'Y_{\rm d},\\
&&Q_{\Lambda}=Y_{\rm c}'\Lambda(\Lambda'Y_{\rm c}Y_{\rm c}'\Lambda)^{-1}\Lambda'Y_{\rm d}=Y_{\rm c}'\Lambda\{U_{(+)}'Y_{\rm c}P_{0}^{\perp}Y_{\rm c}'U_{(+)}\}^{-1}\Lambda'Y_{\rm d}.
\end{eqnarray*}
Note that the $Y_{\rm d}$ in $Q_{\Lambda}$ can be replaced by $\xi_{\rm d}$ [because $U_{(k)}'1=0$ and $U_{(k)}'M=0$ for every $1\leq k\leq m-s$; see below (\ref{eq:P_expression}) and (\ref{eq:equ4})]. Thus, we have
\begin{eqnarray*}
&&P_{1^{\perp}}Y_{\rm c}Q_{\Lambda}Q_{\Lambda}'Y_{\rm c}'P_{1^{\perp}}\\
&=&P_{1^{\perp}}Y_{\rm c}Y_{\rm c}'\Lambda(\Lambda'Y_{\rm c}Y_{\rm c}'\Lambda)^{-1}\Lambda'\xi_{\rm
d}\xi_{\rm d}'\Lambda(\Lambda'Y_{\rm c}Y_{\rm c}'\Lambda)^{-1}\Lambda'Y_{\rm c}Y_{\rm c}'P_{1^{\perp}}\\
&=&P_{1^{\perp}}\Pi_{\rm c}\Lambda'(\xi_{\rm c}\xi_{\rm c}')^{1/2}(\xi_{\rm c}\xi_{\rm c}')^{-1/2}\xi_{\rm d}\xi_{\rm d}'
(\xi_{\rm c}\xi_{\rm c}')^{-1/2}(\xi_{\rm c}\xi_{\rm c}')^{1/2}\Lambda\Pi_{\rm c}'P_{1^{\perp}}\\
&\leq&\|(\xi_{\rm c}\xi_{\rm c}')^{-1/2}\xi_{\rm d}\xi_{\rm d}'(\xi_{\rm c}\xi_{\rm c}')^{-1/2}\|P_{1^{\perp}}\Pi_{\rm c}
\Lambda'\xi_{\rm c}\xi_{\rm c}'\Lambda\Pi_{\rm c}'P_{1^{\perp}}\\
&=&O_{\rm P}(1)P_{1^{\perp}}Y_{\rm c}P_{\Lambda}Y_{\rm c}'P_{1^{\perp}}\\
&=&L_{2},\label{eq:res5}
\end{eqnarray*}
by (\ref{eq:add_cond}) and a previously established result, where $\Pi_{\rm c}=Y_{\rm c}Y_{\rm c}'\Lambda(\Lambda'Y_{\rm c}Y_{\rm c}'\Lambda)^{-1}$. It follows that $P_{1^{\perp}}Y_{\rm c}Q_{\Lambda}=L_{1}$.

On the other hand, write $\tilde{W}=W/\sqrt{\lambda_{m}}$, $\tilde{\alpha}_{\rm c}=\alpha_{\rm c}/\sqrt{n_{\rm c}}$, and
$\tilde{\alpha}_{\rm d}=\alpha_{\rm d}/\sqrt{n_{\rm d}}$. According to the previous arguments, we have
\begin{eqnarray}
\left\{\frac{U_{(k_{0})}'Y_{\rm c}Y_{\rm c}'U_{(k_{0})}}{\lambda_{m}n_{\rm c}}\right\}^{-1}&=&\{U_{(k_{0})}'\tilde{W}\tilde{\alpha}_{\rm c}\tilde{\alpha}_{\rm c}'\tilde{W}'U_{(k_{0})}\}^{-1}
+O_{\rm P}\left(\frac{1}{\sqrt{\lambda_{m}}}+\frac{1}{\sqrt{n}}\right),
\nonumber\\
\frac{Y_{\rm c}Y_{\rm c}'U_{(k_{0})}}{\lambda_{m}n_{\rm c}}&=&\tilde{W}\tilde{\alpha}_{\rm c}\tilde{\alpha}_{\rm c}'\tilde{W}'
U_{(k_{0})}+O_{\rm P}\left(\frac{1}{\sqrt{\lambda_{m}}}+\frac{1}{\sqrt{n}}\right),\nonumber\\
\frac{U_{(k_{0})}'Y_{\rm d}}{\sqrt{\lambda_{m}n_{\rm d}}}&=&U_{(k_{0})}'\tilde{W}\tilde{\alpha}_{\rm d}+O_{\rm P}\left(\frac
{1}{\sqrt{\lambda_{m}}}\right),\label{eq:res6}
\end{eqnarray}
and all of the none $O_{\rm P}$ expressions in (\ref{eq:res6}) are $O_{\rm P}(1)$. Thus, we have
\begin{eqnarray}
P_{1^{\perp}}Y_{\rm c}Q_{0}&=&\sqrt{\lambda_{m}n_{\rm d}}P_{1^{\perp}}\frac{Y_{\rm c}Y_{\rm c}'U_{(k_{0})}}{\lambda_{m}n_{\rm c}}
\left\{\frac{U_{(k_{0})}'Y_{\rm c}Y_{\rm c}'U_{(k_{0})}}{\lambda_{m}n_{\rm c}}\right\}^{-1}\frac{U_{(k_{0})}'Y_{\rm
d}}{\sqrt{\lambda_{m}n_{\rm d}}}\nonumber\\
&=&\sqrt{\lambda_{m}n_{\rm d}}P_{1^{\perp}}\left\{\tilde{W}\tilde{\alpha}_{\rm c}\tilde{\alpha}_{\rm c}'\tilde{W}'U_{(k_{0})}
+O_{\rm P}\left(\frac{1}{\sqrt{\lambda_{m}}}+\frac{1}{\sqrt{n}}\right)\right\}\nonumber\\
&&\times\left[\{U_{(k_{0})}'\tilde{W}\tilde{\alpha}_{\rm c}\tilde{\alpha}_{\rm c}'\tilde{W}'U_{(k_{0})}\}^{-1}+O_{\rm
P}\left(\frac{1}{\sqrt{\lambda_{m}}}+\frac{1}{\sqrt{n}}\right)\right]\nonumber\\
&&\times\left\{U_{(k_{0})}'\tilde{W}\tilde{\alpha}_{\rm d}+O_{\rm P}\left(\frac{1}{\sqrt{m}}\right)\right\}\nonumber\\
&=&\sqrt{\lambda_{m}n_{\rm d}}P_{1^{\perp}}\tilde{W}\tilde{\alpha}_{\rm c}\tilde{\alpha}_{\rm c}'\tilde{W}'U_{(k_{0})}\{U_{(k_{0})}'
\tilde{W}\tilde{\alpha}_{\rm c}\tilde{\alpha}_{\rm c}'\tilde{W}'U_{(k_{0})}\}^{-1}U_{(k_{0})}'\tilde{W}\tilde{\alpha}_{\rm d}\nonumber\\
&&+O_{\rm P}(\sqrt{\lambda_{m}}+\sqrt{n_{\rm d}})\nonumber\\
&=&W\tilde{\alpha}_{\rm c}\tilde{\alpha}_{\rm c}'\tilde{W}'U_{(k_{0})}\{\tilde{W}'U_{(k_{0})}\}^{-1}(\tilde{\alpha}_{\rm
c}\tilde{\alpha}_{\rm c}')^{-1}\{U_{(k_{0})}'\tilde{W}\}^{-1}U_{(k_{0})}'\tilde{W}\alpha_{\rm d}\nonumber\\
&&+O_{\rm P}(\sqrt{\lambda_{m}}+\sqrt{n_{\rm d}})\nonumber\\
&=&W\alpha_{d}+L_{1}.\label{eq:res7}
\end{eqnarray}
Combining the above results, we have
\begin{eqnarray}
P_{1^{\perp}}Y_{\rm c}Q=W\alpha_{\rm d}+L_{1}.\label{eq:Q_res}
\end{eqnarray}

Combining (\ref{eq:P_res}) and (\ref{eq:Q_res}), we have, for any $k>k_{0}$,
\begin{eqnarray}
&&\hat{W}_{(k)}\hat{\alpha}_{(k)}=P_{1^{\perp}}Y_{\rm c}(P\;Q)=[P_{1^{\perp}}Y_{\rm c}P\;\;P_{1^{\perp}}Y_{\rm c}Q]
\nonumber\\
&&=[W\alpha_{\rm c}+L_{1}\;\;W\alpha_{\rm d}+L_{1}]=W\alpha+L_{1}.\label{eq:Wah_ap_Wa}
\end{eqnarray}

Once again, it is seen that the above arguments apply to any fixed $k>k_{0}$, including a fixed sequence $k=k_{N}$ such that $k_{N}\leq m-s$ and $k_{N}\rightarrow\infty$, including $k_{N}=m-s$.
This completes the proof.
\subsubsection*{Proof of Theorem 5}

Let us begin with expression (\ref{eq:Mr}). It can be derived that, in this case, we have $P_{M}={\rm diag}(P_{M_{\rm
r}},I_{m_{\rm nr}})$ and $P_{M^{\perp}}={\rm diag}(P_{M_{\rm r}^{\perp}},0)$, and
\begin{eqnarray}
P_{M^{\perp}}YY'P_{M^{\perp}}={\rm diag}(P_{M_{\rm r}^{\perp}}Y_{\rm r}Y_{\rm r}'P_{M_{\rm r}^{\perp}},0).
\label{eq:prps_ed}
\end{eqnarray}
Let $S_{\rm r}=P_{M_{\rm r}^{\perp}}Y_{\rm r}Y_{\rm r}'P_{M_{\rm r}^{\perp}}=U_{\rm r}D_{\rm r}U_{\rm r}'$ be the
eigenvalue decomposition of $S_{\rm r}$, where $U_{\rm r}$ is an $m_{\rm r}\times m_{\rm r}$ orthogonal matrix, and
$D_{\rm r}$ the diagonal matrix of eigenvalues of $S_{\rm r}$ such that the eigenvalues are arranged in decreasing
order on the diagonal. Let $U={\rm diag}(U_{\rm r}, I_{s_{\rm nr}})$ and $D={\rm diag}(D_{\rm r},0)$. With these notation, we can
express (\ref{eq:prps_ed}) as
\begin{eqnarray}
P_{M^{\perp}}YY'P_{M^{\perp}}={\rm diag}(U_{\rm r}D_{\rm r}U_{\rm r},0)=UDU'.\label{eq:prps_ed1}
\end{eqnarray}
Note that (\ref{eq:prps_ed1}) is in exactly the same form of (\ref{eq:eigen_decom}). In practice, one would
choose $K$ such that $K\leq m_{\rm r}$ (recall $K$ is the maximum $k$ under consideration), which holds
with probability tending to one because $m_{\rm r}\rightarrow\infty$. Then, we have
$$U_{(k)}=\left[\begin{array}{c}U_{{\rm r}(k)}\\
0\end{array}\right],$$
for any $1\leq k\leq K$, where $U_{{\rm r}(k)}$ is the first $k$ columns of $U_{\rm r}$. Now write
$$Y=\left(\begin{array}{c}Y_{\rm r}\\
Y_{\rm nr}\end{array}\right)=(Y_{\rm c}\;Y_{\rm d})=\left(\begin{array}{cc}Y_{\rm rc}&Y_{\rm rd}\\
Y_{\rm nrc}&Y_{\rm nrd}\end{array}\right).$$
Then, it can be shown that
\begin{eqnarray}
\hat{W}_{(k)}\hat{\alpha}_{(k)}=P_{1^{\perp}}Y_{\rm c}(P_{\rm r}\;Q_{\rm r}),\label{eq:prps_Wahat}
\end{eqnarray}
where the $P_{\rm r}$ and $Q_{\rm r}$ have the same expressions as the $P$ and $Q$ in (\ref{eq:Wah_ap_Wa})
[defined below (\ref{eq:proj_decom}) and (\ref{eq:Q_decom}), respectively], with the $Y$ replaced by $Y_{\rm
r}=(Y_{\rm rc}\;Y_{\rm rd})$, and $U_{(k)}$ by $U_{{\rm r}(k)}$, that is,
$P_{\rm r}=Y_{\rm rc}'U_{{\rm r}(k)}\{U_{{\rm r}(k)}'Y_{\rm rc}Y_{\rm rc}'U_{{\rm r}(k)}\}^{-1}U_{{\rm r}(k)}'Y_{\rm
rc}$,
$Q_{\rm r}=Y_{\rm rc}'U_{{\rm r}(k)}\{U_{{\rm r}(k)}'Y_{\rm rc}Y_{\rm rc}'U_{{\rm r}(k)}\}^{-1}U_{{\rm r}(k)}'Y_{\rm
rd}$.

The proof essentially extends the previous results to PRPS. Two notable differences are that the $W$ and
$\epsilon$ are different from those in (\ref{eq:LMM}) in that $W, \epsilon$ now have some special forms [see
(\ref{eq:newMXW})]. For the $W$ this difference is relatively minor. For the $\epsilon$, note that, with the form
in (\ref{eq:newMXW}), the elements of $\epsilon$ are not necessarily independent. However, upon checking the
previous proofs, it turns out that what we need for $\epsilon$ is really an upper-bound order of $\|\epsilon\|$. Note
that
\begin{eqnarray}
\|\epsilon\|\leq\left\|\left(\begin{array}{c}A_{\rm a}\\
I_{m_{0}}\end{array}\right)\right\|\cdot\|\epsilon_{0}\|\leq(\|A_{\rm a}\|+1)\|\epsilon_{0}\|=(\|A_{\rm a}\|+1)\sqrt{n}
O_{\rm P}(1).\label{eq:ub_eps}
\end{eqnarray}
Note that the elements of $A_{\rm a}$ can be expressed by
$$a_{gi}=I_{gi}/n_{g},\; 1\leq g\leq m_{\rm pa}, 1\leq i\leq
m_{0},$$
where $I_{gi}=1$ if assay $i$ is involved in average $g$, and $I_{gi}=0$ otherwise, and $n_{g}>1$ is the number of original assays involved in average $g$. Then, for any vector $v=(v_{i})_{1\leq i\leq m_{0}}\in R^{m_{0}}$, we have, by assumptions (a), (b) for $A_{\rm a}$,
\begin{eqnarray*}
&&|A_{\rm a}v|^{2}=\sum_{g=1}^{m_{\rm pa}}\left(\sum_{i=1}^{m_{0}}\frac{I_{gi}}{n_{g}}v_{i}\right)^{2}\leq
\sum_{g=1}^{m_{\rm pa}}\left(\sum_{i=1}^{m_{0}}I_{gi}\right)\left(\sum_{i=1}^{m_{0}}\frac{I_{gi}}{n_{g}^{2}}v_{i}^{2}
\right)\\
&&\leq b_{1}\sum_{g=1}^{m_{\rm pa}}\sum_{i=1}^{m_{0}}\frac{I_{gi}}{n_{g}^{2}}v_{i}^{2}
=b_{1}\sum_{i=1}^{m_{0}}\left(\sum_{g=1}^{m_{\rm pa}}\frac{I_{gi}}{n_{g}^{2}}\right)v_{i}^{2}\\
&&\leq\frac{b_{1}b_{2}}{4}\sum_{i=1}^{m_{0}}v_{i}^{2}=\left(\frac{b_{1}b_{2}}{4}\right)|v|^{2}.
\end{eqnarray*}
Thus, we have come up with the following inequality
\begin{eqnarray}
\|A_{\rm a}\|=\sqrt{\lambda_{\max}(A_{\rm a}'A_{\rm a})}=\sqrt{\|A_{\rm a}'A_{\rm a}\|}\leq\frac{\sqrt{b_{1}
b_{2}}}{2}.\label{eq:ub_Aa}
\end{eqnarray}
Combining (\ref{eq:ub_Aa}) with (\ref{eq:ub_eps}), it follows that $\|\epsilon\|=\sqrt{n}O_{\rm P}(1)$; in other
words, the order of the upper-bound of $\|\epsilon\|$ has not changed, and this is all we need technically from $\epsilon$.

The rest of the proof can be completed by checking step by step the previous proofs in making sure that the conclusions of the current theorem hold. The details are omitted.

\bibliographystyle{apalike}
\bibliography{references}

\section{Figures}
\FloatBarrier

\FloatBarrier
\newpage

\section{Supplementary Material}

\setcounter{figure}{0}
\renewcommand{\thefigure}{S\arabic{figure}}

\FloatBarrier

\begin{figure}
    \centering
    \begin{tabular}{ccc}
    & $K=1$ & $K=2$ \\
    \rotatebox{90}{\hspace{0.4cm}Number of samples increasing} &
    \includegraphics[scale=0.25]{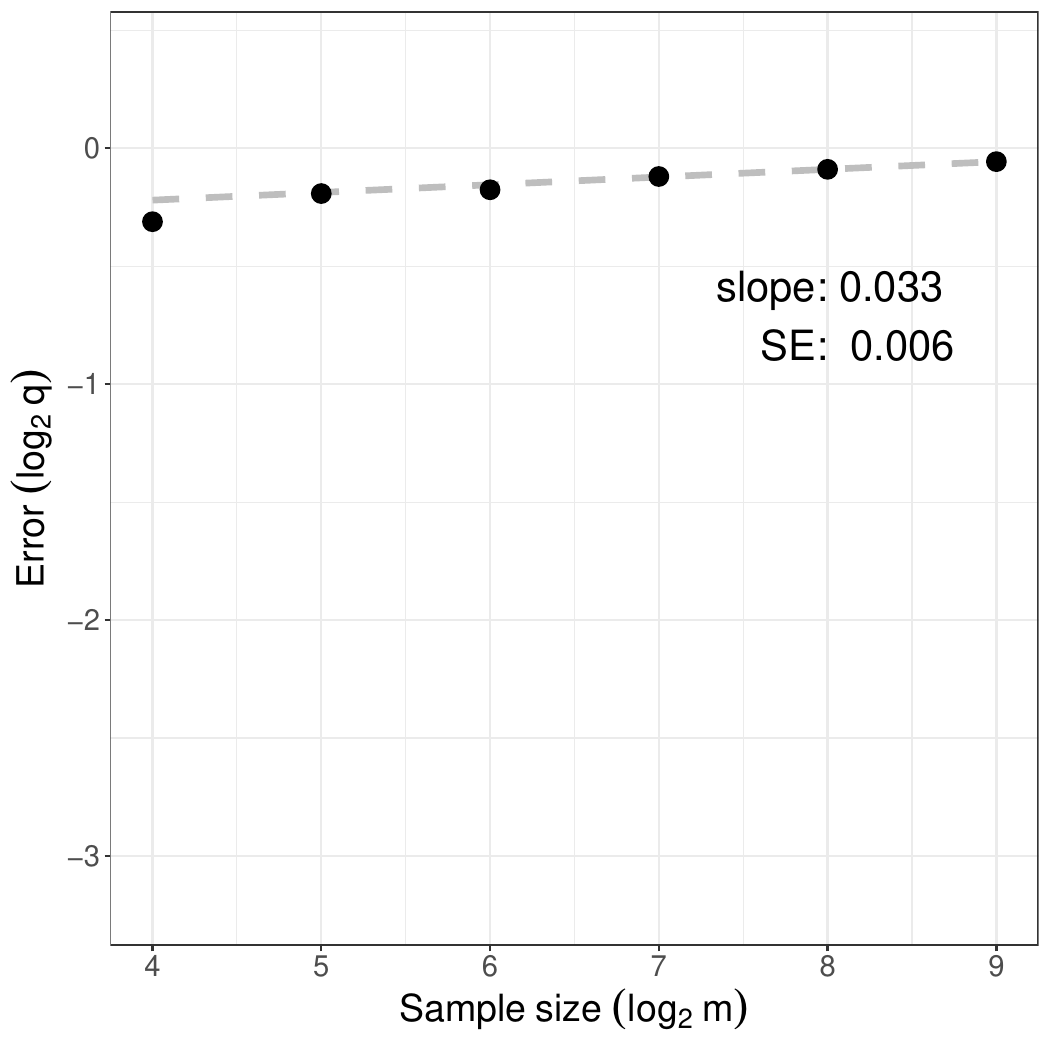} &
    \includegraphics[scale=0.25]{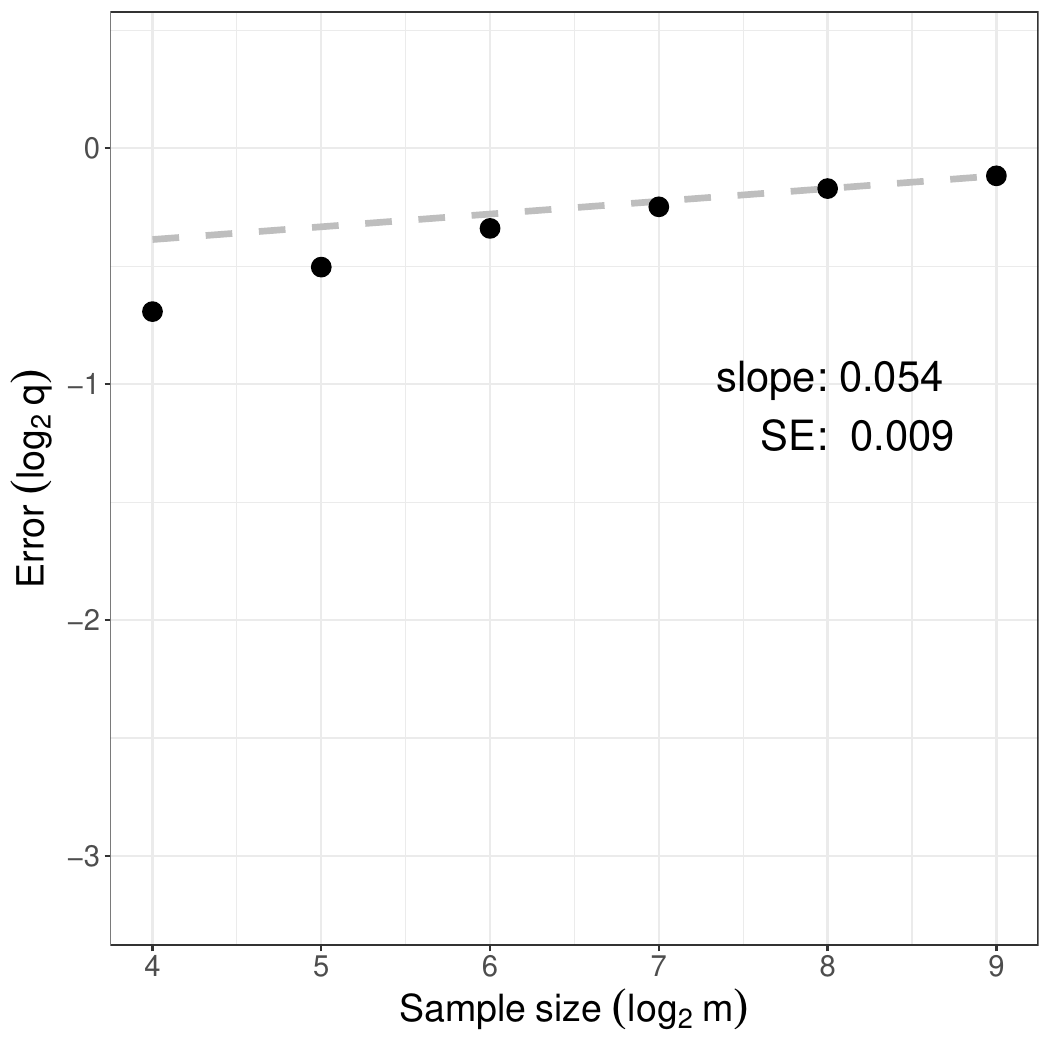} \\    
    \rotatebox{90}{\hspace{0.4cm}Number of replicates increasing} &
    \includegraphics[scale=0.25]{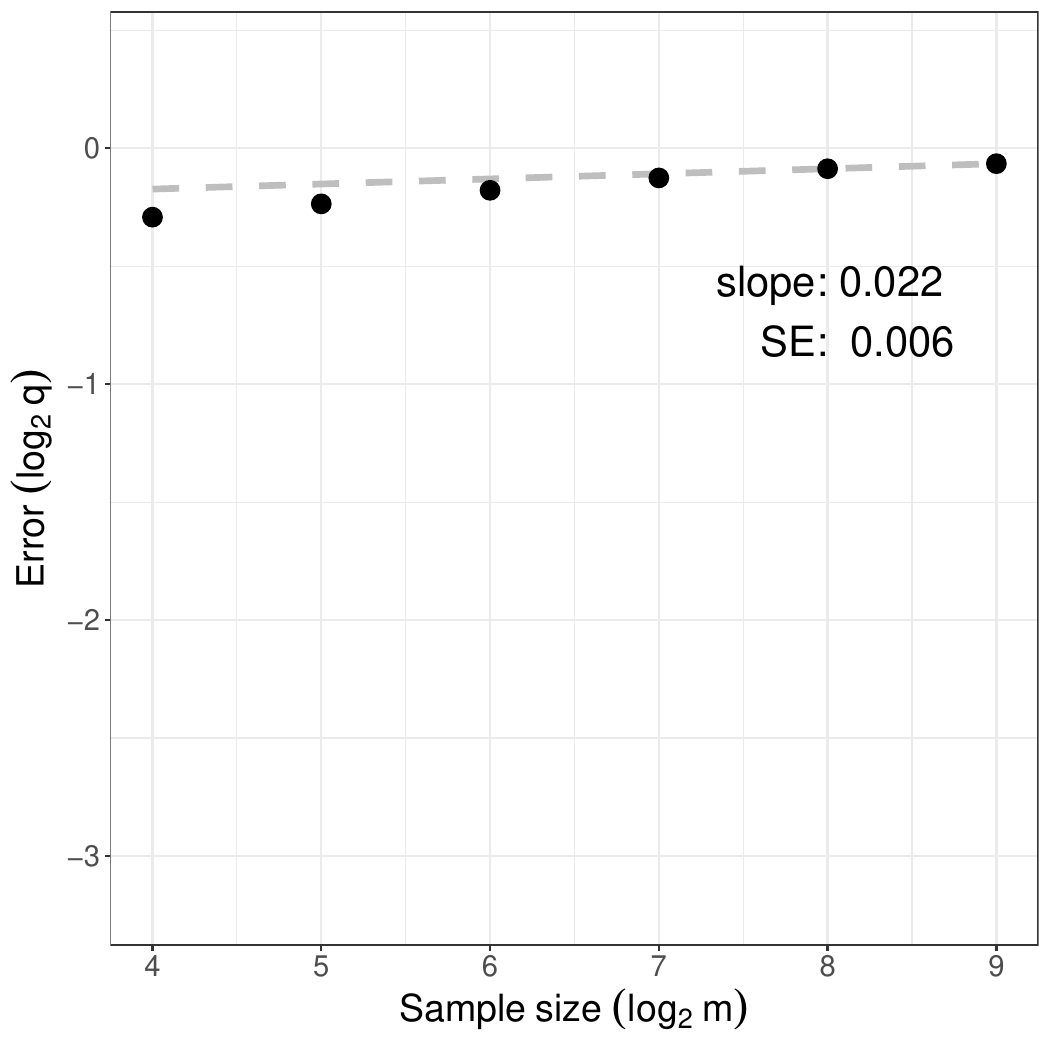} &
    \includegraphics[scale=0.25]{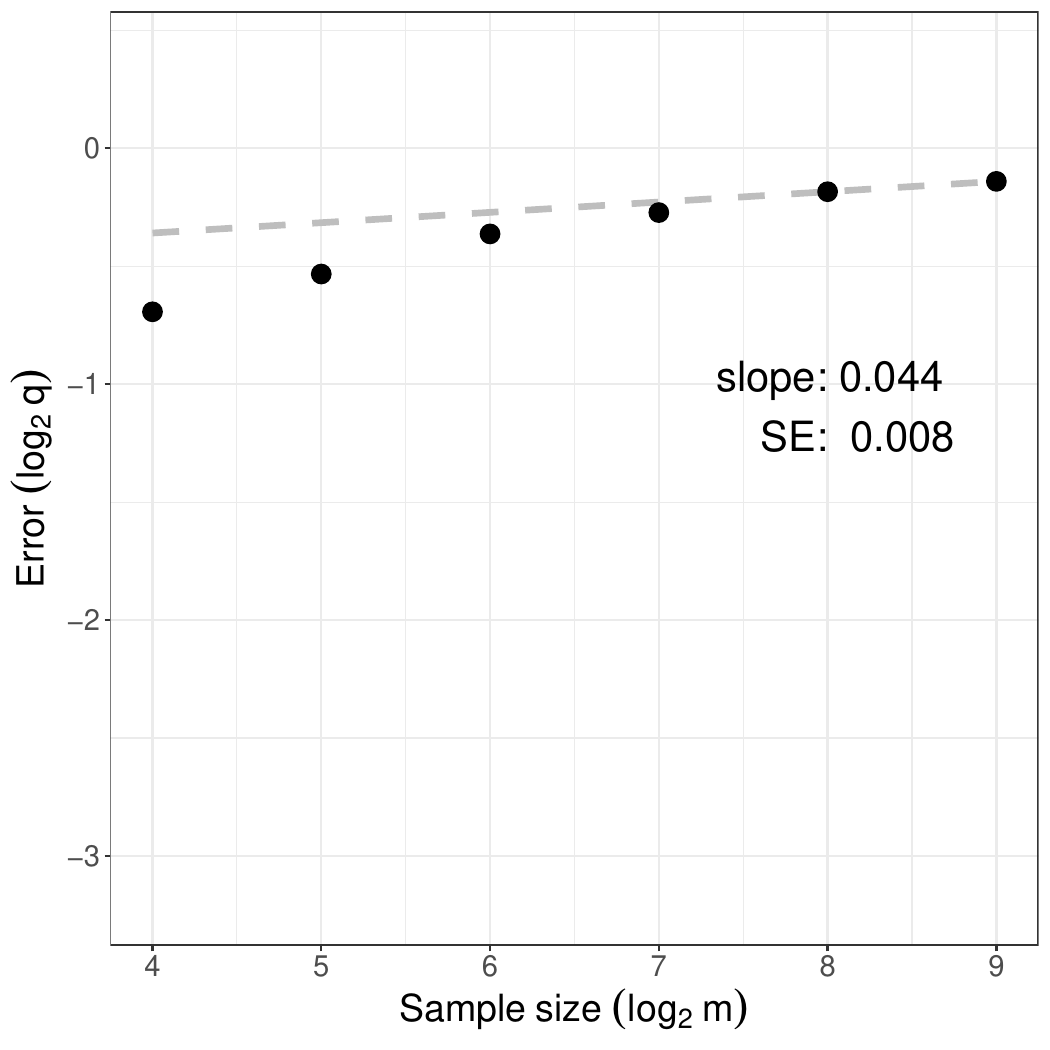}     
    \end{tabular}
    \caption{Like Figure 1, but for $K=1$ and $K=2$.  Note that the vertical scale differs from Figure 1.}
    \label{fig:sim_k12}
\end{figure}

\begin{figure}[ht]
    \centering
    \begin{tabular}{cccc}
    & $K=3$ & $K=10$ & $K=K_{\mathrm{max}}$ \\
    \rotatebox{90}{\hspace{0.4cm}Number of samples increasing} &
    \includegraphics[scale=0.25]{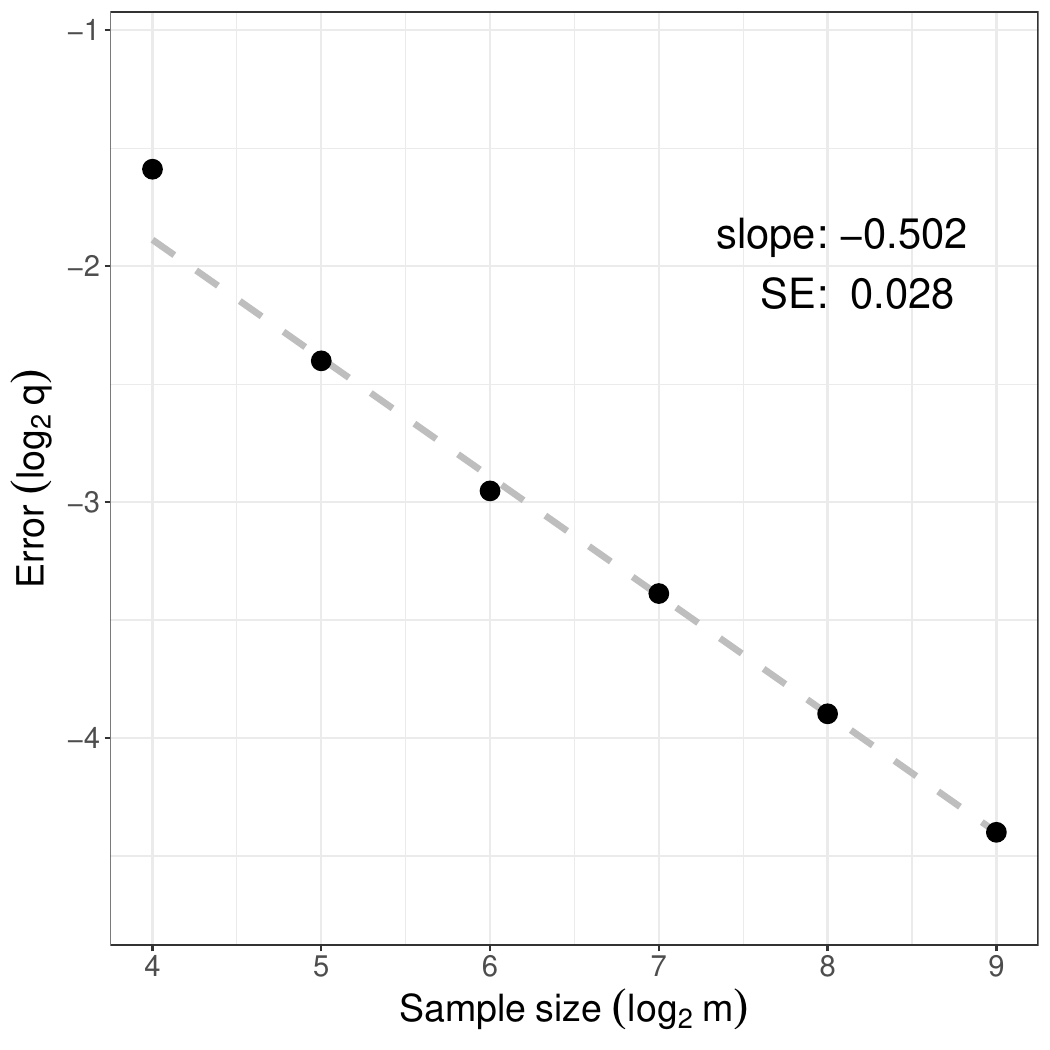} &
    \includegraphics[scale=0.25]{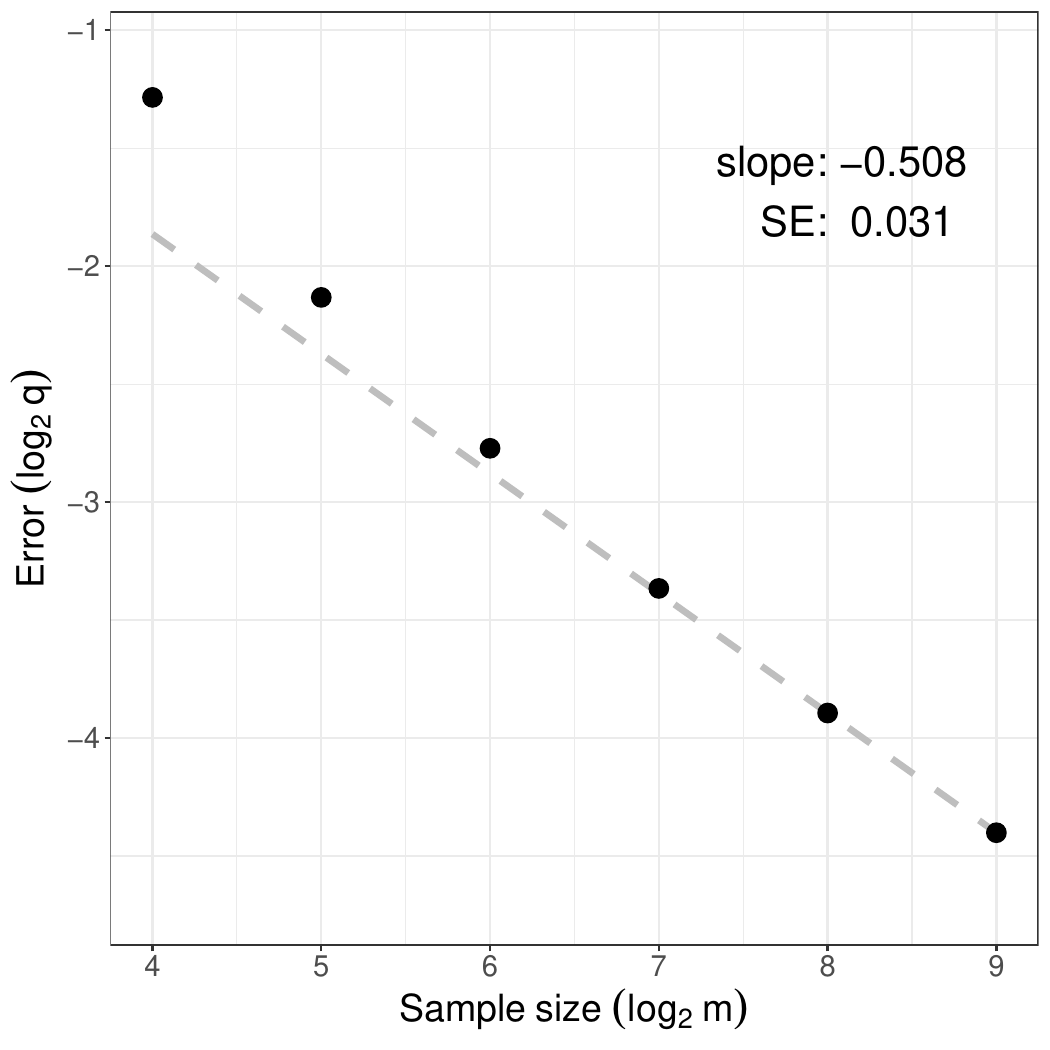} &
    \includegraphics[scale=0.25]{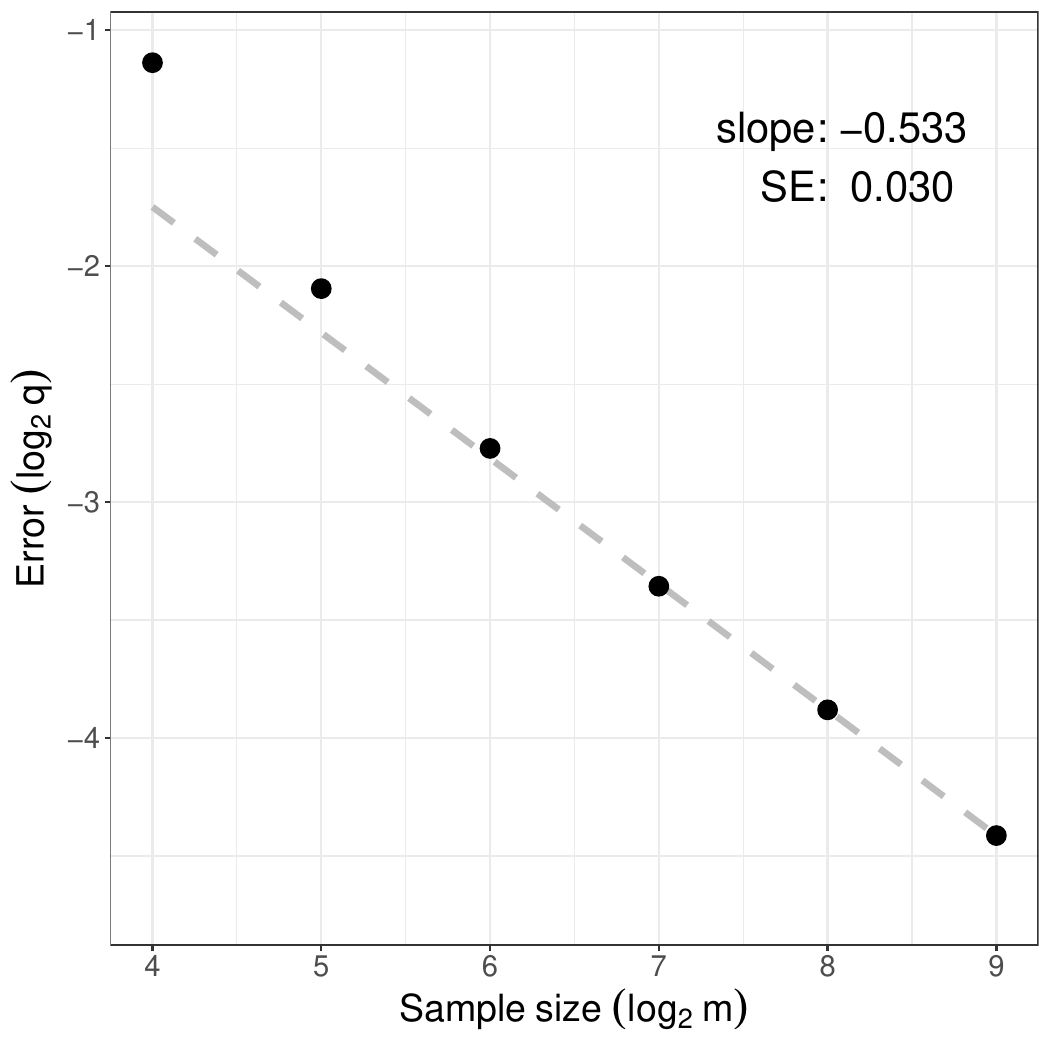} \\    
    \rotatebox{90}{\hspace{0.4cm}Number of replicates increasing} &
    \includegraphics[scale=0.25]{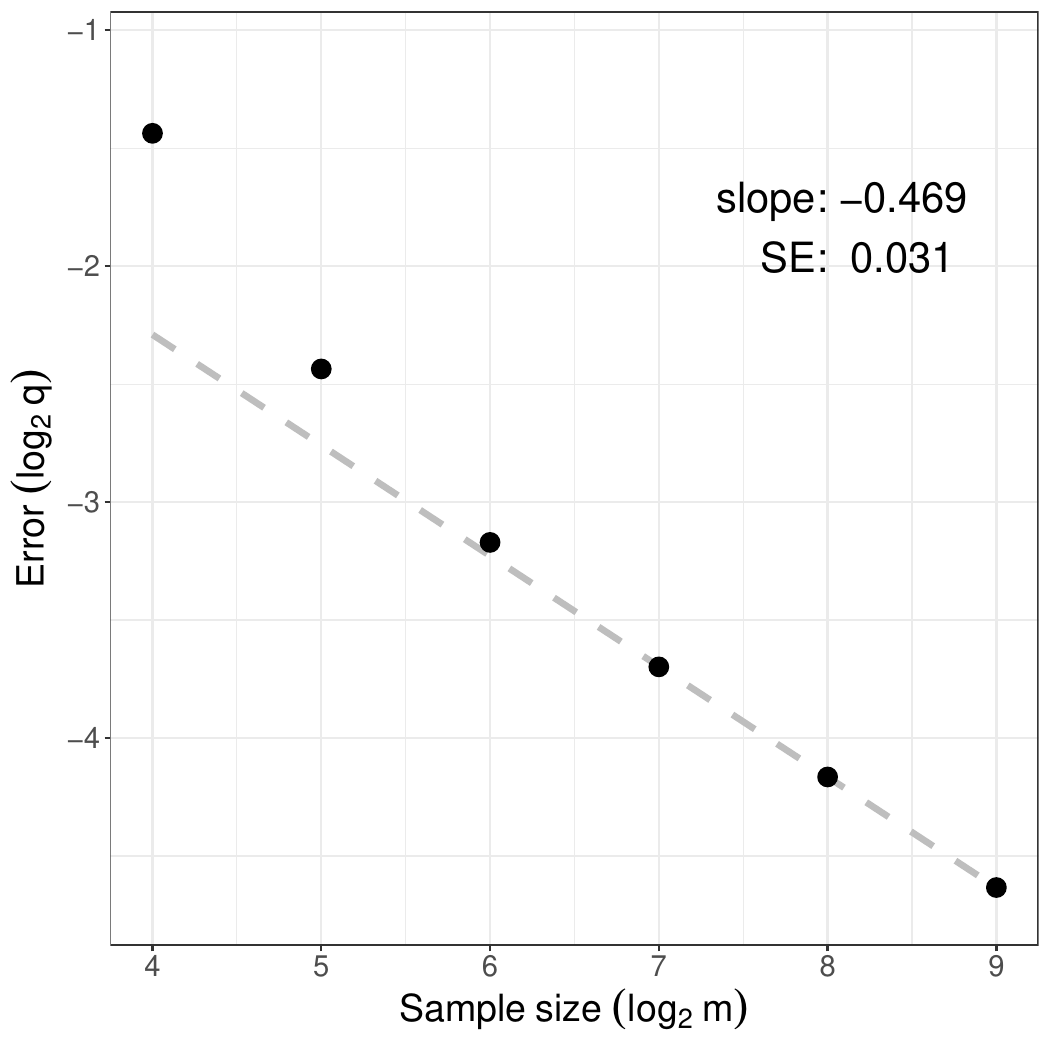} &
    \includegraphics[scale=0.25]{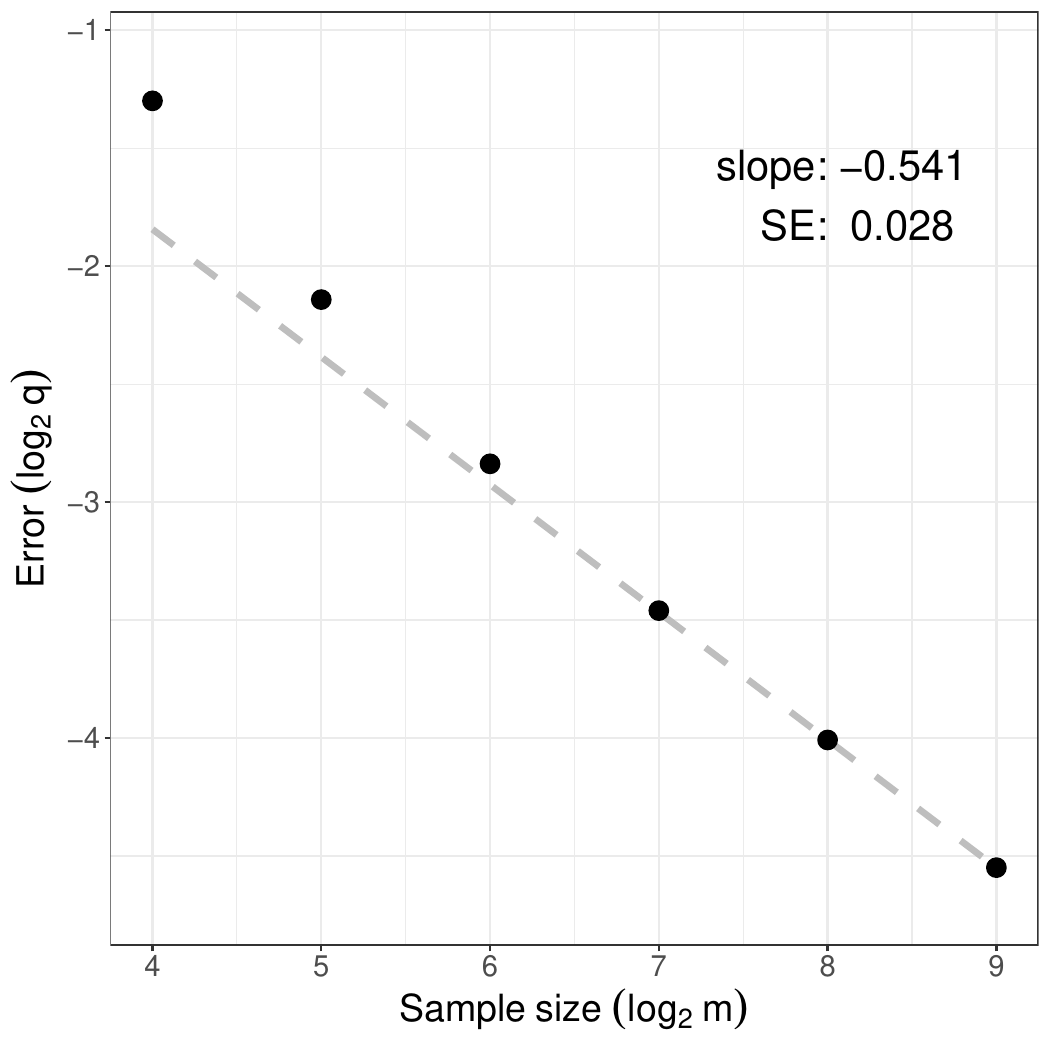} &
    \includegraphics[scale=0.25]{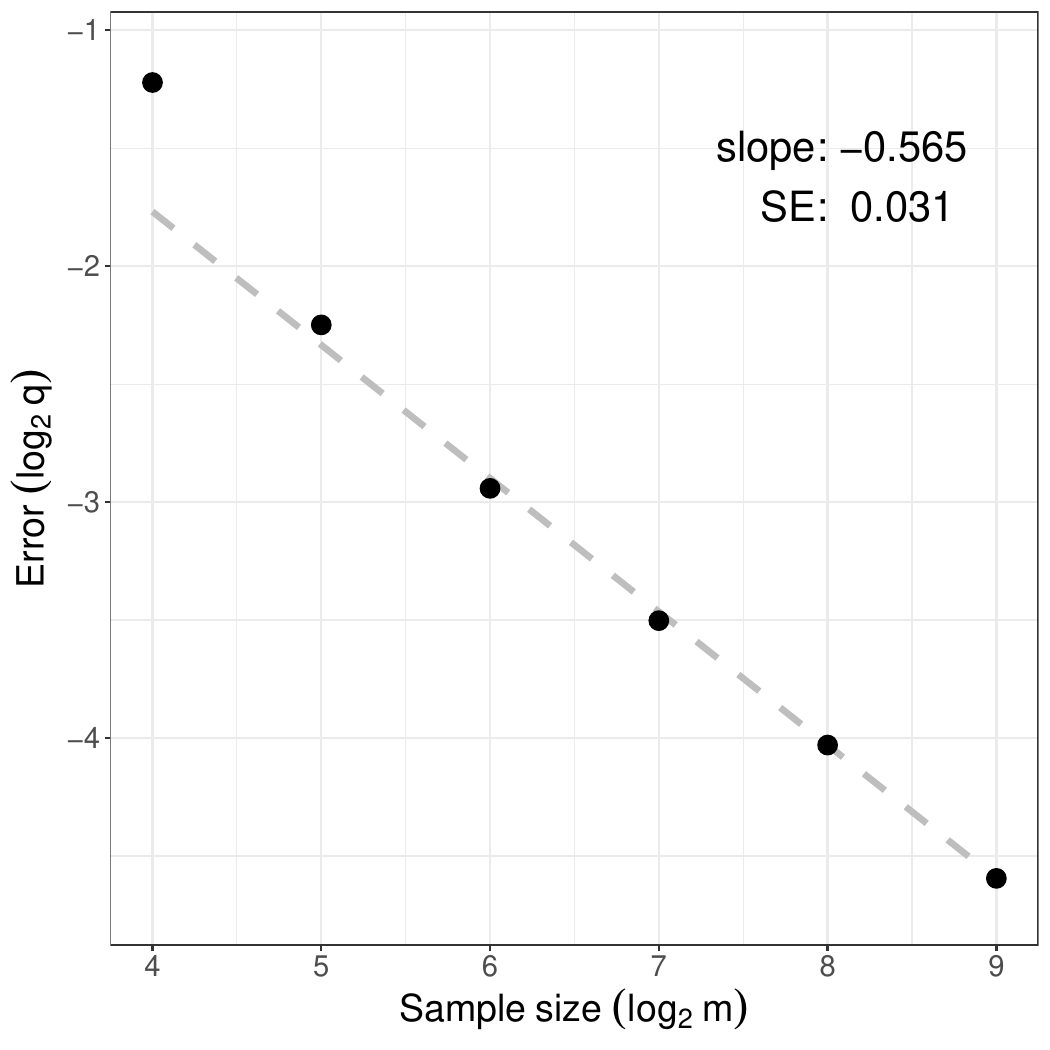}     
    \end{tabular}
    \caption{Like Figure 1, but with Pareto distributed data.}
    \label{fig:simbasicpareto}
\end{figure}

\begin{figure}[ht]
    \centering
    \begin{tabular}{ccc}
    $K=3$ & $K=10$ & $K=K_{\mathrm{max}}$ \\
    \includegraphics[scale=0.25]{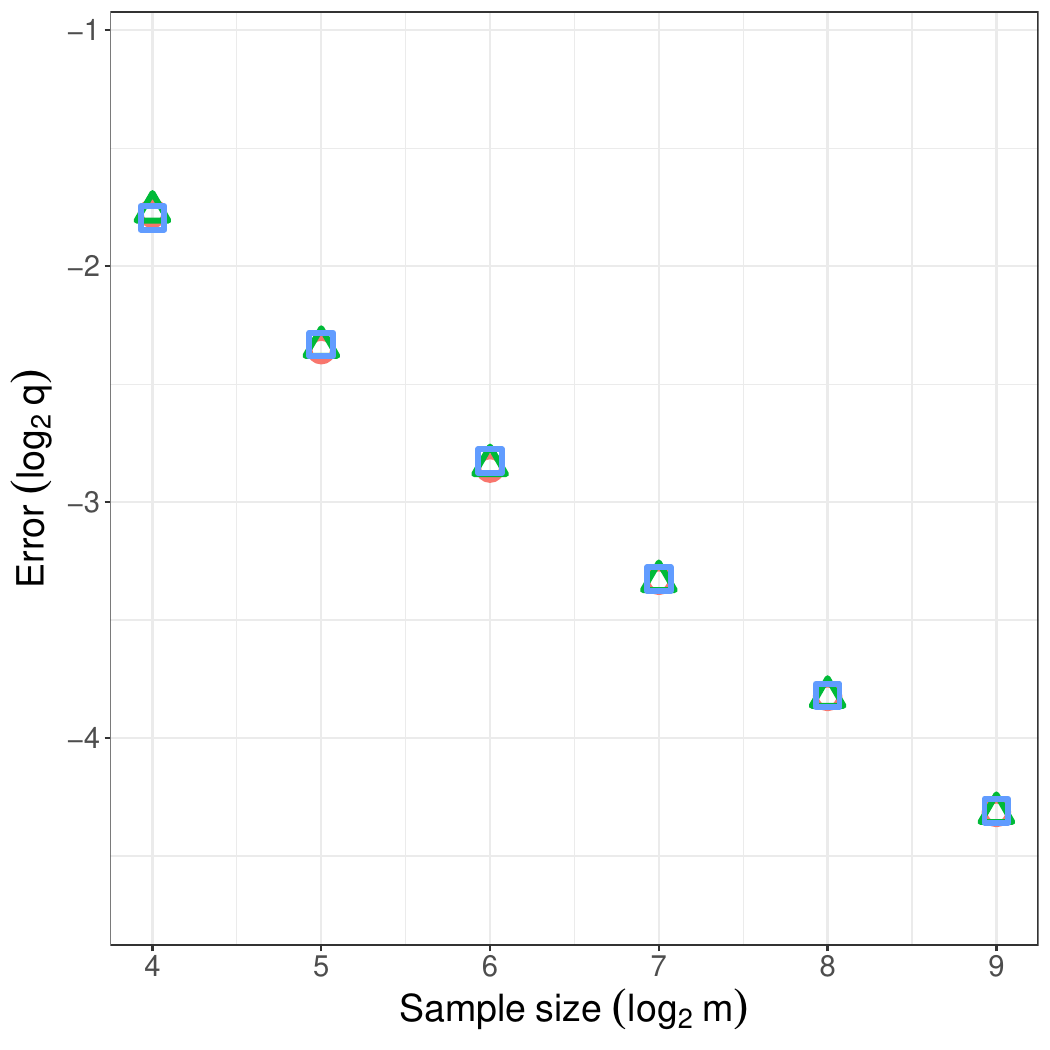} &
    \includegraphics[scale=0.25]{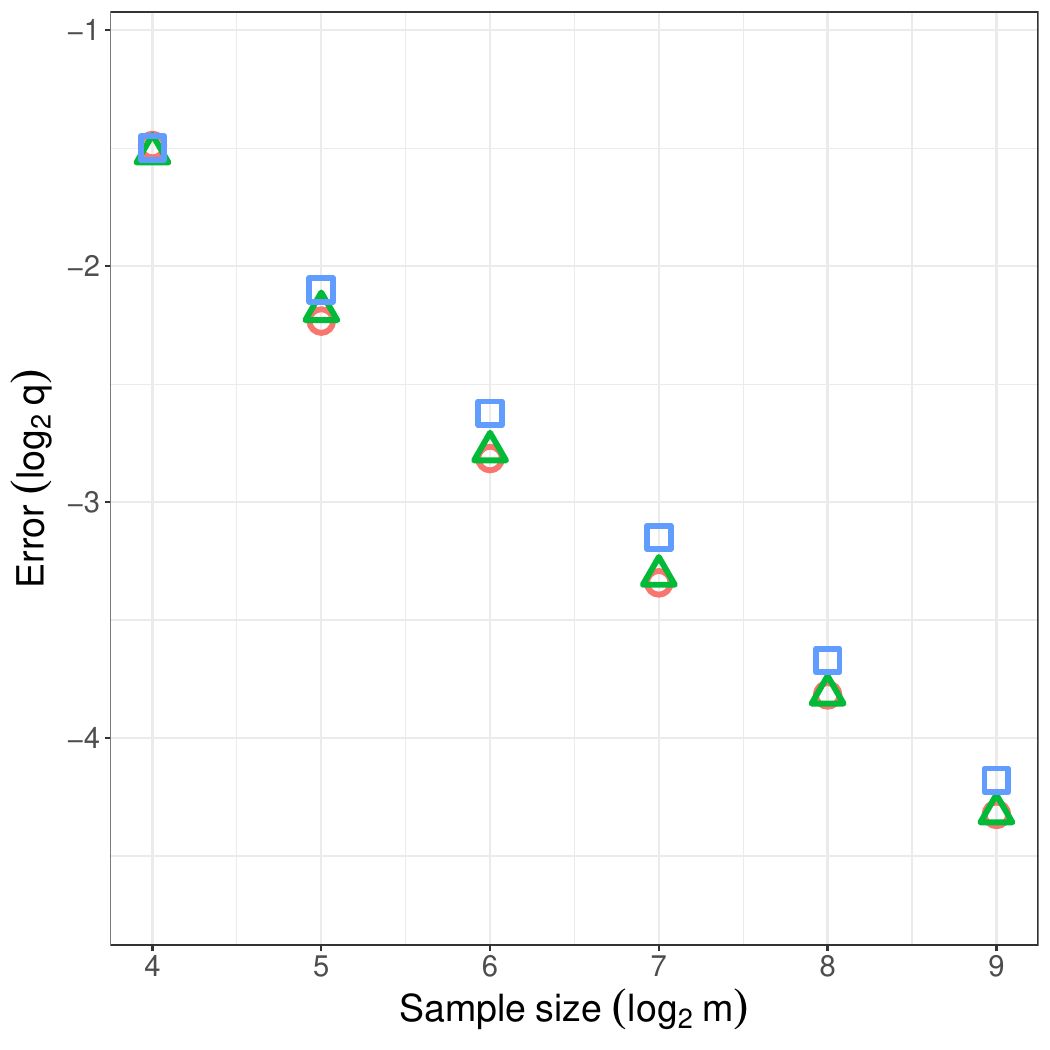} &
    \includegraphics[scale=0.25]{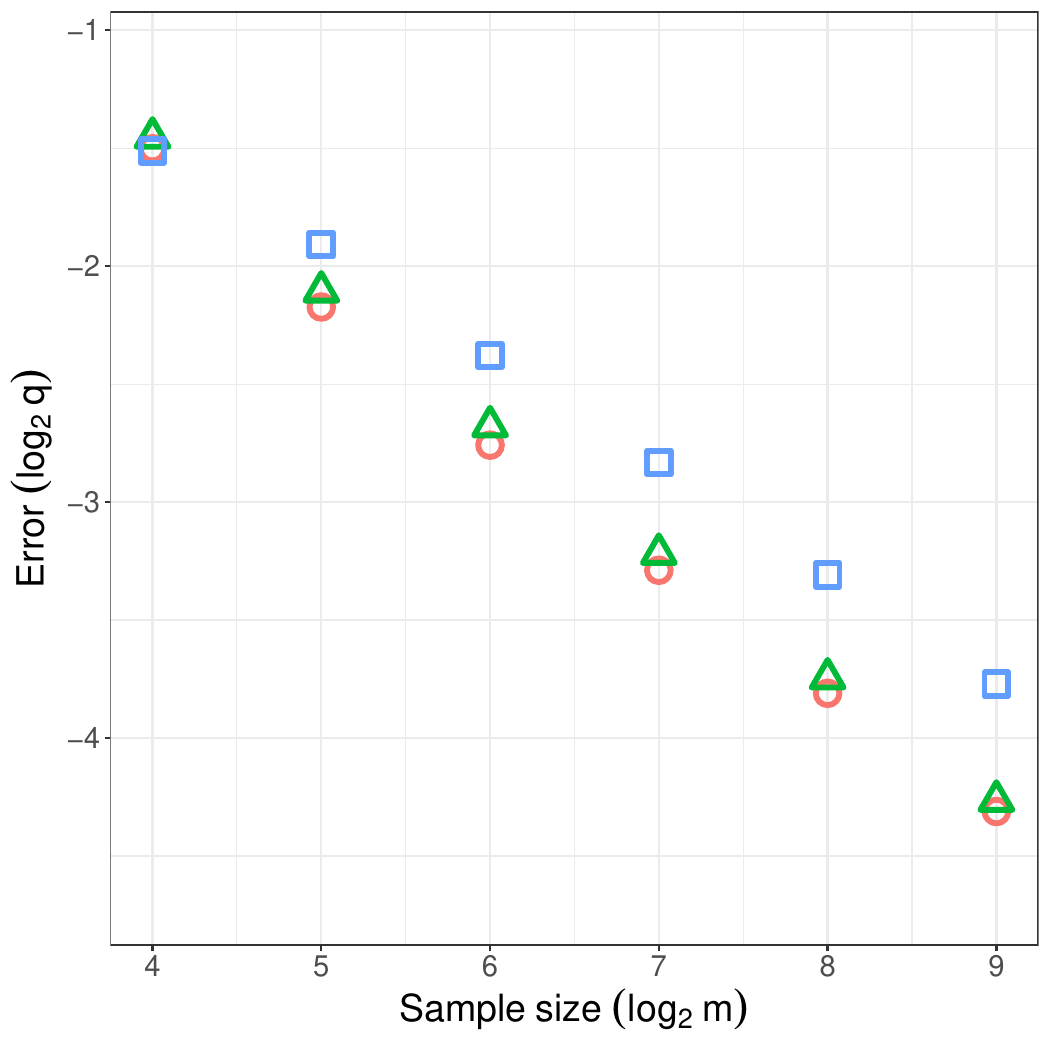} 
    \end{tabular}
    \caption{Like Figure 2, but with Pareto distributed data.}
    \label{fig:simcomparison1pareto}
\end{figure}

\begin{figure}[ht]
    \centering
    \begin{tabular}{ccc}
    $n_c = \frac{1}{8}m^2$ & $n_c = \frac{1}{2}m^\frac{3}{2}$ & $n_c = 2m$ \\
    \includegraphics[scale=0.25]{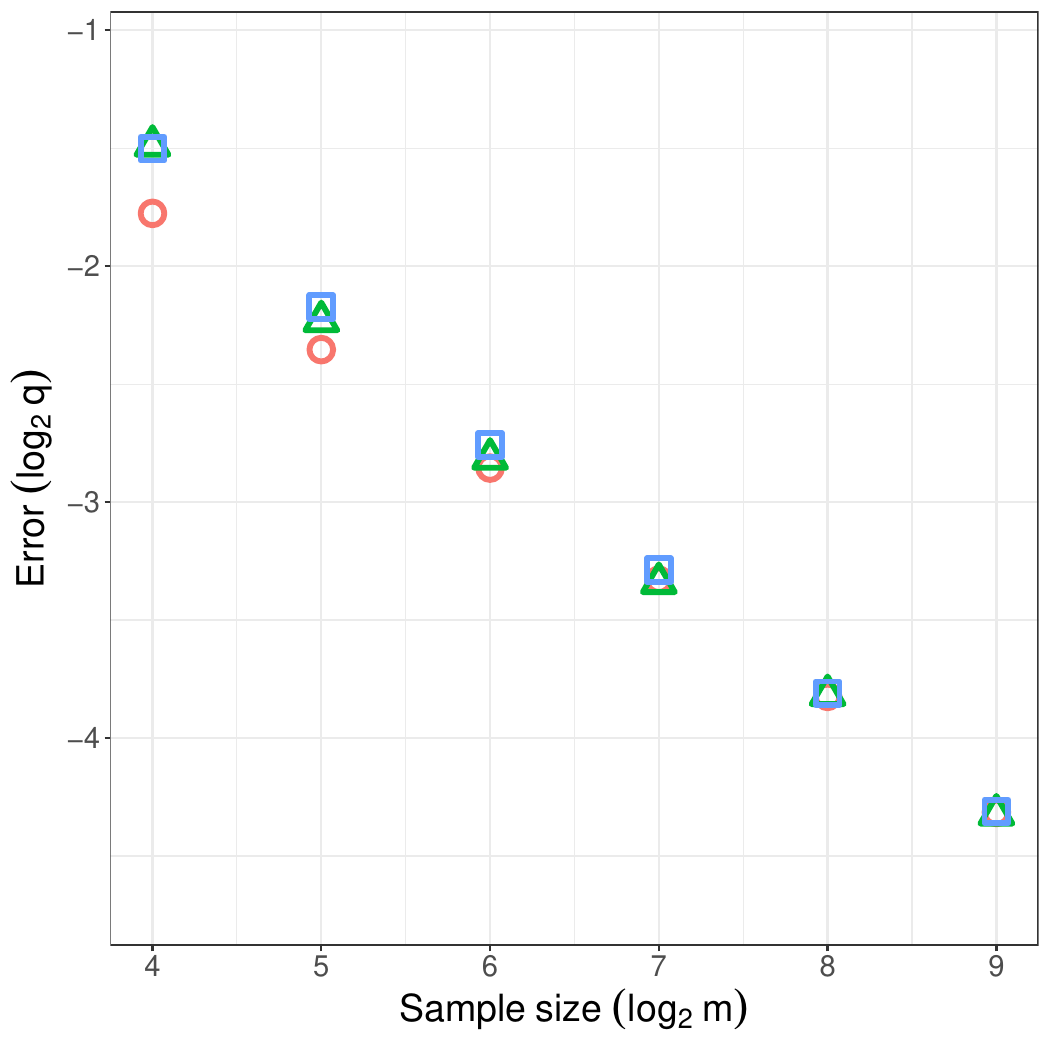} &
    \includegraphics[scale=0.25]{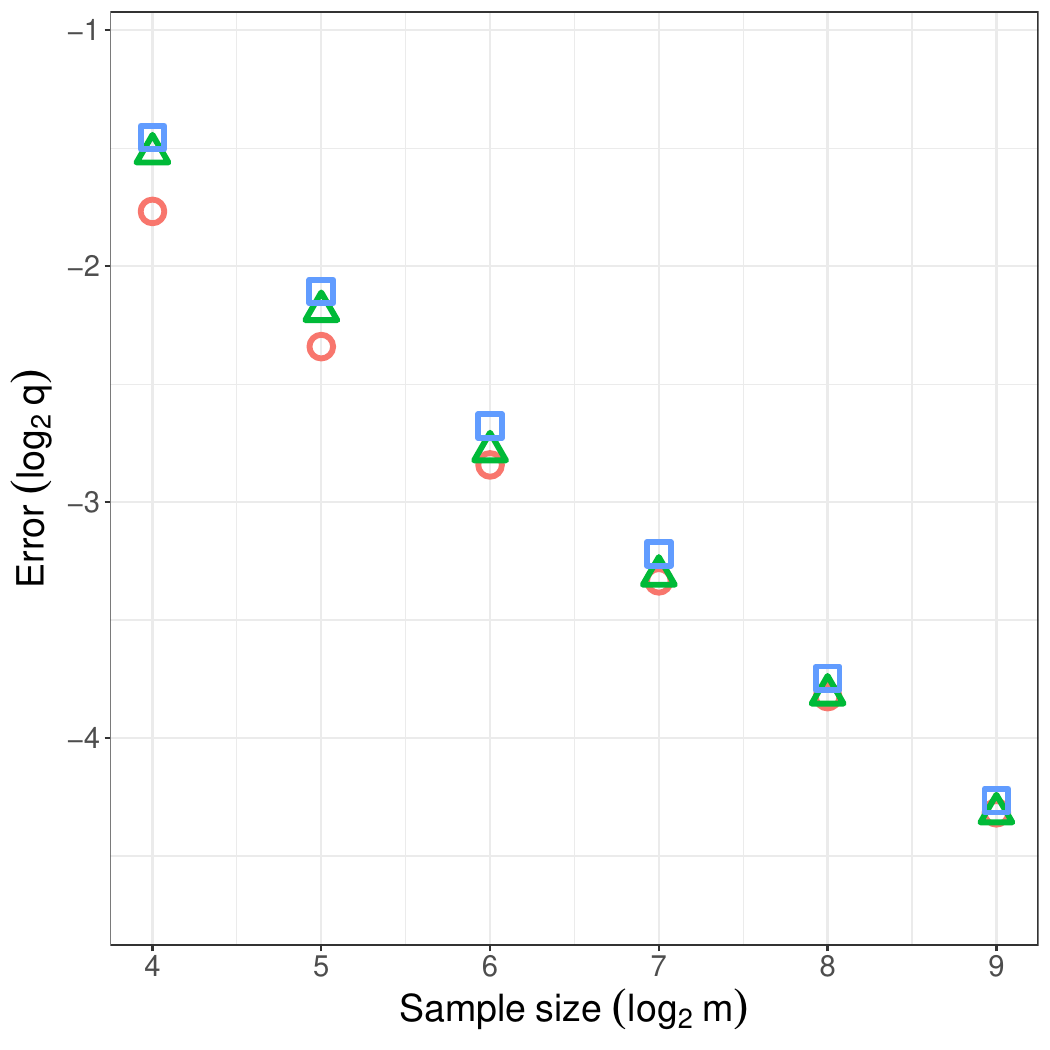} &
    \includegraphics[scale=0.25]{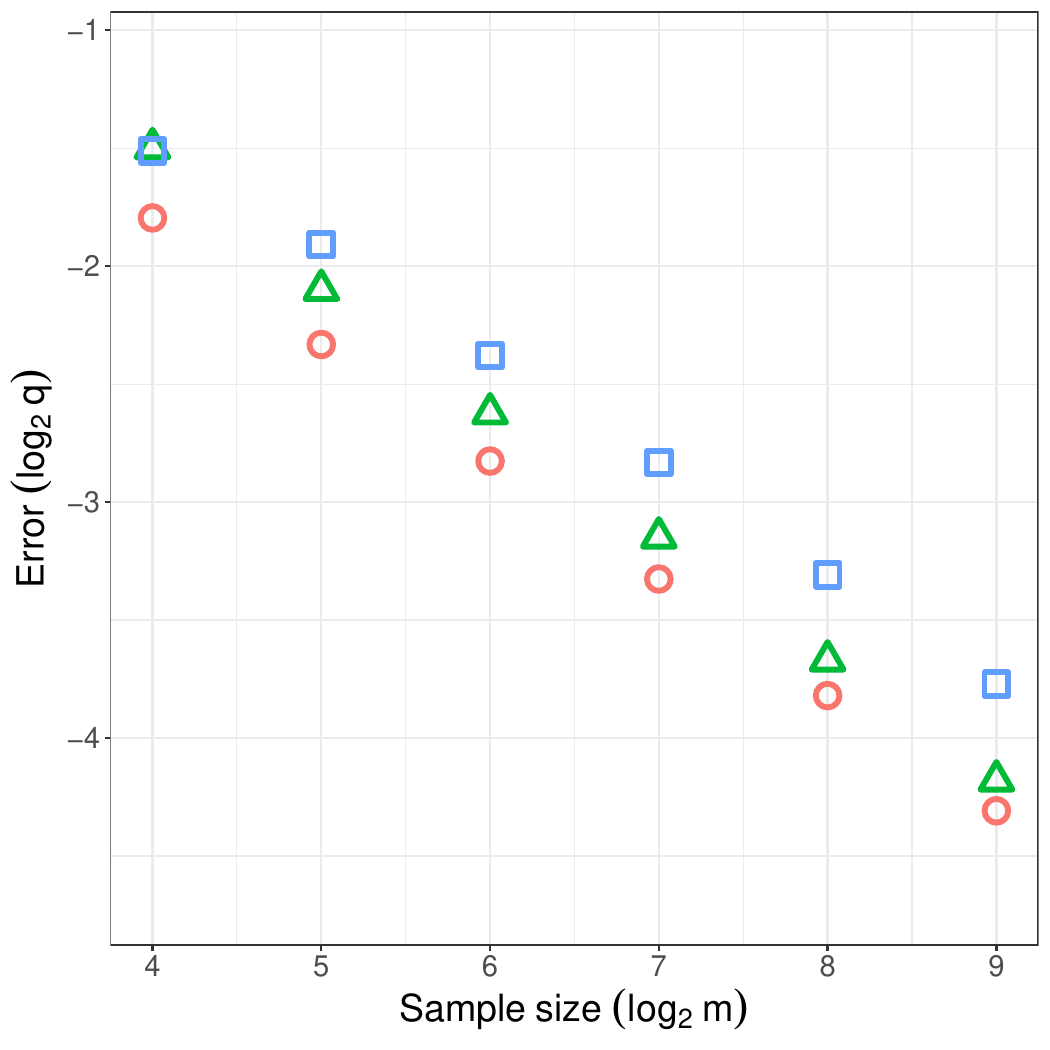} 
    \end{tabular}
    \caption{Like Figure 3, but with Pareto distributed data.}
    \label{fig:simcomparison2pareto}
\end{figure}

\begin{figure}[ht]
    \centering
        \includegraphics[scale=0.25]{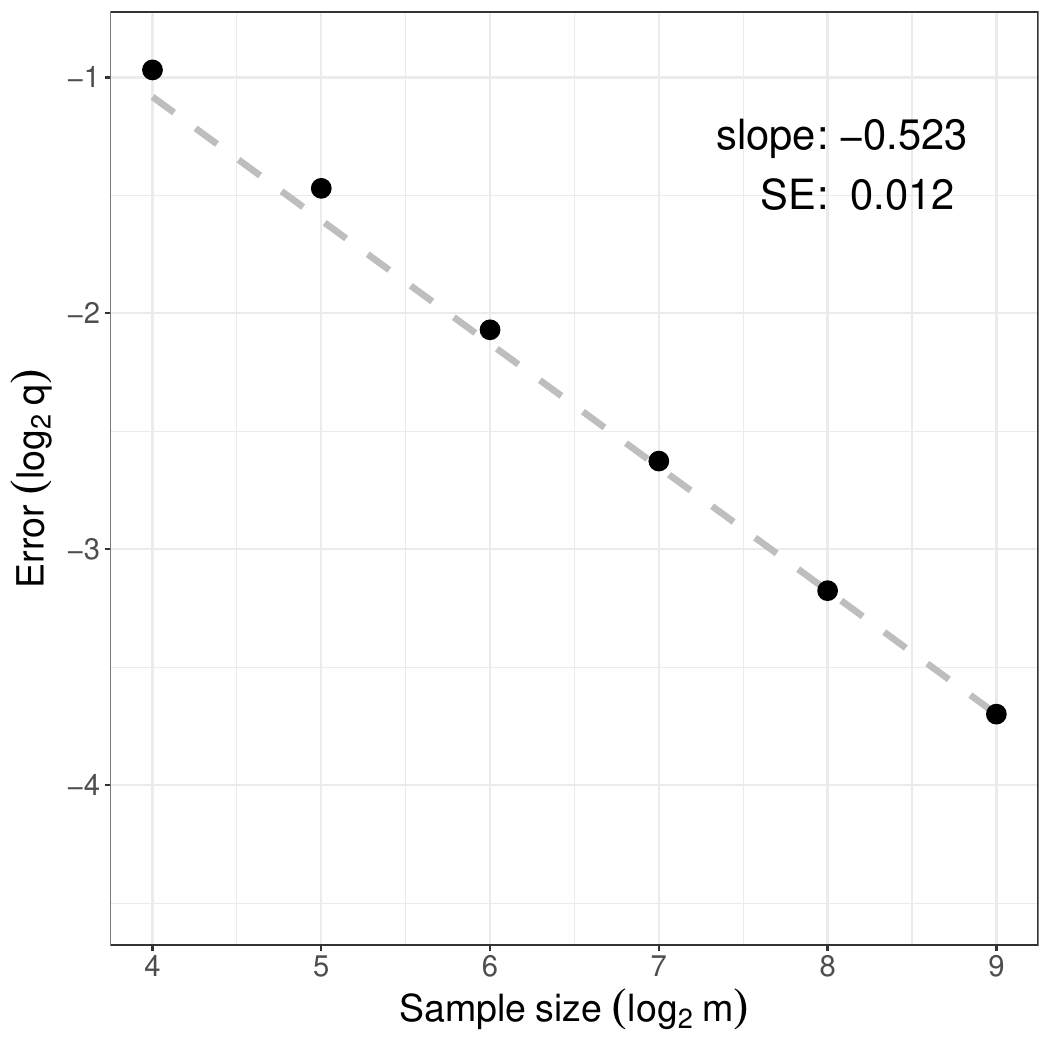} 
    \caption{Results for PRPS.  Note $n_c = \frac{1}{8m^2}$ as in Figure 1.  Note also that the vertical scale differs from Figure 1.  $K=3$.}
    \label{fig:sim_prps}
\end{figure}

\begin{figure}
    \centering
    \begin{tabular}{cc}
    \hspace{-1cm} No Adjustment & 
    \hspace{-1cm} ``Technical'' Adjustment \\
    \includegraphics[scale=0.3]{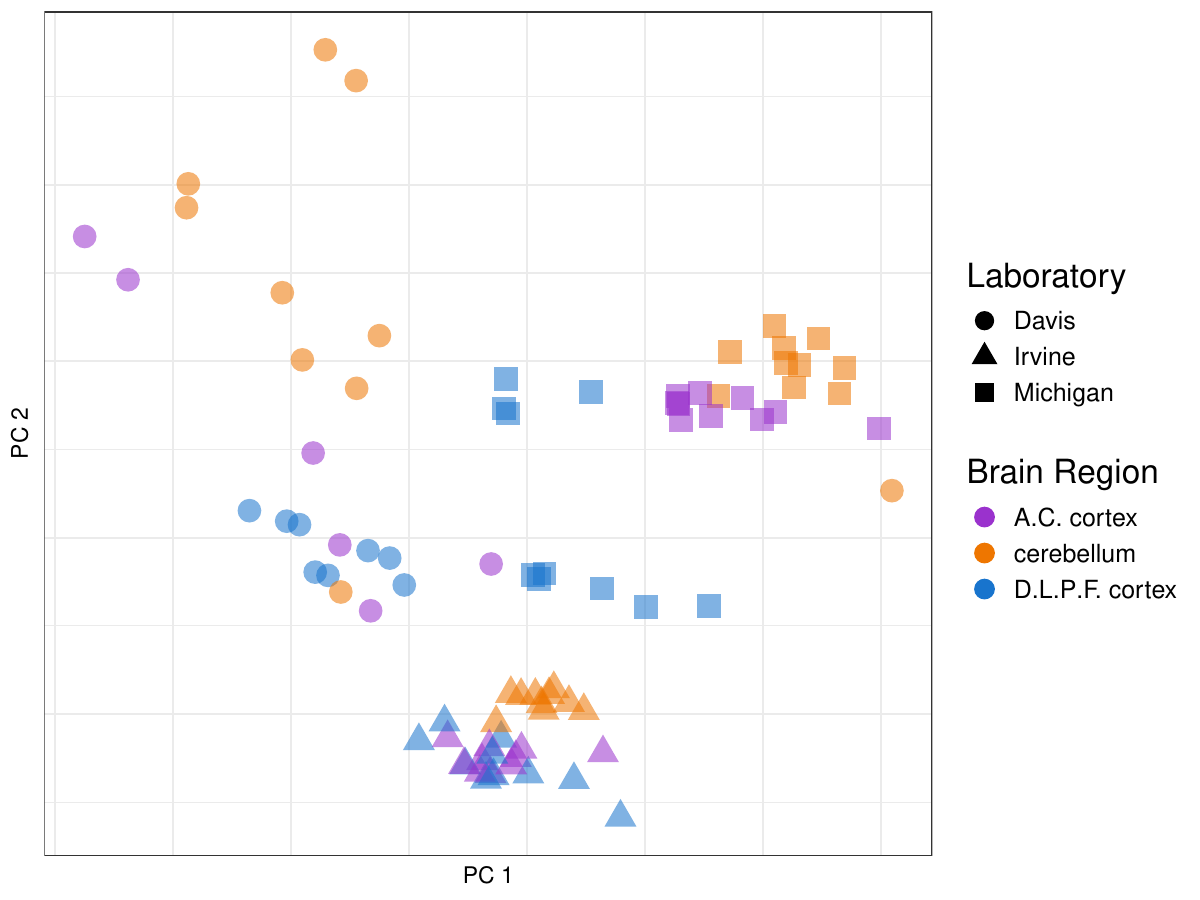} &
    \includegraphics[scale=0.3]{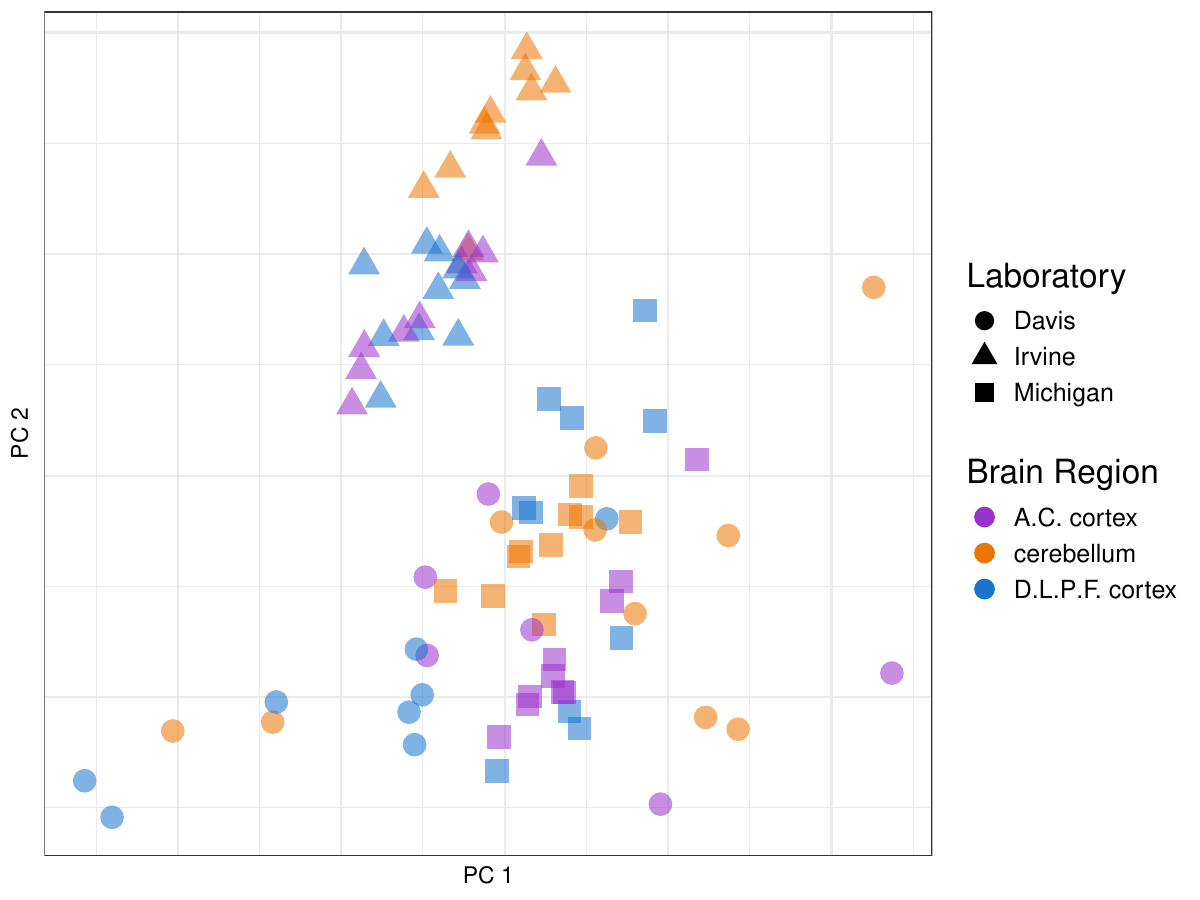} \\    
    \hspace{-1cm} ``Bio'' Adjustment & 
    \hspace{-1cm} ``Bio'' Adjustment (colored by patient) \\
    \includegraphics[scale=0.3]{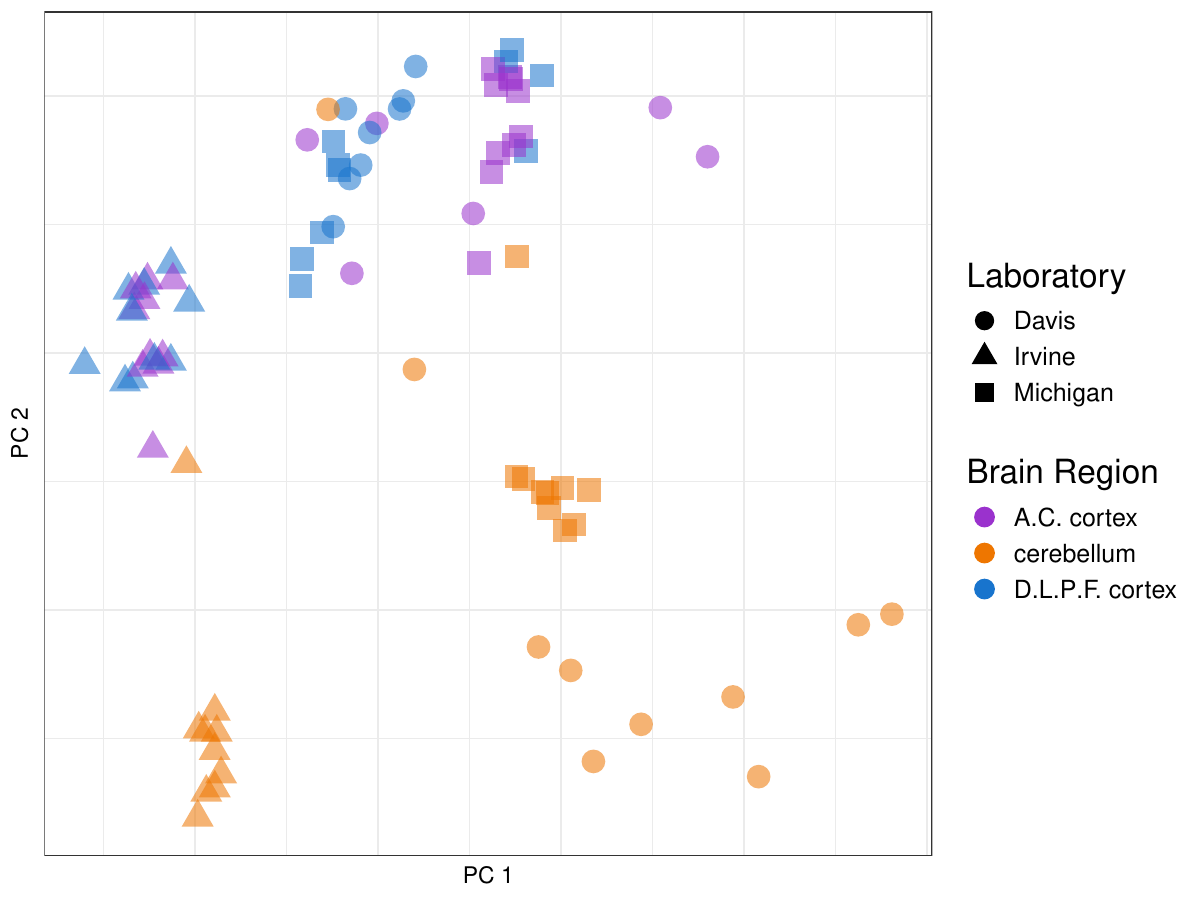}  &
    \includegraphics[scale=0.3]{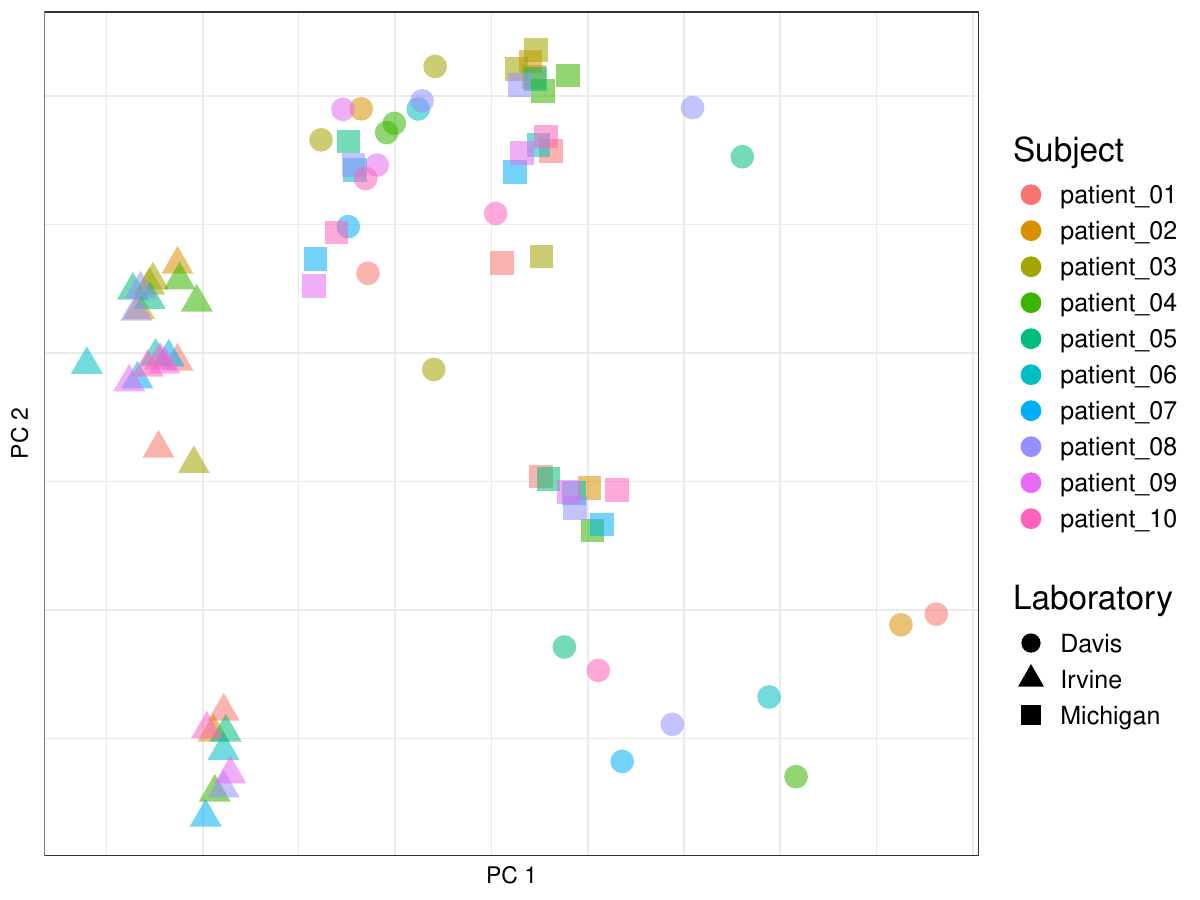}     
    \end{tabular}
    \caption{Like Figure 4, but with $K=1$.}
    \label{fig:K1_svd.all}
\end{figure}

\begin{figure}
    \centering
    \begin{tabular}{cc}
    \hspace{-1cm} No Adjustment & 
    \hspace{-1cm} ``Technical'' Adjustment \\
    \includegraphics[scale=0.3]{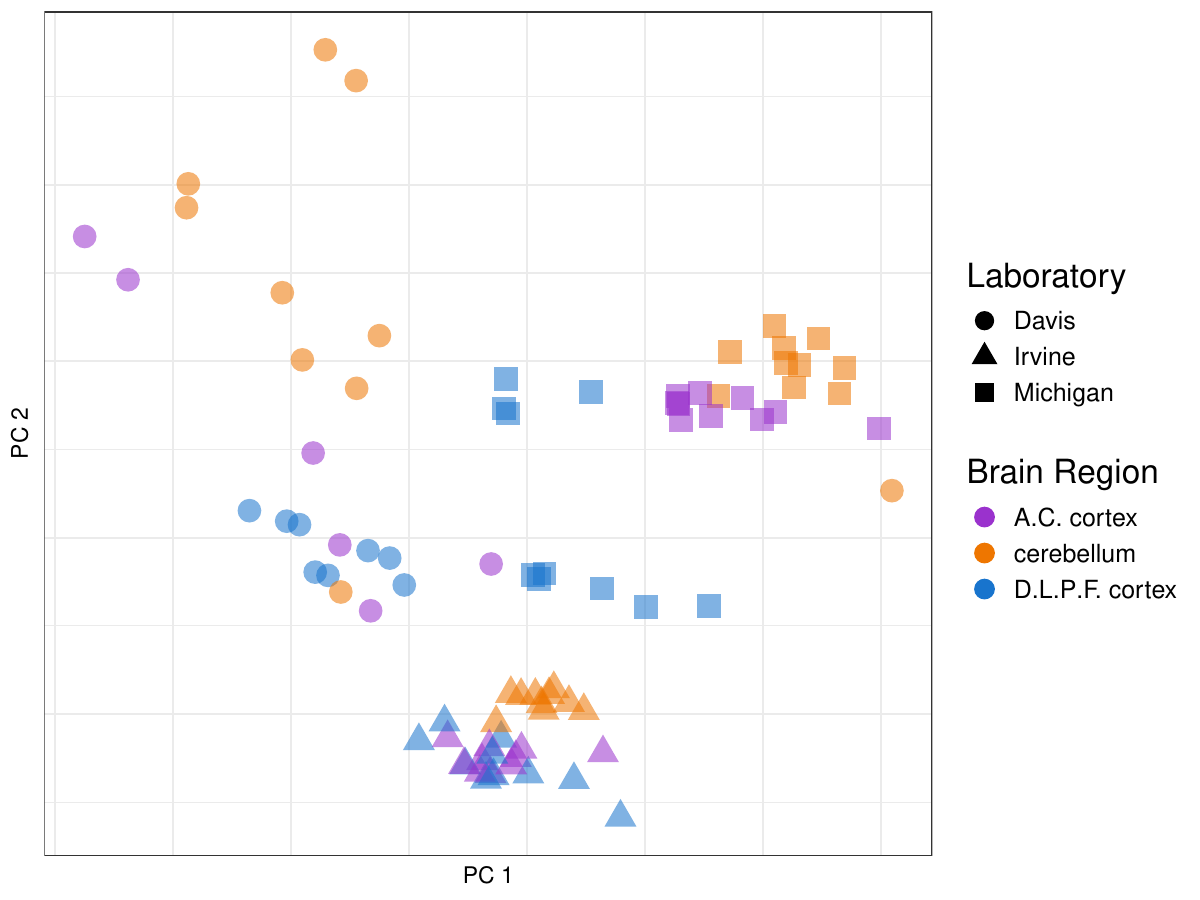} &
    \includegraphics[scale=0.3]{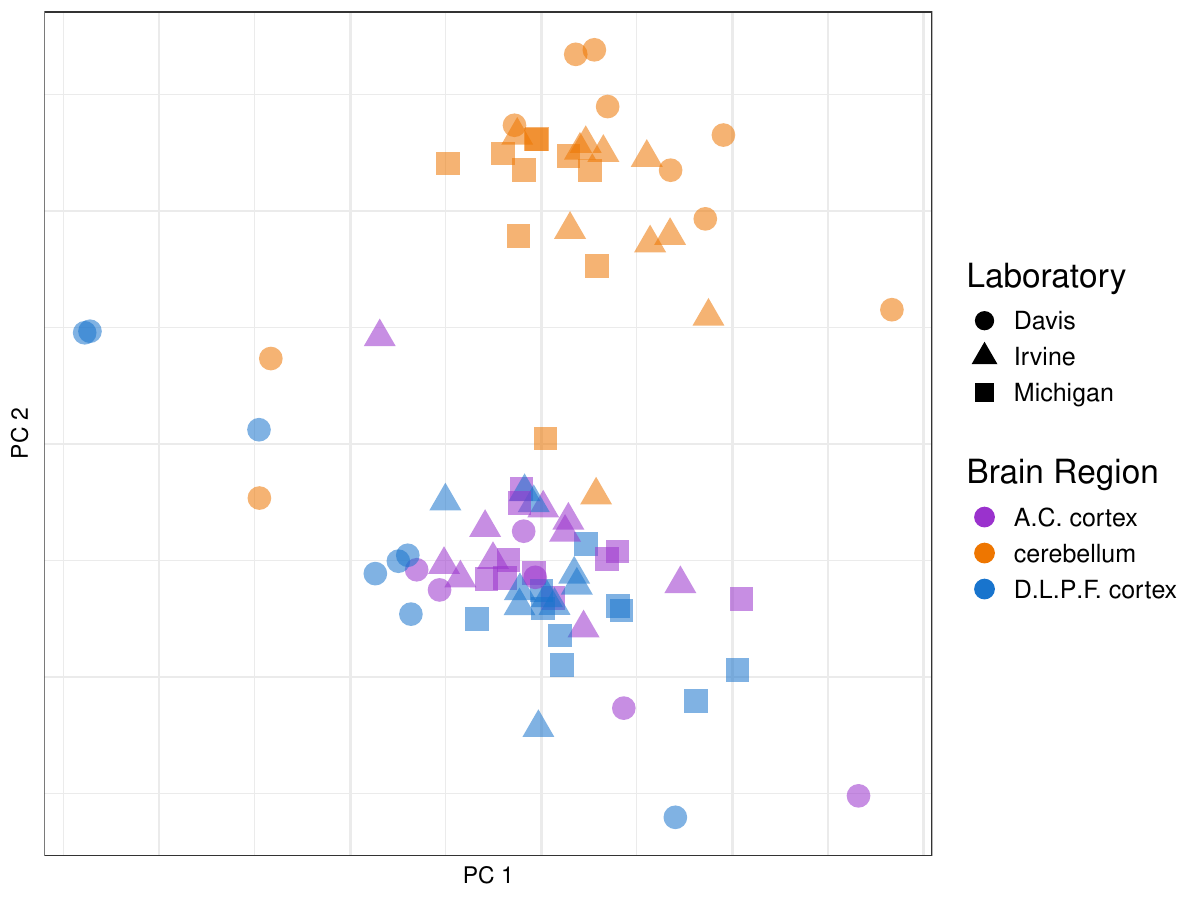} \\    
    \hspace{-1cm} ``Bio'' Adjustment & 
    \hspace{-1cm} ``Bio'' Adjustment (colored by patient) \\
    \includegraphics[scale=0.3]{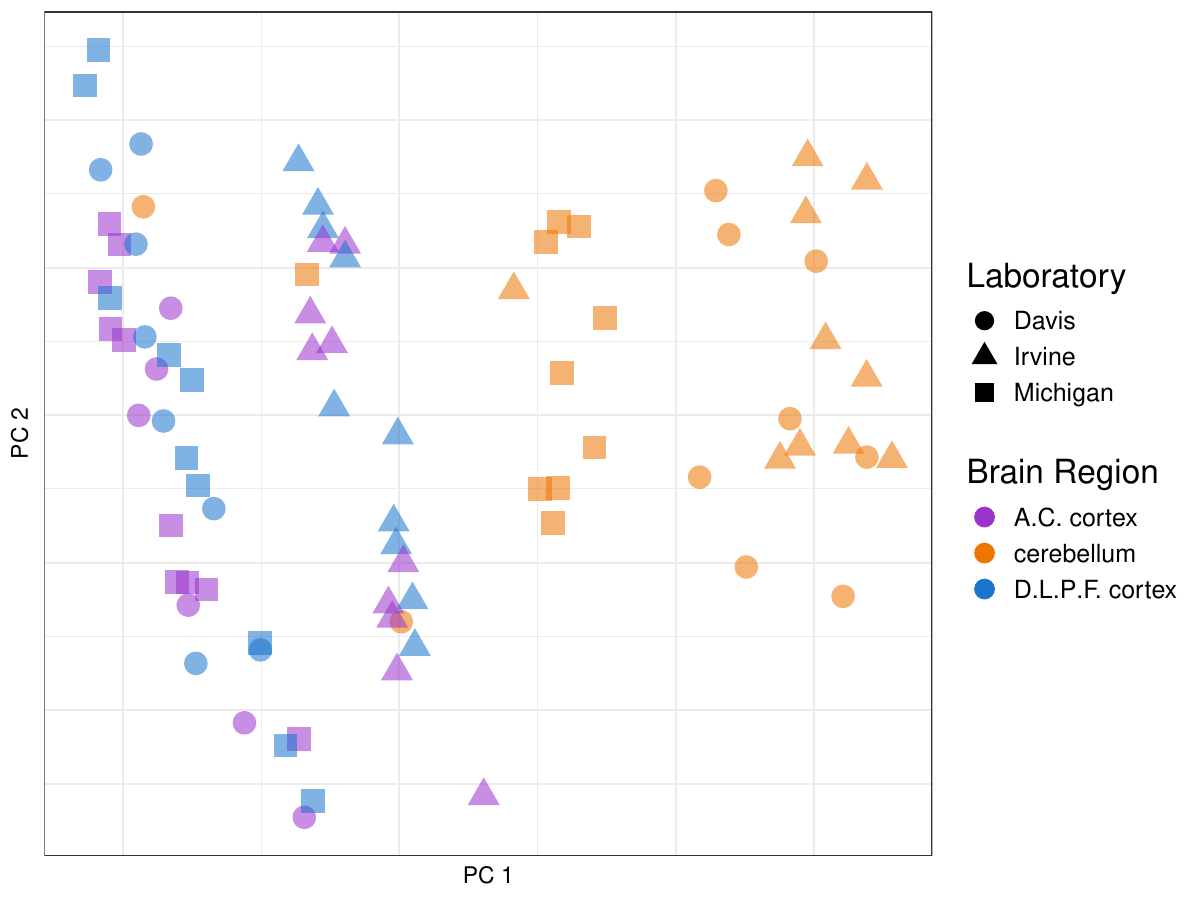}  &
    \includegraphics[scale=0.3]{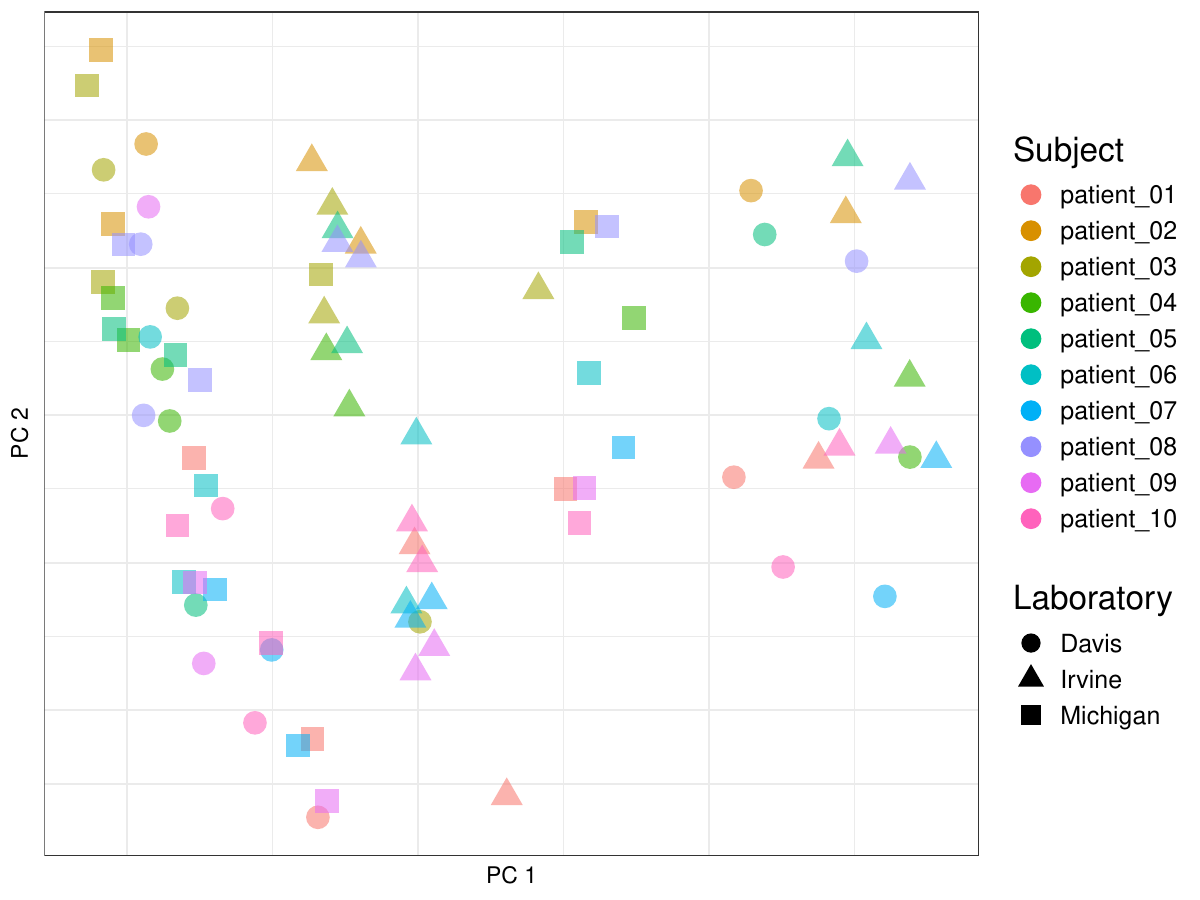}     
    \end{tabular}
    \caption{Like Figure 4, but with $K=2$.}
    \label{fig:K2_svd.all}
\end{figure}

\begin{figure}
    \centering
    \begin{tabular}{cc}
    \hspace{-1cm} No Adjustment & 
    \hspace{-1cm} ``Technical'' Adjustment \\
    \includegraphics[scale=0.3]{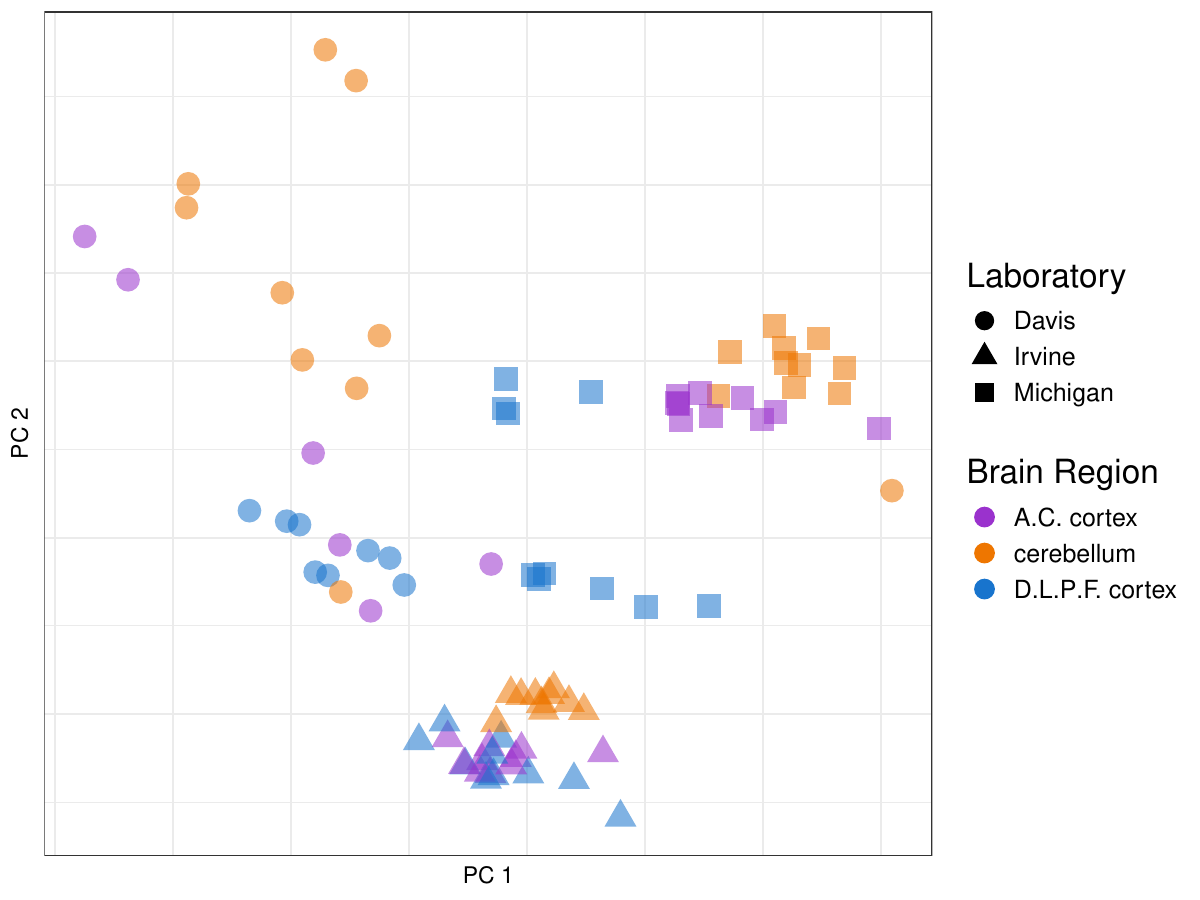} &
    \includegraphics[scale=0.3]{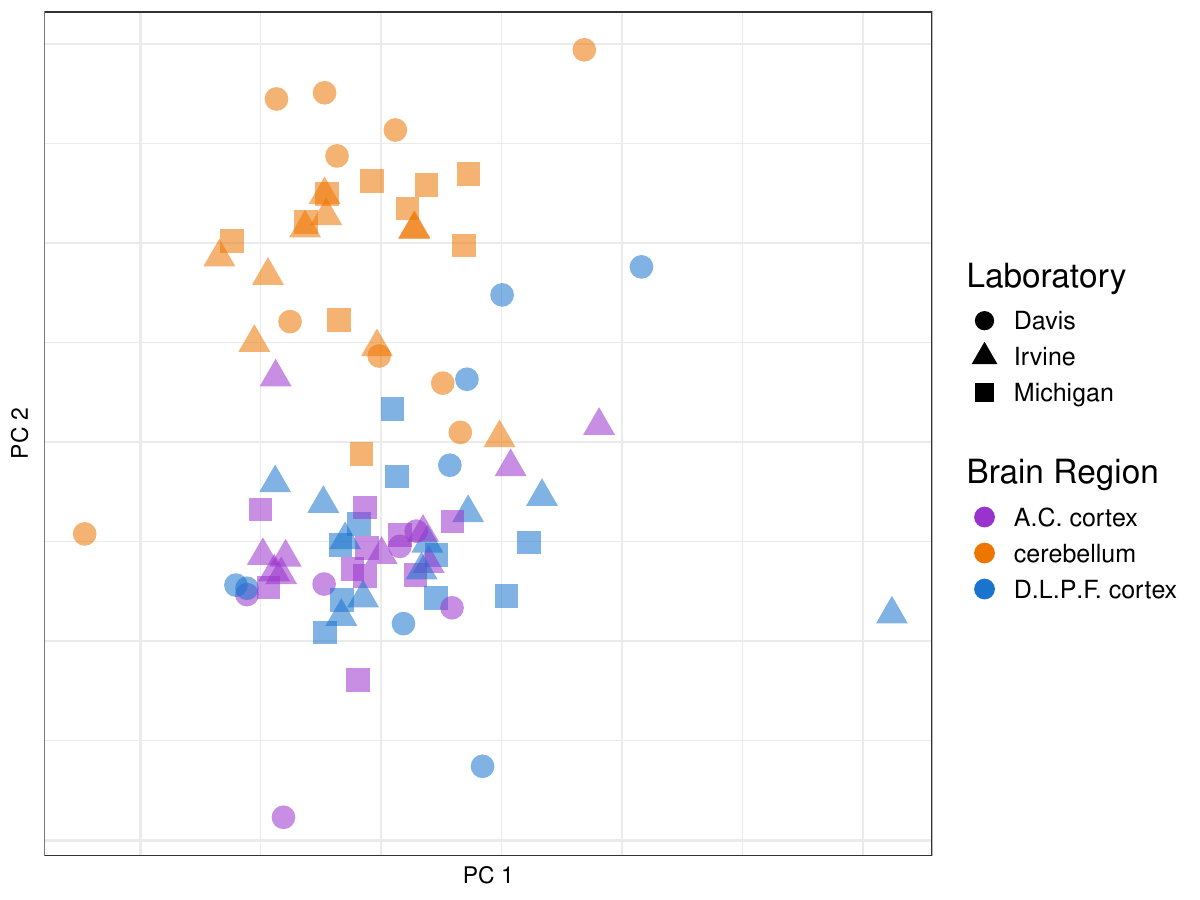} \\    
    \hspace{-1cm} ``Bio'' Adjustment & 
    \hspace{-1cm} ``Bio'' Adjustment (colored by patient) \\
    \includegraphics[scale=0.3]{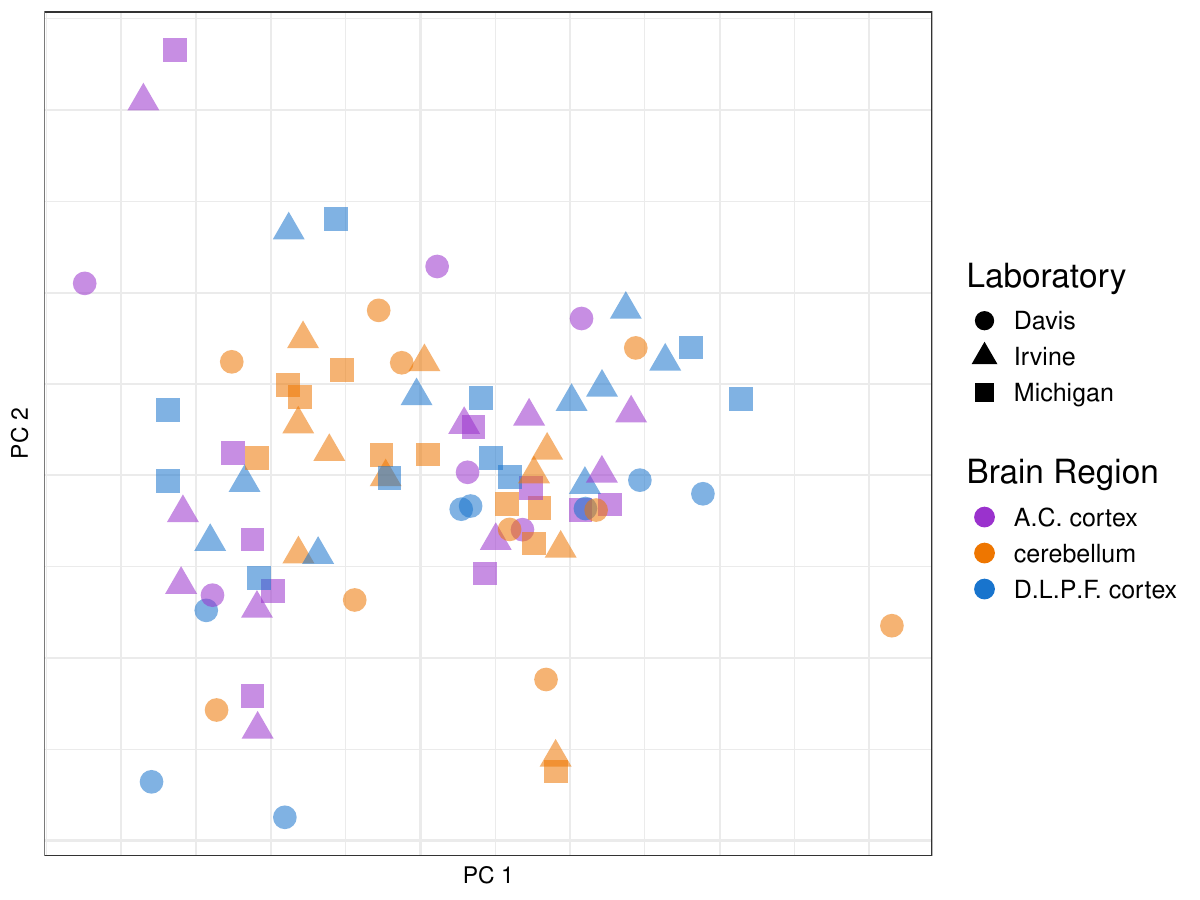}  &
    \includegraphics[scale=0.3]{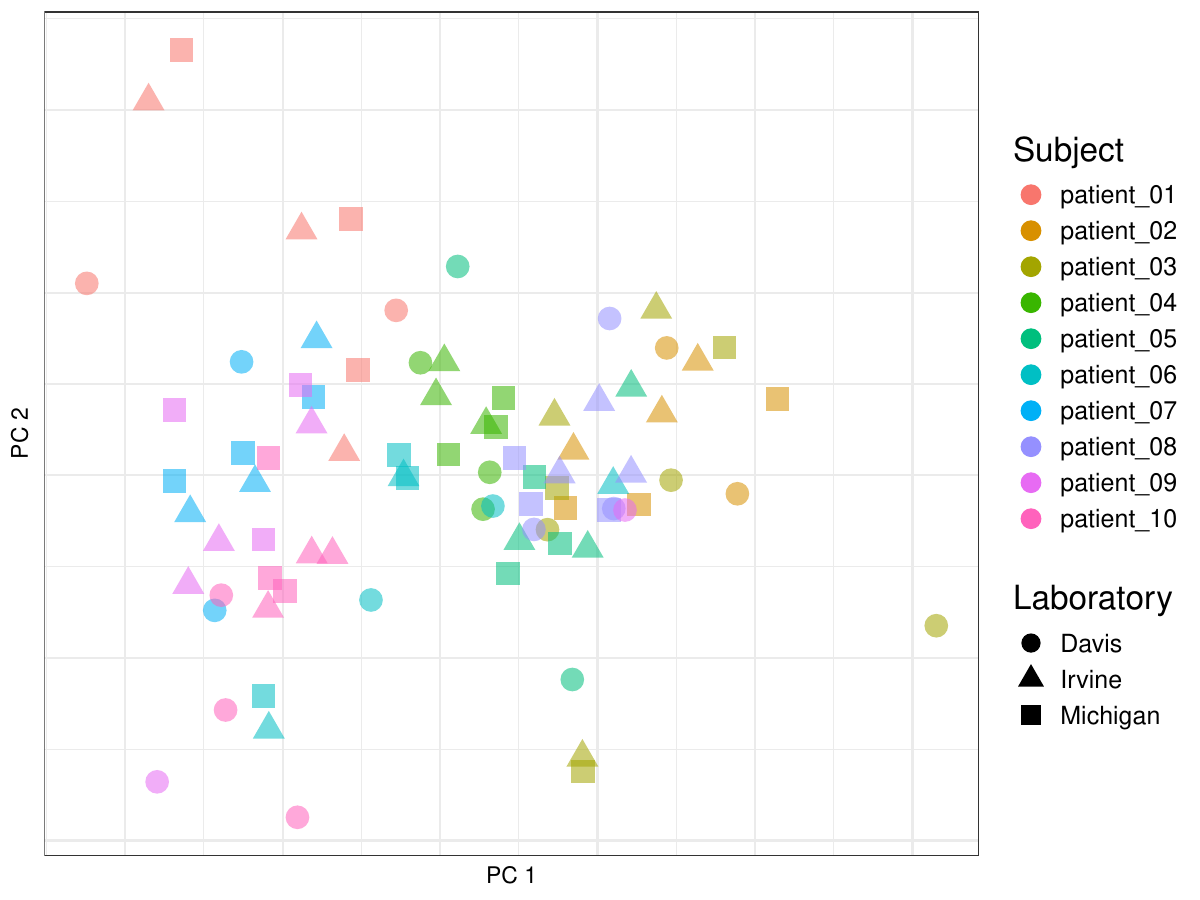}     
    \end{tabular}
    \caption{Like Figure 4, but with $K=5$.}
    \label{fig:K5_svd.all}
\end{figure}

\begin{figure}
    \centering
    \begin{tabular}{cc}
    \hspace{-1cm} No Adjustment & 
    \hspace{-1cm} ``Technical'' Adjustment \\
    \includegraphics[scale=0.3]{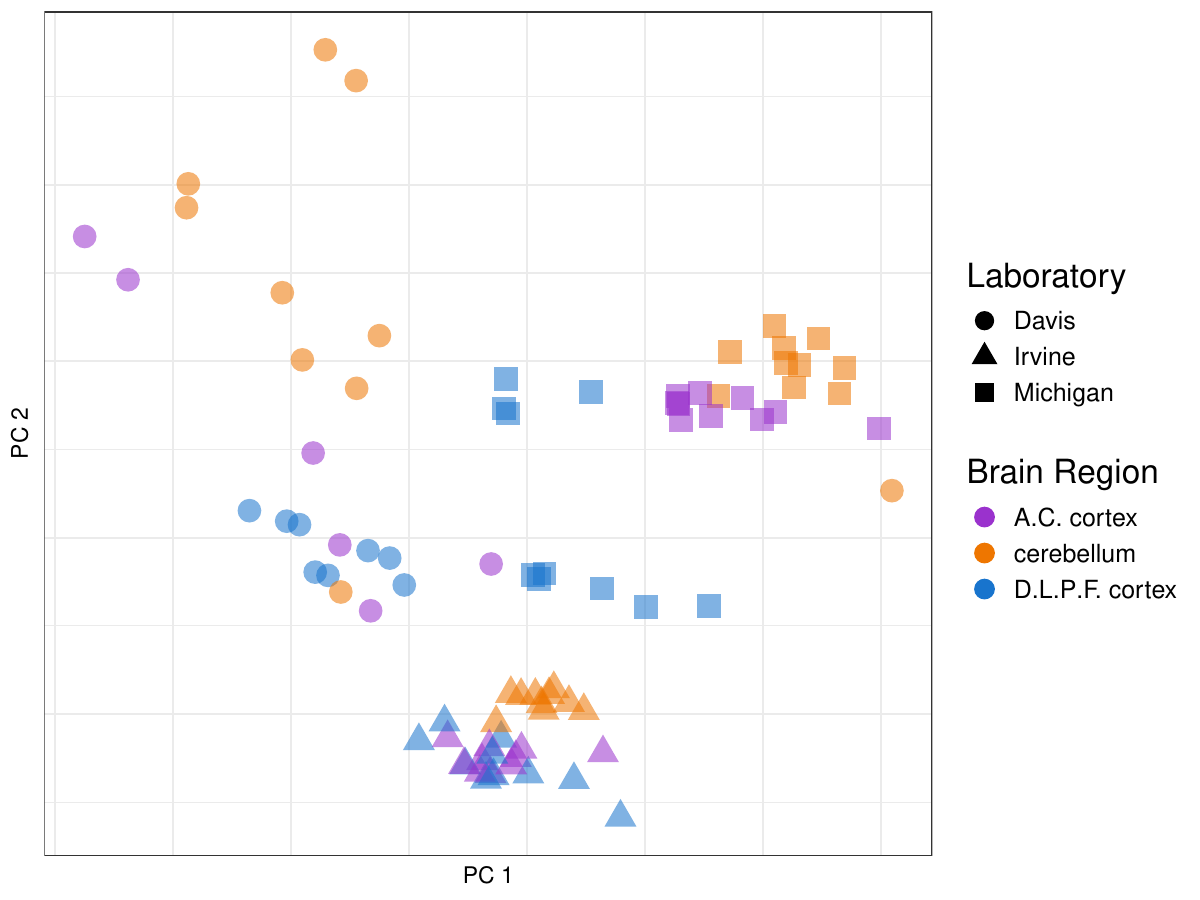} &
    \includegraphics[scale=0.3]{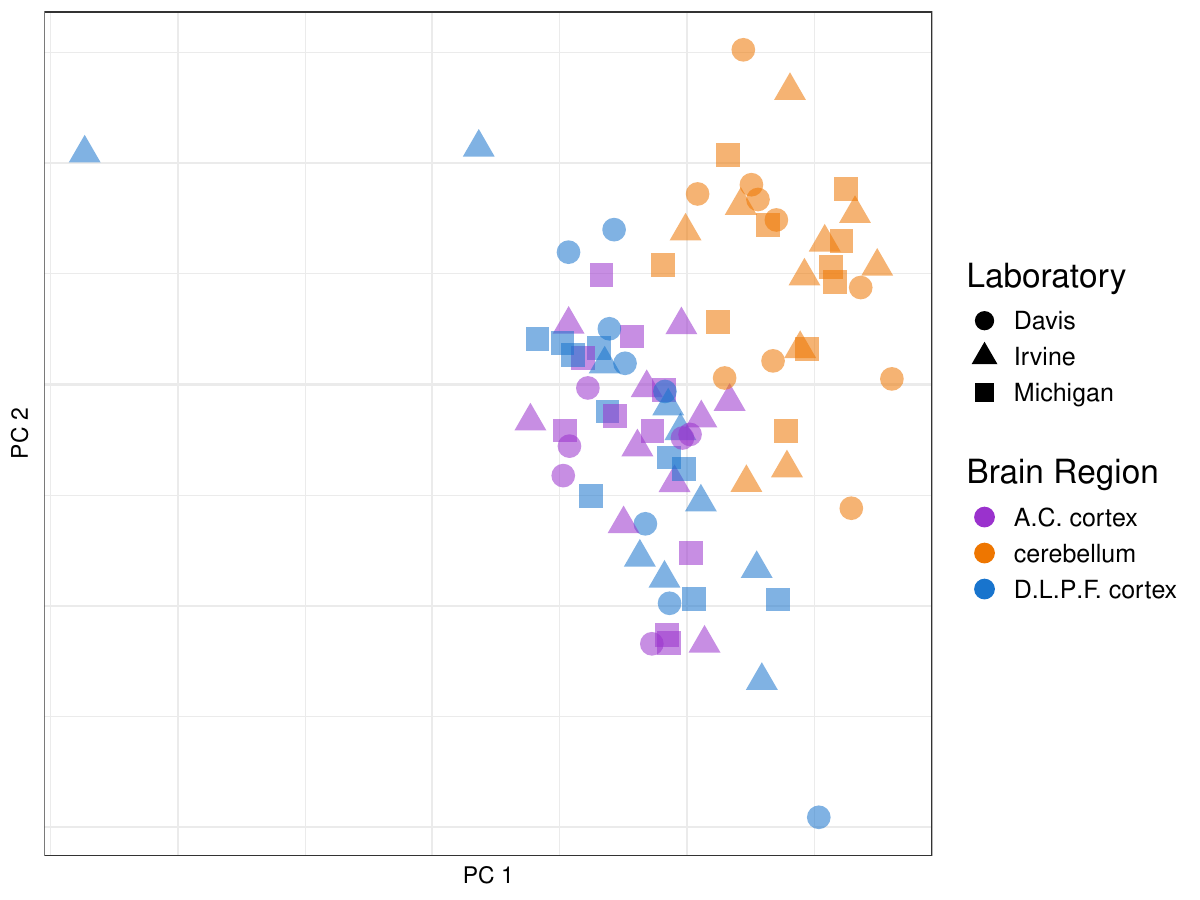} \\    
    \hspace{-1cm} ``Bio'' Adjustment & 
    \hspace{-1cm} ``Bio'' Adjustment (colored by patient) \\
    \includegraphics[scale=0.3]{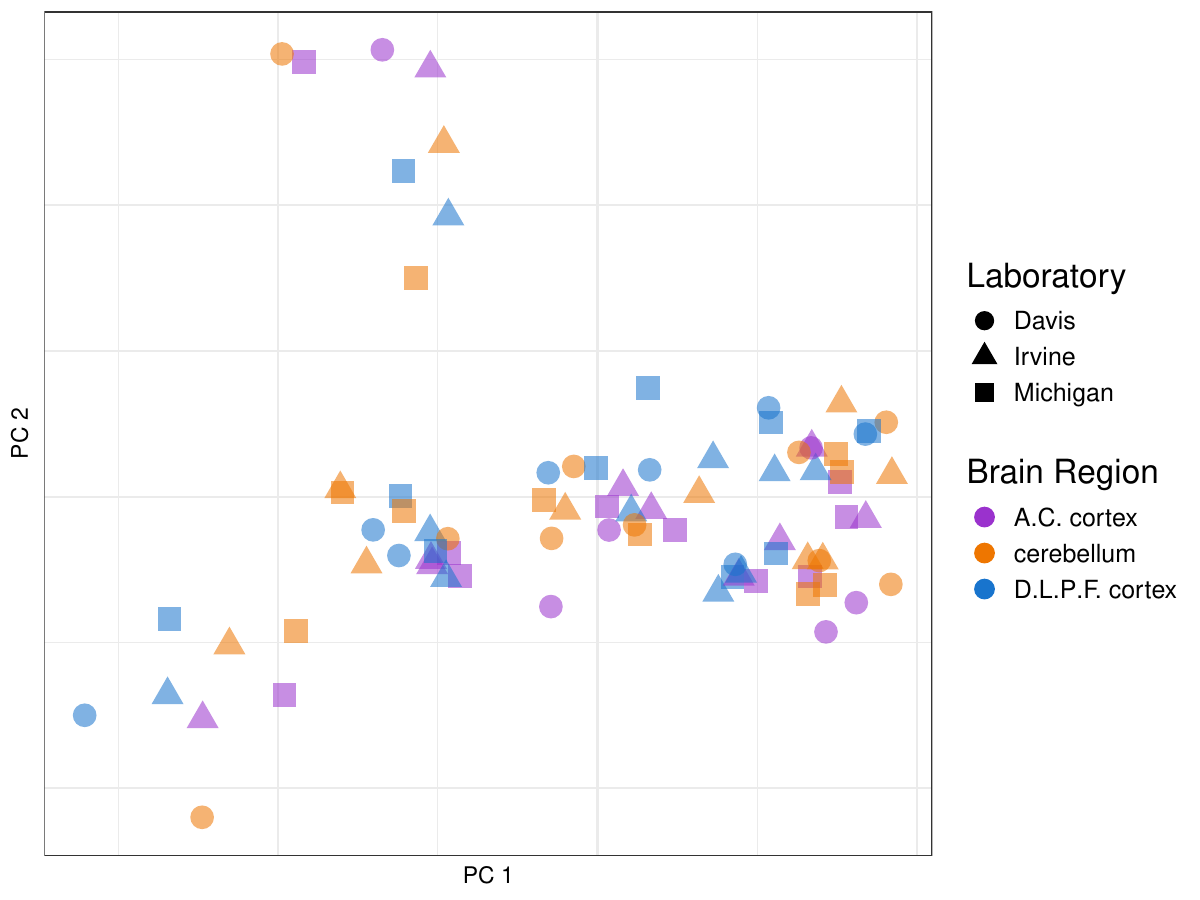}  &
    \includegraphics[scale=0.3]{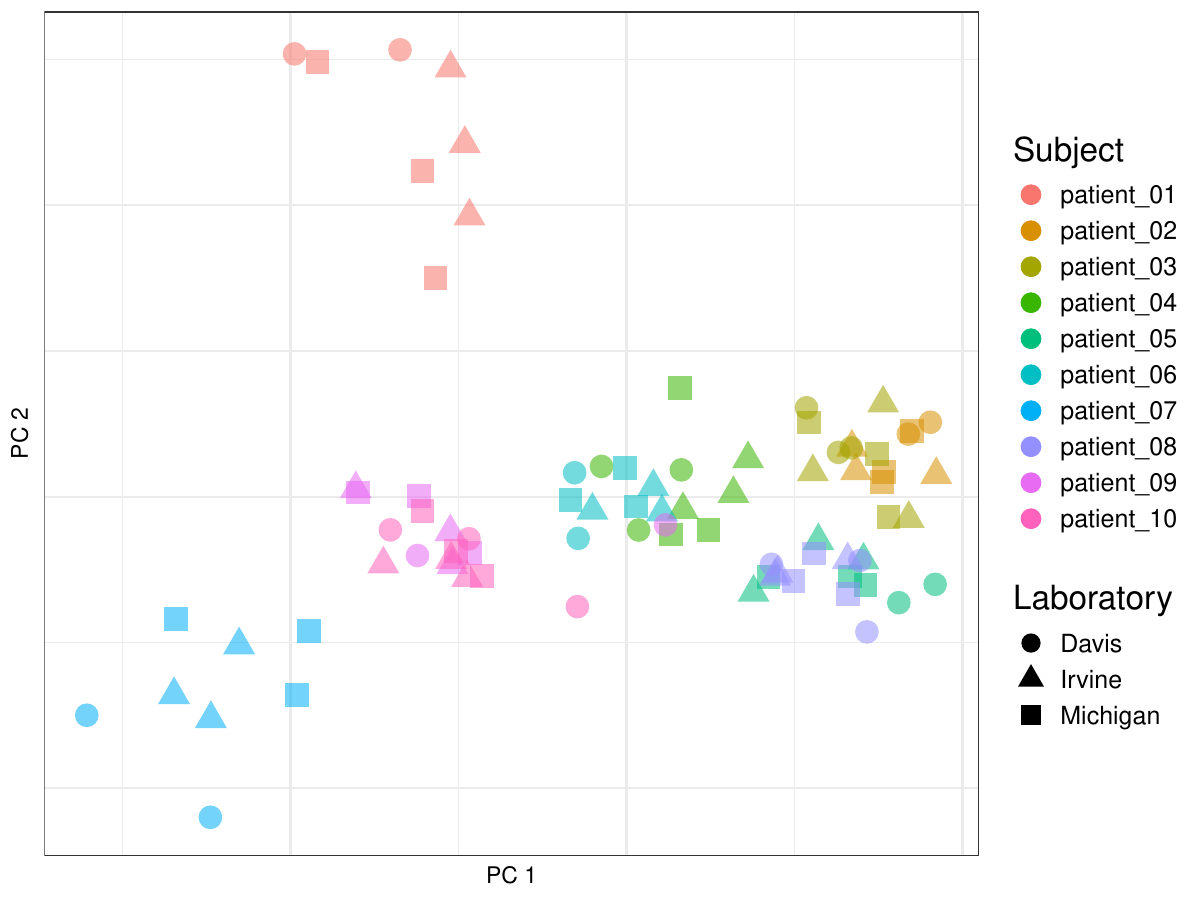}     
    \end{tabular}
    \caption{Like Figure 4, but with $K=15$.}
    \label{fig:K15_svd.all}
\end{figure}

\begin{figure}
    \centering
    \begin{tabular}{cc}
    \hspace{-1cm} No Adjustment & 
    \hspace{-1cm} ``Technical'' Adjustment \\
    \includegraphics[scale=0.3]{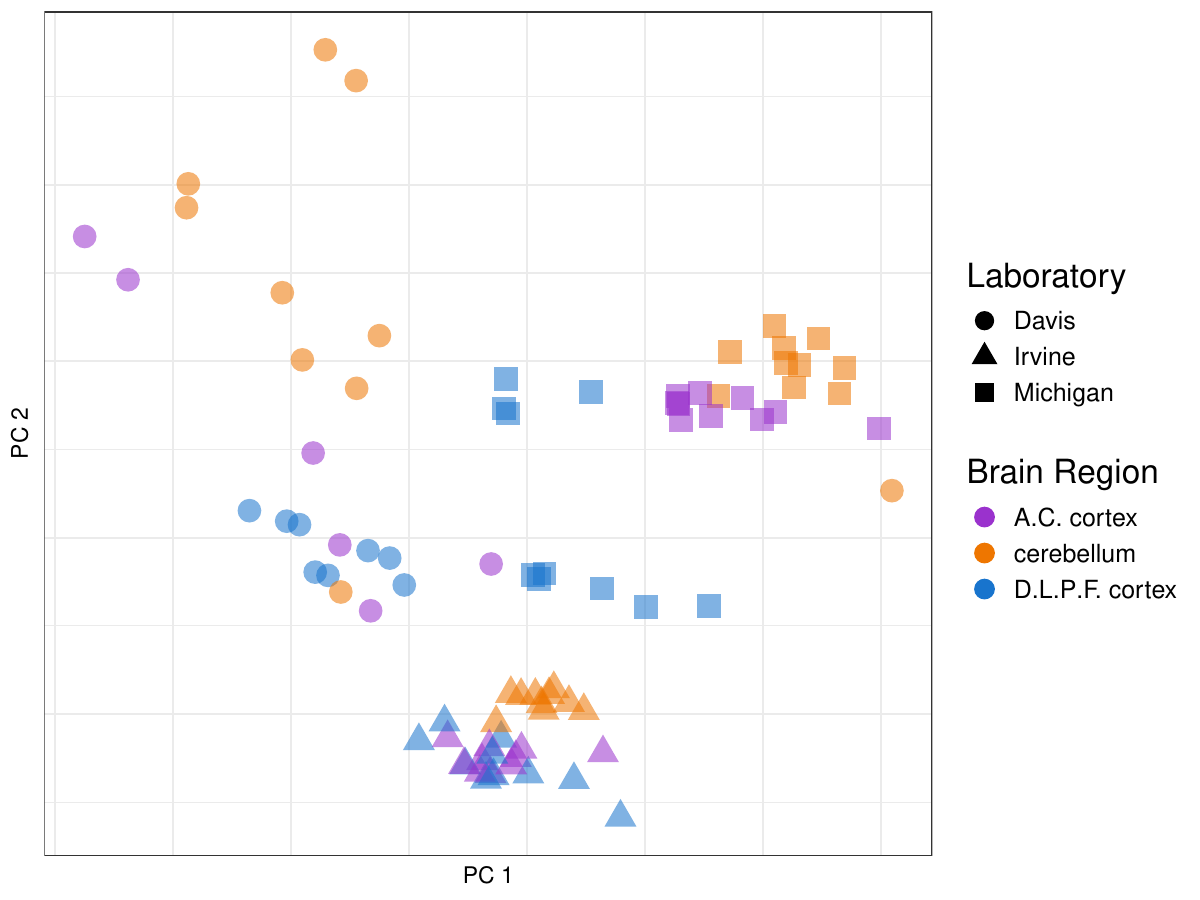} &
    \includegraphics[scale=0.3]{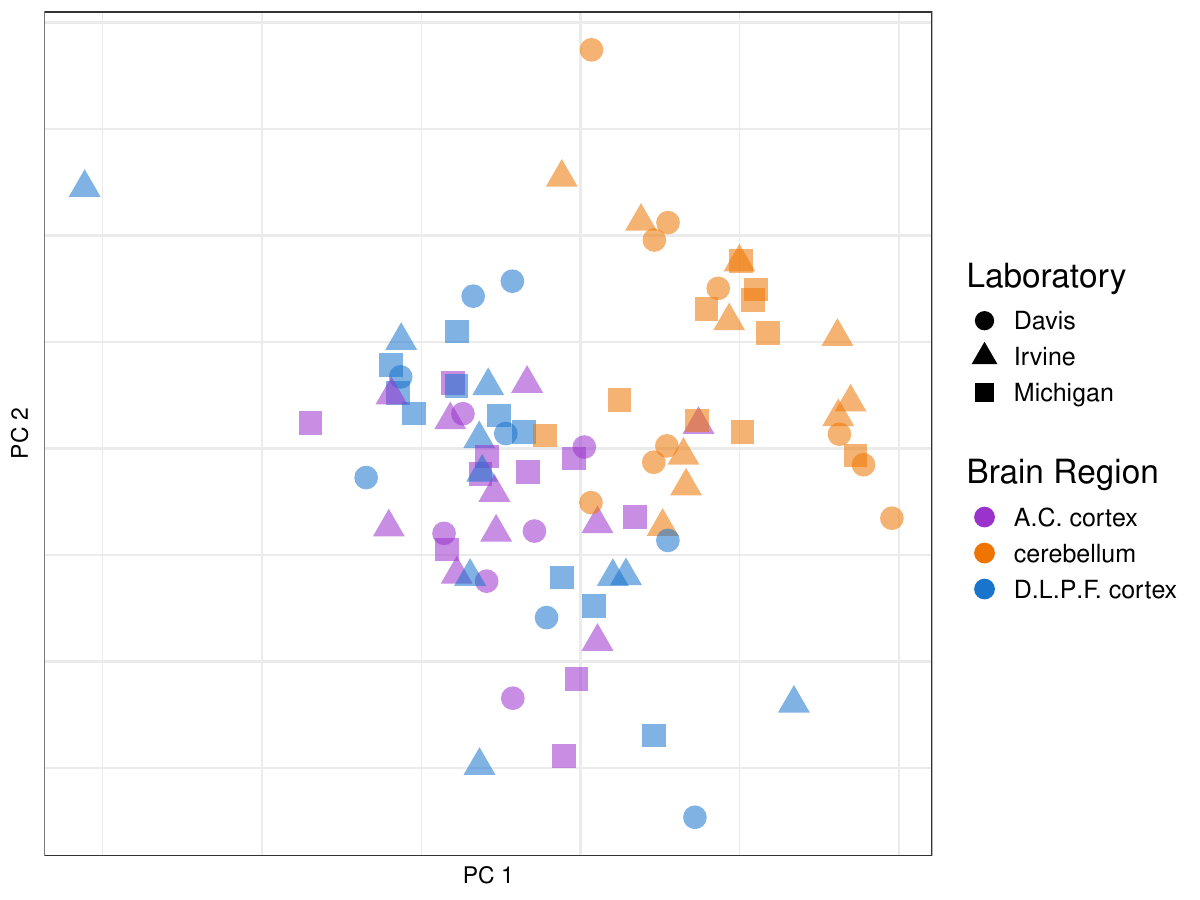} \\    
    \hspace{-1cm} ``Bio'' Adjustment & 
    \hspace{-1cm} ``Bio'' Adjustment (colored by patient) \\
    \includegraphics[scale=0.3]{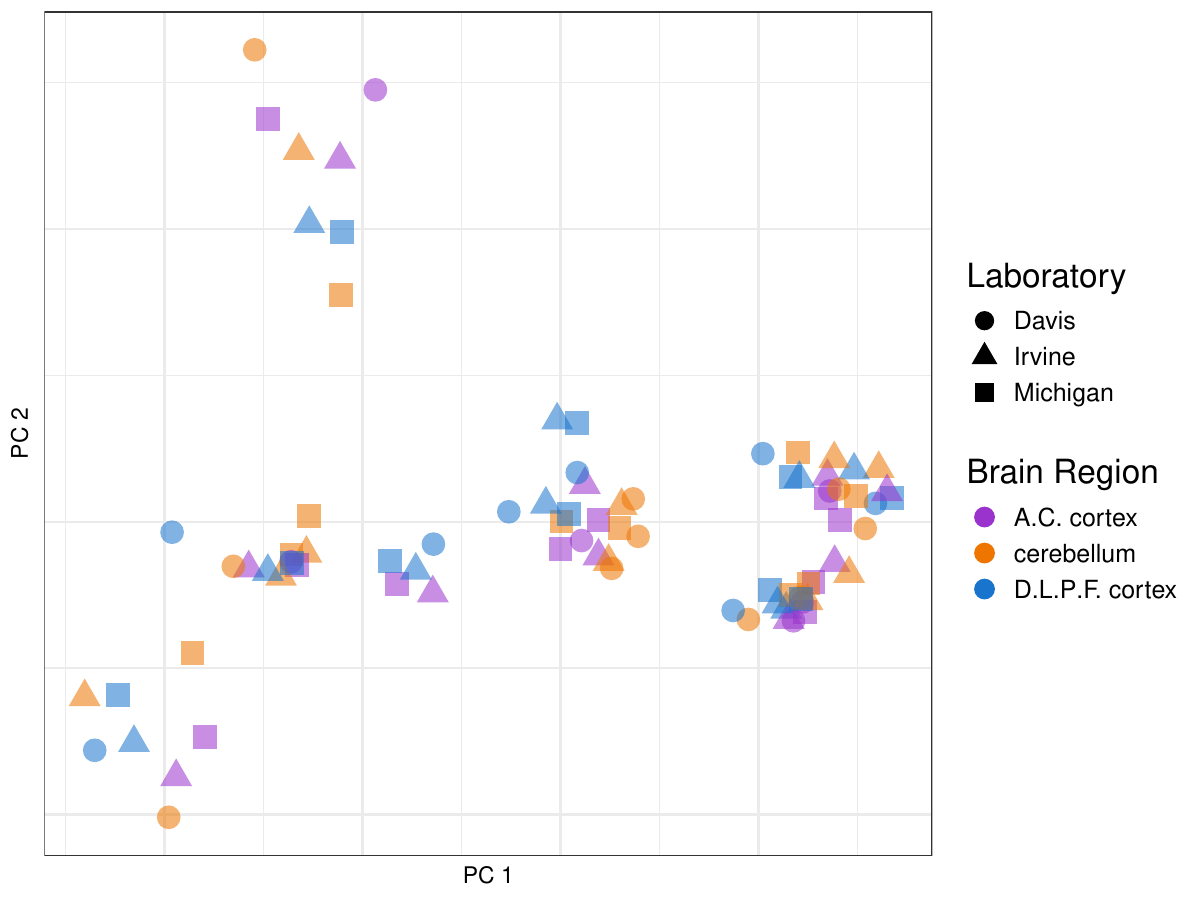}  &
    \includegraphics[scale=0.3]{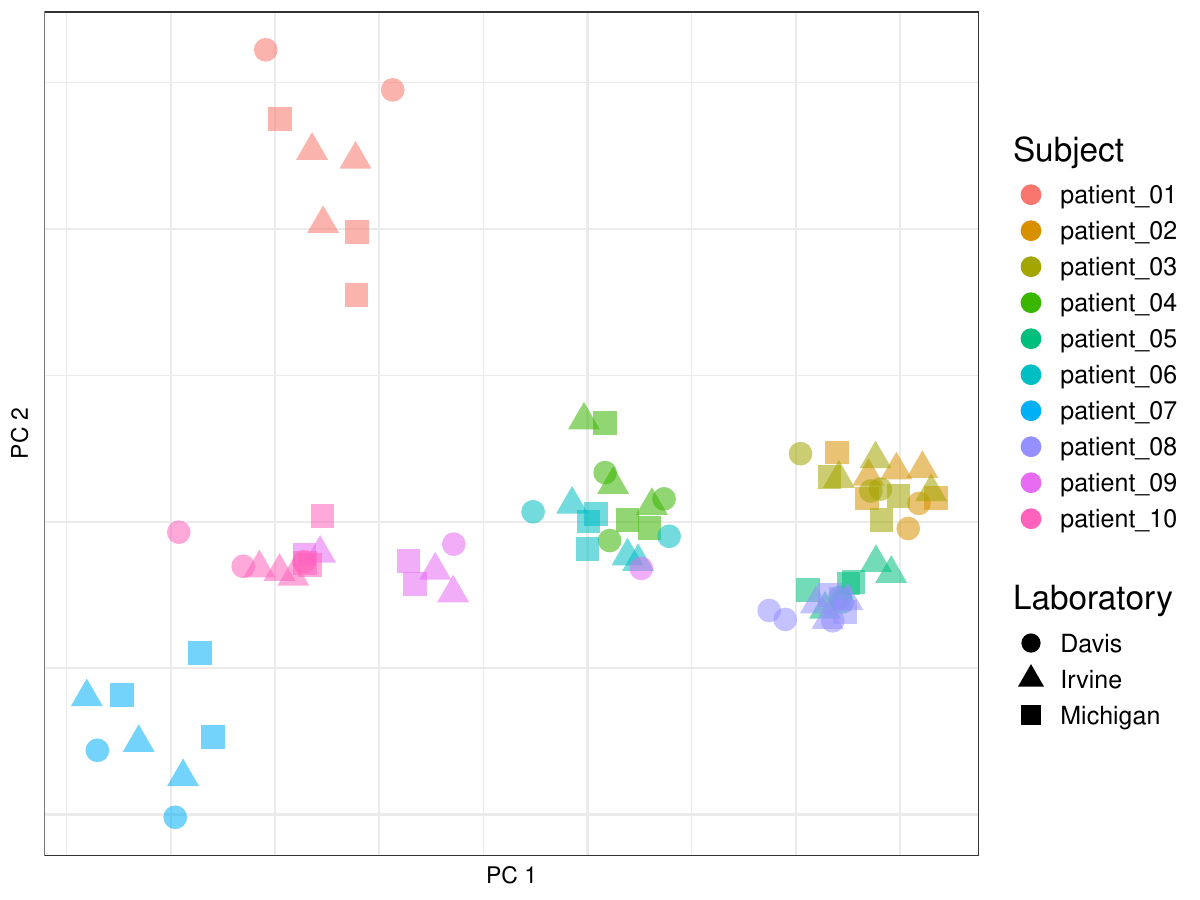}     
    \end{tabular}
    \caption{Like Figure 4, but with $K=20$.}
    \label{fig:K20_svd.all}
\end{figure}

\begin{figure}
    \centering
    \begin{tabular}{cc}
    \hspace{-1cm} No Adjustment & 
    \hspace{-1cm} ``Technical'' Adjustment \\
    \includegraphics[scale=0.3]{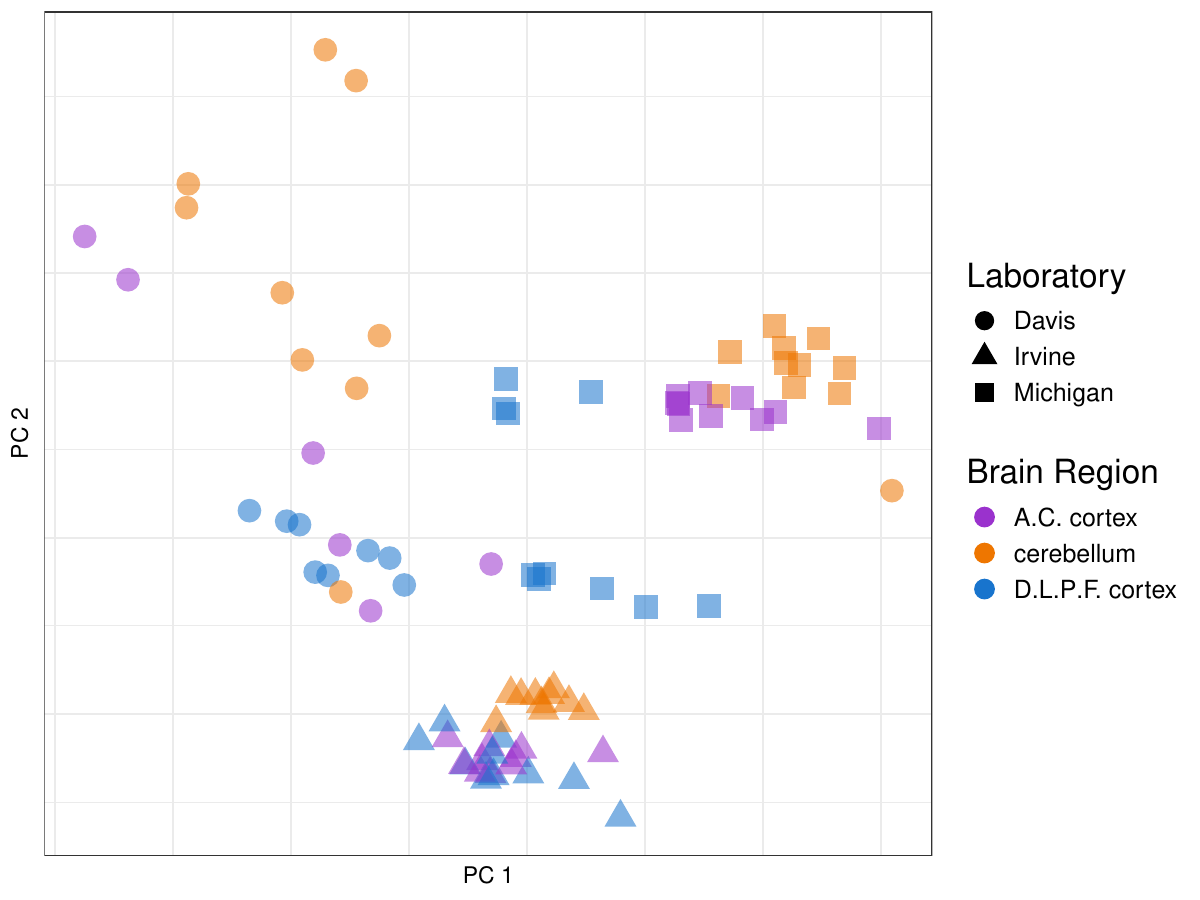} &
    \includegraphics[scale=0.3]{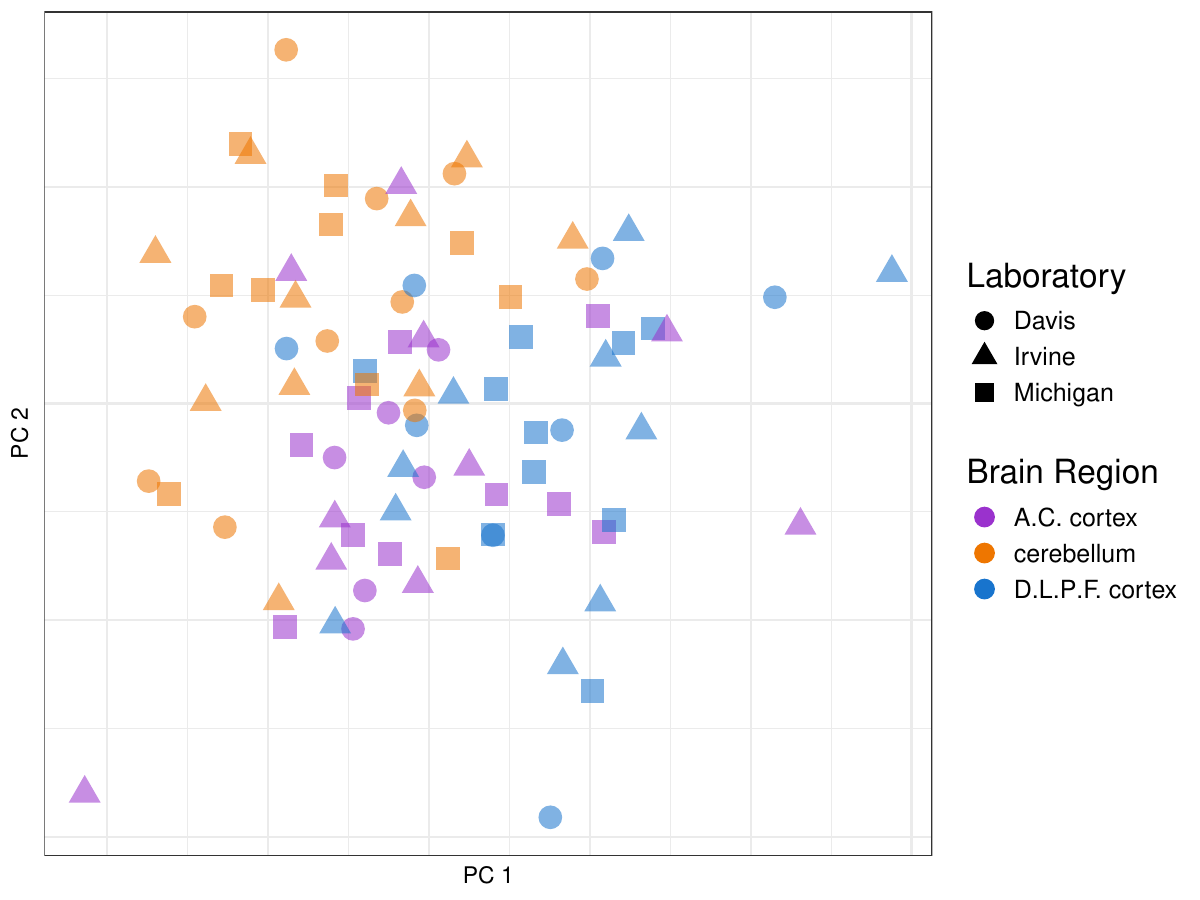} \\    
    \hspace{-1cm} ``Bio'' Adjustment & 
    \hspace{-1cm} ``Bio'' Adjustment (colored by patient) \\
    \includegraphics[scale=0.3]{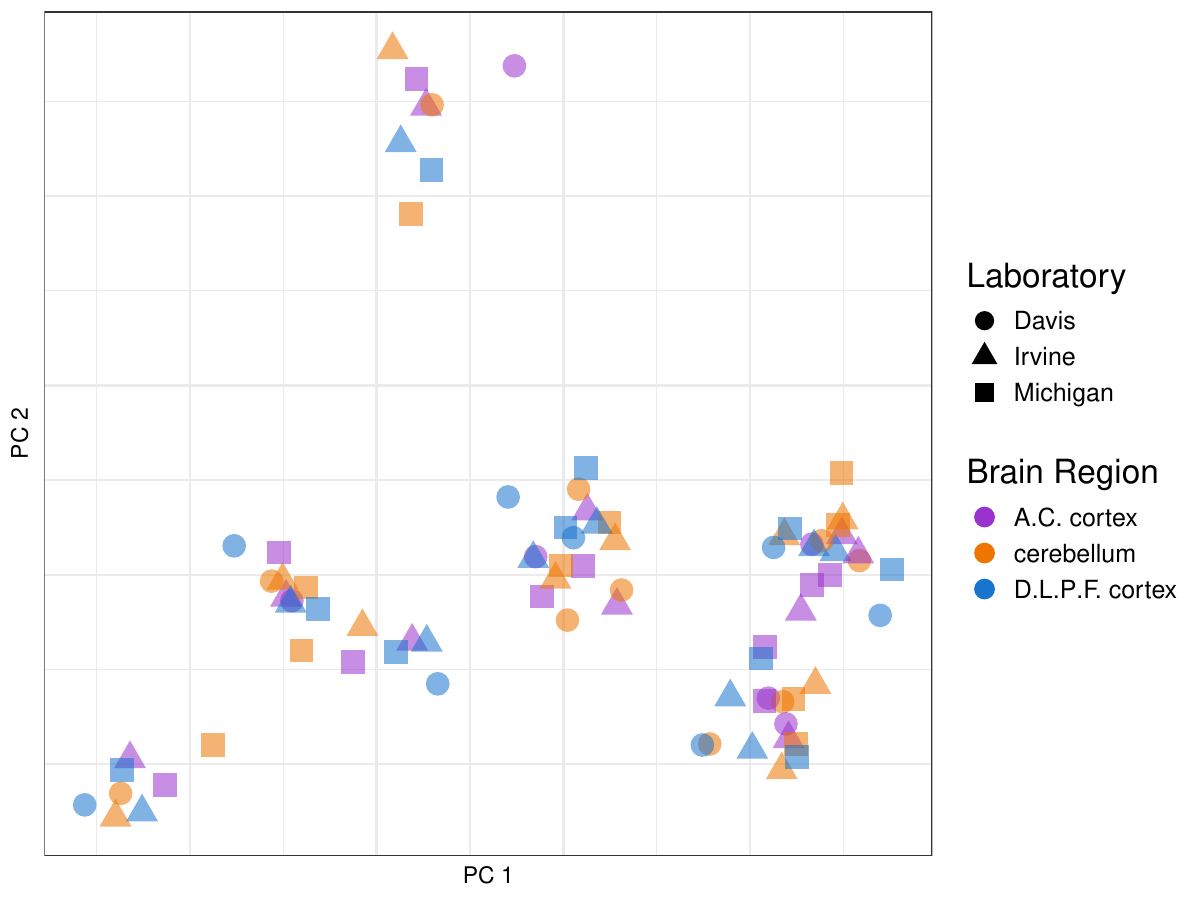}  &
    \includegraphics[scale=0.3]{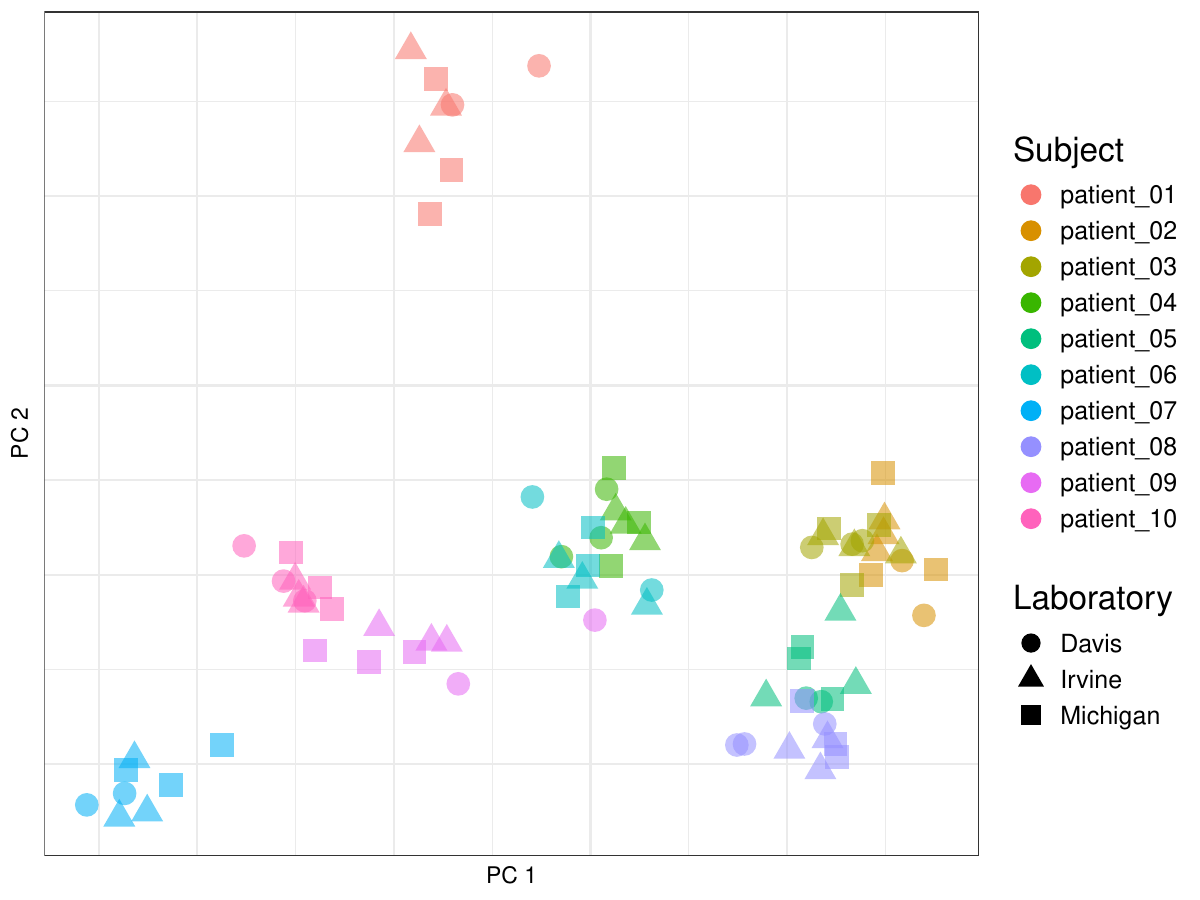}     
    \end{tabular}
    \caption{Like Figure 4, but with $K=30$.}
    \label{fig:K30_svd.all}
\end{figure}

\begin{figure}[h]
    \centering
    \begin{tabular}{cc}
    \multicolumn{2}{c}{No Adjustment}\\
    \includegraphics[scale=0.3]{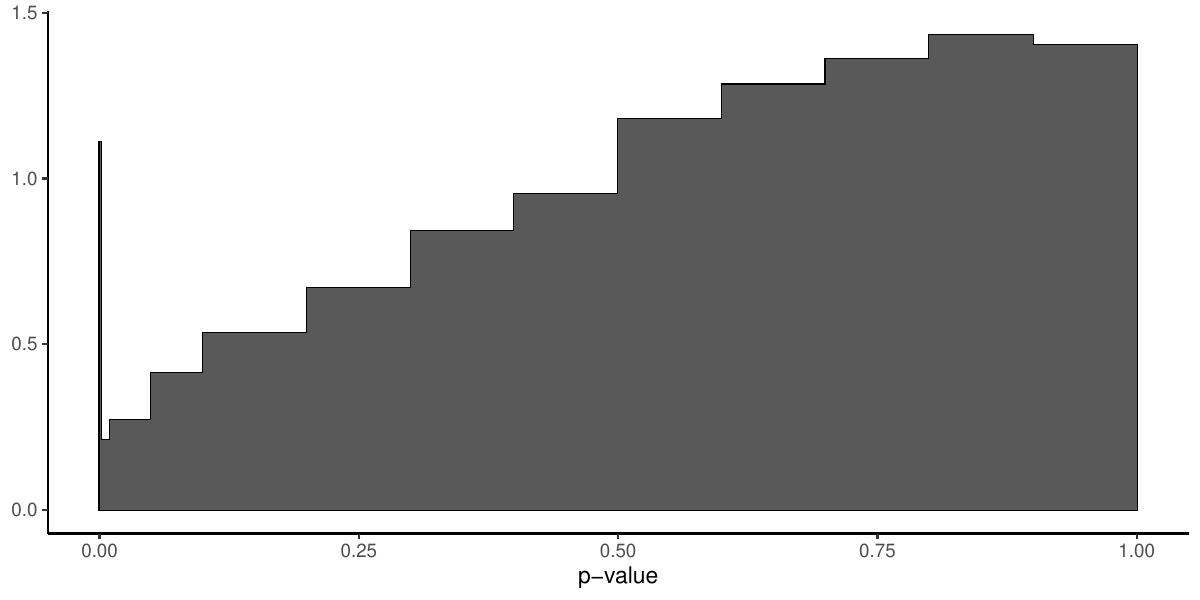} &
    \includegraphics[scale=0.3]{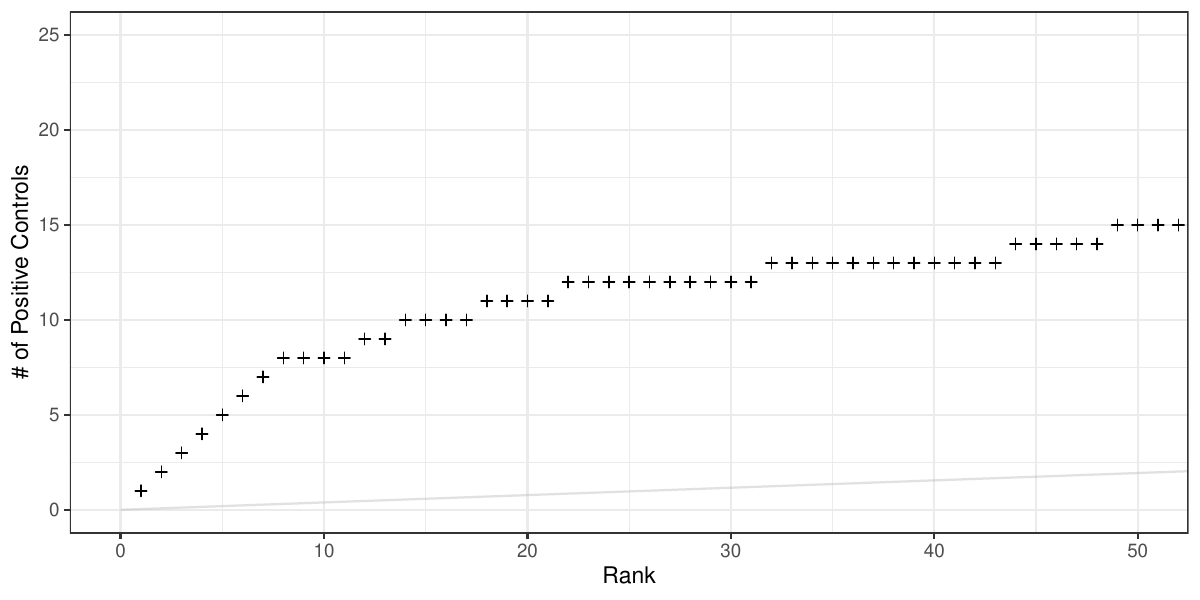} \\    
    \multicolumn{2}{c}{``Technical'' Adjustment}\\
    \includegraphics[scale=0.3]{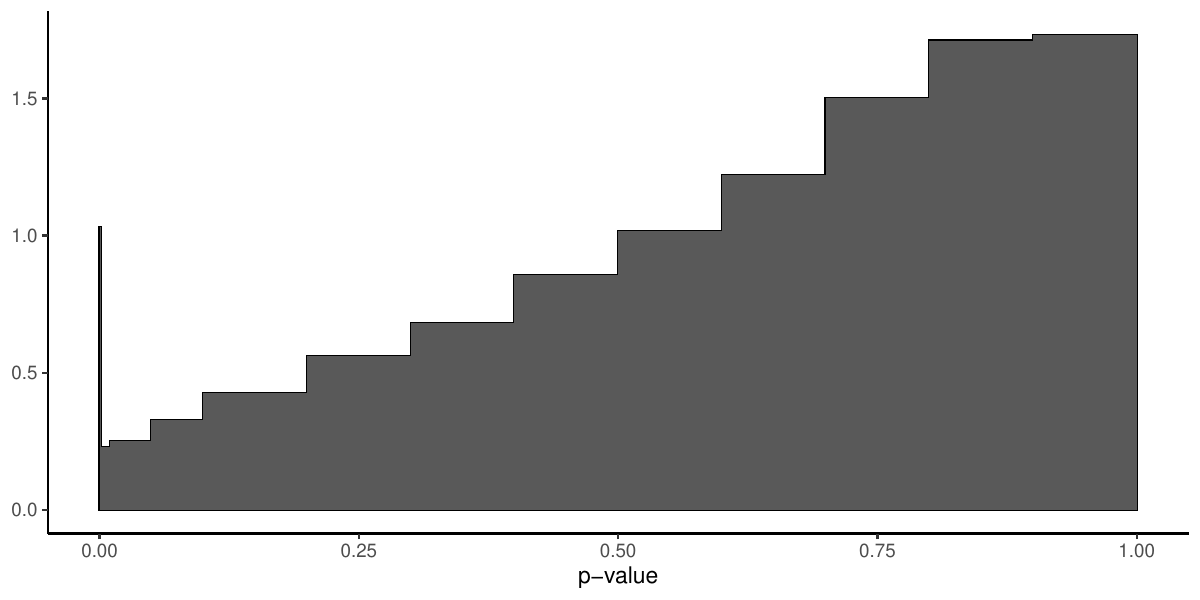} &
    \includegraphics[scale=0.3]{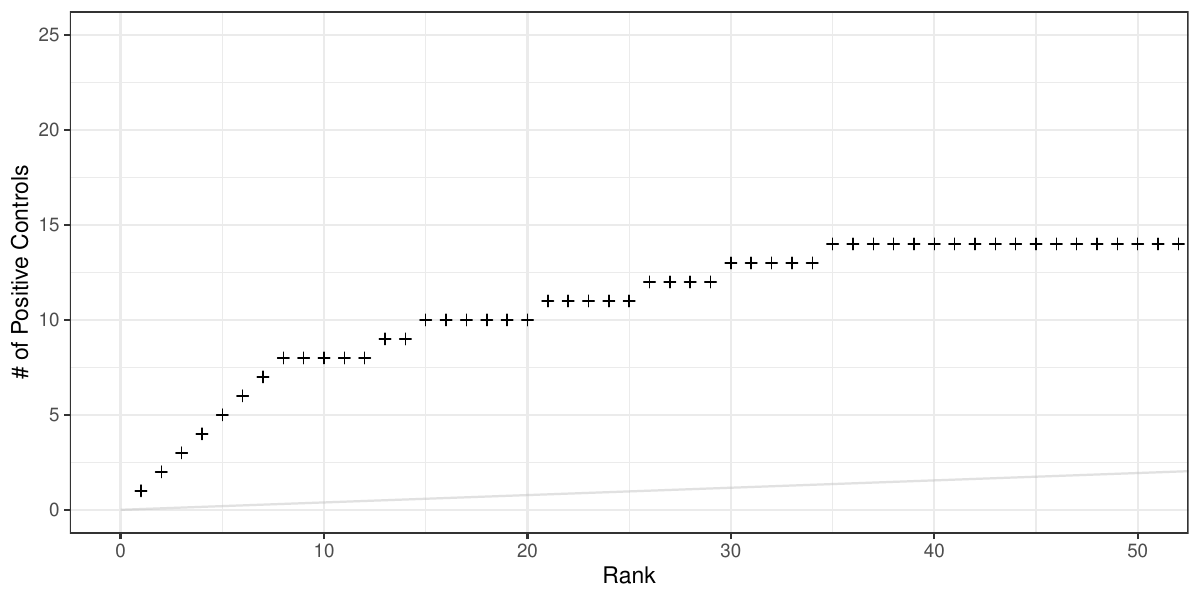} \\    
    \multicolumn{2}{c}{``Bio'' Adjustment}\\
    \includegraphics[scale=0.3]{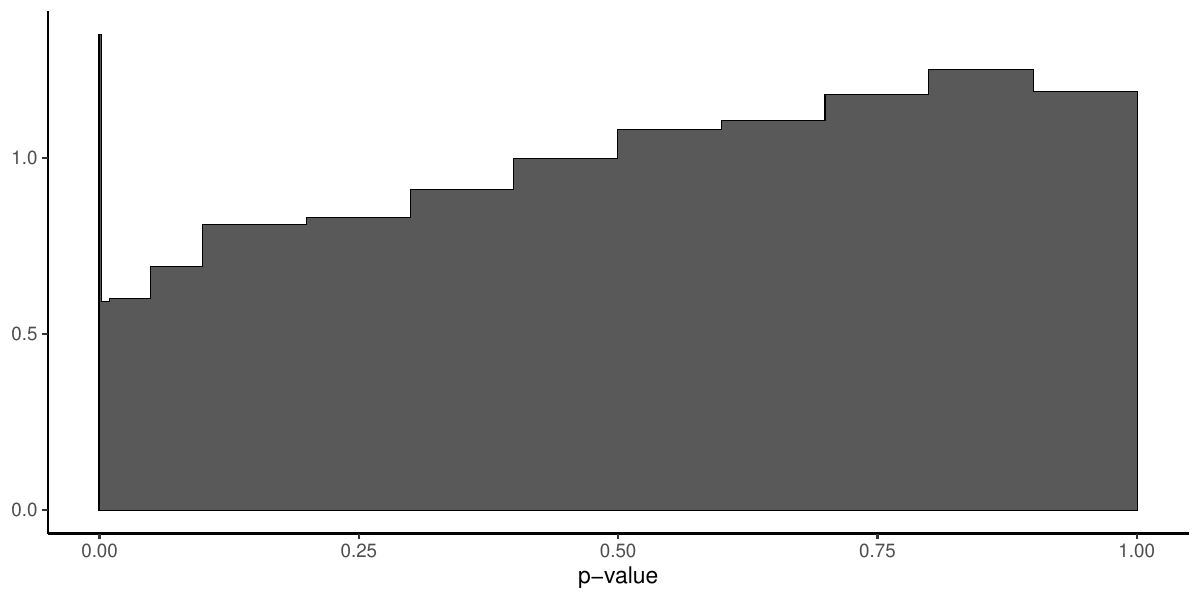} &
    \includegraphics[scale=0.3]{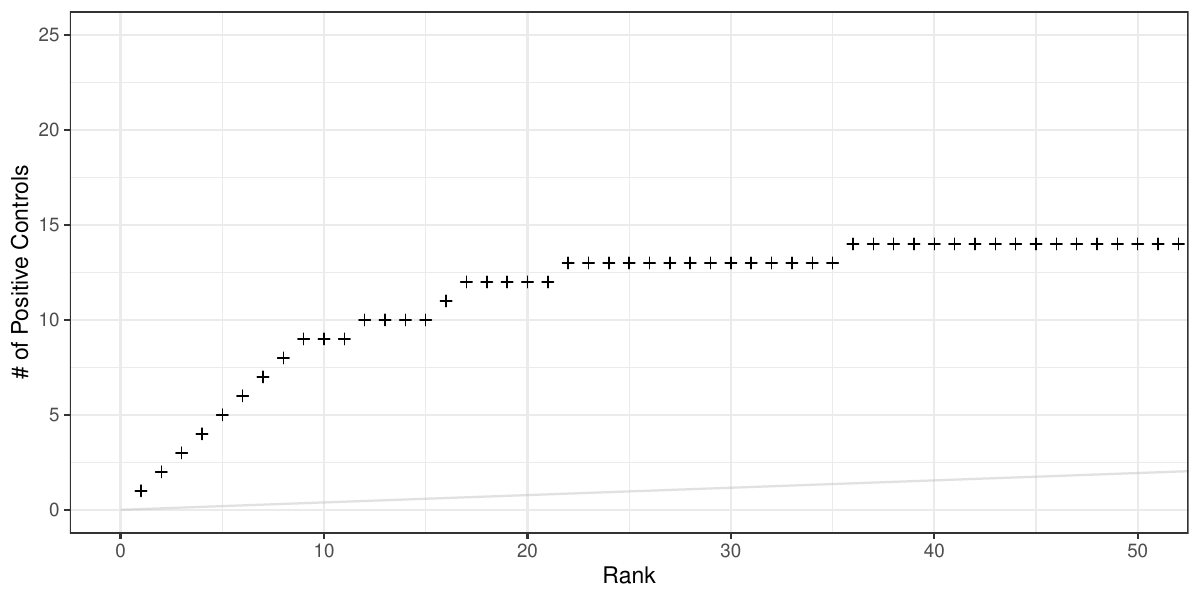} \\    
    \end{tabular}
    \caption{Like Figure 5, but with $K=1$.}
    \label{fig:K1_pvalrank.all}
\end{figure}

\begin{figure}[h]
    \centering
    \begin{tabular}{cc}
    \multicolumn{2}{c}{No Adjustment}\\
    \includegraphics[scale=0.3]{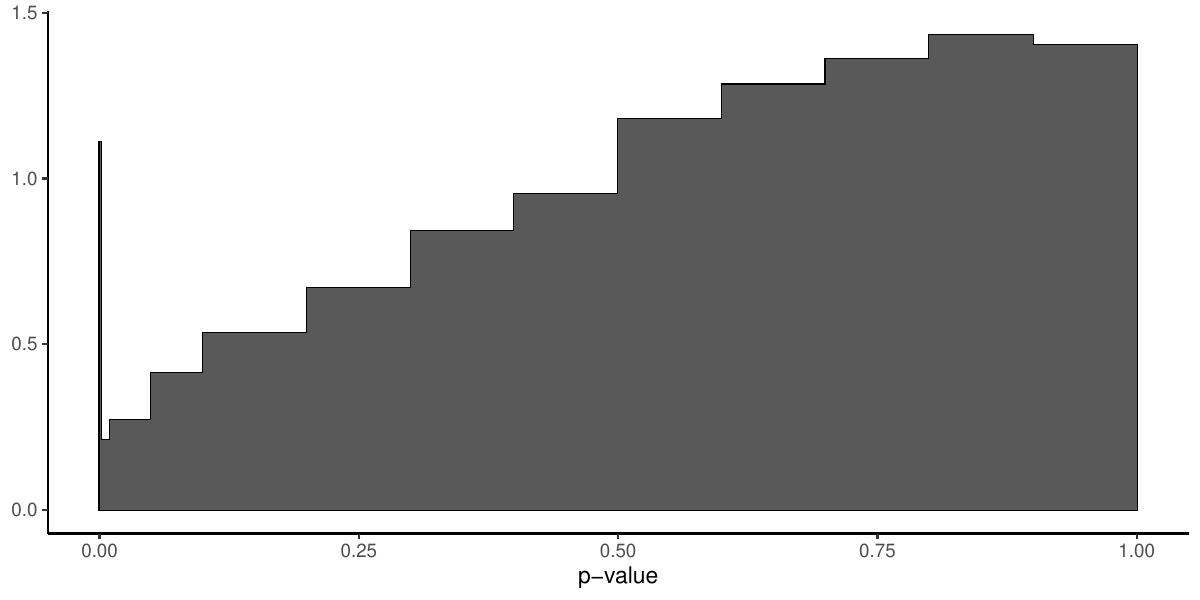} &
    \includegraphics[scale=0.3]{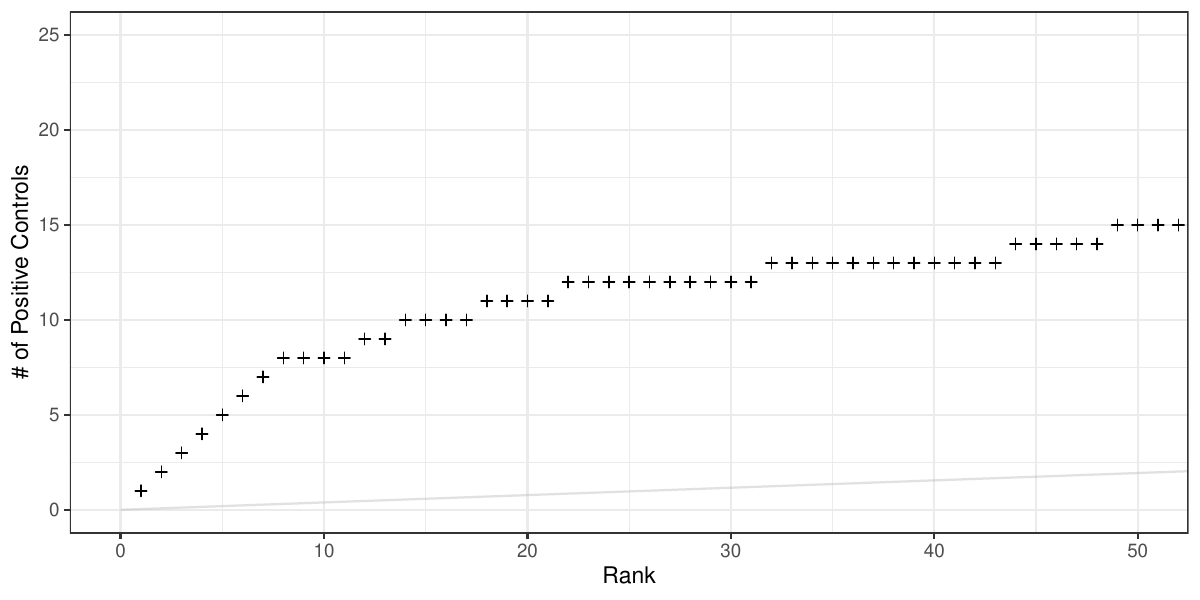} \\    
    \multicolumn{2}{c}{``Technical'' Adjustment}\\
    \includegraphics[scale=0.3]{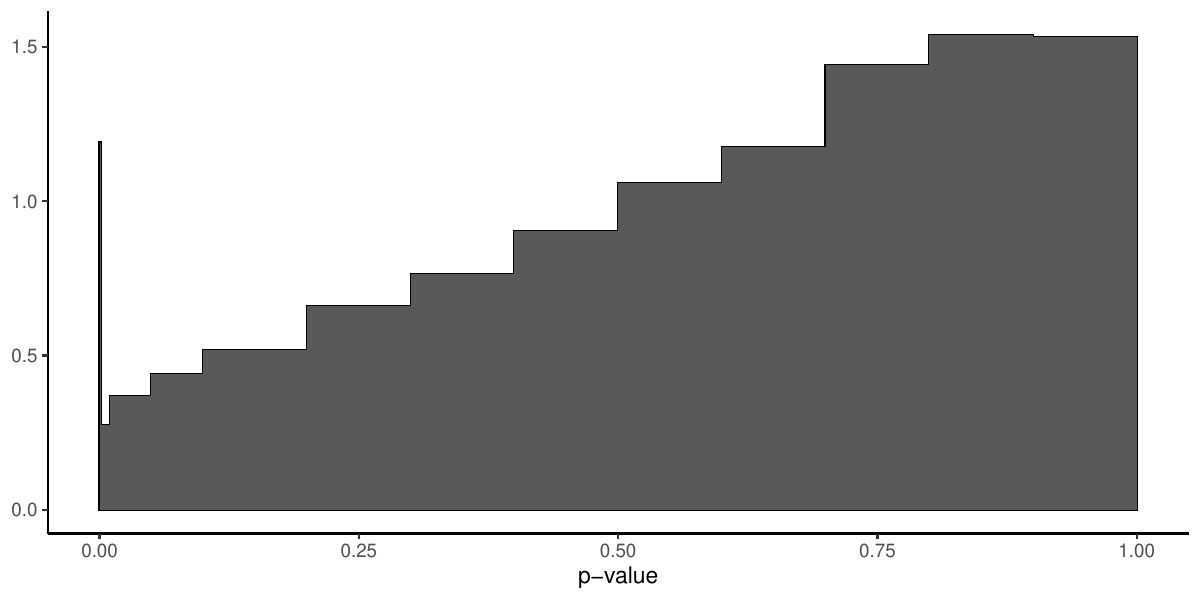} &
    \includegraphics[scale=0.3]{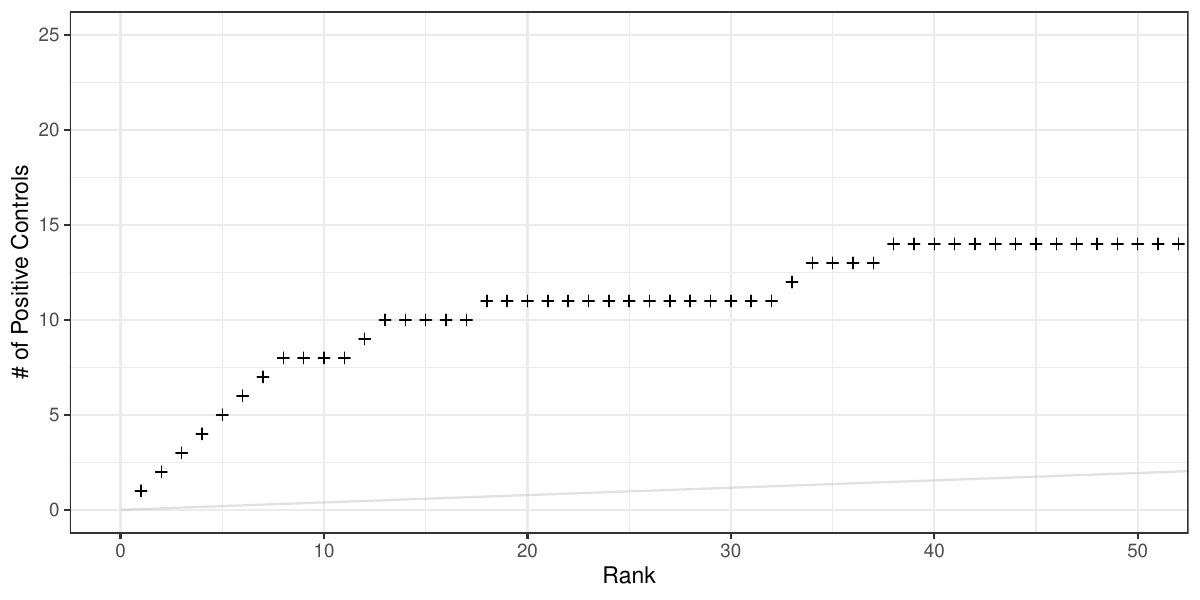} \\    
    \multicolumn{2}{c}{``Bio'' Adjustment}\\
    \includegraphics[scale=0.3]{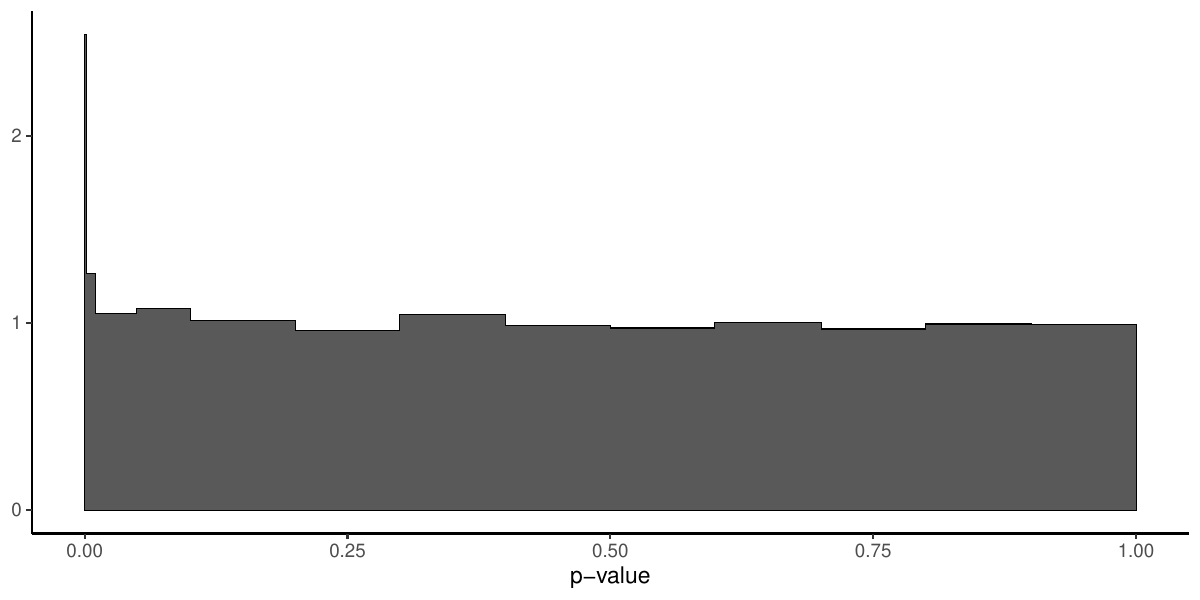} &
    \includegraphics[scale=0.3]{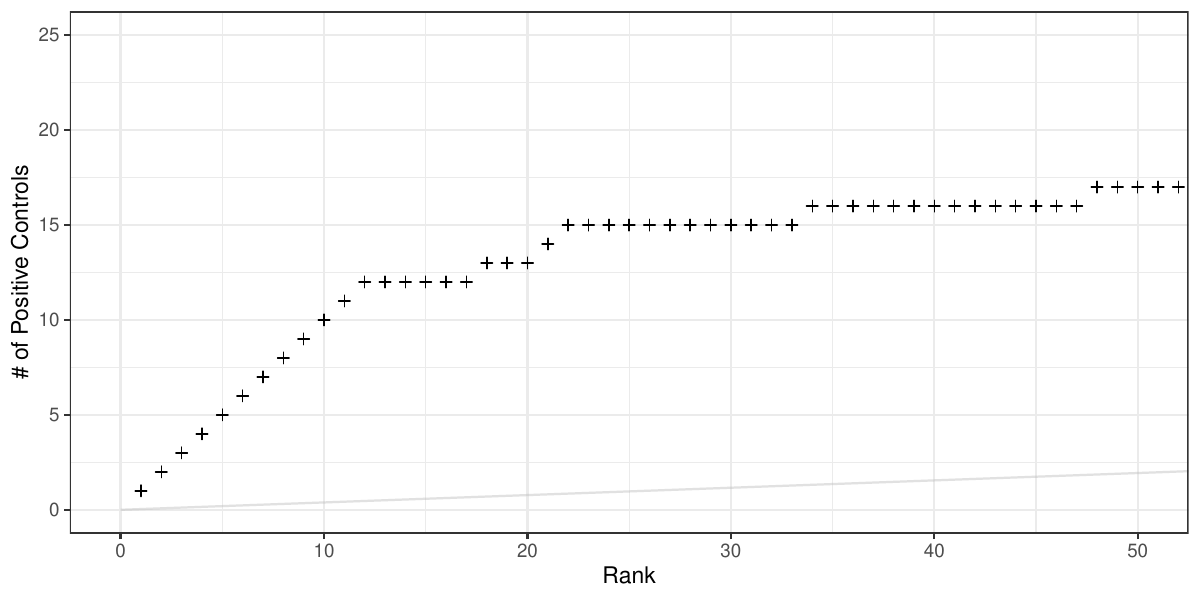} \\    
    \end{tabular}
    \caption{Like Figure 5, but with $K=2$.}
    \label{fig:K2_pvalrank.all}
\end{figure}

\begin{figure}[h]
    \centering
    \begin{tabular}{cc}
    \multicolumn{2}{c}{No Adjustment}\\
    \includegraphics[scale=0.3]{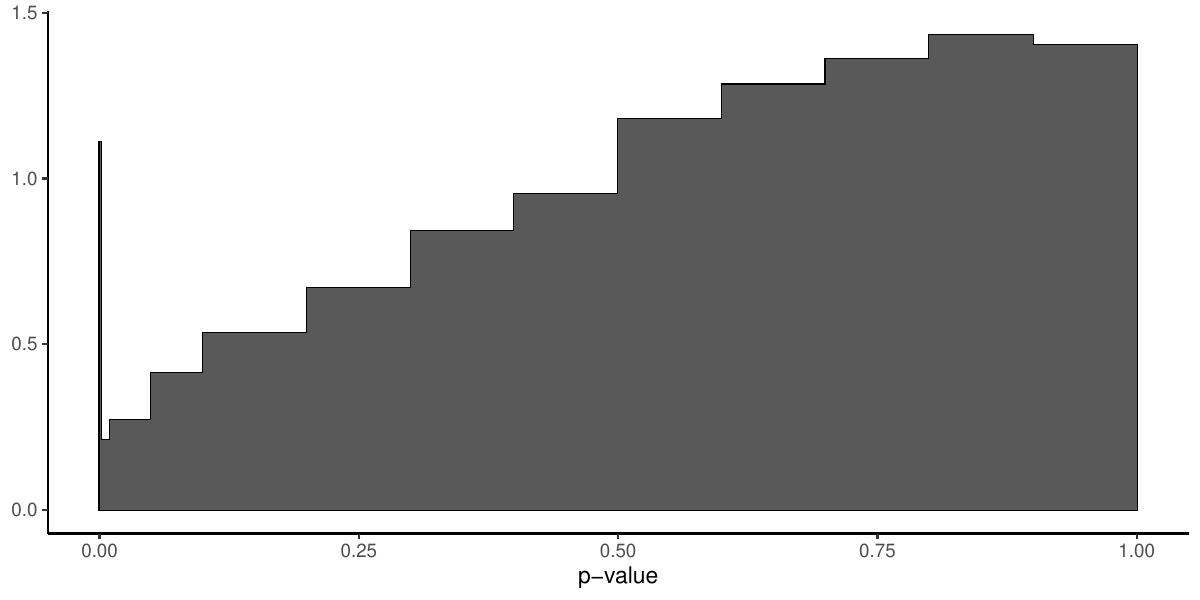} &
    \includegraphics[scale=0.3]{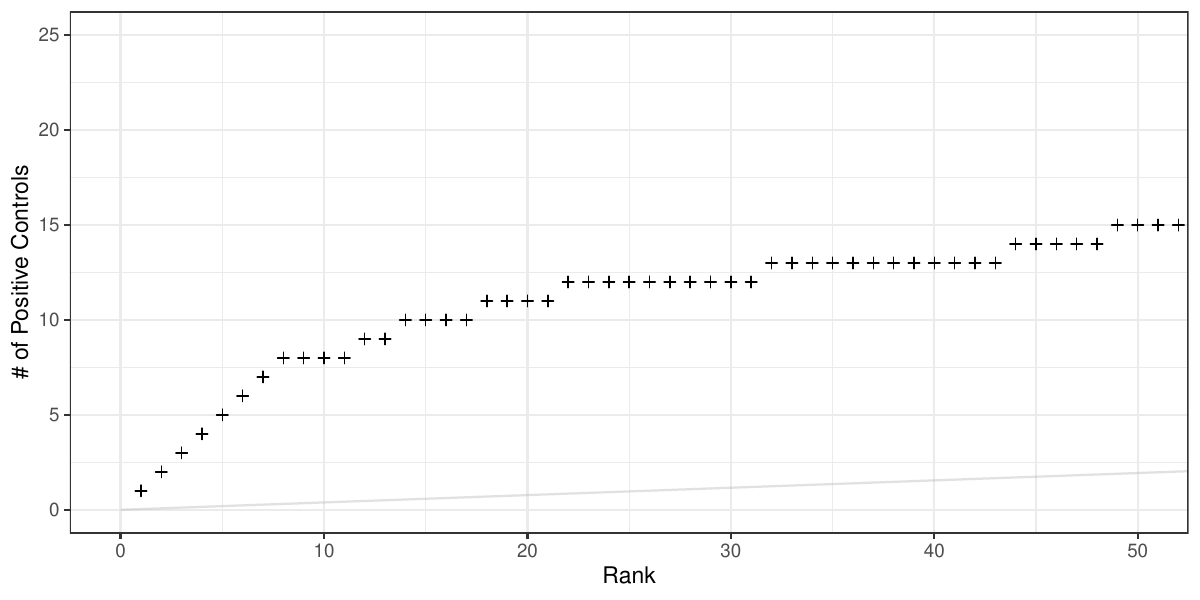} \\    
    \multicolumn{2}{c}{``Technical'' Adjustment}\\
    \includegraphics[scale=0.3]{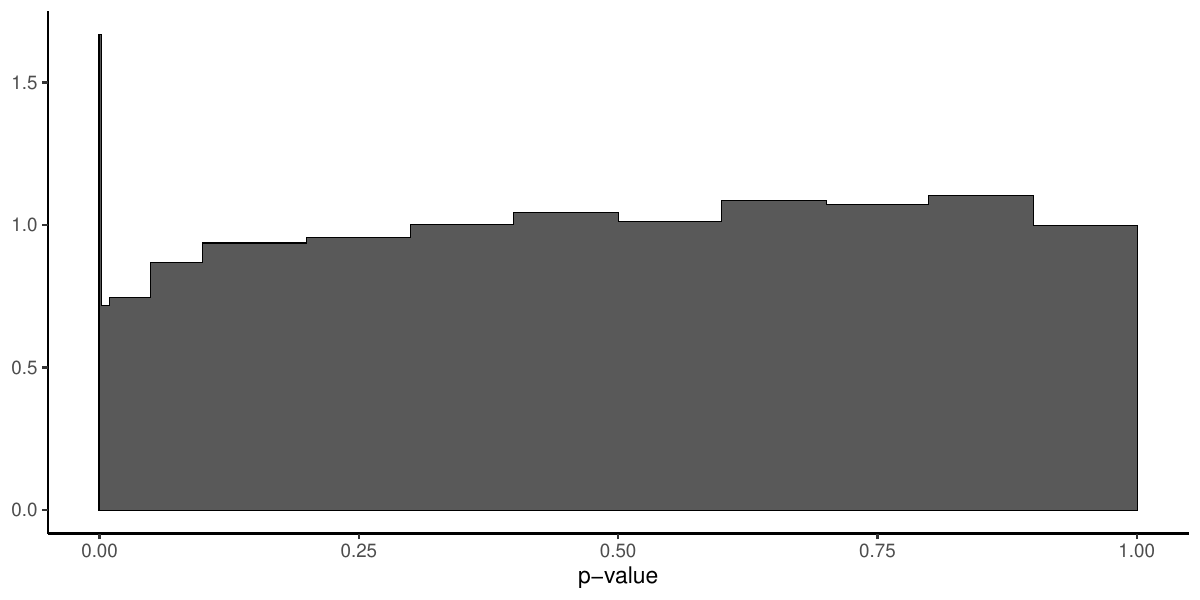} &
    \includegraphics[scale=0.3]{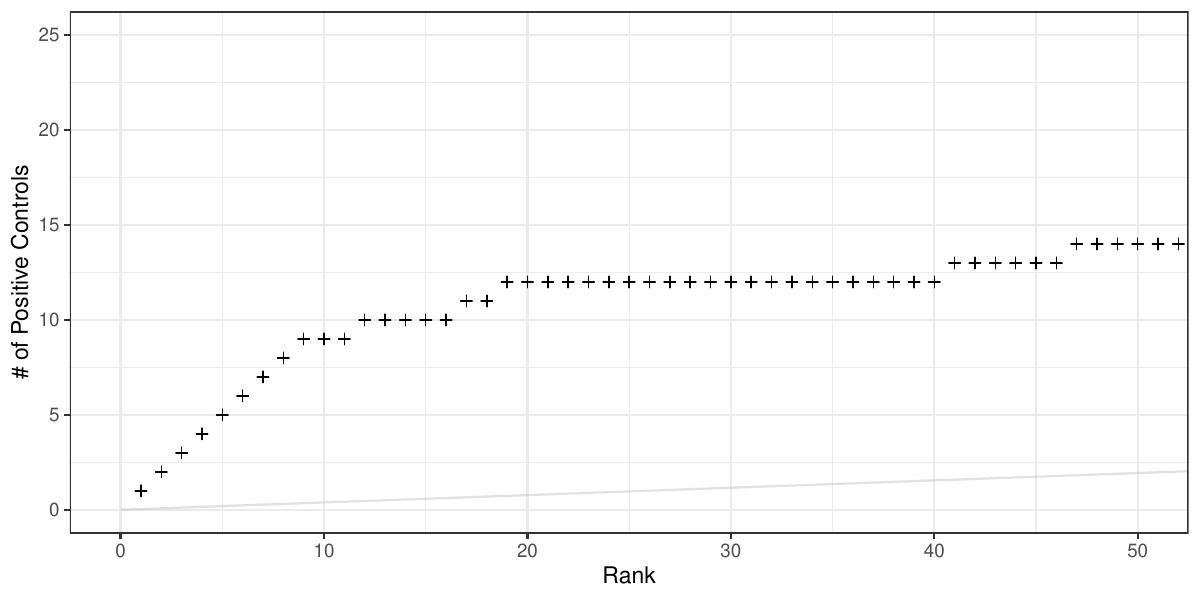} \\    
    \multicolumn{2}{c}{``Bio'' Adjustment}\\
    \includegraphics[scale=0.3]{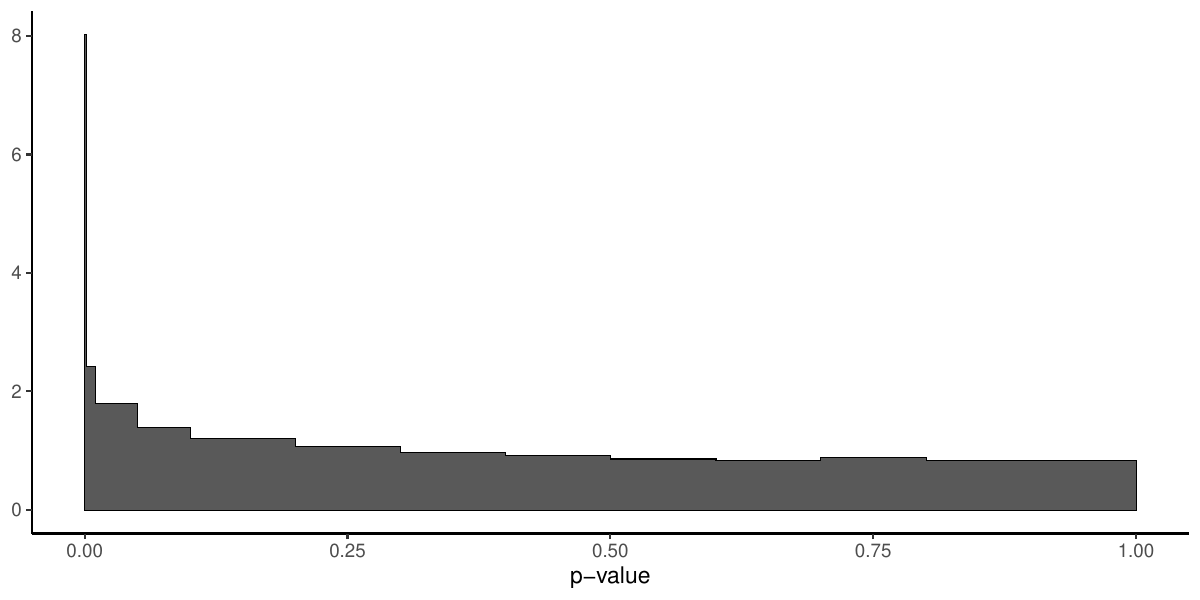} &
    \includegraphics[scale=0.3]{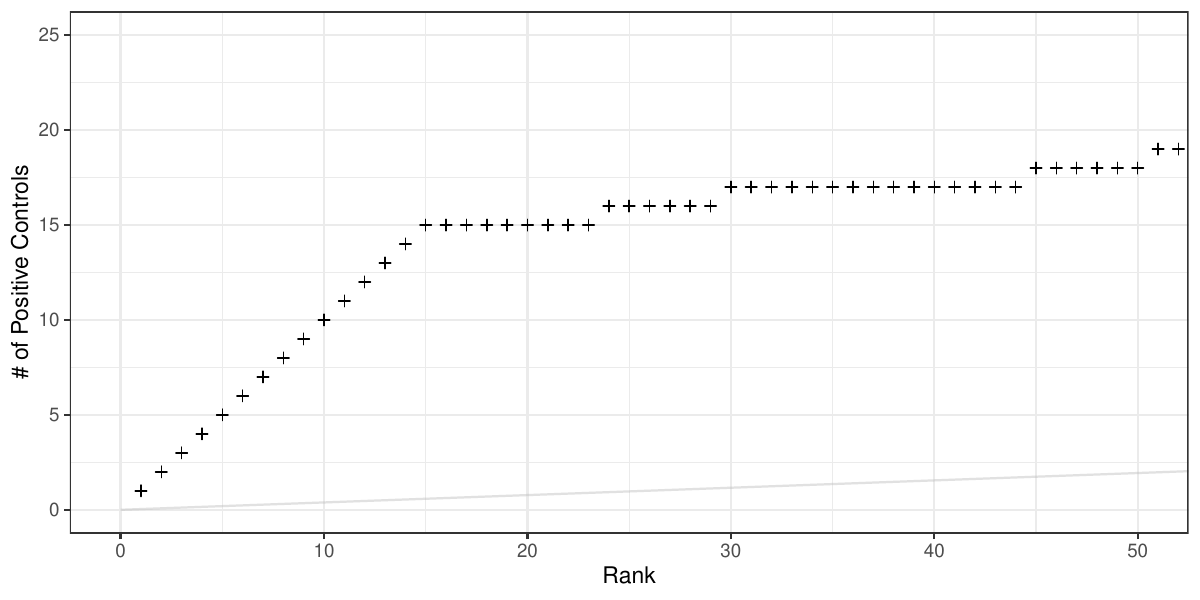} \\    
    \end{tabular}
    \caption{Like Figure 5, but with $K=5$.}
    \label{fig:K5_pvalrank.all}
\end{figure}

\begin{figure}[h]
    \centering
    \begin{tabular}{cc}
    \multicolumn{2}{c}{No Adjustment}\\
    \includegraphics[scale=0.3]{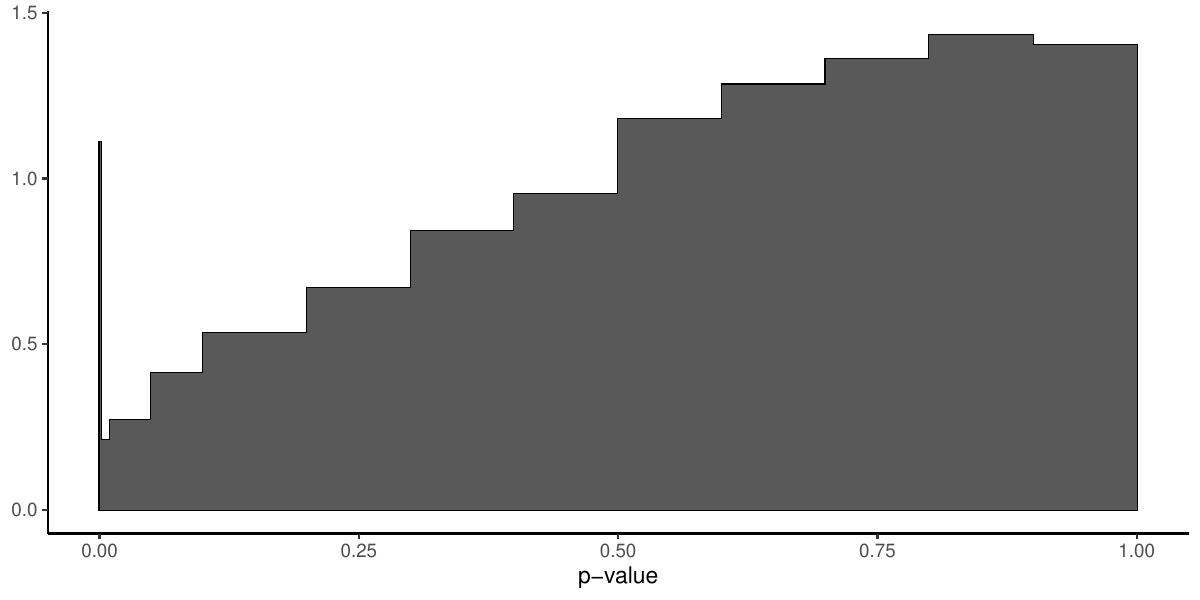} &
    \includegraphics[scale=0.3]{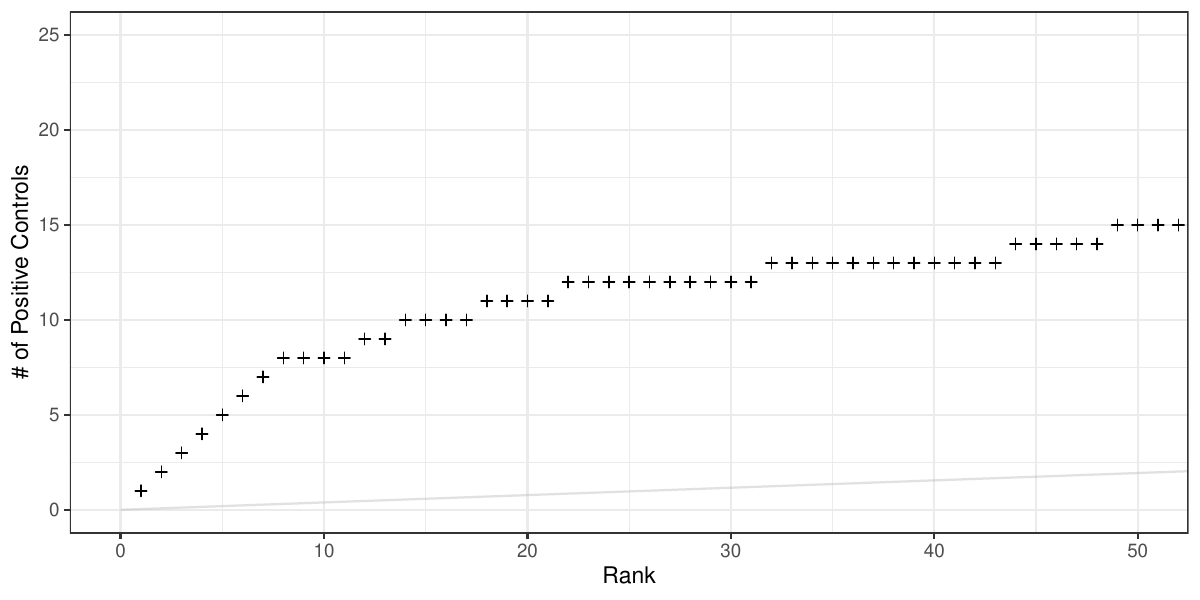} \\    
    \multicolumn{2}{c}{``Technical'' Adjustment}\\
    \includegraphics[scale=0.3]{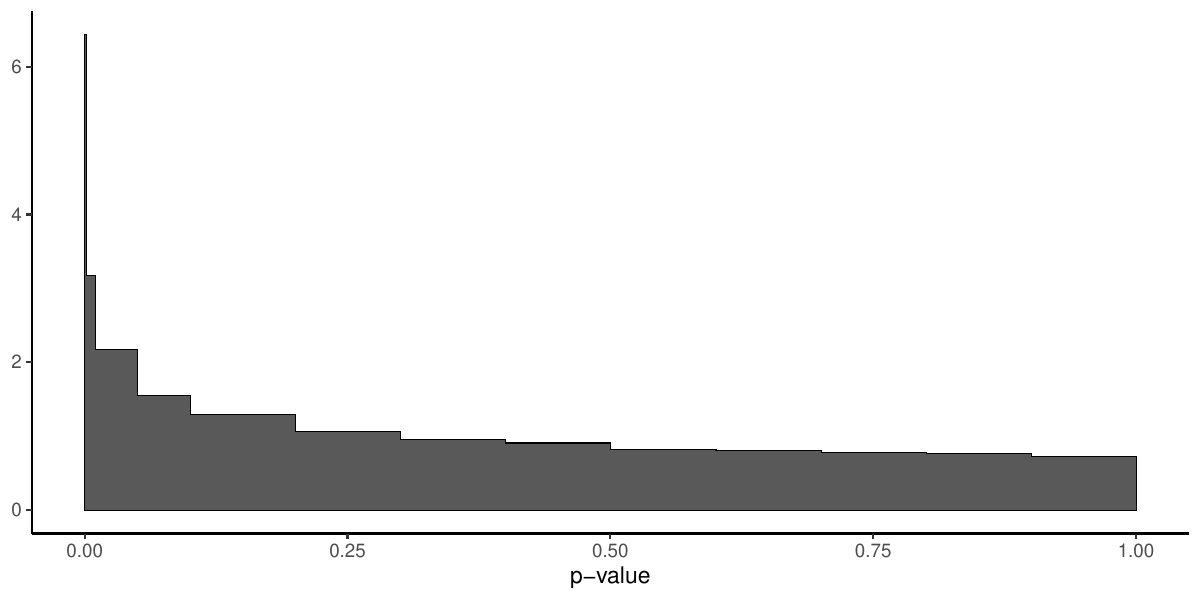} &
    \includegraphics[scale=0.3]{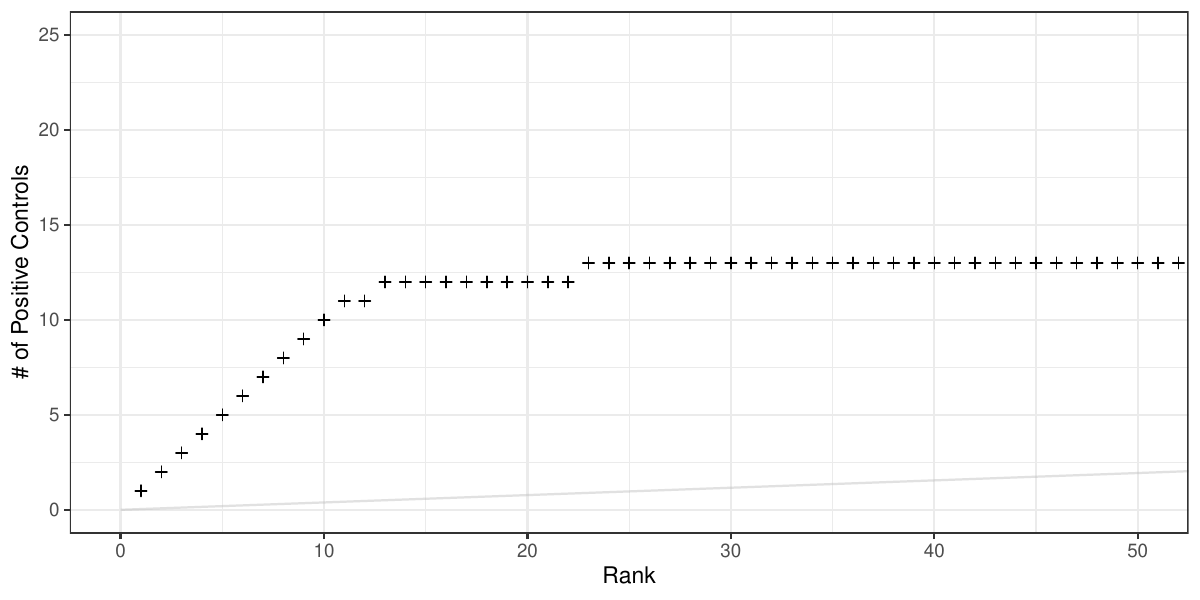} \\    
    \multicolumn{2}{c}{``Bio'' Adjustment}\\
    \includegraphics[scale=0.3]{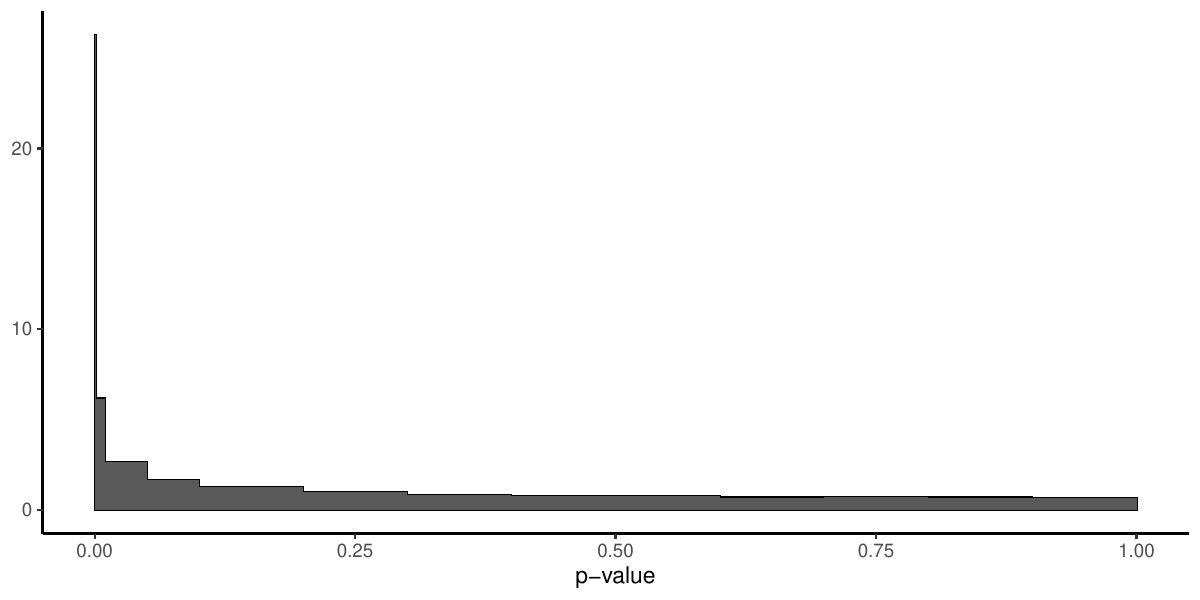} &
    \includegraphics[scale=0.3]{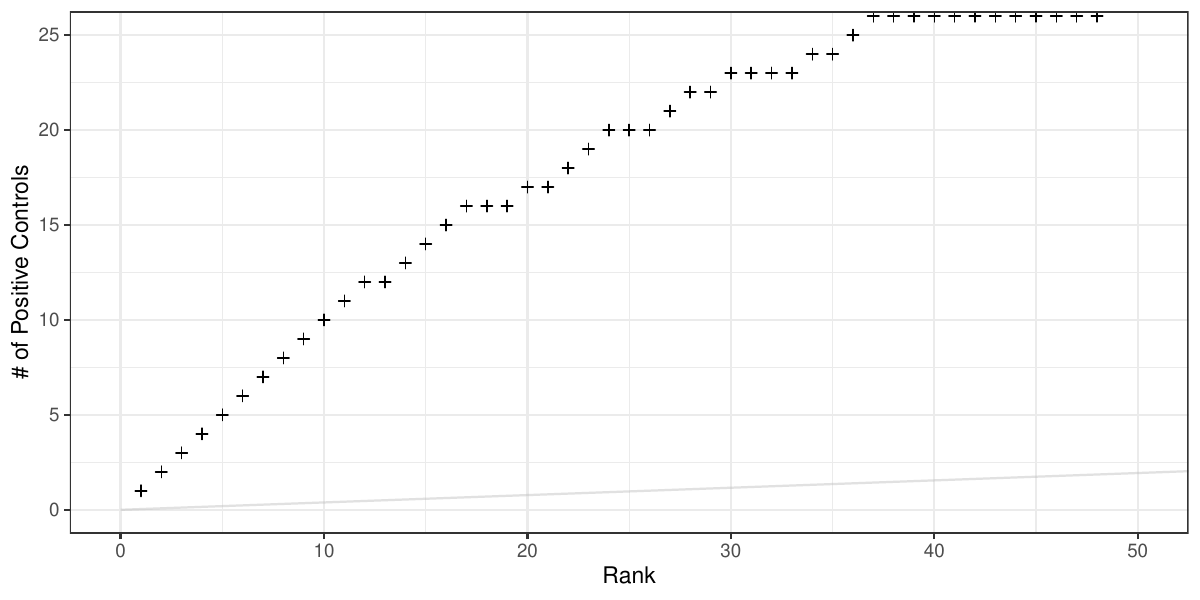} \\    
    \end{tabular}
    \caption{Like Figure 5, but with $K=15$.}
    \label{fig:K15_pvalrank.all}
\end{figure}

\begin{figure}[h]
    \centering
    \begin{tabular}{cc}
    \multicolumn{2}{c}{No Adjustment}\\
    \includegraphics[scale=0.3]{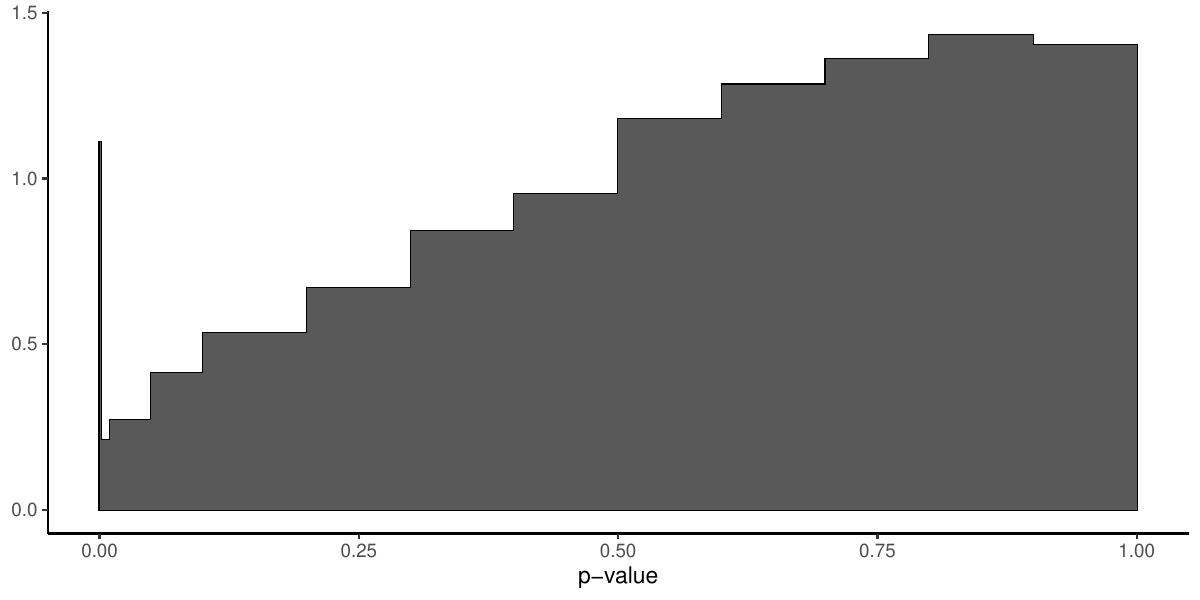} &
    \includegraphics[scale=0.3]{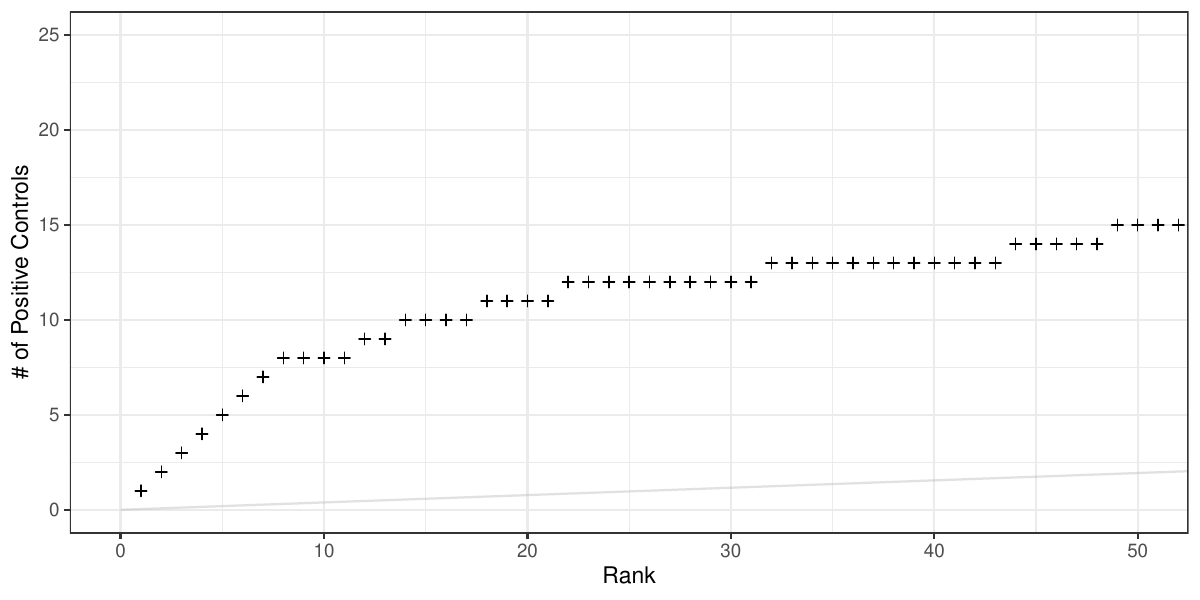} \\    
    \multicolumn{2}{c}{``Technical'' Adjustment}\\
    \includegraphics[scale=0.3]{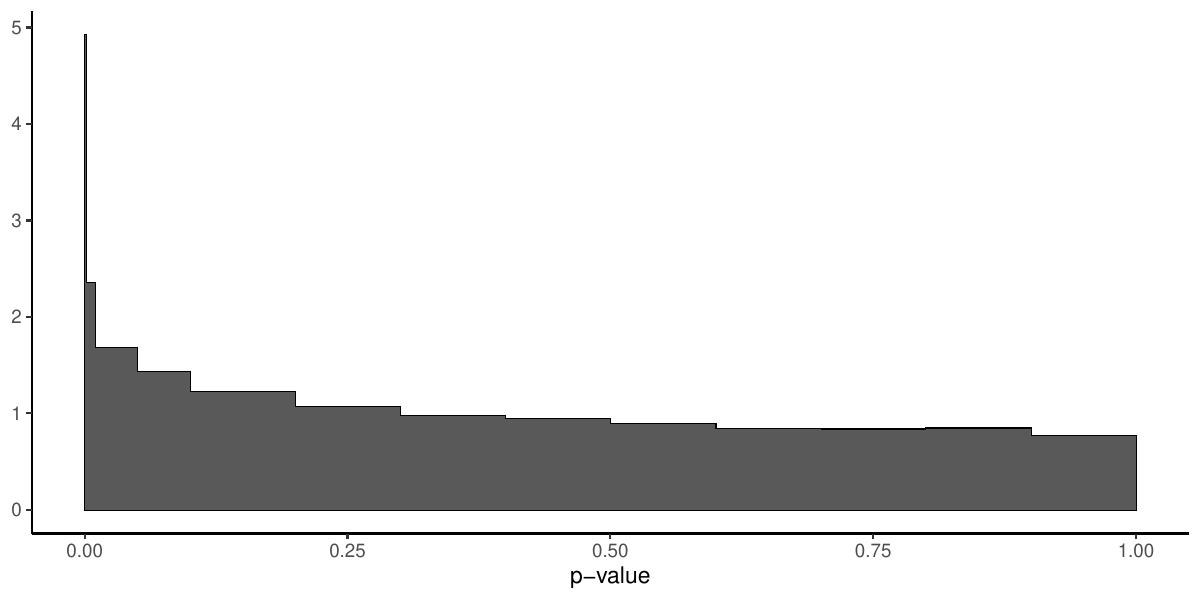} &
    \includegraphics[scale=0.3]{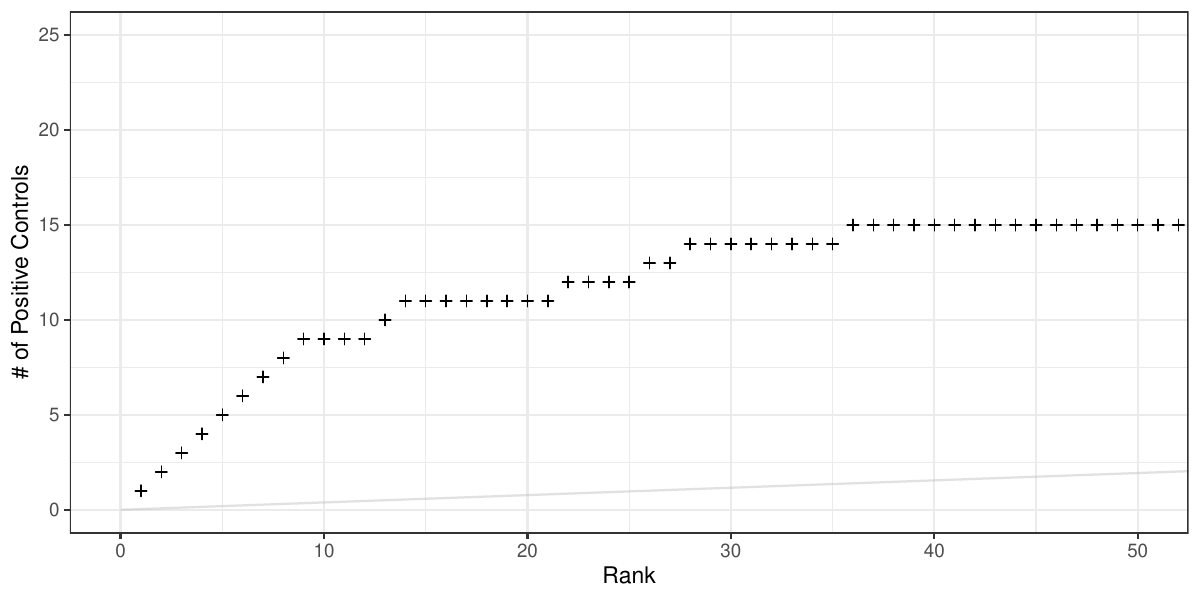} \\    
    \multicolumn{2}{c}{``Bio'' Adjustment}\\
    \includegraphics[scale=0.3]{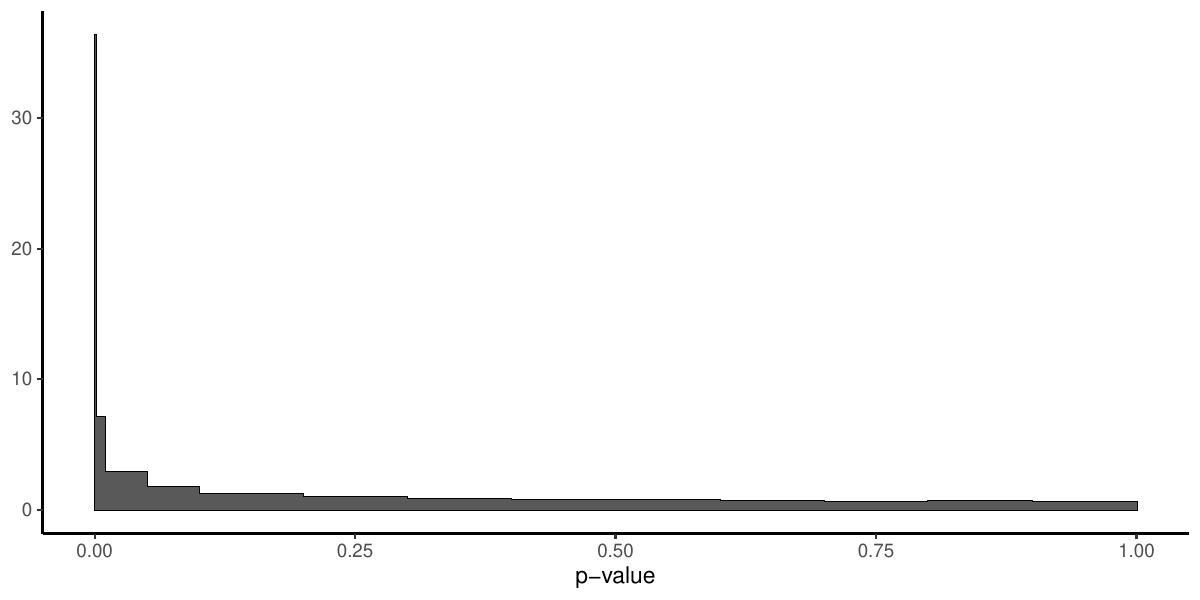} &
    \includegraphics[scale=0.3]{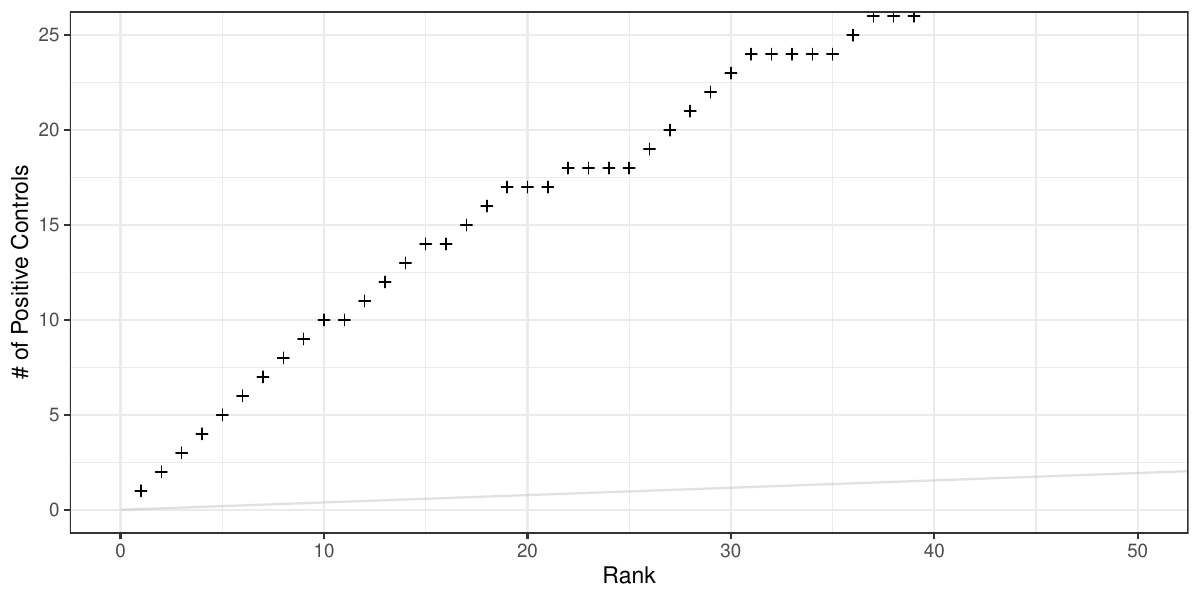} \\    
    \end{tabular}
    \caption{Like Figure 5, but with $K=20$.}
    \label{fig:K20_pvalrank.all}
\end{figure}

\begin{figure}[h]
    \centering
    \begin{tabular}{cc}
    \multicolumn{2}{c}{No Adjustment}\\
    \includegraphics[scale=0.3]{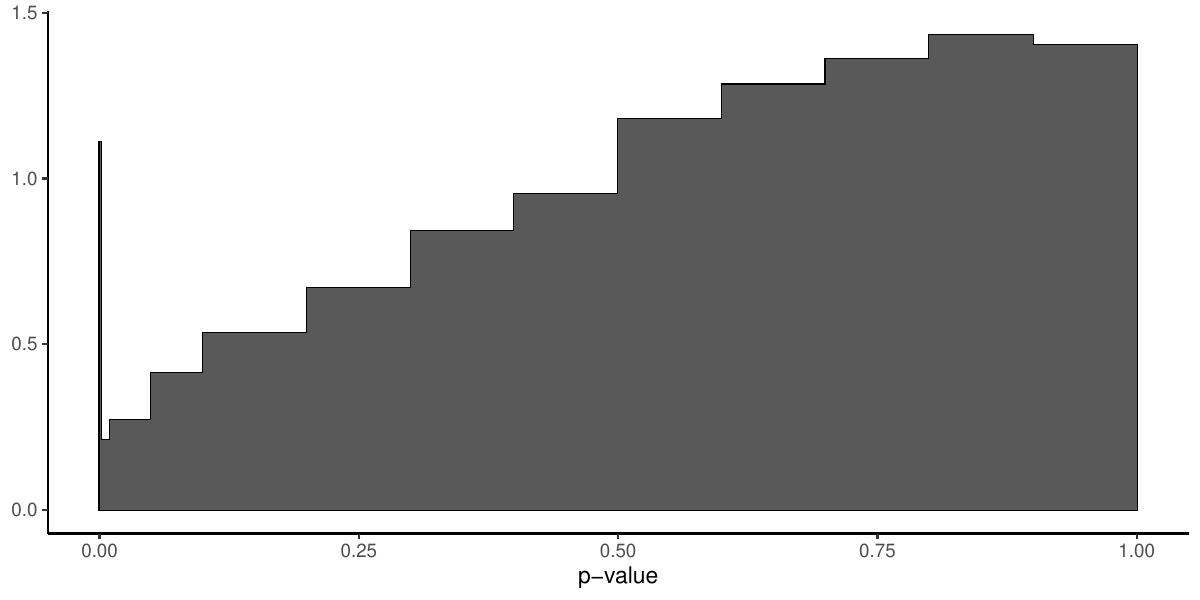} &
    \includegraphics[scale=0.3]{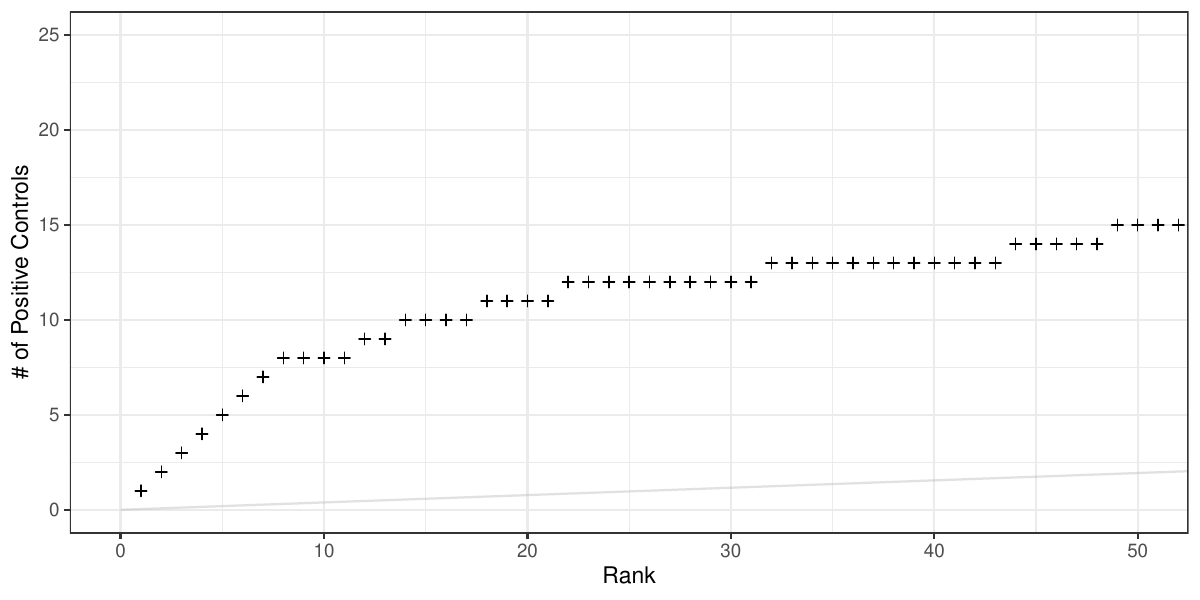} \\    
    \multicolumn{2}{c}{``Technical'' Adjustment}\\
    \includegraphics[scale=0.3]{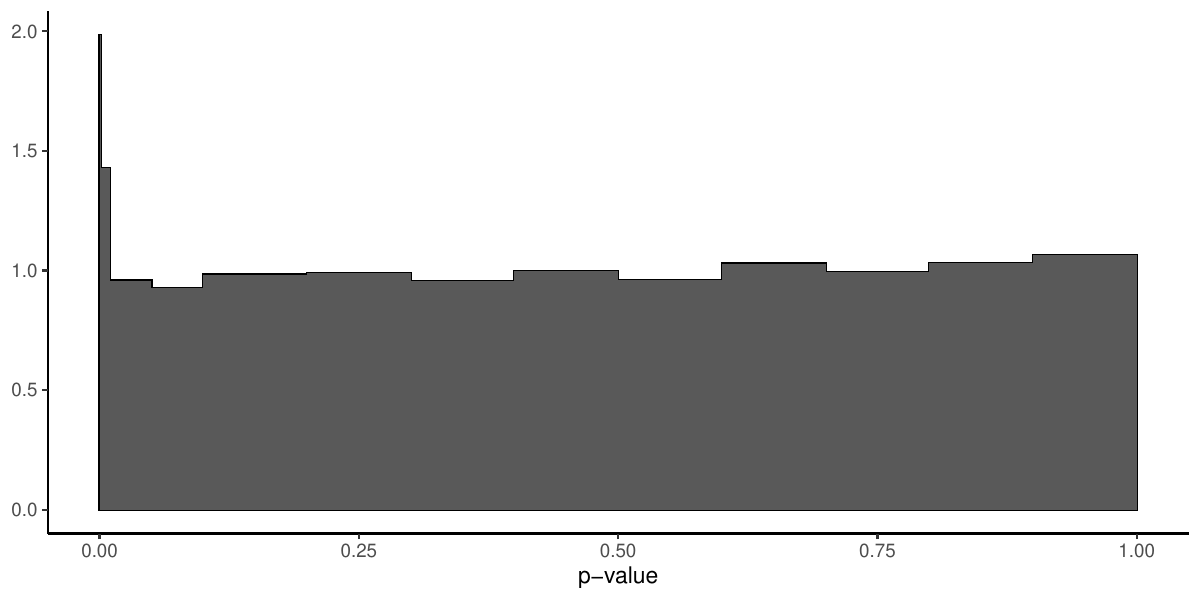} &
    \includegraphics[scale=0.3]{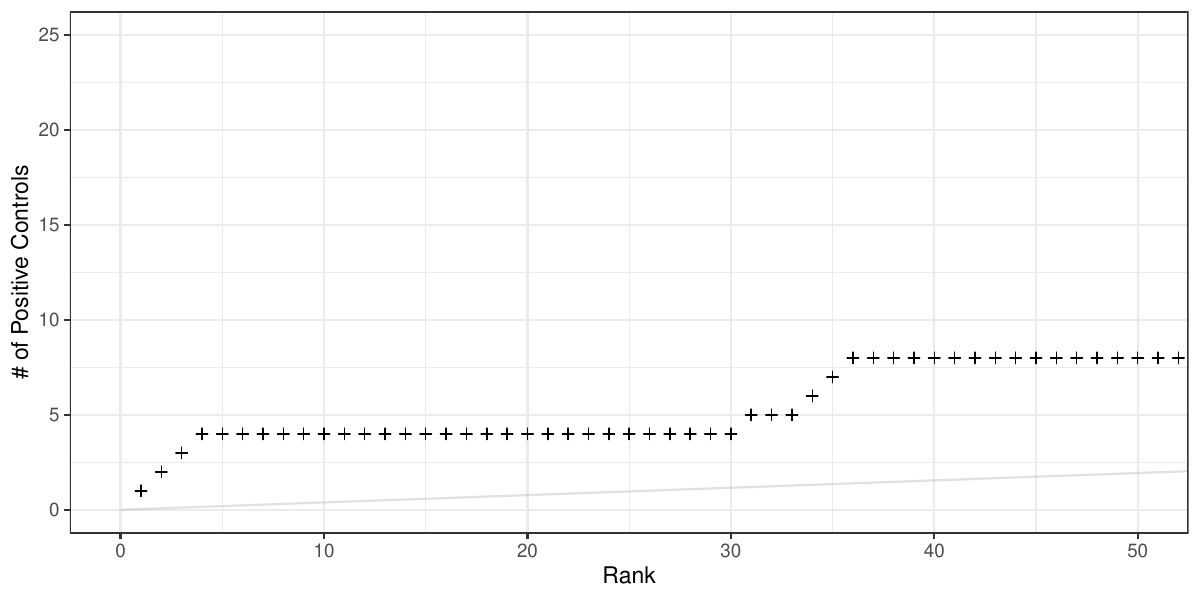} \\    
    \multicolumn{2}{c}{``Bio'' Adjustment}\\
    \includegraphics[scale=0.3]{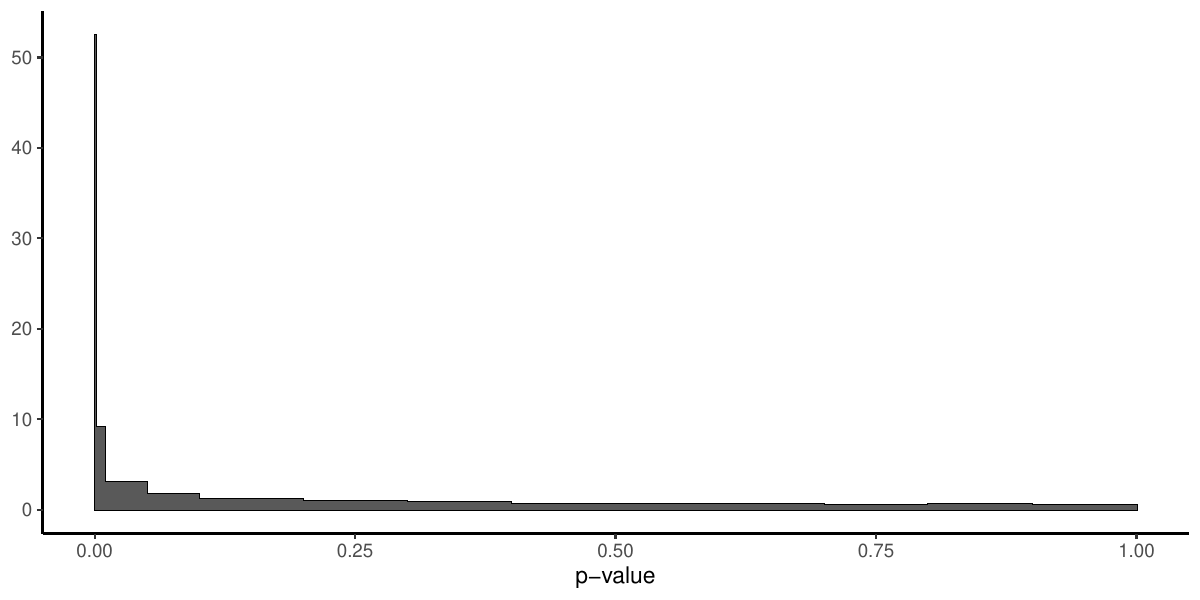} &
    \includegraphics[scale=0.3]{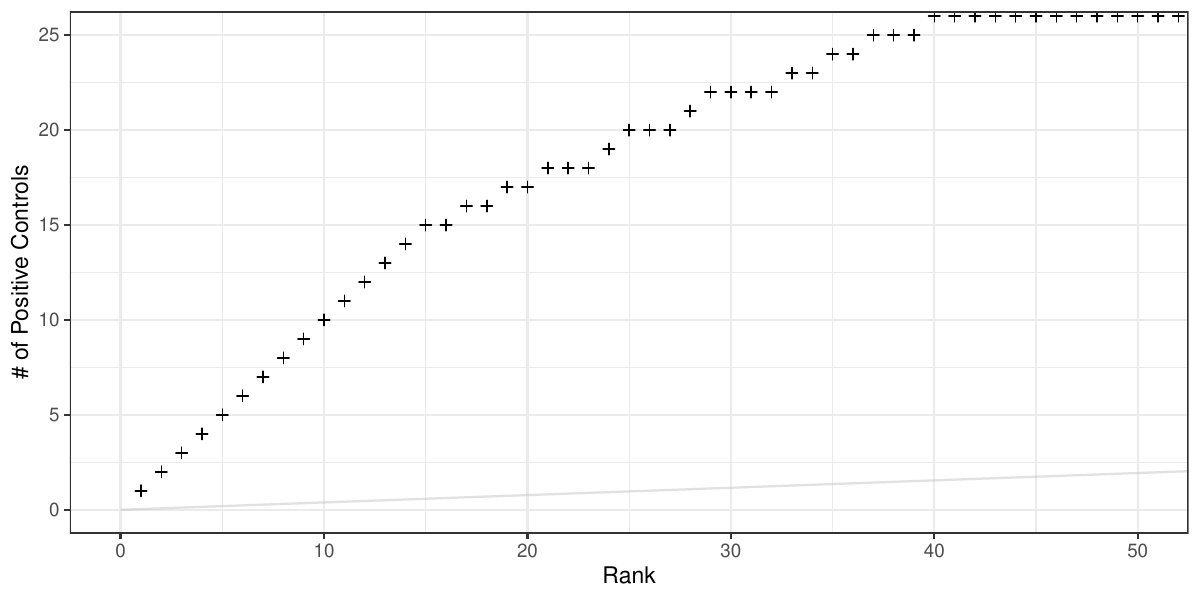} \\    
    \end{tabular}
    \caption{Like Figure 5, but with $K=30$.}
    \label{fig:K30_pvalrank.all}
\end{figure}

\begin{figure}[h]
    \centering
    \begin{tabular}{cc}
    \hspace{-1cm} (\textbf{a}) & 
    \hspace{-1cm} (\textbf{b}) \\
    \includegraphics[scale=0.3]{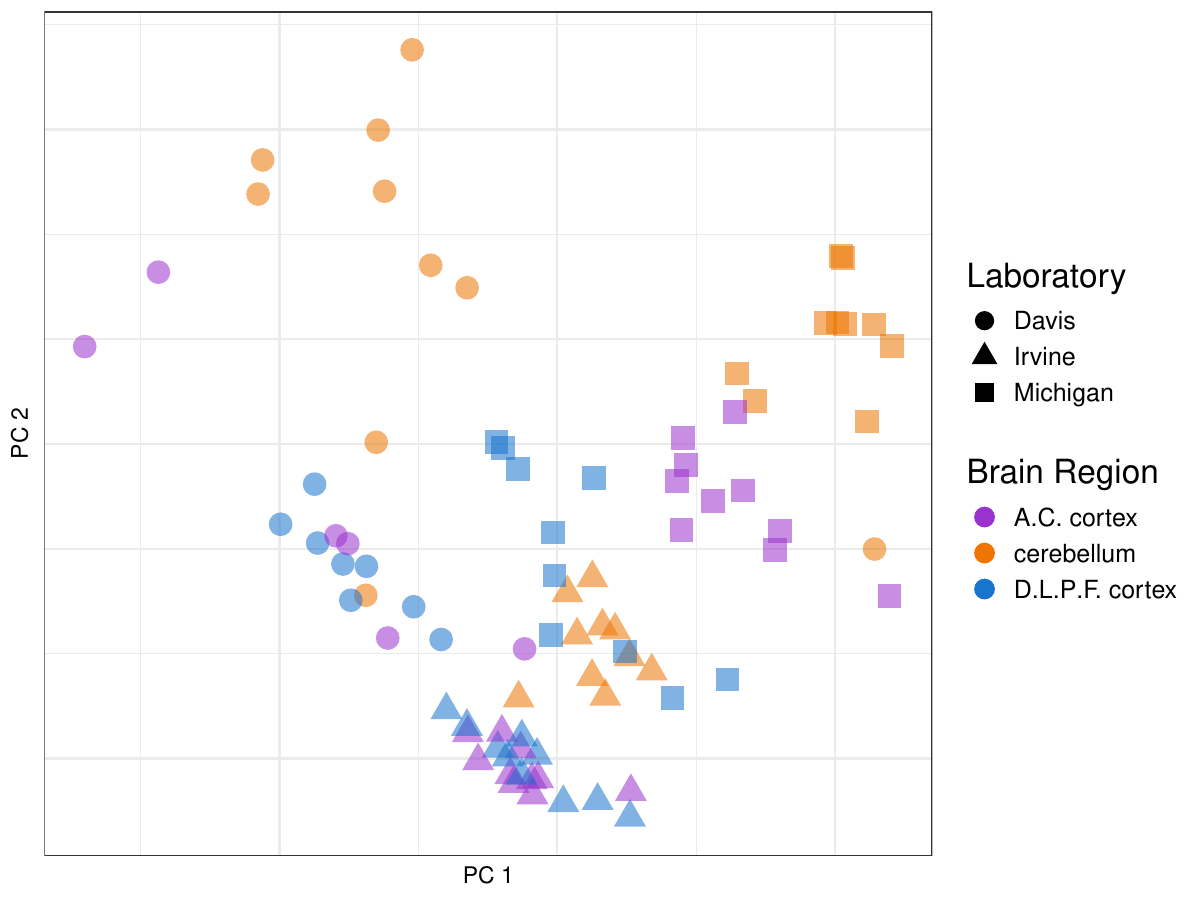} &
    \includegraphics[scale=0.3]{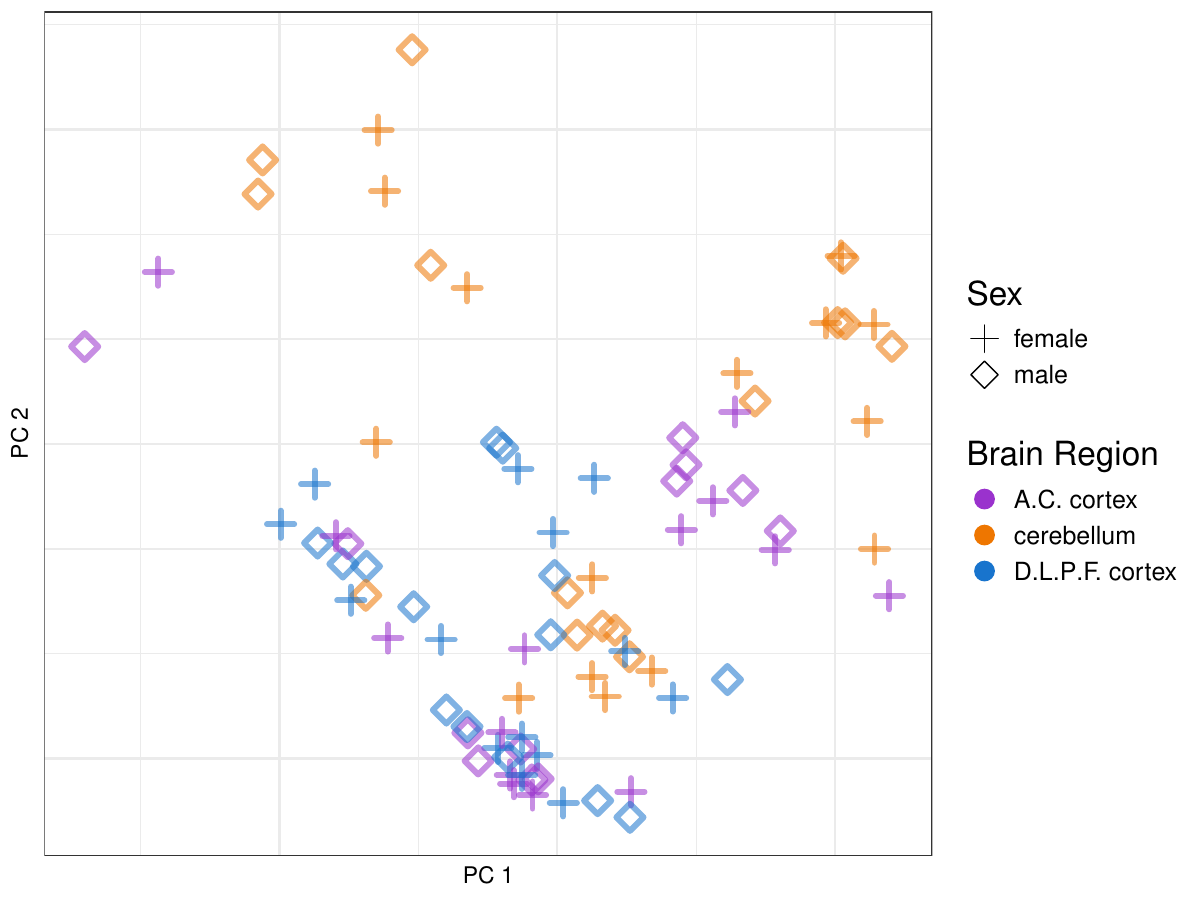} \\    
    \hspace{-1cm} (\textbf{c}) & 
    \hspace{-1cm} (\textbf{d}) \\
    \includegraphics[scale=0.3]{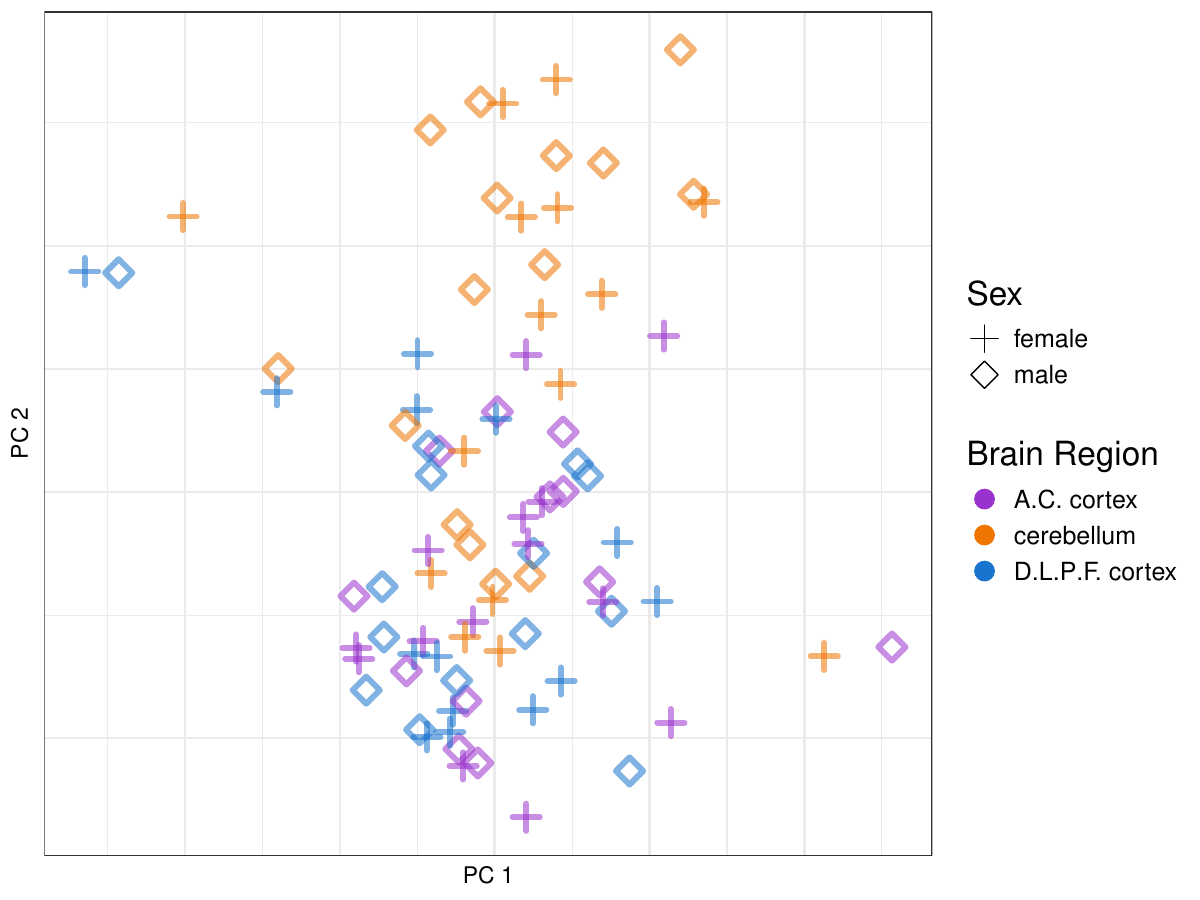}  &
    \includegraphics[scale=0.3]{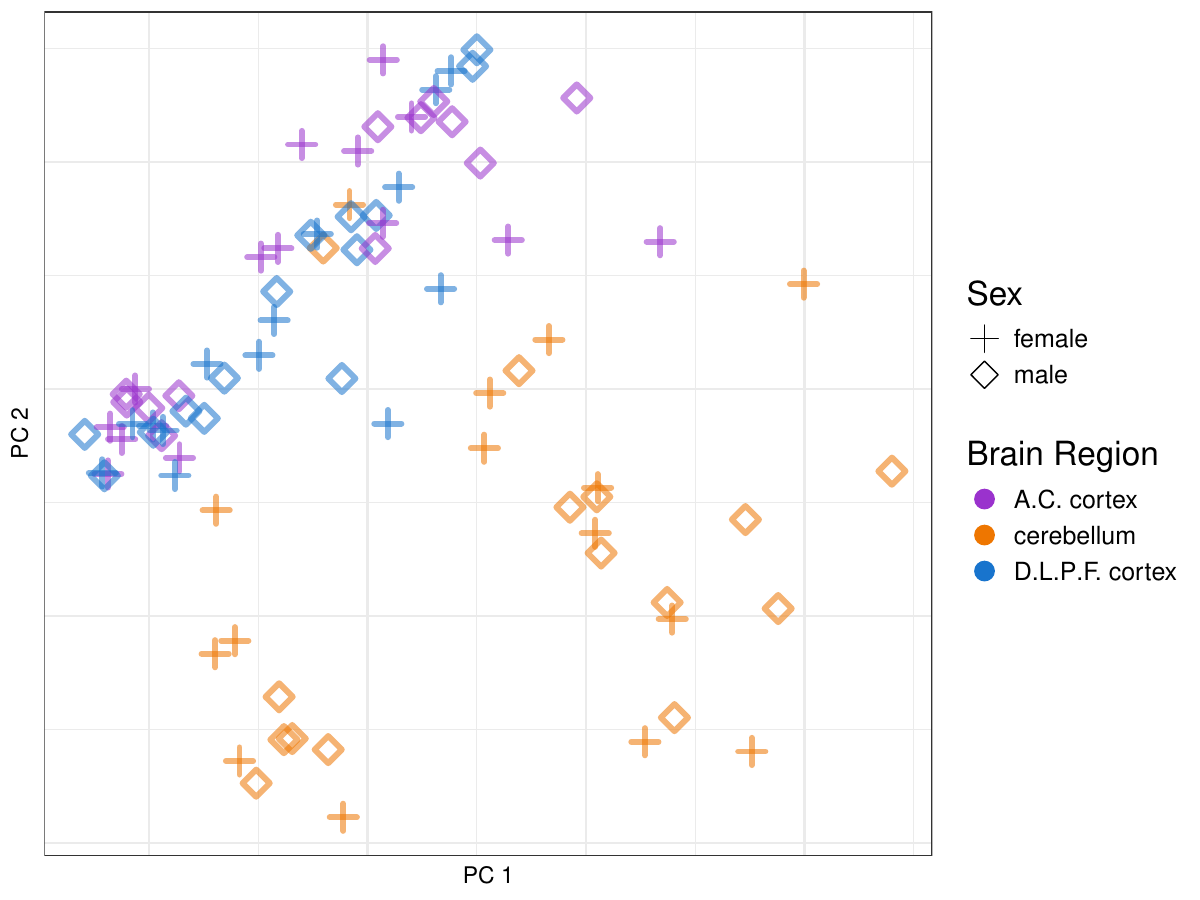}     
    \end{tabular}
    \caption{Like Figure 6, but with $K=1$.}
    \label{fig:K1_svd.xy}
\end{figure}

\begin{figure}[h]
    \centering
    \begin{tabular}{cc}
    \hspace{-1cm} (\textbf{a}) & 
    \hspace{-1cm} (\textbf{b}) \\
    \includegraphics[scale=0.3]{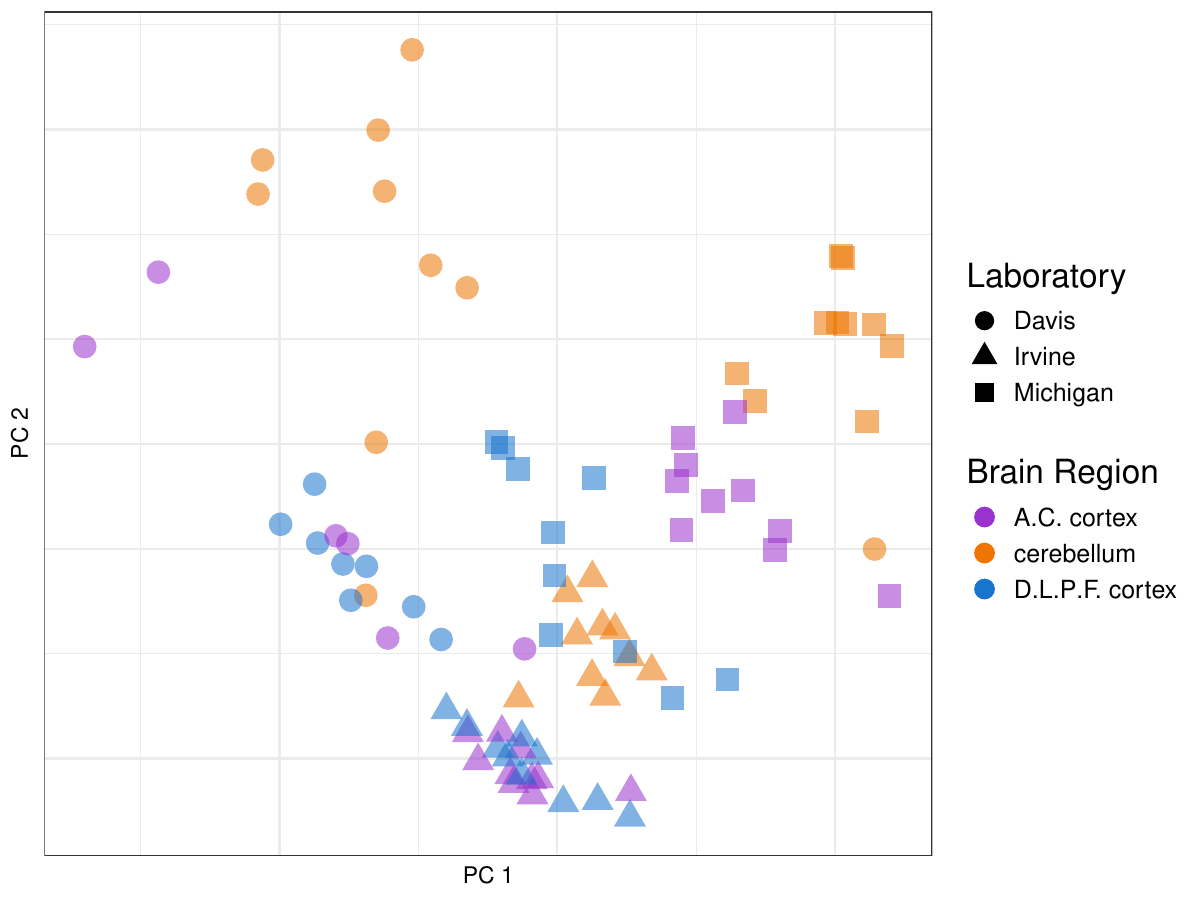} &
    \includegraphics[scale=0.3]{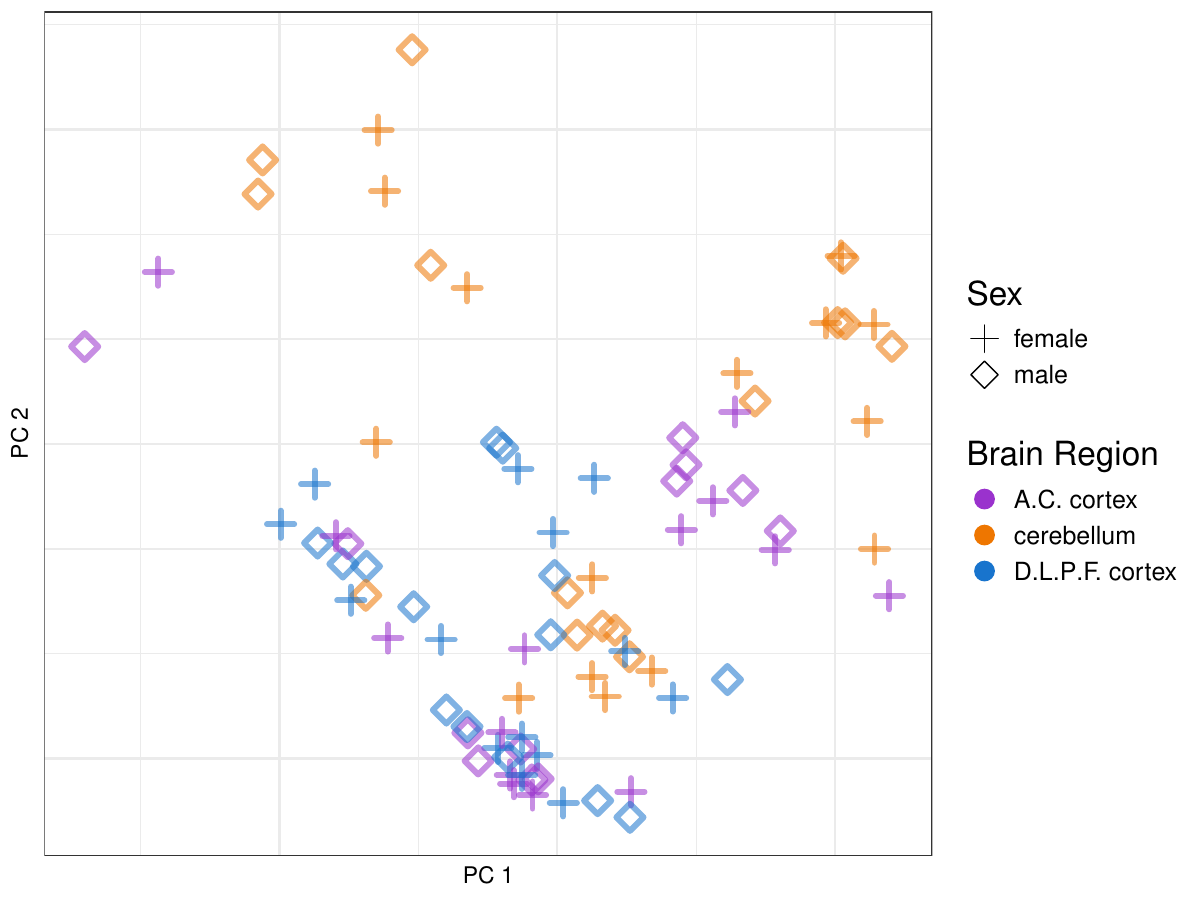} \\    
    \hspace{-1cm} (\textbf{c}) & 
    \hspace{-1cm} (\textbf{d}) \\
    \includegraphics[scale=0.3]{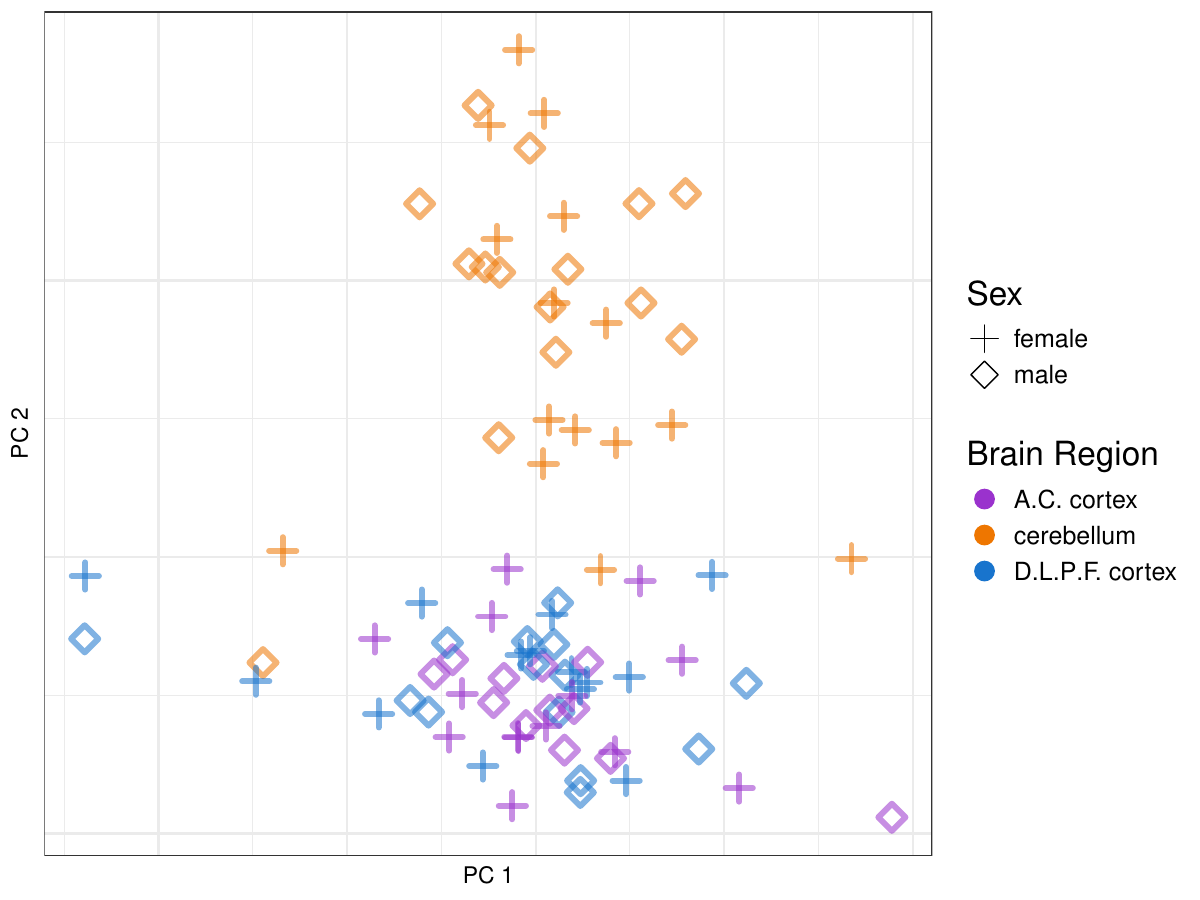}  &
    \includegraphics[scale=0.3]{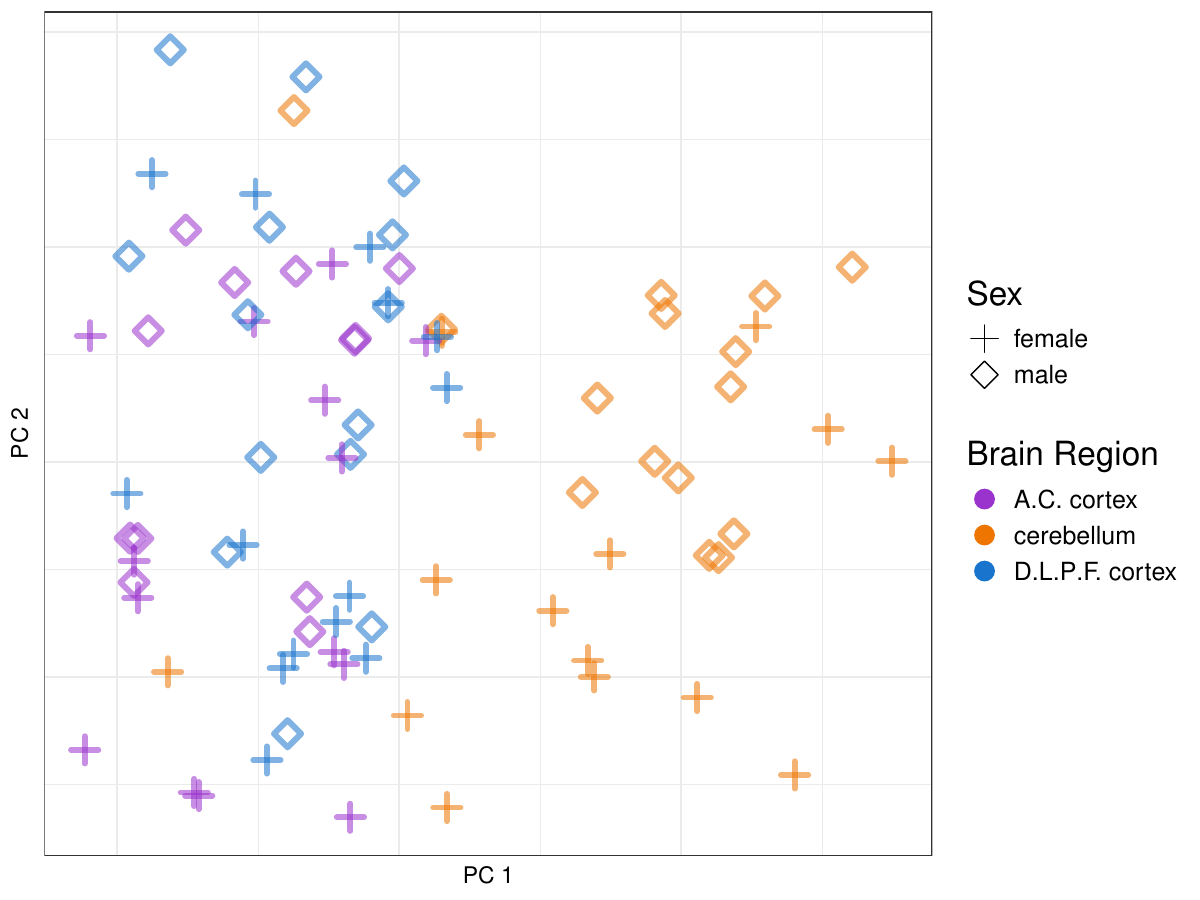}     
    \end{tabular}
    \caption{Like Figure 6, but with $K=2$.}
    \label{fig:K2_svd.xy}
\end{figure}

\begin{figure}[h]
    \centering
    \begin{tabular}{cc}
    \hspace{-1cm} (\textbf{a}) & 
    \hspace{-1cm} (\textbf{b}) \\
    \includegraphics[scale=0.3]{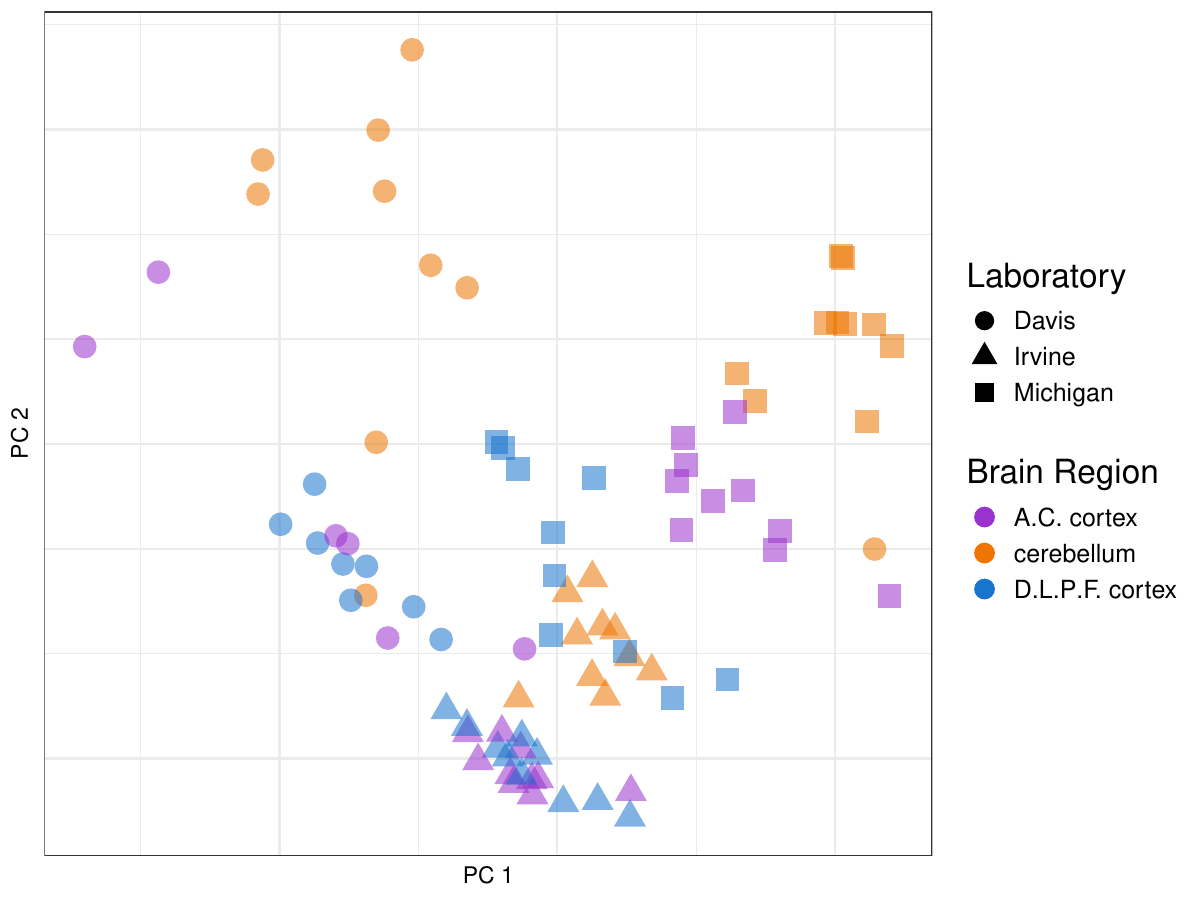} &
    \includegraphics[scale=0.3]{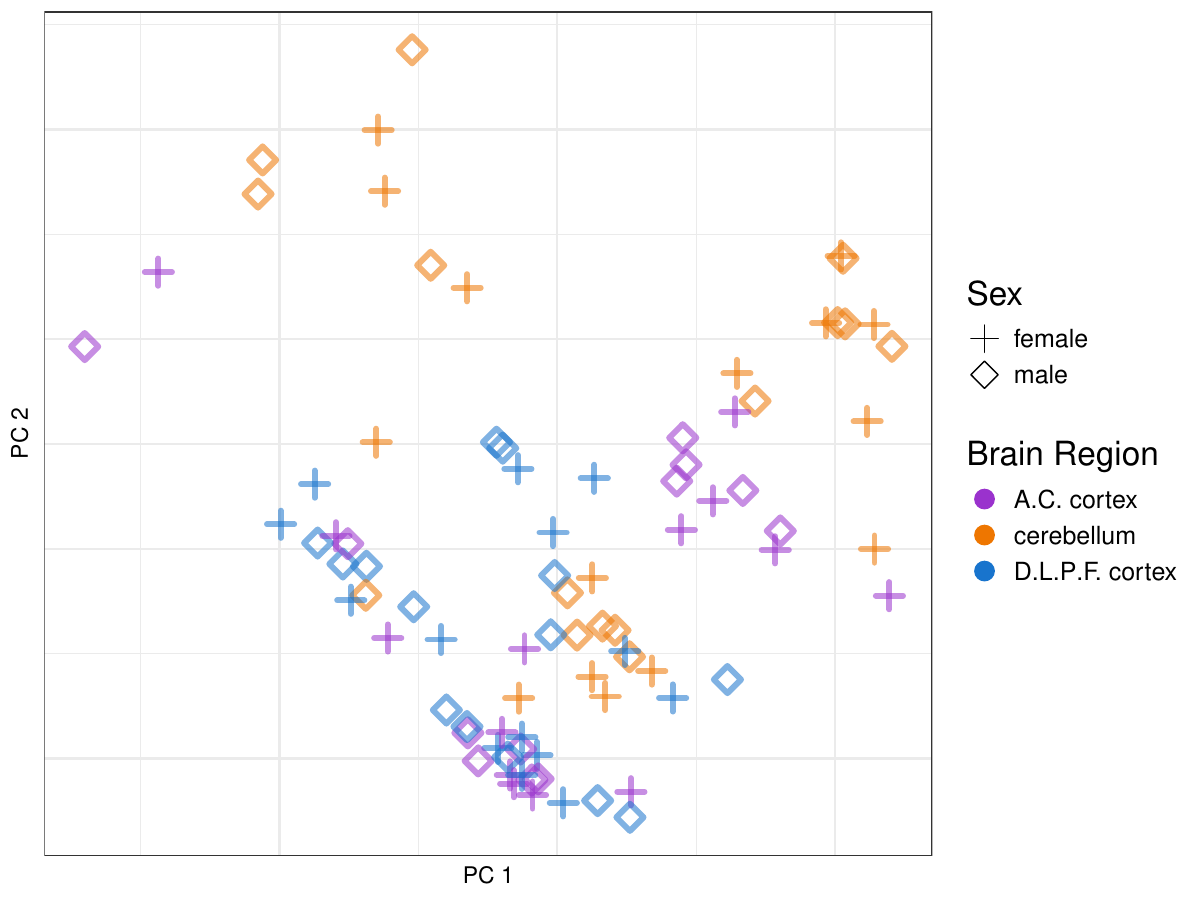} \\    
    \hspace{-1cm} (\textbf{c}) & 
    \hspace{-1cm} (\textbf{d}) \\
    \includegraphics[scale=0.3]{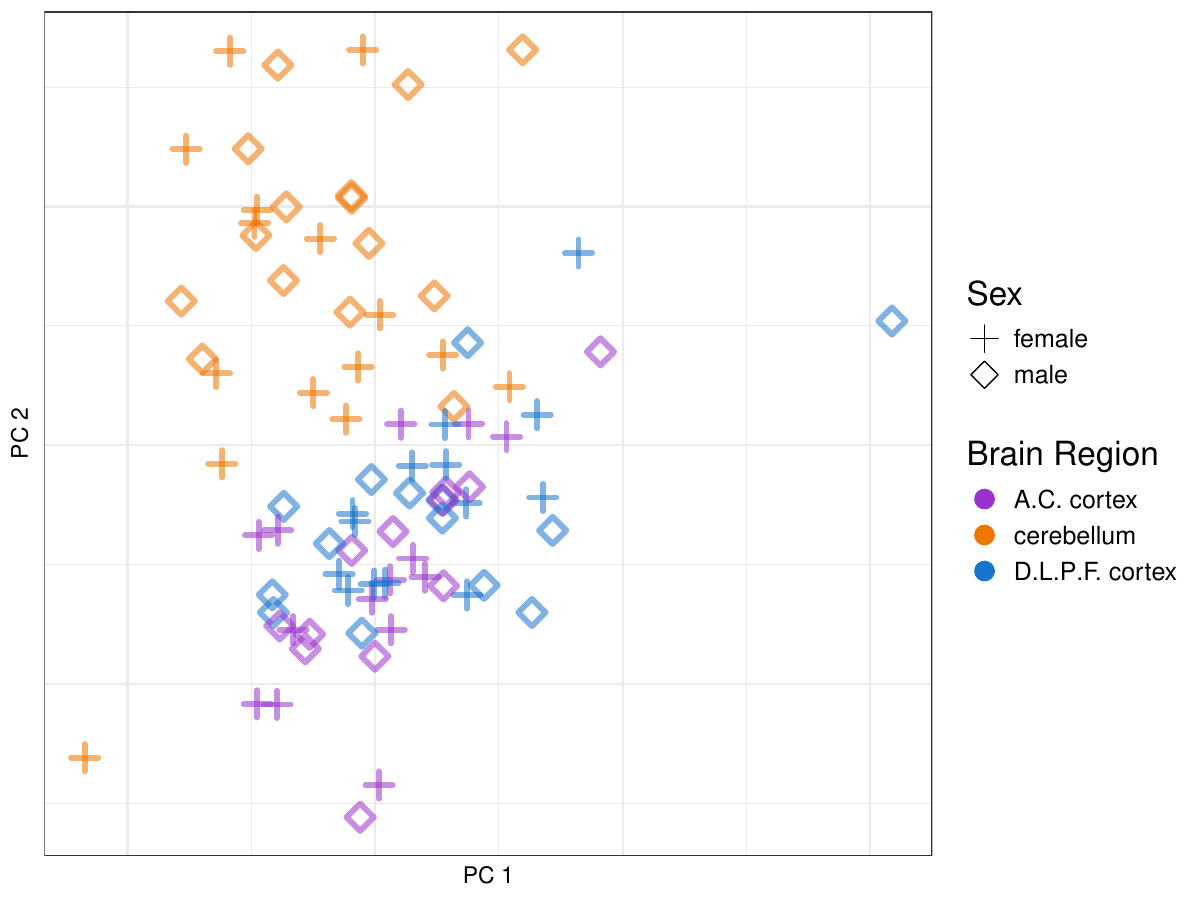}  &
    \includegraphics[scale=0.3]{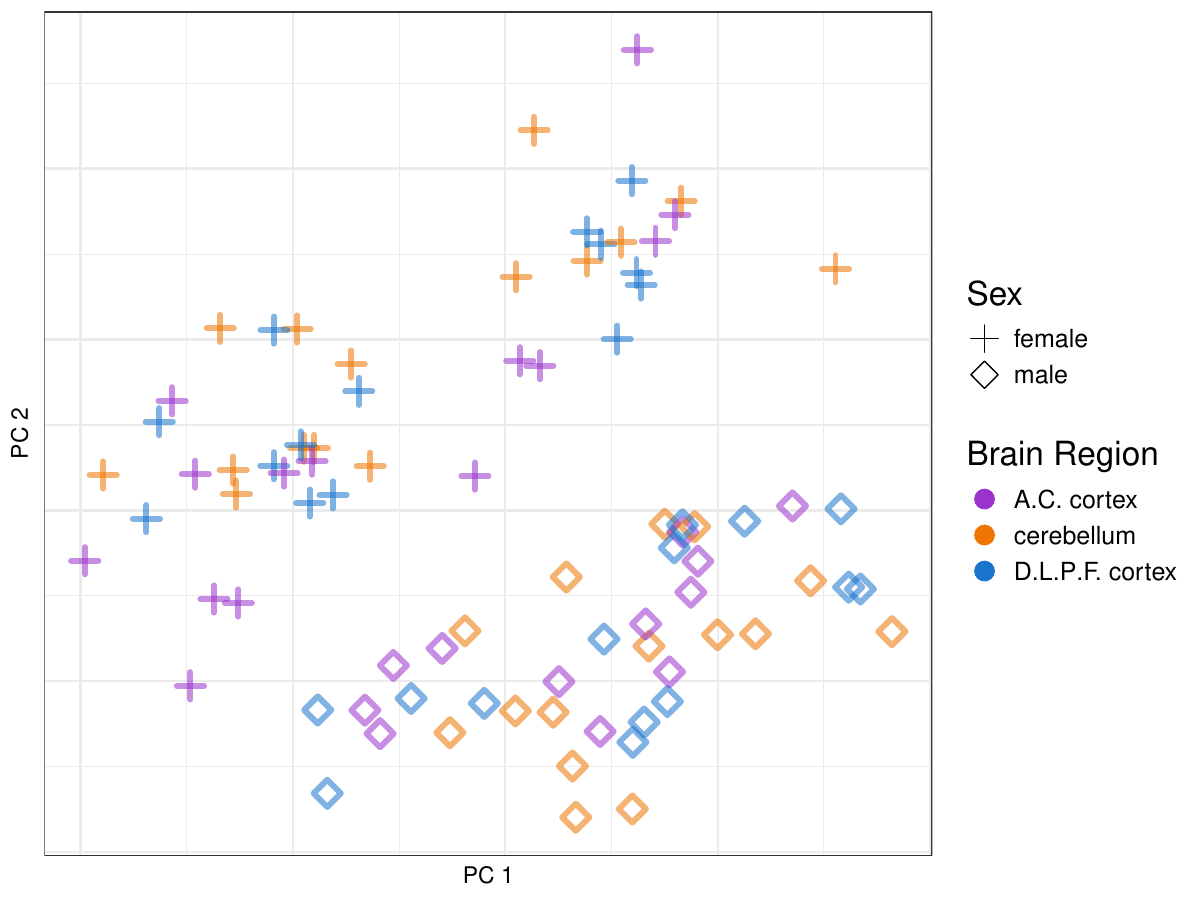}     
    \end{tabular}
    \caption{Like Figure 6, but with $K=5$.}
    \label{fig:K5_svd.xy}
\end{figure}

\begin{figure}[h]
    \centering
    \begin{tabular}{cc}
    \hspace{-1cm} (\textbf{a}) & 
    \hspace{-1cm} (\textbf{b}) \\
    \includegraphics[scale=0.3]{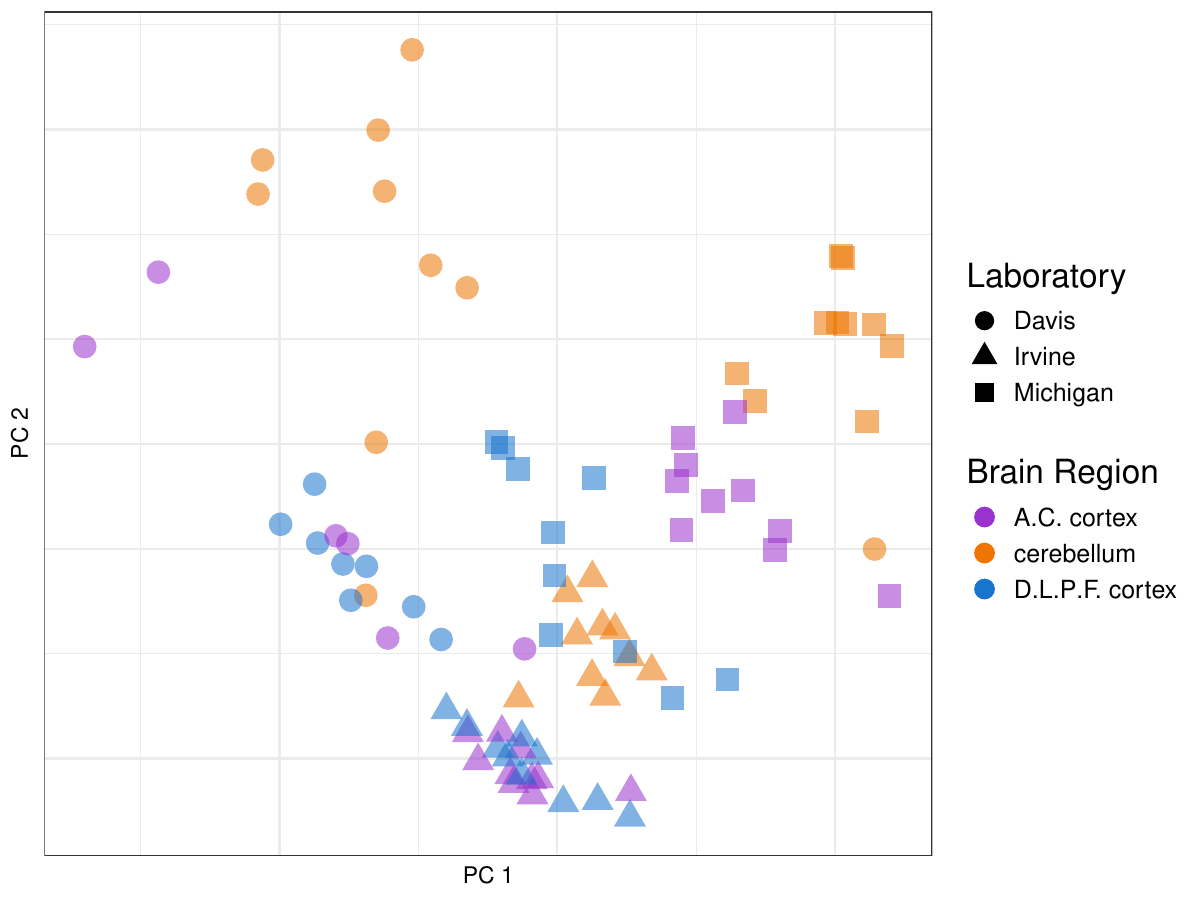} &
    \includegraphics[scale=0.3]{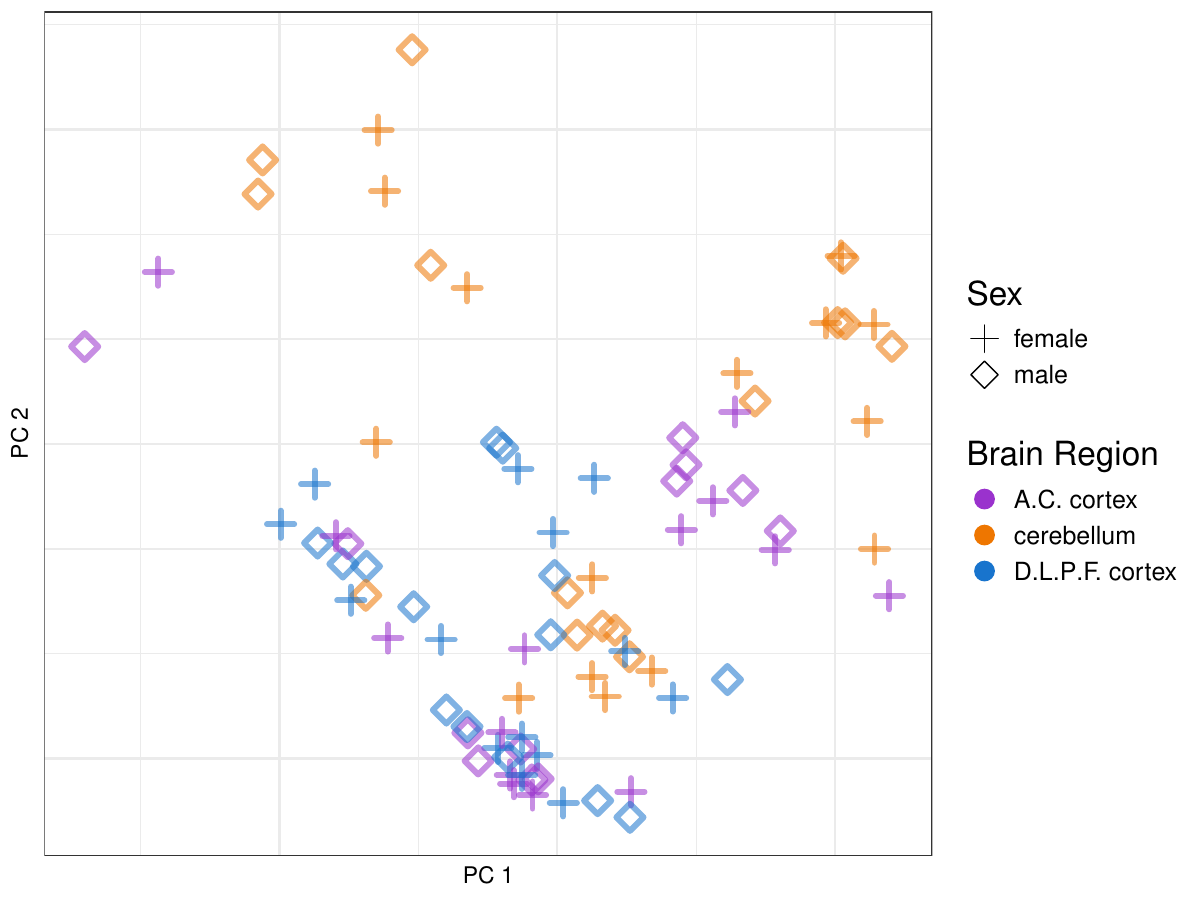} \\    
    \hspace{-1cm} (\textbf{c}) & 
    \hspace{-1cm} (\textbf{d}) \\
    \includegraphics[scale=0.3]{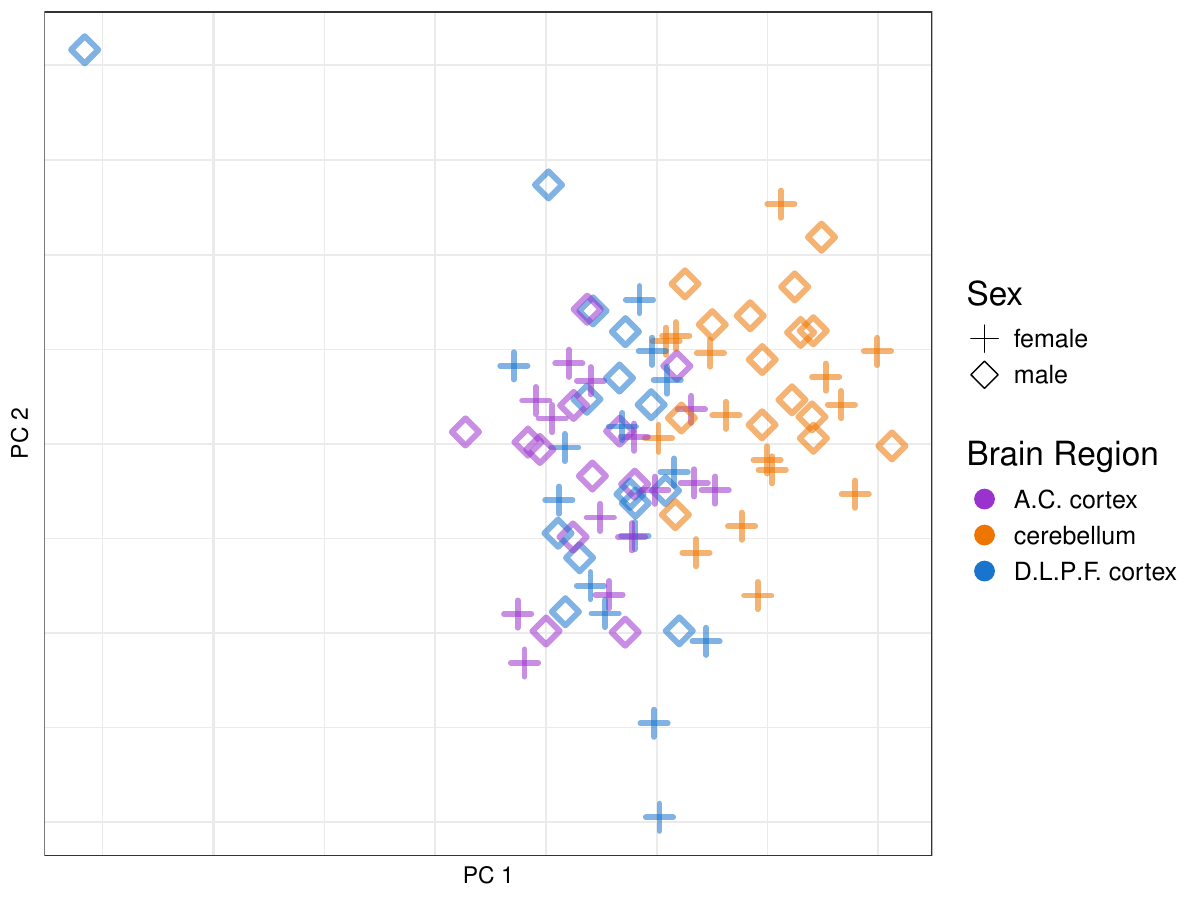}  &
    \includegraphics[scale=0.3]{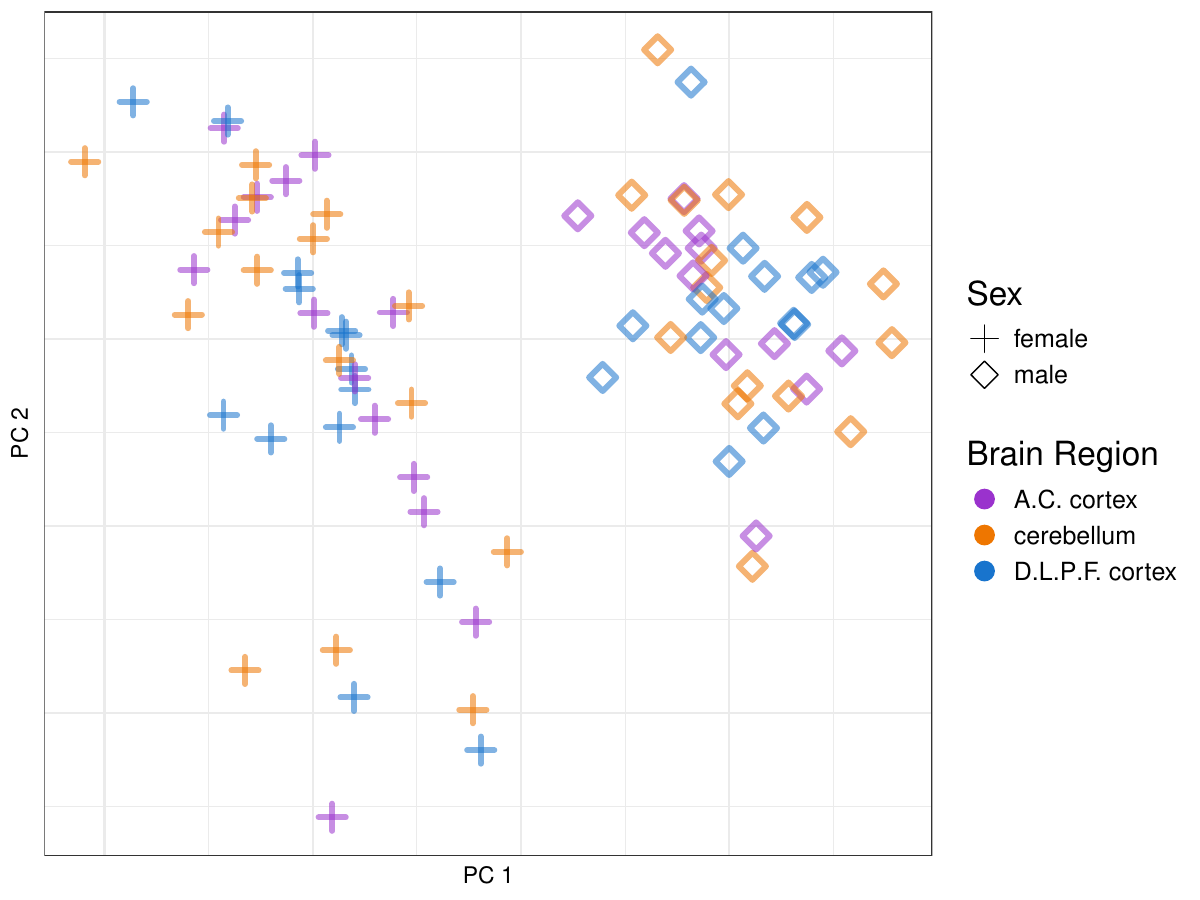}     
    \end{tabular}
    \caption{Like Figure 6, but with $K=15$.}
    \label{fig:K15_svd.xy}
\end{figure}

\begin{figure}[h]
    \centering
    \begin{tabular}{cc}
    \hspace{-1cm} (\textbf{a}) & 
    \hspace{-1cm} (\textbf{b}) \\
    \includegraphics[scale=0.3]{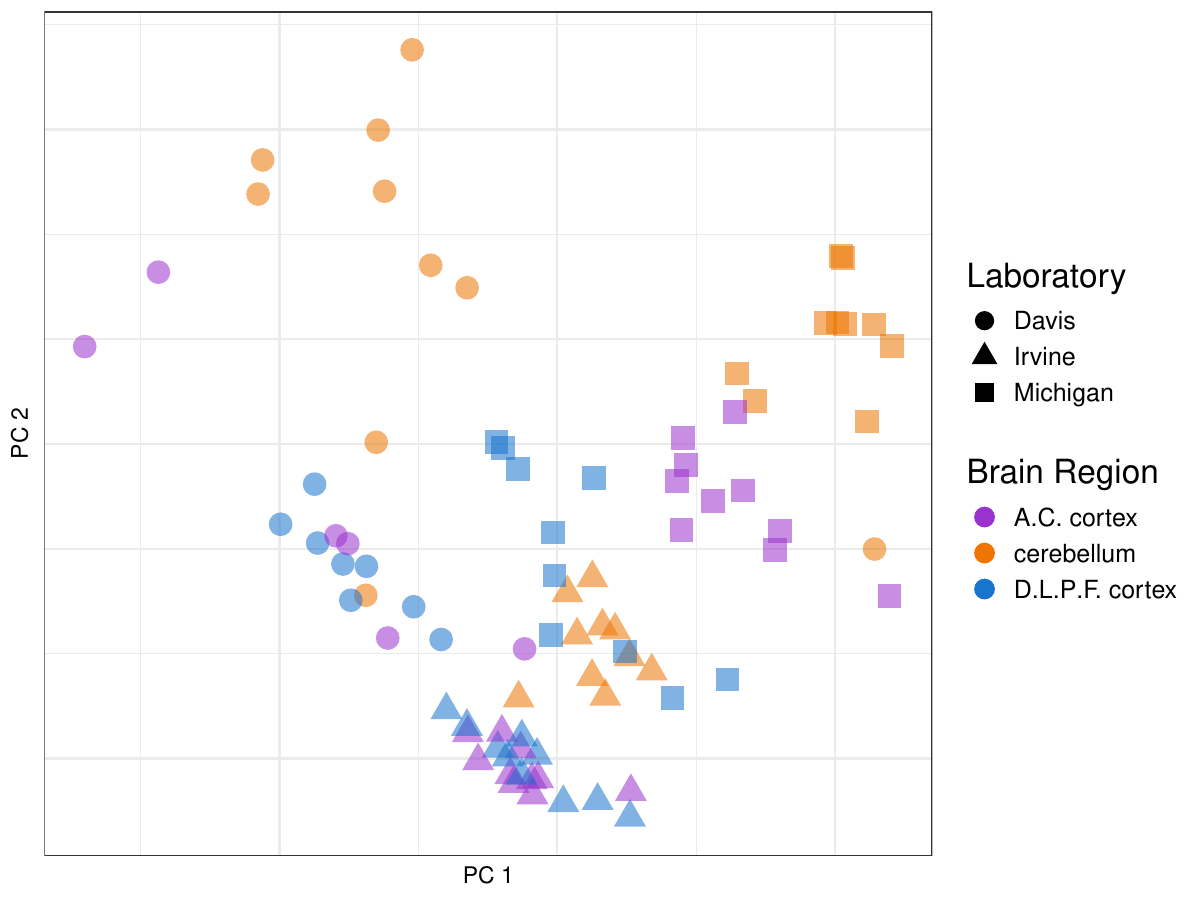} &
    \includegraphics[scale=0.3]{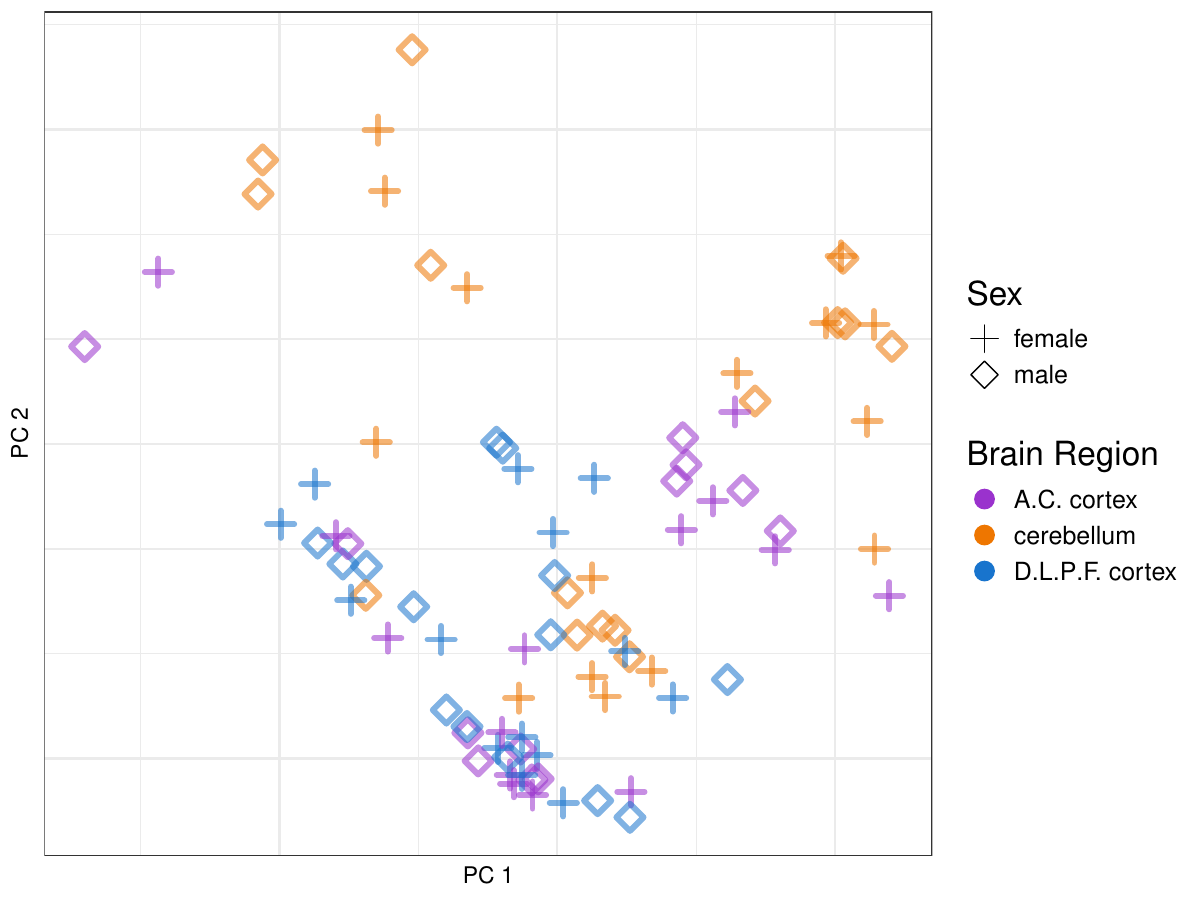} \\    
    \hspace{-1cm} (\textbf{c}) & 
    \hspace{-1cm} (\textbf{d}) \\
    \includegraphics[scale=0.3]{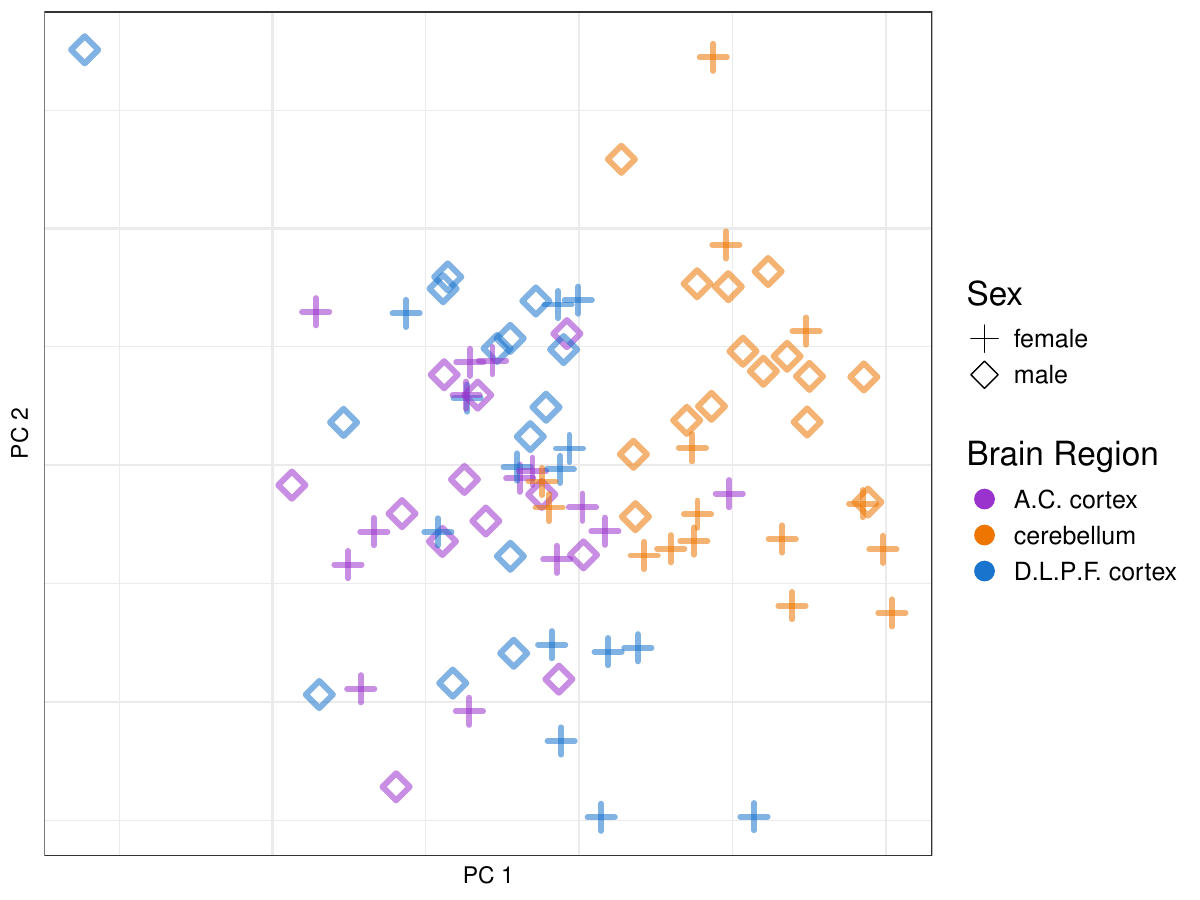}  &
    \includegraphics[scale=0.3]{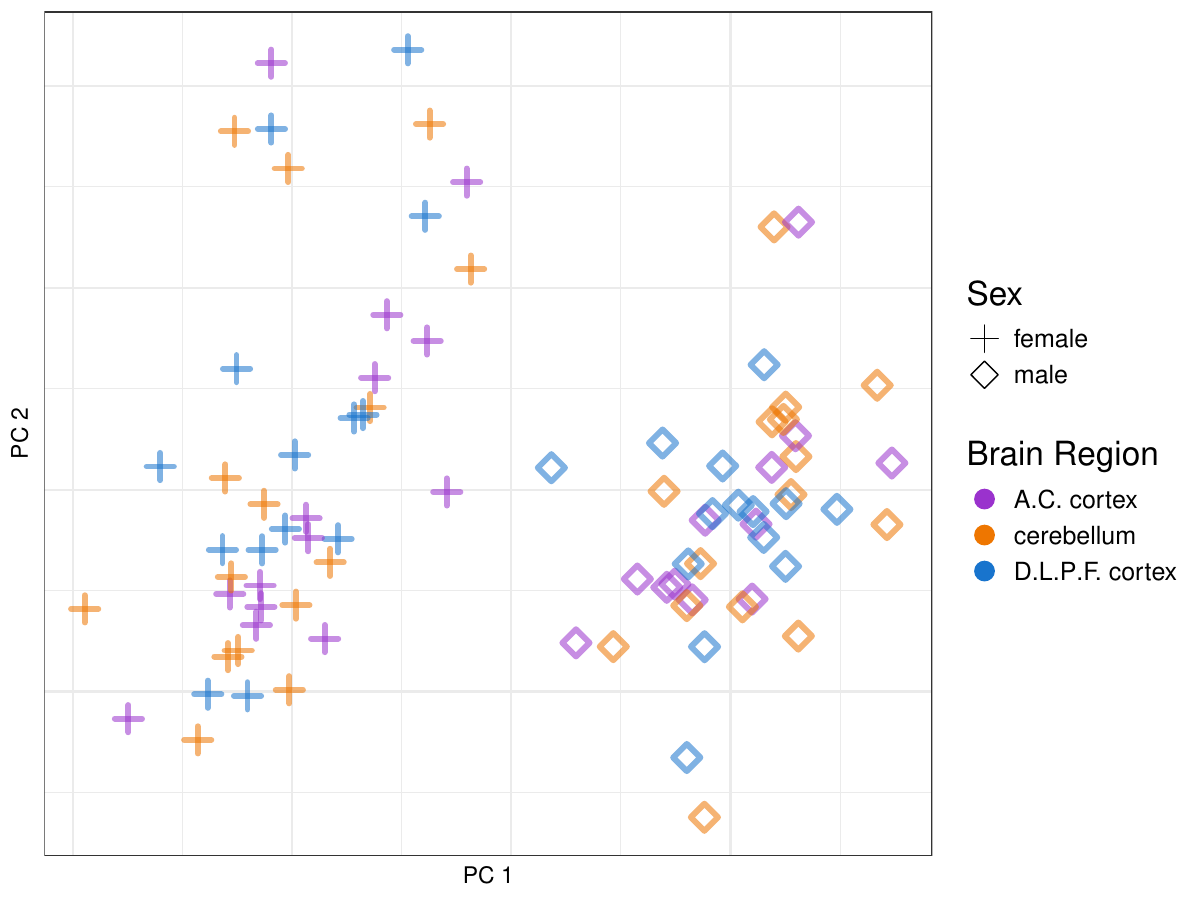}     
    \end{tabular}
    \caption{Like Figure 6, but with $K=20$.}
    \label{fig:K20_svd.xy}
\end{figure}

\begin{figure}[h]
    \centering
    \begin{tabular}{cc}
    \hspace{-1cm} (\textbf{a}) & 
    \hspace{-1cm} (\textbf{b}) \\
    \includegraphics[scale=0.3]{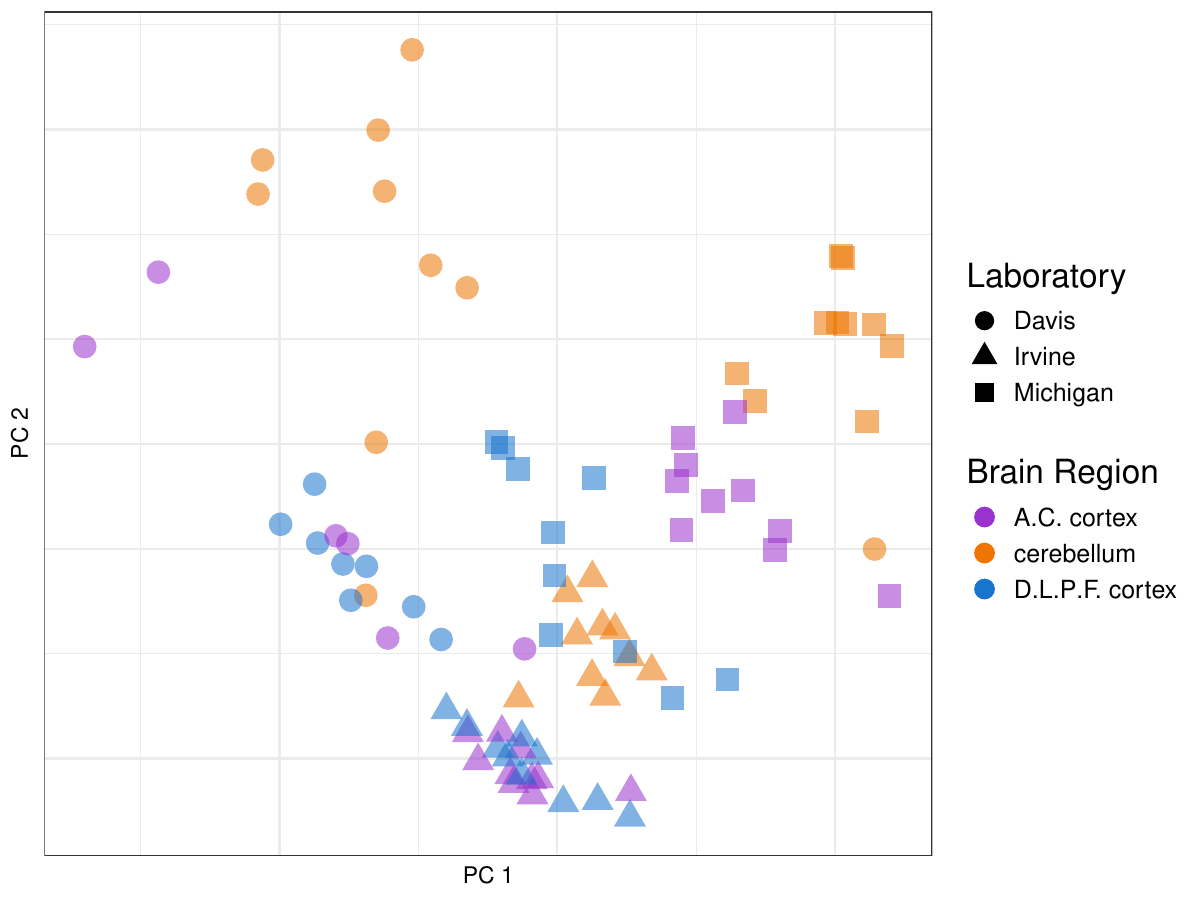} &
    \includegraphics[scale=0.3]{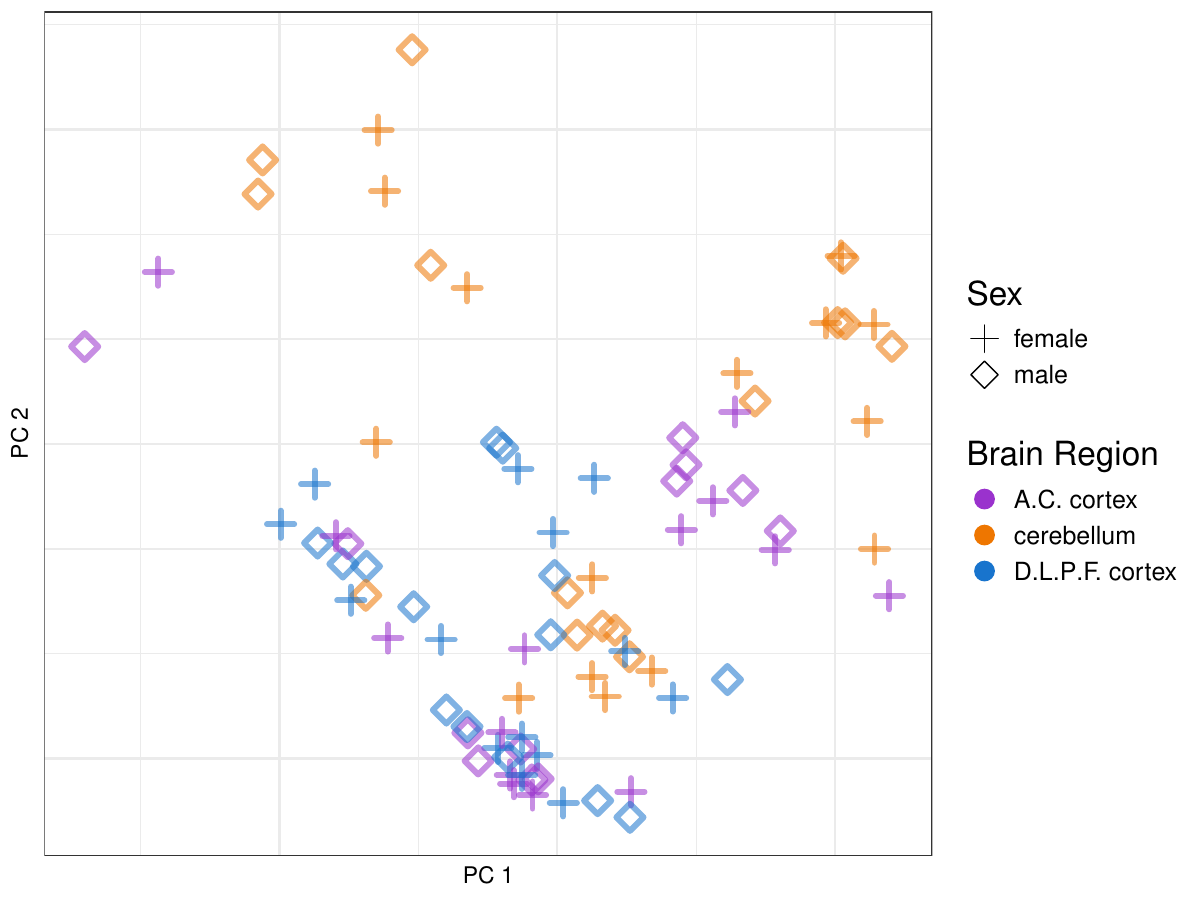} \\    
    \hspace{-1cm} (\textbf{c}) & 
    \hspace{-1cm} (\textbf{d}) \\
    \includegraphics[scale=0.3]{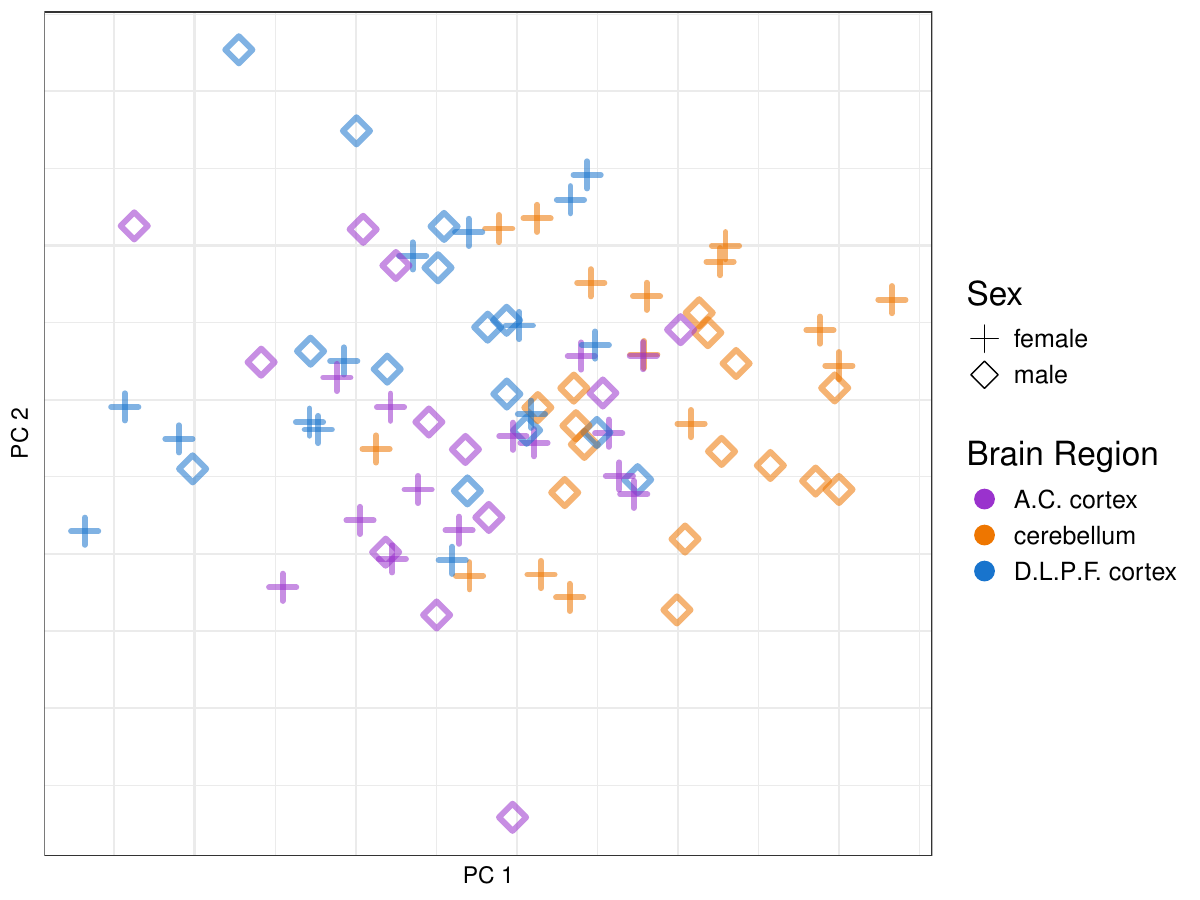}  &
    \includegraphics[scale=0.3]{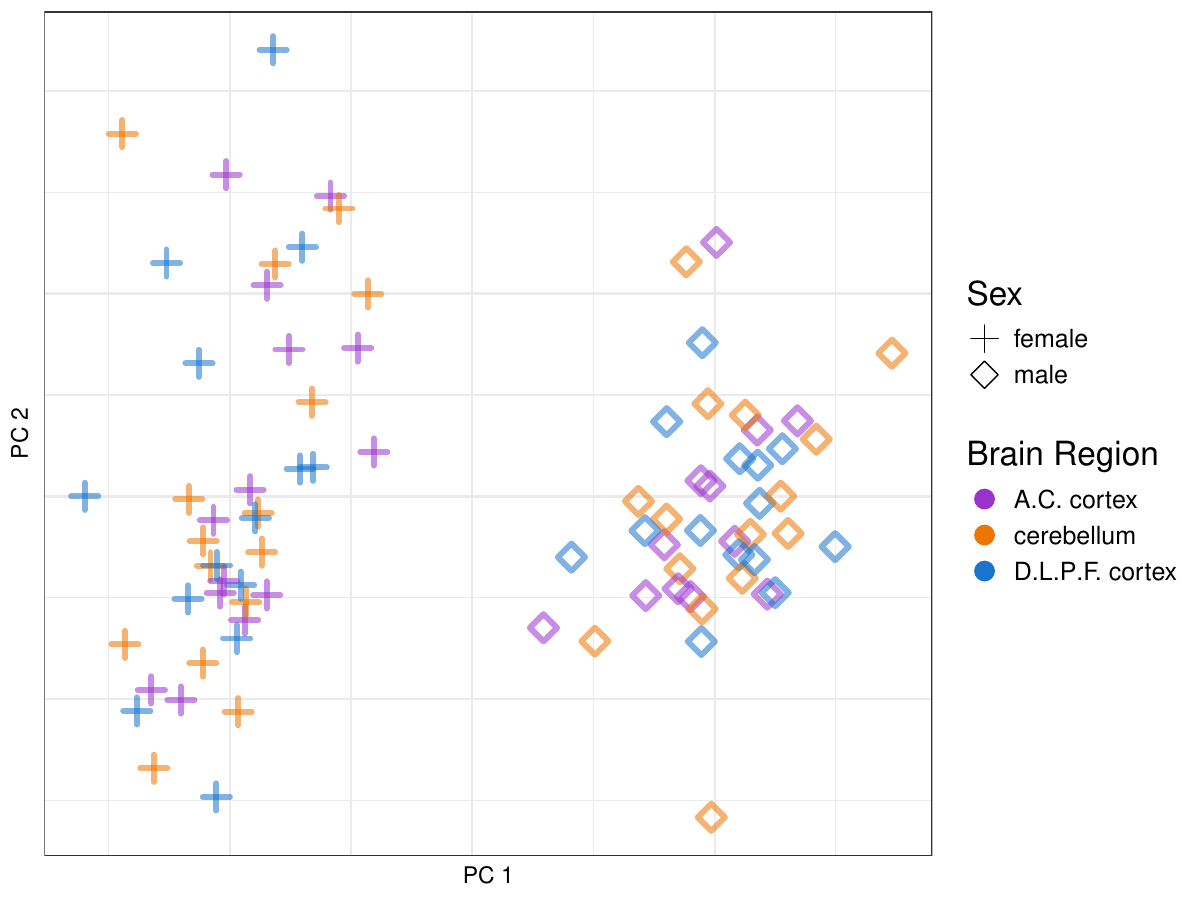}     
    \end{tabular}
    \caption{Like Figure 6, but with $K=30$.}
    \label{fig:K30_svd.xy}
\end{figure}


\end{document}